\documentclass[10pt,a4paper]{amsart}

\usepackage{microtype}
\usepackage{mathtools}     
\usepackage{amssymb}       
\usepackage{amsthm}        
\usepackage{tensor}
\usepackage{bbm}           

\usepackage{graphicx}
\graphicspath{{./images/}}
\usepackage{tikz}
\usetikzlibrary{calc, patterns} 

\usetikzlibrary{matrix}
\usepackage{tikz-cd}

\usepackage[shortlabels]{enumitem}

\usepackage[hidelinks]{hyperref}
\usepackage[nameinlink,noabbrev]{cleveref}
\usepackage{autonum}

\usepackage{float}

\theoremstyle{plain}
\newtheorem{theorem}{Theorem}[section]
\newtheorem{lemma}[theorem]{Lemma}
\newtheorem{prop}[theorem]{Proposition}
\newtheorem{corollary}[theorem]{Corollary}
\newtheorem{question}[theorem]{Question}
\newtheorem{conjecture}[theorem]{Conjecture}

\theoremstyle{definition}
\newtheorem{definition}[theorem]{Definition}

\theoremstyle{remark}
\newtheorem{remark}[theorem]{Remark}
\newtheorem{example}[theorem]{Example}

\newcommand{\bb}[1]{\mathbb{#1}}
\newcommand{\cc}[1]{\mathcal{#1}}
\newcommand{\vp}{\varphi}

\hypersetup{hypertexnames=false}
\usepackage{array}
\usepackage{booktabs}
\usepackage{longtable}

\title{Toward the abelian category of positive Lagrangians}
\author{Asaf Kapota and Jake Solomon} 
\date{\today}

\begin{document}

\begin{abstract}
    We introduce a geometric construction whose terminal
    configuration determines the kernel and cokernel of a
    morphism between positive Lagrangian branes in the
    derived Fukaya category of a Calabi--Yau manifold.
    Thus, if the construction always terminates, the
    full subcategory of positive Lagrangian branes is
    abelian, as conjectured by Joyce.
    We show that this abelian category admits
    Harder--Narasimhan filtrations with respect to the
    central charge defined by the holomorphic volume form.
    Under appropriate split-generation assumptions, we
    construct a Bridgeland stability condition whose stable
    objects are isomorphic to special Lagrangian branes,
    as predicted by the
    Thomas--Yau conjecture.

    The geometric construction is based on a one-parameter
    family of special Lagrangian submanifolds with boundary.
    One part of the boundary is contained in the Lagrangian
    underlying the domain of the morphism and the other in the Lagrangian underlying the codomain.
    Surgeries occur as the family encounters singularities
    or obstruction walls.
    We prove that the construction always terminates
    for the two-torus and provide higher dimensional examples.
\end{abstract}

\maketitle

\tableofcontents

\section{Introduction}

Let $(X,\omega,J,\Omega)$ be a Calabi--Yau manifold.
An oriented immersed Lagrangian $L$ is \emph{positive}
if $\operatorname{Re}\Omega|_L>0$, and \emph{special
    of phase zero} if $\operatorname{Im}\Omega|_L=0$
and $\operatorname{Re}\Omega|_L>0$.
We call $L$ \emph{imaginary special} if
$\operatorname{Re}\Omega|_L=0$ and
$\operatorname{Im}\Omega|_L>0$. Lagrangians are implicitly equipped with a spin structure and a
\emph{grading}: a continuous function
$\theta_L:L\to\mathbb R$ such that
$e^{-i\theta_L}\Omega|_L$ is a positive real volume form. The notation $L[1]$ denotes the shift that
adds $\pi$ to the grading and reverses the orientation.
A \emph{Lagrangian brane} $\mathcal L$ is a Lagrangian $L$
together with a rank one local system on $L$ and a bounding cochain in the Fukaya $A_\infty$ algebra of $L$; see Definition~\ref{def:immersed-lagrangian-brane}. We write $\mathcal L[1]$ for the
brane whose underlying graded Lagrangian is $L[1]$ equipped with the same bounding cochain as $\mathcal L$.
A brane is \emph{positive} if its underlying Lagrangian is positive.

Thomas~\cite{Thomas-moment-maps}
and Thomas--Yau~\cite{ThomasYau}
conjectured that an appropriate notion of stability should characterize the
Hamiltonian isotopy classes containing a unique
special Lagrangian representative.
They proposed Lagrangian mean curvature flow as a
geometric mechanism for constructing such representatives.
In its categorical formulation, the Thomas--Yau
conjecture relates Bridgeland stability~\cite{Bridgeland-stability} to the
existence of special Lagrangian
representatives within isomorphism classes of the
derived Fukaya category
\cite[Conjecture~3.2]{Joyce-stability}.
As part of this program,
Joyce~\cite[Conjecture~3.5]{Joyce-stability}
conjectures that, under suitable hypotheses,
positive Lagrangian branes form an abelian subcategory
of the derived Fukaya category.

In this paper, we propose a geometric construction of kernels and
cokernels for morphisms between positive immersed Lagrangian branes. The construction is based on evolving an immersed brane representing the mapping cone of a morphism, preserving the
isomorphism class of the brane.
The geometric evolution uses families of imaginary-special
Lagrangians with boundary, building on the cylinder
constructions of the second author and
Yuval~\cite{SY-cylinders,SY-webbing}.
If this evolution reaches a terminal configuration,
the terminal brane determines the kernel and cokernel. Under our continuation and termination conjectures,
we establish the abelianity predicted by Joyce and
prove that this category admits Harder--Narasimhan
filtrations for the geometric central charge. Assuming these conjectures for every phase rotation of
$\Omega$ and suitable split-generation hypotheses,
we further construct a Bridgeland stability condition
on the derived Fukaya category whose semistable objects
admit special Lagrangian brane representatives.
For the two-torus, we establish the continuation
and termination conjectures, proving that every morphism admits such an evolution terminating in finite time.

\subsection{Web evolution and abelianity}
Let $\mathcal D$ denote the derived immersed Fukaya category of $X$,
and let $\mathcal A\subset\mathcal D$ be the full
subcategory whose objects are isomorphic to finite
direct sums of unobstructed immersed branes with
positive or imaginary-special components; see
Definition~\ref{def:positive-brane-subcategory}.
Let $\mathcal L_0,\mathcal L_1\in\mathcal A$ be
represented by branes supported on $L_0,L_1$,
and consider a morphism
\[
    f:\mathcal L_0\longrightarrow\mathcal L_1.
\]
A closed degree-zero Floer cochain representing $f$,
together with the brane structures of
$\mathcal L_0[1]$ and $\mathcal L_1$,
equips the graded
immersed union
\[
    C_0=L_0[1]\sqcup L_1
\]
with a brane structure.
The resulting brane $\mathcal C_0$ represents
the categorical mapping cone $\operatorname{Cone}(f)$.
A criterion of Beilinson--Bernstein--Deligne
\cite[Propositions~1.2.2 and~1.2.4]{BBD}
shows that if $\mathcal K,\mathcal K^*\in\mathcal A$
fit into an exact triangle
\begin{equation}
    \mathcal K[1]
    \longrightarrow
    \operatorname{Cone}(f)
    \longrightarrow
    \mathcal K^*
    \longrightarrow
    \mathcal K[2],
    \phantomsection
    \label{eq:terminal-web-triangle}
\end{equation}
then the induced maps
\[
    \mathcal K\longrightarrow\mathcal L_0,
    \qquad
    \mathcal L_1\longrightarrow\mathcal K^*
\]
are a kernel and a cokernel of $f$ in $\mathcal A$,
respectively.
If every morphism in $\mathcal A$ admits such a
triangle, then
$\mathcal A$ is abelian.

We define an evolution for $\mathcal C_0$,
preserving its isomorphism class. If this evolution continutes until it reaches a terminal configuration, its terminal brane determines an exact triangle of the
form~\eqref{eq:terminal-web-triangle}.
We describe the initialization, continuation,
and termination of this evolution below.

\subsubsection{Initialization and regular evolution}
\begin{figure}[htbp]
    \centering
    \includegraphics[scale=0.12]{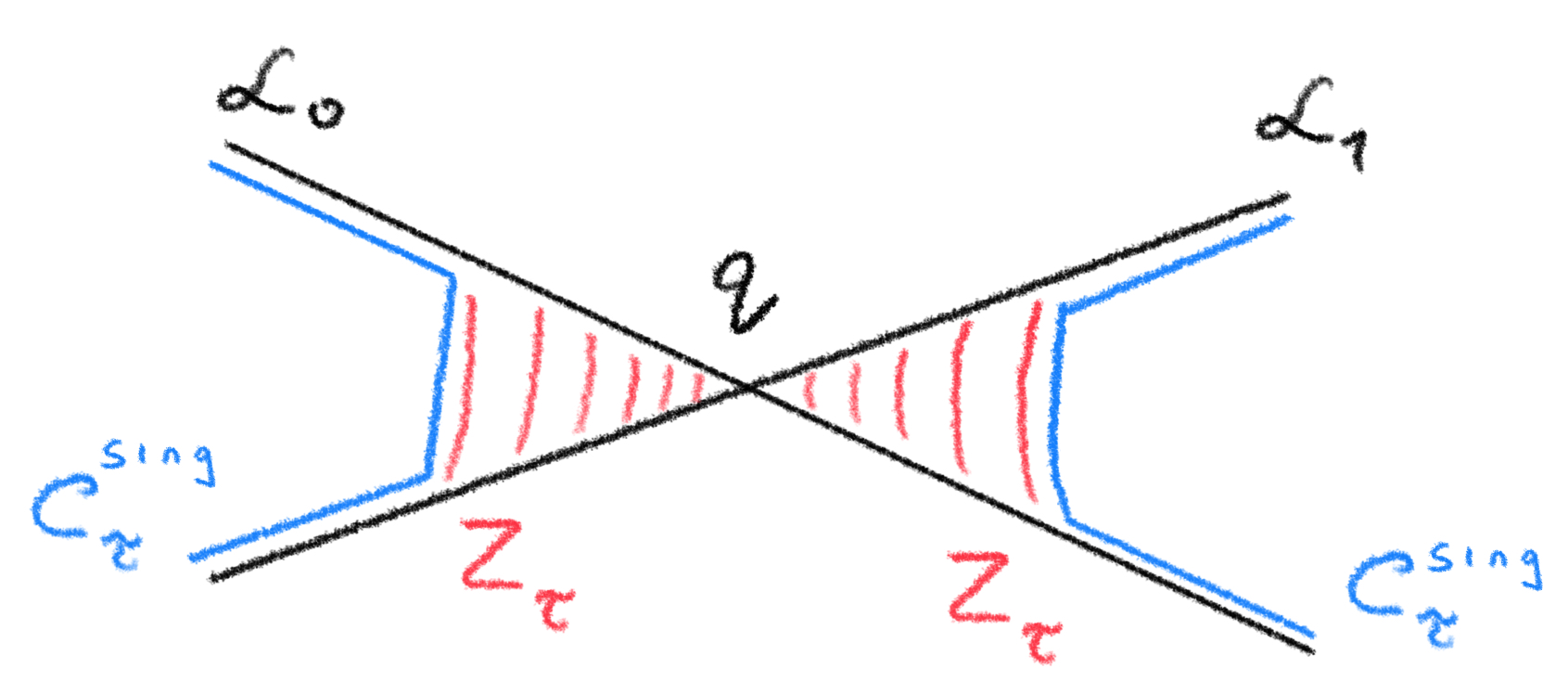}
    \caption{Schematic picture of the web-cone construction near
        an intersection $q$ of Maslov index zero.
        The Lagrangians $L_0,L_1$ are shown in black and
        successive imaginary-special cylinders $Z_\tau$ in red.
        Gluing $Z_\tau$ to the unswept regions
        $R_0(\tau)[1]$ and $R_1(\tau)$ gives the singular web cone (blue), drawn slightly displaced for clarity.
        Rounding its corners gives the web cone $C_\tau$.}
    \label{fig:web-cone}
\end{figure}

The starting point is the special Lagrangian
webbing construction of \cite[Theorem~1.1(a)]{SY-webbing}.
Suppose that $L_0$ and $L_1$ intersect transversely at a point $q$
of Maslov index zero. The construction then gives a one-parameter family of imaginary-special Lagrangian
cylinders $\{Z_\tau\}_{\tau \in (0,\varepsilon)}$, each diffeomorphic to
$S^{n-1}\times[0,1]$, where $n=\dim_{\mathbb C}X$.
The boundary components $S^{n-1}\times\{0\}$ and
$S^{n-1}\times\{1\}$ lie on $L_0$ and $L_1$,
respectively, and the cylinders collapse to $q$
as $\tau\to0$.
As the cylinders move away from $q$, their boundaries sweep out
regions in $L_0$ and $L_1$. Define the \emph{residues} $R_i(\tau)$ to be the closure of
the complement in $L_i$ of the region swept out up to time $\tau$,
for $i=0,1$.
Gluing the residues $R_i(\tau)$ to $Z_\tau$ along the corresponding boundary
components gives the piecewise smooth graded Lagrangian
\[
    C_\tau^{\mathrm{sing}}
    =
    R_0(\tau)[1]\cup Z_\tau\cup R_1(\tau).
\]
Here $R_i(\tau)$ carries the grading inherited from $L_i$,
while $Z_\tau$ carries its natural grading of phase $\pi/2$.
We call the resulting immersion the
\emph{singular web cone}; see Figure~\ref{fig:web-cone}.
Rounding its corners gives an immersed Lagrangian $C_\tau$,
the \emph{web cone}. The roundings can be chosen coherently so
that the resulting family is an exact deformation.
As $\tau\to0$, the cylinder collapses to $q$ and the singular
web cone converges to $C_0$.

We call the moving imaginary-special Lagrangians with boundary
\emph{web fronts}.
They form the moving pieces of the singular web cone, while
their boundary motion determines the residues
$R_0(\tau)$ and $R_1(\tau)$.
Thus, evolving the fronts determines the geometric evolution
of the assembled web cone.

The brane structure of $\mathcal C_0$ supplies the initial
brane data for the web cone.
These data must be transported along the evolution of the
web cones so that the resulting branes $\mathcal C_\tau$ satisfy
\[
    \mathcal C_\tau\simeq \mathcal C_0 \simeq \operatorname{Cone}(f).
\]
We call such an initial evolution, together with
its brane data and comparison, an
\emph{initialization of the web evolution of $f$};
see Definition~\ref{def:web-initialization}.
\begin{conjecture}[Initialization]
    \label{conj:intro-initialization}
    Every closed degree-zero Floer cochain $f$
    representing a morphism in $\mathcal A$
    admits an initialization of its web evolution.
\end{conjecture}

The existence of the web front underlying the initialization is given by~\cite[Theorem~1.1(a)]{SY-webbing}. So, the content of Conjecture~\ref{conj:intro-initialization} is purely Floer theoretic. We expect
the required brane structures and comparison
with $\operatorname{Cone}(f)$ to follow by
arguments similar to those of
Palmer--Woodward~\cite{PalmerWoodwardSurgery}
and Fukaya--Oh--Ohta--Ono~\cite{FOOOChapter10}.
Alternatively, it should be possible to use the cobordism approach to surgery and
mapping cones developed by
Biran--Cornea~\cite{BiranCorneaCobordismI,BiranCorneaFukaya}

\subsubsection{Geometric degenerations and obstruction walls}
Continuing the evolution of the web cone and its brane structure
requires addressing two kinds of events.

First, a web front may develop a singularity.
We describe local surgery models for these
degenerations in Section~\ref{subsec:geometric-degenerations}.
Surgery may change the topology of the web cone.
The surgery gives rise to a \emph{compatible continuation} if there exists a matching of the
pieces of the incoming and outgoing web cones unaffected by the surgery that provides the necessary geometric data to continue the web evolution;
see Definition~\ref{def:surgery-compatible}.
The brane structure must be transferred from the incoming to the outgoing
web cone of a compatible continuation while preserving the represented Fukaya object; see Definition~\ref{def:brane-surgery-transfer}.

Second, even if the underlying special Lagrangian of the web front can be smoothly evolved, it may become impossible to further transport the bounding cochain because a valuation positivity condition is violated. Following Joyce's discussion
\cite[Section~3.4]{Joyce-stability}, we call such an event an
\emph{obstruction wall}.
Compatible continuation past an obstruction wall requires modifying the web cone by surgery
or algebraic reduction as described in Subsection~\ref{subsec:obstruction-walls}. Again, the brane structure must be transferred to the outgoing web cone.

After surgery, the web front need no longer be a cylinder, so the moduli space of nearby Lagrangians with boundary need no longer be one dimensional. Nonetheless, the motion of the web front is determined uniquely as follows.
For a smooth front $Z$, the nearby special Lagrangians
with boundary on $L_0$ and $L_1$ form a smooth moduli
space, as discussed in~\cite[Remark~4.9]{SY-cylinders}
and proved by
Pidaparthy~\cite[Theorem~2.2]{Pidaparthy-deformations}.
The tangent space to the moduli space at $Z$ is the space of harmonic
one-forms on $Z$ whose pullbacks to $\partial Z$ vanish,
which is isomorphic to
$H^1(Z,\partial Z;\mathbb R)$.
The motion of a front $Z$ is determined
by the harmonic $1$-form given by the differential of the unique harmonic function $\sigma: Z \to \mathbb R$,
taking values $0$ and $1$ on the components of $\partial Z$ contained in
$L_0$ and $L_1$, respectively. 

An immersed Lagrangian that cannot be written as the union of non-empty immersed Lagrangians is called irreducible.
If the boundary of an irreducible component $Z'$ of $Z$ is contained entirely in only one of
$L_0$ or $L_1$, then $\sigma|_{Z'}$ is constant and the
motion of $Z'$ halts.
We retain such components, together with any closed
components, as \emph{inactive web fragments}.

An evolution extending a web initialization of $f$,
with regular brane transport between events and
compatible continuation and brane transfer at events,
is called a \emph{special Lagrangian web for $f$};
see Definition~\ref{def:f-admissible-web-evolution}.

\begin{conjecture}[Continuation]
    \label{conj:intro-continuation}
    For each morphism $f:\mathcal L_0\to\mathcal L_1$
    in $\mathcal A$, after generic
    Hamiltonian perturbations of the strictly positive
    components of $L_0$ and $L_1$, with the brane
    structures and $f$ transported accordingly,
    an initialization as in
    Conjecture~\ref{conj:intro-initialization}
    can be continued through isolated geometric
    degenerations and obstruction walls by
    compatible continuation and brane transfer.
\end{conjecture}

For arbitrary representatives of $L_0$ and $L_1$,
we do not expect to control all singularities of
the web evolution.
We expect generic Hamiltonian perturbations to
avoid singularities that prevent compatible continuation
by surgery, while those that occur generically
should admit such a continuation.
This does not require the local models
in Section~\ref{sec:singular-events} to exhaust the possible generic degenerations.

\begin{theorem}
    \label{thm:intro-elliptic-continuation}
    Conjectures~\ref{conj:intro-initialization}
    and~\ref{conj:intro-continuation}
    hold on the two-torus.
\end{theorem}

\subsubsection{Termination and abelianity}
As the evolution continues through singular events,
active fronts may become inactive fragments. If no active fronts remain, we call the web
\emph{geometrically exhausted}; its cone is
then stationary. However, its bounding cochain
may continue to evolve so that the web evolution subsequently reaches another obstruction wall, at which surgery may
initiate new active fronts. A similar phenomenon was observed by Joyce for bounding cochains of Lagrangians evolving by mean curvature flow
\cite[compare Section~3.4]{Joyce-stability}.
We call the web evolution \emph{terminal from time $\tau_*$}
if, after compatible continuation at $\tau_*$,
the web has no active fronts and stationary brane
continuation encounters no further obstruction wall;
see Definition~\ref{def:brane-terminal-web}.

\begin{conjecture}[Termination]
    \label{conj:intro-termination}
    \label{conj:web-cone-termination}
    For each morphism $f:\mathcal L_0\to\mathcal L_1$
    in $\mathcal A$, after generic Hamiltonian
    perturbations as in
    Conjecture~\ref{conj:intro-continuation},
    a special Lagrangian web for $f$ becomes terminal
    in finite time, after finitely many geometric
    singularities and obstruction walls.
\end{conjecture}

\begin{theorem}
    \label{thm:intro-elliptic-termination}
    Conjecture~\ref{conj:intro-termination}
    holds on the two-torus.
\end{theorem}

The proofs of
Theorems~\ref{thm:intro-elliptic-continuation}
and~\ref{thm:intro-elliptic-termination}
are given in Section~\ref{sec:elliptic-curve},
using results from the forthcoming companion
papers~\cite{KapotaWebGeometry,KapotaWebFloer};
see Subsection~\ref{subsec:intro-elliptic-curve}.
Let a special Lagrangian web for $f$ be terminal
from time $\tau_*$.
We show that smoothing an exact deformation
of its terminal configuration gives an unobstructed
brane $\mathcal C_{\tau_*}\simeq\operatorname{Cone}(f)$
with underlying immersed Lagrangian
\[
    C_{\tau_*} = K[1] \sqcup K^*,
\]
where each connected component of $K$ or $K^*$ is a positive
Lagrangian or a closed imaginary-special Lagrangian.
These Lagrangians are constructed from the terminal
web configuration as follows.

\begin{figure}[htbp]
    \centering
    \includegraphics[scale=0.12]{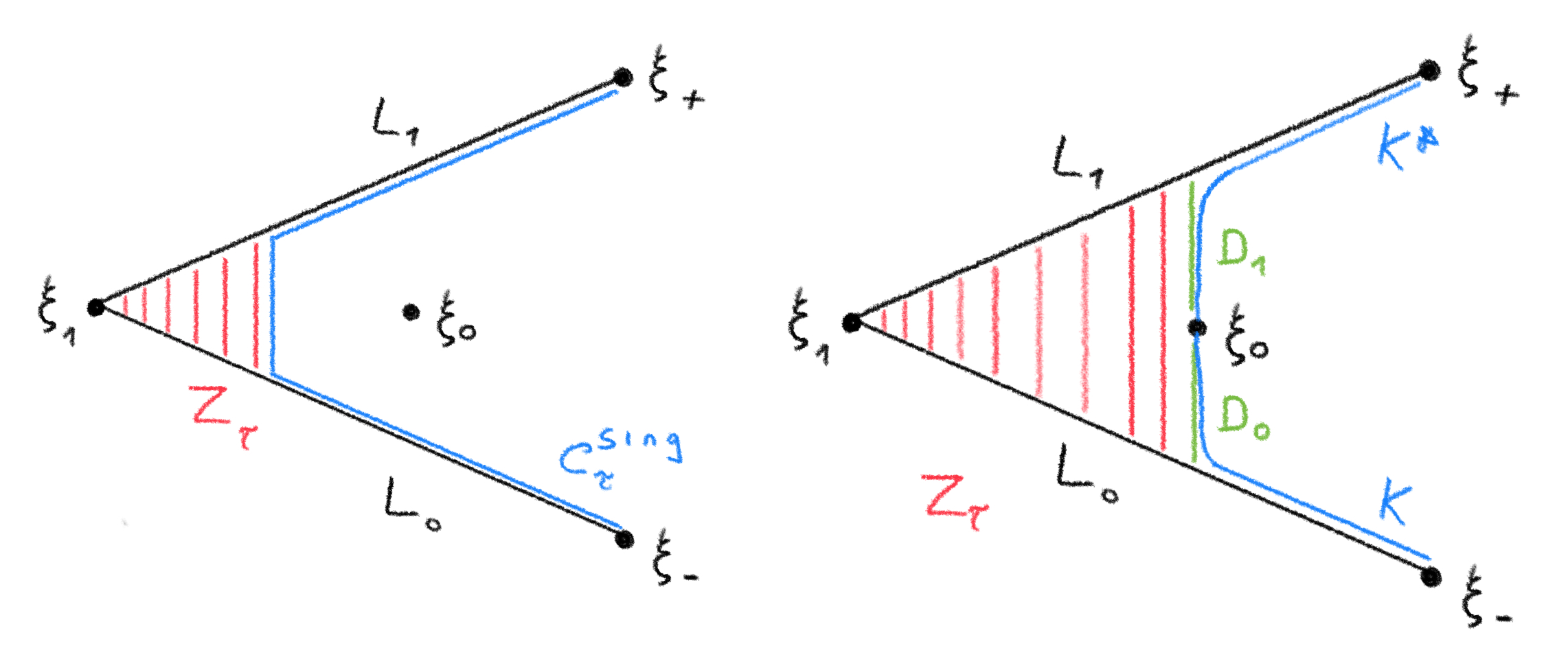}
    \par\medskip
    \includegraphics[scale=0.12]{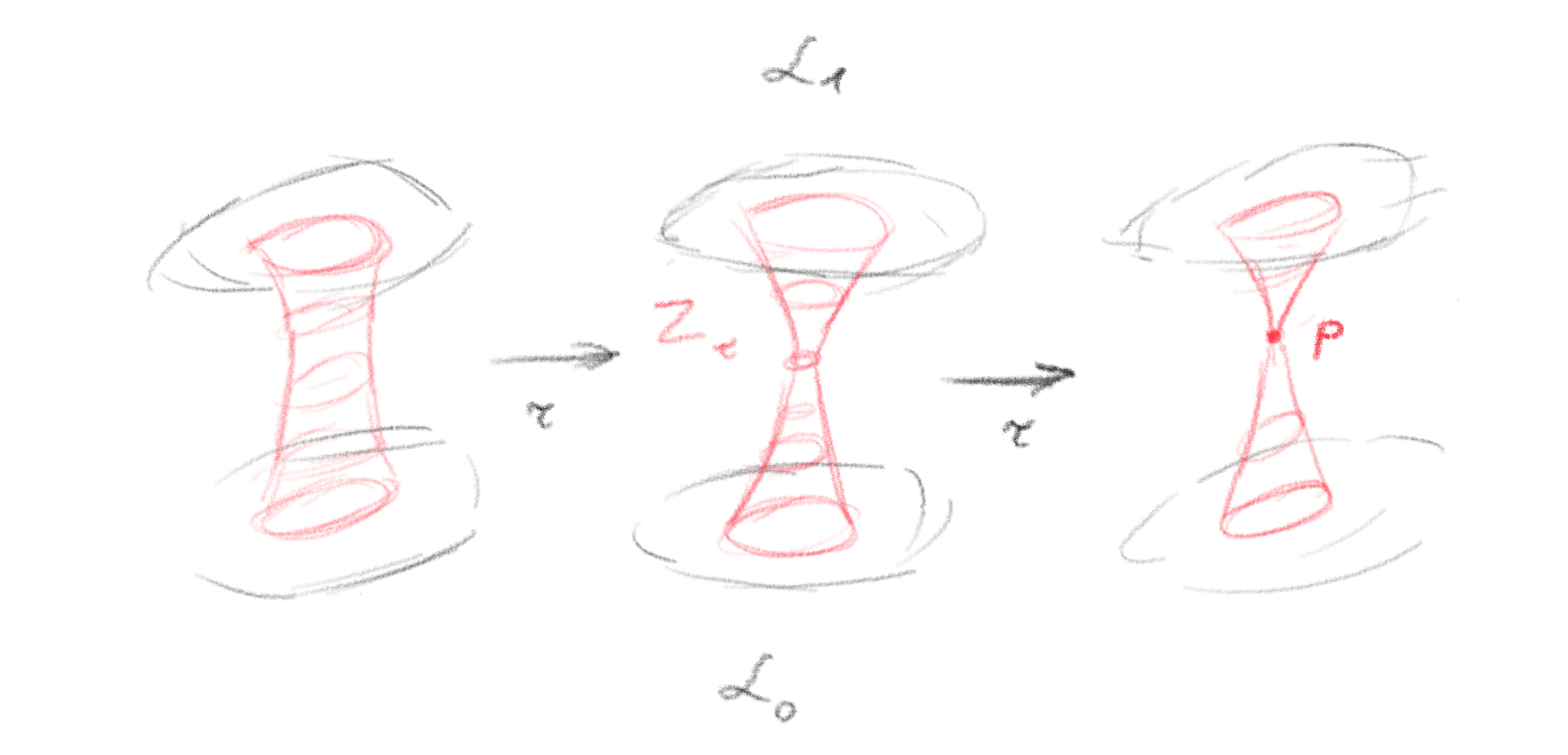}
    \caption{The web evolution in
        Example~\ref{ex:intro-milnor-fiber}.
        \emph{Top:} Projections under $\pi$.
        On the left, successive web fronts $Z_\tau$ are shown
        in red and the assembled web cone $C_\tau$ in blue.
        On the right, the inactive discs $D_0,D_1$ are shown
        in green; gluing them to the corresponding unswept
        regions and smoothing gives the spheres $K,K^*$
        shown in blue.
        \emph{Bottom:} The interior degeneration of the web
        front. As the projected segment reaches $\xi_0$,
        an interior circle of the cylinder collapses to
        $p=(0,0,\xi_0)$, leaving two discs meeting at $p$,
        with boundary on $L_0$ and $L_1$, respectively.}
    \label{fig:milnor-web}
\end{figure}

We glue the closures of the terminal residues
$R_0(\tau_*)$ and $R_1(\tau_*)$ to the inactive fragments with
boundary on $L_0$ and $L_1$, respectively.
The corners are smoothed and then an exact deformation is applied so that the resulting Lagrangians are positive; see Figure~\ref{fig:milnor-web} and Section~\ref{subsec:terminal-smoothing}.
Together with the retained closed fragments, these deformed smoothings give $K$
and $K^*$, respectively.
The following proposition describes the categorical structure
of this terminal configuration; its proof is given in Subsection~\ref{subsec:terminal-triangle}.

\begin{prop}[The terminal triangle]
    \label{prop:terminal-web-triangle}
Consider a special Lagrangian web for $f$
that is terminal from time $\tau_*$.
Its terminal brane data induce unobstructed
immersed branes $\mathcal K$ and $\mathcal K^*$
on the Lagrangians $K$ and $K^*$ described above.
    Each connected component of $K$ or $K^*$ is a
    positive Lagrangian or a closed imaginary-special
    Lagrangian.
    These branes fit into an exact triangle of the
    form~\eqref{eq:terminal-web-triangle}
    in $\mathcal D$.
\end{prop}

The terminality condition that no further obstruction walls occur forces
the part of the terminal bounding cochain of $\mathcal C$ directed from $\mathcal K[1]$
to $\mathcal K^*$ to vanish. The remaining parts of the bounding cochain determine the extension
expressed by the triangle~\eqref{eq:terminal-web-triangle}.
The triangle and the categorical criterion
of~\cite[Propositions~1.2.2 and~1.2.4]{BBD} give the following conditional result.
Its proof is given in
Section~\ref{subsec:positive-branes-abelian}.

\begin{theorem}[Abelianity]
    \label{thm:intro-abelianity}
    Assume Conjectures~\ref{conj:intro-initialization},~\ref{conj:intro-continuation}
    and~\ref{conj:intro-termination}.
    Then $\mathcal A$ is an abelian subcategory
    of $\mathcal D$.
    For each morphism $f$ in $\mathcal A$, its kernel
    and cokernel are represented by the terminal branes
    $\mathcal K$ and $\mathcal K^*$, respectively.
    Every short exact sequence in $\mathcal A$ extends
    to a distinguished triangle in $\mathcal D$.
\end{theorem}

The following example illustrates the geometric
evolution and the construction of $K$ and $K^*$ in a simple case.

\begin{example}[A Milnor fiber]
    \label{ex:intro-milnor-fiber}
    Let $P(z)=z(z^3+1)$ and consider the $A_3$ Milnor
    fiber of complex dimension two,
    \[
        X=\bigl\{(u,v,z)\in\mathbb C^3:
        u^2+v^2=P(z)\bigr\},
    \]
    equipped with the induced Kähler form and the
    holomorphic volume form $\Omega$, given on
    $\{v\neq0\}$ by
    \[
        \Omega=\frac{du\wedge dz}{2v}.
    \]
    Let $\pi:X\to\mathbb C$ be the projection
    $\pi(u,v,z)=z$.
    The zeros of $P$ are
    \begin{equation}
        \xi_0 = 0, \quad
        \xi_1 = -1, \quad
        \xi_\pm = \frac12 \pm \frac{\sqrt{3}}{2}i.
    \end{equation}
    The group $SO(2)$ acts on $X$ by real rotations
    of $(u,v)$, with moment map
    \[
        \mu(u,v,z)=\operatorname{Im}(u\overline v).
    \]
    For the straight segments
    \[
        \gamma_0=[\xi_1,\xi_-],
        \qquad
        \gamma_1=[\xi_1,\xi_+],
    \]
    define
    \[
        L_i:=\mu^{-1}(0)\cap\pi^{-1}(\gamma_i),
        \qquad i=0,1.
    \]
    See Figure~\ref{fig:milnor-web}. Over each interior point of $\gamma_i$, the fibre
    of $\pi|_{L_i}$ is a single circle orbit.
    These circles collapse to the fixed points
    $(0,0,\xi)$ at the endpoints $\xi$ of $\gamma_i$.
    The resulting smooth Lagrangian two-spheres
    are called \emph{matching spheres}.

    With their positive orientations, $L_0$ and $L_1$
    have constant phases $-\pi/6$ and $\pi/6$,
    respectively.
    They intersect transversely at the unique point
    $q=(0,0,\xi_1)$, and their natural positive gradings
    give $m(q;L_0,L_1)=0$.
    We initiate a webbing at $q$. The imaginary-special
    Lagrangian cylinders project under $\pi$ to straight
    vertical segments moving to the right.
    During smooth evolution, we assemble the web cone $C_\tau$
    at each time $\tau$ as described above.

    The cylinders evolve smoothly until the projected
    segment reaches the root $\xi_0$, at time $\tau_*$.
    At this time, a geometric singularity occurs due to the collapse of a family of circles in the cylinders to the point $p=(0,0,\xi_0)$.
    The limiting web front consists of two discs $D_0$ and
    $D_1$ meeting at $p$, with boundary on $L_0$ and $L_1$,
    respectively.
    Compatible continuation through this singularity retains both
    discs as inactive fragments.
    Neither has boundary on both Lagrangians, so no active
    fronts remain and the geometric evolution stops.

    Gluing these discs to the corresponding unswept regions,
    then smoothing the corners and applying the exact
    deformations described above, gives
    \[
        C_{\tau_*}=K[1]\cup K^*.
    \]
    Here $K$ and $K^*$ are positive Lagrangian matching
    spheres whose projected paths join $\xi_0$ to $\xi_-$
    and $\xi_0$ to $\xi_+$, respectively; see
    Figure~\ref{fig:milnor-web}.
\end{example}

\subsection{Harder--Narasimhan filtrations and stability}

Assume that Conjectures~\ref{conj:intro-initialization}, \ref{conj:intro-continuation}
and~\ref{conj:intro-termination} hold.
Then $\mathcal A$ is abelian by
Theorem~\ref{thm:intro-abelianity}.
We study the stability properties of $\mathcal A$.
The existence of Harder--Narasimhan filtrations for
the geometric central charge is part of Joyce's
conjectural stability picture
\cite[Conjectures~3.2 and~3.5]{Joyce-stability}.
We use the geometric central charge
\[
    Z(\mathcal L)=\int_L\Omega.
\]
The open--closed construction shows that these periods
are invariant under isomorphism and additive on short
exact sequences; see
Remark~\ref{rem:homological-compatibility}.
The phase conditions defining $\mathcal A$ ensure
that every nonzero object has a well-defined phase
\[
    \phi_Z(\mathcal L)
    :=\arg Z(\mathcal L)\in(-\pi/2,\pi/2].
\]
Thus $Z$ is a stability function on $\mathcal A$.
We recall the relevant definitions, following
\cite[Definitions~2.2 and~2.3]{Bridgeland-stability}.

\begin{definition}
    \label{def:intro-semistability-HN}
    A nonzero object $\mathcal L\in\mathcal A$ is
    \emph{semistable} if
    $\phi_Z(\mathcal M)\leq\phi_Z(\mathcal L)$
    for every nonzero subobject
    $\mathcal M\subset\mathcal L$.
    A subobject violating this inequality is
    \emph{destabilising}; a nonzero quotient
    $\mathcal L\twoheadrightarrow\mathcal Q$
    is destabilising if
    $\phi_Z(\mathcal Q)<\phi_Z(\mathcal L)$.

    A \emph{Harder--Narasimhan filtration} of $\mathcal L$
    is a finite filtration by subobjects in $\mathcal A$,
    \[
        0=\mathcal L_0\subset\mathcal L_1
        \subset\cdots\subset\mathcal L_m=\mathcal L,
    \]
    whose factors
    $\mathcal F_j=\mathcal L_j/\mathcal L_{j-1}$
    are nonzero and semistable, with
    $\phi_Z(\mathcal F_1)>\cdots>\phi_Z(\mathcal F_m)$.
    We say that $Z$ has the \emph{Harder--Narasimhan property}
    if every nonzero object of $\mathcal A$ admits such a
    filtration.
\end{definition}

If $X$ is noncompact, we impose the geometric
bounds specified in
Subsection~\ref{subsec:harder-narasimhan}.
These bounds are automatic when $X$ is compact.
Our proof of the Harder--Narasimhan property
uses the terminal representatives of
cokernels.
Fix a geometric representative of a nonzero object
$\mathcal L\in\mathcal A$.
Given a destabilising quotient
$\mathcal L\twoheadrightarrow\mathcal Q$, apply the
web evolution to its kernel inclusion
$\mathcal M\hookrightarrow\mathcal L$.
The terminal brane $\mathcal K^*$ represents
$\mathcal Q$.
Its geometry, together with the destabilising phase
inequality, gives a bound on its total volume depending
only on the fixed representative of $\mathcal L$. In the noncompact case, the monotonicity lemma
applied to the special fragments also confines
these representatives to a compact subset
depending only on the fixed representative
of $\mathcal L$.
This yields the following proposition, proved in
Subsection~\ref{subsec:harder-narasimhan}.

\begin{prop}[Bounded destabilising quotients]
    \label{prop:intro-destabilising-volume-bound}
    For every nonzero $\mathcal L\in\mathcal A$,
    there exist a compact subset
    $Y_{\mathcal L}\subset X$ and a constant
    $C_{\mathcal L}>0$ such that every destabilising
    quotient $\mathcal L\twoheadrightarrow\mathcal Q$
    admits a geometric representative
    $\mathcal K^*\simeq\mathcal Q$ satisfying
    \[
        K^*\subset Y_{\mathcal L},
        \qquad
        \operatorname{Vol}_g(K^*)\leq C_{\mathcal L}.
    \]
\end{prop}

For a fixed $\mathcal L$, the representatives
supplied by
Proposition~\ref{prop:intro-destabilising-volume-bound}
lie in a common compact subset and have uniformly
bounded volume.
Their integral homology classes therefore range
over a finite set.
Consequently, the central charges of the
destabilising quotients of $\mathcal L$
also form a finite set.
For a destabilising subobject
$\mathcal M\subset\mathcal L$, the quotient
$\mathcal Q=\mathcal L/\mathcal M$ is destabilising,
and additivity gives
\[
    Z(\mathcal M)=Z(\mathcal L)-Z(\mathcal Q).
\]
Thus the central charges, and hence the phases, of
the destabilising subobjects also form a finite set.
This excludes the infinite chains appearing in
Bridgeland's criterion
\cite[Proposition~2.4]{Bridgeland-stability}
and gives the following theorem.
The proof is given in
Subsection~\ref{subsec:harder-narasimhan}.

\begin{theorem}
    \label{thm:intro-positive-branes-HN}
    The central charge $Z$ has the
    Harder--Narasimhan property on $\mathcal A$.
\end{theorem}

Assume Conjectures~\ref{conj:intro-initialization},~\ref{conj:intro-continuation}
    and~\ref{conj:intro-termination}, so $\mathcal A$ is an
admissible abelian subcategory of $\mathcal D$,
stable under extensions.
By \cite[Proposition~1.3.13]{BBD}, it is therefore
the heart of a bounded $t$-structure on the full
triangulated subcategory
$\cc{D}^b_{\mathcal A}\subset\mathcal D$
consisting of finite extensions of shifts of
objects of $\mathcal A$.
This subcategory is also closed under direct
summands; see Proposition~\ref{prop:BBD-generated-heart}.
The Harder--Narasimhan property and
\cite[Proposition~5.3]{Bridgeland-stability}
therefore give a stability condition on
$\cc{D}^b_{\mathcal A}$ with heart $\mathcal A$
and central charge $\widehat Z=iZ$.
A split-generation hypothesis gives a stability
condition on the whole Fukaya category.

\begin{theorem}[Stability from split-generation]
    \label{prop:intro-stability-from-generation}
    Assume Conjectures~\ref{conj:intro-initialization},~\ref{conj:intro-continuation}
    and~\ref{conj:intro-termination}.
    Suppose that $\mathcal A$ split-generates
    $\mathcal D$.
    Then $\mathcal D$ admits a stability condition
    in the sense of
    \cite[Definition~1.1]{Bridgeland-stability},
    with heart $\mathcal A$ and central charge
    $\widehat Z=iZ$.
\end{theorem}
The proof is given in
Section~\ref{subsubsec:stability-generation}.
These constructions also provide a route to the
correspondence between stability and special
Lagrangian representatives underlying the
Thomas--Yau program
\cite{ThomasYau,Joyce-stability}.
The additional ingredient is to apply the
construction at different phases.
For $c\in\mathbb R$, put
\[
    \Omega_c=e^{i(\pi/2-c)}\Omega,
\]
and let $\mathcal A_c$ be the positive-brane
subcategory defined using $\Omega_c$.
Its defining representatives have source
components with phases in $(c-\pi,c)$, together
with closed special components of phase $c$.
Thus $\mathcal A_{\pi/2}=\mathcal A$.
The continuation and termination conjectures
for $\Omega_c$ make $\mathcal A_c$ an abelian
category, and the same BBD construction identifies
it as a bounded heart in
$\cc{D}^b_{\mathcal A_c}$.
Comparing these hearts gives the following
consequence. We express its conclusion using geometric phases:
a Bridgeland phase $\varphi$ corresponds to the
geometric phase $\pi\varphi-\pi/2$.

\begin{theorem}[Special representatives]
    \label{thm:intro-special-representatives}
    Suppose that
    Conjectures~\ref{conj:intro-initialization},
    \ref{conj:intro-continuation}
    and~\ref{conj:intro-termination}
    hold for every $\Omega_c$.
    Assume that each $\mathcal A_c$ split-generates
    $\mathcal D$. 
    Then, for every $c\in\mathbb R$, a nonzero
    object of $\mathcal D$ is semistable of
    geometric phase $c$ with respect to the
    stability condition of
    Theorem~\ref{prop:intro-stability-from-generation}
    if and only if it admits a special brane
    representative of phase $c$.
\end{theorem}

Here the special brane representative is isomorphic
to the given object in $\mathcal D$, as in Joyce's
categorical formulation~\cite{Joyce-stability}.
Its underlying Lagrangian may be immersed and
disconnected, with all components of the same phase.
In particular, the theorem gives special
representatives of stable objects.

The hypothesis of Theorem~\ref{thm:intro-special-representatives}
holds, in
particular, if $\mathcal D$ is split-generated
by special Lagrangian branes; see
Remark~\ref{rem:special-generators}.

\subsection{The two-torus}
\label{subsec:intro-elliptic-curve}
We now outline the proofs of
Theorems~\ref{thm:intro-elliptic-continuation}
and~\ref{thm:intro-elliptic-termination};
the details are given in Section~\ref{sec:elliptic-curve}. The arguments use geometric smoothing and perturbation
results whose proofs are deferred to the forthcoming
companion paper~\cite{KapotaWebGeometry}, together with
Floer-theoretic results on brane transport, transfer
across surgery, and closed-component reduction whose
proofs are deferred to~\cite{KapotaWebFloer}.
Let
\[
    X=\mathbb C/\Gamma,
    \qquad
    \omega=dx\wedge dy,
    \qquad
    \Omega=dz,
\]
where $\Gamma\subset\mathbb C$ is a lattice and $z=x+iy$.
In the universal cover, positive Lagrangians are
graphs over the $x$-axis, oriented towards
increasing $x$, while imaginary-special web fronts
are vertical segments.
The evolution moves these segments horizontally,
with their endpoints constrained to the positive Lagrangian boundary conditions.

To prove Theorem~\ref{thm:intro-elliptic-termination}, we count crossings in suitable
perturbations of the web cone and add the total covering
multiplicity of the closed label-zero components.
The resulting nonnegative integer does not increase
during regular evolution or geometric singularities,
and decreases with each obstruction wall.
Thus only finitely many such walls can occur.
Geometric singularities remove labeled endpoints of active fronts,
while each obstruction wall introduces only
a uniformly bounded number.
An endpoint count therefore also bounds the number
of geometric singularities.

The crossing count gives a uniform upper bound on
the lengths of active fronts.
Since each front moves horizontally at speed the reciprocal
of its length, this yields a positive lower bound
on the rate at which a nonempty front sweeps out
the positive Lagrangian boundary conditions, measured by $dx$.
The bounded total real period of the swept regions
then bounds the total time during which fronts are active.

Whenever the geometry becomes stationary before
terminality is attained, the brane evolution reaches
another obstruction wall after a finite waiting time.
Since only finitely many such walls occur, the total
stationary waiting time is finite as well. These bounds show that the evolution undergoes only
finitely many events and reaches a terminal configuration
in finite time, proving
Theorem~\ref{thm:intro-elliptic-termination}.

The following example shows how the coefficients
of a morphism determine the timing of obstruction
walls and the resulting terminal configuration.
The boundary curves remain fixed, while the
valuation of one coefficient varies.

\begin{example}
    \label{ex:intro-degree-two}
    Consider the square torus and the positive Lagrangians
    \[
        X=\mathbb R^2/\mathbb Z^2,
        \qquad
        L_0=\{y=0\},
        \qquad
        L_1=\{y=2x\},
    \]
    with the standard Calabi--Yau structure on $X$. We use the fundamental domain $[-1/2,1/2)^2$,
    with opposite sides identified.
    Equip the curves with their positive gradings,
    trivial rank-one local systems over the complex
    Novikov field, and zero bounding cochains.
    Denote the resulting branes by $\mathcal L_0$
    and $\mathcal L_1$.
    Their two intersections
    \[
        p_0=(0,0),
        \qquad
        p_1=(-1/2,0)
    \]
    have Maslov index zero.
    For $\lambda>0$, consider the morphism
    \[
        f_\lambda=[p_0]+T^\lambda [p_1]
        \in HF^0(\mathcal L_0,\mathcal L_1),
    \]
    where $T$ is the Novikov parameter and
    $\operatorname{val}(T^\lambda)=\lambda$.
    Write $F_c=\{x=c\}$ for the vertical circle at $x=c$.
    We compute the terminal configuration of an $f_\lambda$-admissible web directly from the web evolution.

    \begin{figure}[htbp]
        \centering
        \includegraphics[width=\textwidth]{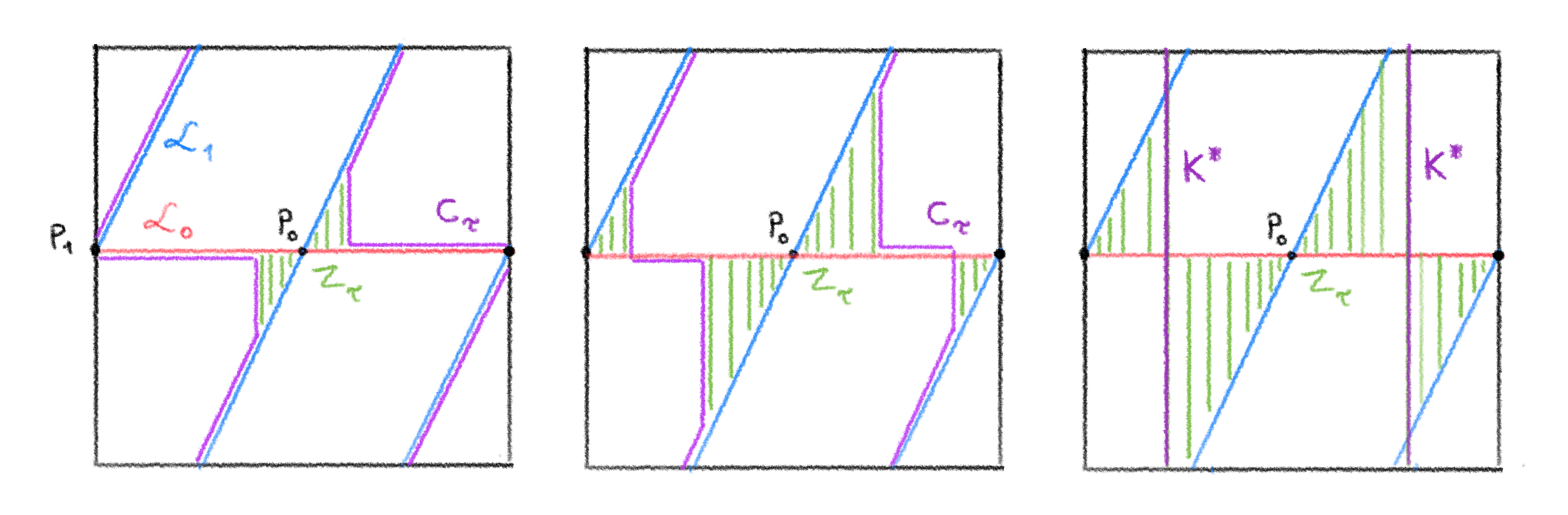}
        \caption{Web evolution for
        $f_\lambda=[p_0]+T^\lambda [p_1]$, with $0<\lambda<1/4$.
        The initial geometric cone $C_0=L_0[1]\sqcup L_1$
        consists of two embedded circles intersecting
        transversely at $p_0$ and $p_1$.
        Left: initiating webbing at $p_0$ joins these into a single
        immersed circle with one self-intersection, at $p_1$.
        Middle: the obstruction wall at $\tau=\lambda$
        initiates a second webbing at $p_1$; the resulting
        surgery resolves this self-intersection, leaving
        two disjoint embedded circles.
        Right: at the terminal time, these become the
        vertical circles $F_{-\frac14-\lambda}$ and
        $F_{\frac14+\lambda}$.
        The curves $L_0$ and $L_1$ are shown in red and blue,
        respectively; successive front positions are green,
        and the web cone is purple.
        Each square is a fundamental domain centered at $p_0$,
        with opposite sides identified.}
        \label{fig:elliptic-web}
    \end{figure}

    On the initial union $L_0\cup L_1$,
    the cone bounding cochain is determined by $f_\lambda$.
    The coefficient $1$ at $p_0$ is resolved by initiating
    a webbing there at time $\tau=0$.
    We use the parametrization in which each front
    sweeps out symplectic area at unit speed; see
    Lemma~\ref{lem:elliptic-harmonic-flux}.
    The area swept out from
    $p_0$ to horizontal coordinate $x$ is
    $\int_0^x 2s\,ds=x^2$.
    Thus, at time $\tau$, the two fronts have coordinates
    \[
        x=\pm\sqrt{\tau}\pmod{\mathbb Z}.
    \]
    While $p_1$ remains an intersection of the unswept
    regions, it remains a generator of the Floer complex
    of the web cone.
    Write $b_{\tau,p_1}$ for its coefficient in the
    transported bounding cochain.
    The valuation-transport formula
    \eqref{eq:web-cone-valuation-transport} gives
    \[
        \operatorname{val}(b_{\tau,p_1})
        =\operatorname{val}(b_{0,p_1})-\tau
        =\lambda-\tau.
    \]
    Suppose first that $0<\lambda<1/4$.
    At time $\tau=\lambda$, the coefficient
    $b_{\tau,p_1}$ reaches valuation zero.
    Further regular transport would make its valuation
    negative, so the web cone encounters an obstruction wall.
    The existing fronts are still away from $p_1$,
    since $\sqrt{\lambda}<1/2$.
    We resolve this wall by initiating a second webbing
    at $p_1$, inducing surgery of the web cone there.
    Using the lifts of $p_1$ at $x=\pm1/2$, the two
    new fronts move into the fundamental domain with
    coordinates
    \[
        x=-\frac12+\sqrt{\tau-\lambda},
        \qquad
        x=\frac12-\sqrt{\tau-\lambda}.
    \]
    Consider the right-moving front from $p_0$ and the
    left-moving front from the lift of $p_1$ at $x=1/2$.
    When they reach the same horizontal coordinate $x=c$,
    their endpoints coincide on each of $L_0$ and $L_1$.
    The areas swept out by these fronts are $c^2$ and
    $(1/2-c)^2$, respectively.
    Since their initiation times differ by $\lambda$,
    \[
        c^2-\left(\frac12-c\right)^2=\lambda,
    \]
    thus $c=\frac14+\lambda.$ At this position, the two fronts are complementary
    arcs of the vertical circle $F_c$, and boundary
    surgery joins them into that circle.
    The other two fronts undergo the same collision
    at $x=-c$.
    The unswept regions disappear, and the terminal
    web cone consists of
    \[
        F_{-\frac14-\lambda}
        \quad\text{and}\quad
        F_{\frac14+\lambda}.
    \]
    Both closed components belong to $K^*$; see Remark~\ref{rem:elliptic-boundary-collision-fragments}. 
    Thus the kernel of $f_\lambda$ is zero, and the
    induced brane structures on these circles represent
    its cokernel.

    If $\lambda>1/4$, the order of events changes.
    The fronts from $p_0$ reach $x=\pm1/2$ at time
    $\tau=1/4$.
    These boundary positions are identified on the torus,
    and the front endpoints collide at $p_1$ while its
    coefficient still has positive valuation.
    No second webbing is initiated.
    The resulting closed curve covers $F_{-1/2}$ twice.
    Its brane structure represents a direct sum of two
    copies of $F_{-1/2}$ carrying distinct rank-one local
    systems.
    In this regime, the terminal geometric support is
    independent of $\lambda$, while the brane data retain
    the dependence on the morphism.
    At $\lambda=1/4$, the obstruction wall and the
    boundary collision occur simultaneously.

    The mirror description recovers these configurations.
    Under the correspondence described in
    \cite[Sections~4 and~5.5]{PolishchukZaslow},
    $f_\lambda$ becomes a nonzero morphism from
    the structure sheaf to a line bundle of degree two.
    It is therefore injective, and its cokernel
    is supported on the zero divisor of the
    corresponding section.
    This divisor consists of two distinct points.
    For $0<\lambda<1/4$, their valuations are
    \[
        -\frac14-\lambda,
        \qquad
        \frac14+\lambda,
    \]
    recovering the two terminal fibers.
    For $\lambda\geq1/4$, both points have valuation
    $-\frac12$ modulo $\mathbb Z$, but remain distinct
    points of the mirror curve.
    They correspond to the two different local systems
    on $F_{-1/2}$.
\end{example}

\subsection{Organization of the paper}
Section~\ref{sec:regular-web-evolution} develops the framework
for regular web evolution.
Subsection~\ref{subsec:reminder-webbing} recalls the geometric
background from \cite{SY-cylinders,SY-webbing}.
Subsections~\ref{subsec:admissible-web-data}
and~\ref{subsec:web-cone} introduce admissible web data
and construct the associated web cone and its smoothing.
Subsection~\ref{subsec:brane-structures} fixes the
Floer-theoretic conventions and describes transport
of brane structures.

Section~\ref{sec:singular-web-evolution} treats
initialization, compatible continuation, and termination.
Subsection~\ref{subsec:initial-cone} specifies the initial
cone data and defines initialization.
Subsection~\ref{subsec:web-continuation} defines
compatible continuation and brane transfer.
Subsection~\ref{subsec:terminal-triangle} defines
terminality and constructs the terminal branes and
the exact triangle~\eqref{eq:terminal-web-triangle}.
Subsection~\ref{subsec:Dirichlet-energy} gives a
criterion for finite total active time and a bound
on the total special volume.

Section~\ref{sec:singular-events} describes the local
models for singular events.
Subsection~\ref{subsec:local-continuation-conventions}
sets out the common notations and convergence conventions.
Subsection~\ref{subsec:geometric-degenerations} treats
geometric degenerations and their continuation by surgery.
Subsection~\ref{subsec:obstruction-walls} describes
obstruction walls and their resolutions by surgery
or algebraic reduction of closed inactive components.

Section~\ref{sec:categorical-consequences} derives the
categorical consequences under the stated conjectures
and hypotheses.
Subsection~\ref{subsec:positive-branes-abelian} uses
the terminal triangle to prove abelianity and identify
geometric kernels and cokernels.
Subsection~\ref{subsec:harder-narasimhan} proves the
Harder--Narasimhan property using geometric bounds
on destabilizing quotients.
Subsection~\ref{subsec:stability-conditions} constructs
stability conditions and compares rotated hearts
to establish the correspondence between semistability
and special Lagrangian representatives.

Section~\ref{sec:elliptic-curve} establishes initialization,
continuation, and finite-time termination for
two-torus.
Subsections~\ref{subsec:elliptic-regular-evolution}--%
\ref{subsec:elliptic-obstruction-walls} describe regular
evolution, geometric degenerations, and obstruction walls,
together with the required brane transport.
Subsection~\ref{subsec:elliptic-counting} develops the
counting arguments used in the proofs of
Subsection~\ref{subsec:elliptic-continuation-termination}.
Subsection~\ref{subsec:elliptic-examples} illustrates
the construction with explicit examples.

Finally, Section~\ref{sec:abelian-surface} considers
an example on an abelian surface.
Subsections~\ref{subsec:abelian-product-morphism}
and~\ref{subsec:abelian-terminal-surface} describe
the morphism and its mirror and derive predictions
for the topology and intersections of the terminal cokernel.
Subsection~\ref{subsec:abelian-web-evolution} presents
a proposed web evolution consistent with these predictions.

\subsection{AI declaration}
ChatGPT suggested the use of minimal counting
diagrams in Subsection~\ref{subsec:elliptic-counting}.
Its assistance also led to simplifications of the
authors' original argument for Harder--Narasimhan
filtrations, including reformulating the relevant
estimates in terms of the terminal representatives
$K$ and $K^*$ and simplifying the use of Stokes'
theorem on chains associated with the web.
ChatGPT was used for exploratory
calculations and checks in examples on the two-torus. The system also suggested useful references and
helped identify connections with existing literature.

During the writing process, author-written paragraphs were submitted to ChatGPT,
which proposed revised text in LaTeX.
These proposals were inspected carefully, often edited
further, and incorporated manually into the source files.

\subsection{Acknowledgements}
The authors would like to thank M.~Abouzaid, O.~Chodosh, Y.~Eliashberg, D.~Kazhdan, E.~Kosloff and D.~Nadler for helpful conversations. The authors were partially supported by ISF grants 1127/22 and 3029/26. J.~S. was partially supported by the Miriam and Julius Vinik Chair in Mathematics.

\section{Webs, cones, and regular brane evolution}
\label{sec:regular-web-evolution}

\subsection{Background on special Lagrangian webbing}
\label{subsec:reminder-webbing}

We recall the special Lagrangian webbing construction of
\cite{SY-cylinders,SY-webbing}.
Let $(X,\omega,J,\Omega)$ be a Calabi-Yau manifold of complex dimension $n$ in the sense of \cite[Definition~2.17]{SY-cylinders}: $(X,\omega)$ is symplectic, $J$ is an $\omega$-compatible integrable complex structure, and $\Omega$ is a nowhere-vanishing holomorphic $(n,0)$-form. Let
\begin{equation}
    g=\omega(\cdot,J\cdot)
\end{equation}
be the K\"ahler metric. Following \cite{SY-cylinders}, let
$\rho:X\to\mathbb R_{>0}$ denote the Calabi-Yau density, characterized by
\begin{equation}
    \rho^2\frac{\omega^n}{n!}
    =
    (-1)^{\frac{n(n-1)}{2}}
    \left(\frac{\sqrt{-1}}{2}\right)^n
    \Omega\wedge\overline{\Omega}.
\end{equation}

\begin{definition}
    \label{def:immersed-submanifold}
    We briefly recall the convention for immersed submanifolds used in
    \cite{SY-cylinders}. Let $M$ be a manifold with or without boundary. An \emph{immersed submanifold $Z$ in $X$ of type $M$} is an equivalence class of immersions
    \begin{equation}
        Z=[f:M\longrightarrow X]
    \end{equation}
    where the immersions $f$, $f'$ are \emph{equivalent} if there exists a diffeomorphism of manifolds with boundary $\vp: M \to M$ such that $f' = f \circ \vp$. We call $M$ the \emph{source} of $Z$.

    Points, functions, and differential forms on $Z$
    are understood on the source; restriction of an
    ambient form means pullback by $f$.
    In particular, distinct source points remain distinct
    even when they have the same image in $X$.
    At a source point $p$, the tangent plane is
    $df_p(T_pM)\subset T_{f(p)}X$.
\end{definition}

\begin{definition}
    Following \cite[Definition~4.1]{SY-cylinders}, let
    $L_0,L_1\subset X$ be embedded Lagrangians and let $N$
    be a closed connected $(n-1)$-manifold.
    A \emph{Lagrangian cylinder} between $L_0$ and $L_1$
    is an immersed Lagrangian
    \[
        Z=[f:N\times[0,1]\longrightarrow X]
    \]
    whose boundary restrictions
    $f|_{N\times\{i\}}:N\to L_i$ are embeddings and satisfy
    \[
        df_x\bigl(T_x(N\times[0,1])\bigr)
        \neq T_{f(x)}L_i,
        \qquad x\in N\times\{i\},\quad i=0,1.
    \]
    We write $C_i\subset Z$ for the boundary component
    corresponding to $N\times\{i\}$.
    Denote the space of such cylinders by
    $\cc{L}C(N;L_0,L_1)$.
    When $L_0,L_1$ are positive, let
    $\cc{S}LC(N;L_0,L_1)\subset\cc{L}C(N;L_0,L_1)$
    be the subspace of imaginary special Lagrangian cylinders.
\end{definition}

We recall the deformation convention for immersed Lagrangians,
which will also be used for fronts of general topology and
for web cones.

\begin{definition}
    \label{def:lagrangian-deformation}
    \label{def:local-Hamiltonian}
    Let $P$ be a manifold, possibly with boundary, and let
    $I\subset\mathbb R$ be an interval.
    A family $(Z_s)_{s\in I}$ of immersed Lagrangians of
    type $P$ is \emph{smooth} if it admits a smooth map
    \[
        \Phi:P\times I\longrightarrow X
    \]
    such that $\Phi_s:=\Phi(\,\cdot\,,s)$ represents $Z_s$.
    We call $\Phi$ a \emph{lifting} of the family.
    The derivative of the family is the closed $1$-form
    \[
        \frac{d}{ds}Z_s
        :=\alpha_s
        :=\Phi_s^*
        \bigl(\iota_{\partial_s\Phi_s}\omega\bigr).
    \]
    Viewed as a form on the immersed Lagrangian $Z_s$,
    it is independent of the lifting
    \cite[Remark~2.14]{SY-cylinders}.
    Its closedness follows from Cartan's formula and
    $\Phi_s^*\omega=0$.

    The family is an \emph{exact deformation}, or
    \emph{local Hamiltonian isotopy}, if there is a smooth
    family of functions $H_s:P\to\mathbb R$ such that
    \[
        \alpha_s=dH_s.
    \]
    We call $H_s$ a \emph{local Hamiltonian}
    \cite[Section~2.6]{Joyce-stability}.
    It is a function on the source and need not be the
    pullback of a function on $X$.
\end{definition}

For a smooth family of cylinders between fixed Lagrangians
$L_0,L_1$, the deformation form $\alpha_s$ restricts to zero
on the boundary. Indeed, both the boundary velocity and
the tangent vectors to the boundary lie in $TL_i$, on
which $\omega$ vanishes.

Since the source of a cylinder is $N\times[0,1]$, every
closed $1$-form whose restriction to the boundary vanishes
is exact \cite[Lemma~4.3]{SY-cylinders}.
Thus, for a smooth family of cylinders $Z_s$ with boundary
components $C_{i,s}$, there is a unique normalized primitive
$h_s$ satisfying
\[
    \frac{d}{ds}Z_s=dh_s,
    \qquad
    h_s|_{C_{0,s}}=0,
    \qquad
    h_s|_{C_{1,s}}\text{ is constant}.
\]
In particular, every such family is a local Hamiltonian
isotopy.

Suppose now that $L_0,L_1$ are positive.
On an imaginary special Lagrangian cylinder $Z$, define
the weighted Laplacian by
\[
    \Delta_\rho u:=*d(\rho *du),
\]
where $\rho$ is the restriction of the Calabi--Yau density
and $*$ is the Hodge star of the induced metric.

\begin{definition}
    Following \cite[Definition~4.13]{SY-cylinders}, the
    \emph{fundamental harmonic} of
    $Z\in\cc{S}LC(N;L_0,L_1)$ is the unique function
    $\sigma_Z\in C^\infty(Z)$ solving
    \[
        \Delta_\rho\sigma_Z=0,
        \qquad
        \sigma_Z|_{C_0}=0,
        \qquad
        \sigma_Z|_{C_1}=1.
    \]
\end{definition}

The fundamental harmonic determines the infinitesimal
deformations of imaginary special Lagrangian cylinders.

\begin{prop}[\cite{SY-cylinders}, Proposition~4.7]
    \label{prop:slc-deformation-space}
    The space $\cc{S}LC(N;L_0,L_1)$ is a smooth
    one-dimensional submanifold of
    $\cc{L}C(N;L_0,L_1)$.
    Under the identification of deformation vectors
    with closed $1$-forms,
    \[
        T_Z\cc{S}LC(N;L_0,L_1)
        =\mathbb R\,d\sigma_Z.
    \]
\end{prop}

We use relative Lagrangian flux to fix the direction and
speed of the cylinder evolution.

\begin{definition}
    \label{def:relative-lagrangian-flux}
    Let $(Z_s)_{s\in[s_0,s_1]}$ be a smooth family of
    Lagrangian cylinders, with boundary components $C_{i,s}$.
    Write
    \[
        \frac{d}{ds}Z_s=dh_s,
        \qquad
        h_s|_{C_{0,s}}=0,
        \qquad
        h_s|_{C_{1,s}}=A_s.
    \]
    Following \cite[Definition~4.11]{SY-cylinders}, the
    \emph{relative Lagrangian flux} of the family is
    \[
        \operatorname{RelFlux}
        \bigl((Z_s)_{s\in[s_0,s_1]}\bigr)
        :=-\int_{s_0}^{s_1}A_s\,ds.
    \]
\end{definition}

Geometrically, relative flux is the signed symplectic area
swept out by a path joining the two boundary components.
Indeed, choose a lifting $\Phi$ of the family and a path
$\gamma:[0,1]\to N\times[0,1]$ with
$\gamma(i)\in N\times\{i\}$.
Setting $\bar\gamma(t,s)=\Phi(\gamma(t),s)$, one has
\[
    \int_{[0,1]\times[s_0,s_1]}\bar\gamma^*\omega
    =
    \operatorname{RelFlux}
    \bigl((Z_s)_{s\in[s_0,s_1]}\bigr),
\]
independently of these choices
\cite[Lemma~4.10 and Remark~4.12]{SY-cylinders}.

For a smooth family of imaginary special Lagrangian
cylinders with nonvanishing derivative,
Proposition~\ref{prop:slc-deformation-space} gives
\[
    \frac{d}{ds}Z_s=a_s\,d\sigma_{Z_s},
    \qquad a_s\neq0.
\]
Thus $A_s=a_s$, so relative flux is strictly monotone.
We choose its negative as the parameter $\tau$,
unique up to an additive constant. In this parameter,
\[
    \frac{d}{d\tau}Z_\tau=d\sigma_{Z_\tau},
    \qquad
    \operatorname{RelFlux}
    \bigl((Z_\tau)_{\tau\in[\tau_0,\tau_1]}\bigr)
    =-(\tau_1-\tau_0).
\]
We use this normalization for the evolution below.

The following theorem supplies the local starting point
for the web evolution. We equip positive Lagrangians with their
natural gradings $\theta_L\in(-\pi/2,\pi/2)$ and use the Maslov-index convention of
\cite[Definition~3.12]{SY-webbing}.

\begin{theorem}[\cite{SY-webbing}, Theorem~1.1(a)]
    \label{thm:SY-start-web}
    Let $L_0,L_1\subset X$ be smoothly embedded positive
    Lagrangians intersecting transversely at $q$, with
    $m(q;L_0,L_1)=0$.
    Then there exist $\epsilon>0$ and a smooth
    one-parameter family
    \[
        (Z_s)_{s\in(0,\epsilon)}
        \subset
        \cc{S}LC(S^{n-1};L_0,L_1)
    \]
    converging regularly to $q$ as $s\to0$.
    We call this family a \emph{special Lagrangian webbing}
    emanating from $q$.
\end{theorem}

\begin{remark}
    We use regular convergence in the sense of
    \cite[Definition~4.17(3)]{SY-cylinders}.
    Roughly, the cylinders shrink to $q$ while retaining
    a smooth, nondegenerate shape after rescaling.
    This controls both the collapsing cylinders and the
    regions swept out by their boundaries on $L_0$ and
    $L_1$, providing the local geometric model for
    initiating the web-cone evolution.
\end{remark}

\begin{remark}
    \label{rem:maslov-n-webbing}
    If $m(q;L_0,L_1)=n$, the theorem applies after
    exchanging $L_0$ and $L_1$, since
    \[
        m(q;L_1,L_0)=n-m(q;L_0,L_1)=0;
    \]
    see \cite[Remark~1.2]{SY-webbing}.
    Exchanging the boundary labels replaces the
    fundamental harmonic $\sigma$ by $1-\sigma$,
    reversing the prescribed deformation direction.
    Thus, for the original ordering $(L_0,L_1)$,
    the evolution traverses this family toward $q$.
    We call it a \emph{special Lagrangian webbing
        terminating at $q$}.
\end{remark}

\subsection{Admissible web data}
\label{subsec:admissible-web-data}

In this subsection we define the geometric data needed to
assemble the web cone on an interval of smooth evolution.
These data include imaginary special Lagrangian fronts of
general topology, stationary inactive fragments, and boundary
traces recording the swept regions. An incoming boundary
trace records regions swept out before the interval begins.
The admissibility conditions below ensure that the boundaries
of the fronts and fragments match those of the unswept regions,
providing the attachments used in the web cone construction.

\subsubsection{Regular front evolution}
\label{subsubsec:regular-front-evolution}

Throughout this subsection, fix an open interval $I=(a,b)$
and strictly positive immersed Lagrangians
\[
    L_i=[\iota_i:M_i\longrightarrow X],
    \qquad i=0,1,
\]
which serve as fixed boundary conditions along $I$.
Their sources are compact and may be empty or disconnected.
Since the boundary conditions are immersed, we must specify
the sheets on which each boundary lies.

\begin{definition}
    \label{def:immersed-boundary-lifts}
    An \emph{immersed submanifold with boundary on $L_0, L_1$},
    is an equivalence class of triples
    \[
        Z=[(F,\beta_0,\beta_1)],
    \]
    where $F:P\to X$ is an immersion, and
    \[
        \partial P=\partial_0P\sqcup\partial_1P
    \]
    is a decomposition of its boundary. We call $\partial_iP$ the \emph{boundary with label $i$}. The maps $\beta_i$ are immersions satisfying
    \[
        \beta_i:\partial_iP\longrightarrow M_i,
        \qquad
        F|_{\partial_iP}=\iota_i\circ\beta_i.
    \]
    We call $\beta_i$ the \emph{chosen boundary lifts}; they are part of the data and specify the sheets of $L_i$ on which the boundary lies. Two triples $(F:P\to X,\beta_0,\beta_1)$ and $(\widetilde F:\widetilde P\to X,\widetilde\beta_0,\widetilde\beta_1)$ are equivalent if there is a
    diffeomorphism $\varphi:\widetilde P\to P$
    such that $\varphi(\partial_i \widetilde P) = \partial_i P$ and
    \[
        \widetilde F=F\circ\varphi,
        \qquad
        \widetilde\beta_i
        =\beta_i\circ
        \varphi|_{\partial_i\widetilde P}.
    \]
    We will usually write $Z=[F]$, suppressing the boundary maps.
    The notation $\partial_iZ$ denotes the immersed submanifold $[\iota_i\circ\beta_i: \partial_i P \to X]$.
\end{definition}

\begin{definition}
    A \emph{special Lagrangian front of type $P$}
    between $L_0$ and $L_1$ is an immersed imaginary
    special Lagrangian
    \[
        Z=[(F:P\longrightarrow X,\beta_0,\beta_1)]
    \]
    with boundary on $L_0,L_1$, whose source $P$ is a
    compact $n$-manifold and whose boundary lifts $\beta_i$
    are embeddings.
    We require every connected component of $P$ to meet
    both $\partial_0P$ and $\partial_1P$.
    The source may be empty or disconnected.

    The connected components of $Z$ are called its
    \emph{front components}.
    For a component $Y$, write $P_Y$ for its source,
    $F_Y=F|_{P_Y}$ for its immersion, and
    $\beta_i^Y=\beta_i|_{\partial_iP_Y}$ for its boundary lifts.
\end{definition}

For fronts of general topology, the space of imaginary
special Lagrangian deformations need not be one-dimensional;
see \cite[Remark~4.9]{SY-cylinders} and~\cite[Theorem~2.2]{Pidaparthy-deformations}.
Nevertheless, the Dirichlet problem defining the fundamental
harmonic selects an exact deformation direction.

\begin{definition}
    Let $Y=[F_Y:P_Y\to X]$ be a connected special
    Lagrangian front.
    Its \emph{fundamental harmonic} is the unique function
    $\sigma_Y\in C^\infty(P_Y)$ satisfying
    \[
        \Delta_{\rho\circ F_Y}\sigma_Y=0,
        \qquad
        \sigma_Y|_{\partial_0P_Y}=0,
        \qquad
        \sigma_Y|_{\partial_1P_Y}=1,
    \]
    where the weighted Laplacian is computed using the
    induced metric $F_Y^*g$.
    For a disconnected front $Z$, we define $\sigma_Z$
    componentwise.
\end{definition}

The form $d\sigma_Y$ restricts to zero on the boundary
and satisfies the linearized imaginary special Lagrangian
equation \cite[Lemma~4.5(a)]{SY-cylinders}.
We use it to prescribe the motion of the front.

\begin{definition}
    \label{def:regular-front-evolution}
    A \emph{regular front evolution} on $I$ is a family
    of special Lagrangian fronts
    \[
        Z_\tau=[(F_\tau,\beta_{0,\tau},\beta_{1,\tau})],
        \qquad \tau\in I,
    \]
    on a fixed source $P$, with the immersions $F_\tau$
    and chosen boundary lifts $\beta_{i,\tau}$ depending
    smoothly on $\tau$, and satisfying
    \begin{equation}
        \label{eq:fundamental-harmonic-evolution}
        \frac{d}{d\tau}Z_\tau=d\sigma_{Z_\tau}.
    \end{equation}
    The derivative is understood as a $1$-form on the
    source, as in Definition~\ref{def:lagrangian-deformation}.
\end{definition}

\begin{remark}
    Regularity is understood using the chosen boundary lifts:
    a boundary point may pass through the image of another
    sheet of $L_i$ while the front and its lifts remain smooth.
\end{remark}

Thus, at each time, the Dirichlet problem determines the
exact deformation direction of the current front.
Together with the local deformation theory of
\cite[Proposition~4.7 and Remark~4.9]{SY-cylinders},
this gives local existence and uniqueness of the
regular front evolution, up to reparameterization
of the source.

\begin{remark}
    The relative-flux formula of
    Definition~\ref{def:relative-lagrangian-flux} extends
    to exact families of fronts
    $\frac{d}{ds}Z_s=dh_s$ whose normalized primitives
    take common values $0$ and $A_s$ on the boundaries
    with labels $0$ and $1$, respectively:
    \[
        \operatorname{RelFlux}
        \bigl((Z_s)_{s\in[s_0,s_1]}\bigr)
        :=-\int_{s_0}^{s_1}A_s\,ds.
    \]
    For a regular front evolution, the primitive
    $\sigma_{Z_\tau}$ has boundary values $0$ and $1$.
    Consequently, every front component accumulates
    relative flux $-(\tau_1-\tau_0)$ during
    $[\tau_0,\tau_1]\subset I$.
\end{remark}

\subsubsection{Inactive fragments}
\label{subsubsec:inactive-fragments}

Alongside the moving fronts, the web data may contain
imaginary special Lagrangians whose boundary lies on only
one side, or which have no boundary.
These pieces remain stationary throughout $I$.

\begin{definition}
    \label{def:terminal-fragment}
    An \emph{inactive fragment} is a connected immersed
    imaginary special Lagrangian
    \[
        Y=[(F_Y:P_Y\longrightarrow X,\beta_0^Y,\beta_1^Y)]
    \]
    with compact source and boundary on $L_0,L_1$, such that
    \[
        \partial_0P_Y=\emptyset
        \qquad\text{or}\qquad
        \partial_1P_Y=\emptyset.
    \]
    In particular, an inactive fragment may be closed.
\end{definition}

\begin{definition}
    \label{def:inactive-fragment-notation}
    Let $\mathcal F$ be a finite collection of inactive
    fragments with fixed immersions and boundary lifts
    throughout $I$.
    For $i=0,1$, let $D_i$ denote the disjoint union of
    the nonclosed fragments whose boundary has label $i$.

    Each closed fragment is equipped with an additional
    label $i\in\{0,1\}$.
    Let $D_{\mathrm{cl},i}$ denote the disjoint union of
    the closed fragments carrying label $i$.
    Thus,
    \[
        \bigsqcup_{Y\in\mathcal F}Y
        =
        D_0\sqcup D_1
        \sqcup D_{\mathrm{cl},0}\sqcup D_{\mathrm{cl},1}.
    \]
    All unions retain the distinct source components,
    even when their images in $X$ coincide.
\end{definition}

\begin{remark}
    For a closed fragment, the label does not refer to a
    boundary component. It will determine the fragment's
    grading and height in the web-cone construction.
\end{remark}

\subsubsection{Boundary traces and admissible web data}

To assemble the web cone, we must record which regions of
the boundary sources have already been swept out.
On the interval $I=(a,b)$, this record consists of an
incoming trace, supplied as initial data, enlarged by
the subsequent motion of the lifted front boundaries.

\begin{definition}
    \label{def:regular-web-data}
    \emph{Web data on $I$}, with boundary conditions
    $L_0,L_1$, consist of:
    \begin{enumerate}
        \item a regular front evolution
              $(Z_\tau)_{\tau\in I}$;
        \item a finite collection $\mathcal F$ of stationary
              inactive fragments, with chosen boundary lifts
              and labels on the closed fragments;
        \item open subsets
              $B_i^{\mathrm{in}}\subset M_i$, for $i=0,1$,
              called the \emph{incoming boundary traces}.
    \end{enumerate}
    We write
    \[
        \mathcal W_I
        =
        \bigl(
        (Z_\tau)_{\tau\in I},
        \mathcal F,
        B_0^{\mathrm{in}},B_1^{\mathrm{in}}
        \bigr).
    \]
    The incoming traces record the regions swept out
    before the interval begins.
\end{definition}

For a local webbing starting at its collapsing end,
the incoming traces are empty.
More generally, we allow nonempty incoming traces
without specifying the preceding evolution.
Their compatibility across singular times is specified in
Definition~\ref{def:surgery-compatible}.

\begin{definition}
    \label{def:boundary-trace}
    For web data $\mathcal W_I$, the \emph{boundary trace}
    at time $\tau\in I$ is
    \[
        B_i(\tau)
        :=
        \operatorname{Int}_{M_i}
        \overline{
            B_i^{\mathrm{in}}
            \,\cup\!
            \bigcup_{a<s<\tau}
            \beta_{i,s}(\partial_iP)
        },
        \qquad i=0,1,
    \]
    where the closure is taken in $M_i$.
\end{definition}

Only the moving front boundaries enlarge the trace;
the inactive fragments remain stationary.
The interior of the closure fills punctures left by
a collapsing initial front while keeping the trace open.
For example, for a local webbing emanating from $q$,
the trace contains the chosen source point
$x_i\in M_i$ with $\iota_i(x_i)=q$.
On the corresponding locally embedded sheet, it agrees
with the unpunctured swept submanifold of
\cite[Definition~5.7]{SY-cylinders}.

\begin{remark}
    \label{rem:boundary-trace-monotone}
    The boundary traces depend only on the chosen boundary
    lifts, not on their source parameterizations.
    They increase with time:
    \[
        B_i(\tau_0)\subset B_i(\tau_1),
        \qquad a<\tau_0<\tau_1<b.
    \]
\end{remark}

The boundary of the swept region must match the boundaries
of all the pieces to be attached to the unswept region,
including the inactive fragments.
This is the purpose of the following admissibility condition.

\begin{definition}
    \label{def:admissible-web}
    For web data $\mathcal W_I$, define the combined
    boundary lift
    \[
        \beta_{i,\tau}^{\mathrm{tot}}
        :=
        \beta_{i,\tau}
        \sqcup
        \bigsqcup_{Y\in\mathcal F}\beta_i^Y
        :
        \partial_iP
        \sqcup
        \bigsqcup_{Y\in\mathcal F}\partial_iP_Y
        \longrightarrow M_i,
    \]
    and write
    \[
        \Gamma_i(\tau)
        :=
        \operatorname{im}\beta_{i,\tau}^{\mathrm{tot}}
    \]
    for the \emph{boundary front}.
    The web data are \emph{admissible} if, for every
    $\tau\in I$ and $i=0,1$,
    \begin{enumerate}
        \item $\beta_{i,\tau}^{\mathrm{tot}}$ is an embedding;
        \item $\partial_{M_i}B_i(\tau)=\Gamma_i(\tau)$.
    \end{enumerate}
    The open subset
    \[
        R_i(\tau)
        :=
        M_i\setminus\overline{B_i(\tau)}
    \]
    is called the \emph{residue}.
\end{definition}

For admissible web data, the boundary source decomposes as
\[
    M_i
    =
    B_i(\tau)\sqcup\Gamma_i(\tau)\sqcup R_i(\tau).
\]
Thus $B_i(\tau)$ is the swept region, $R_i(\tau)$ is the
unswept region, and $\Gamma_i(\tau)$ is their common boundary.
The combined boundary lift identifies the boundaries of
the fronts and inactive fragments with the boundary of
$\overline{R_i(\tau)}$.
These are the identifications used to assemble the
singular web cone.
All these conditions are imposed in the sources $M_i$.
In particular, points of $B_i(\tau)$ and $R_i(\tau)$
may have the same image under the immersion $\iota_i$.

\subsection{The web cone}
\label{subsec:web-cone}

We now assemble the fronts, inactive fragments, and closed
residues into a piecewise smooth Lagrangian.
The fundamental harmonics extend to a continuous height
function on its source, which will determine the local
Hamiltonian of the smoothed evolution.

Fix admissible web data $\mathcal W_I$ on $I=(a,b)$,
with the notation of the preceding subsection. Write $\overline R_i(\tau):=\overline{R_i(\tau)}\subset M_i.$

\begin{definition}
    \label{def:singular-web-cone}
    For $\tau\in I$, the source of the
    \emph{singular web cone} is
    \[
        N_\tau^{\mathrm{sing}}
        :=
        \left(
        \overline R_0(\tau)
        \sqcup P
        \sqcup\bigsqcup_{Y\in\mathcal F}P_Y
        \sqcup\overline R_1(\tau)
        \right)\big/\sim,
    \]
    where the equivalence relation is generated by
    \[
        x\sim\beta_{i,\tau}^{\mathrm{tot}}(x),
        \qquad
        x\in
        \partial_iP
        \sqcup\bigsqcup_{Y\in\mathcal F}\partial_iP_Y,
        \qquad i=0,1.
    \]
    The immersions of the pieces agree along these
    identifications and therefore induce a continuous,
    piecewise smooth map
    \[
        \iota_\tau^{\mathrm{sing}}:
        N_\tau^{\mathrm{sing}}\longrightarrow X.
    \]
    The \emph{singular web cone} is the piecewise smooth
    immersed Lagrangian
    \[
        C_\tau^{\mathrm{sing}}
        :=
        [\iota_\tau^{\mathrm{sing}}].
    \]
    The pieces and their gradings are recorded by
    \begin{equation}
        \label{eq:singular-web-cone-pieces}
        C_\tau^{\mathrm{sing}}
        =
        \overline R_0(\tau)[1]
        \cup D_0
        \cup D_{\mathrm{cl},0}[1]
        \cup Z_\tau
        \cup D_{\mathrm{cl},1}
        \cup D_1
        \cup\overline R_1(\tau),
    \end{equation}
    where $D_i$ and $D_{\mathrm{cl},i}$ are as in
    Definition~\ref{def:inactive-fragment-notation}.
    Here $\overline R_i(\tau)$ is regarded as a
    Lagrangian piece via the restriction of $\iota_i$,
    with the grading inherited from $L_i$.
    Before the indicated shifts, $Z_\tau$ and all
    inactive fragments carry the constant grading
    $\pi/2$.
\end{definition}

The unions in \eqref{eq:singular-web-cone-pieces} denote
gluing along the chosen boundary lifts.
No further source points are identified merely because
they have the same image in $X$.
The fundamental harmonic on each active front takes
the boundary value $i$ along its attachment to
$\overline R_i(\tau)$.
We extend it over the stationary pieces by these
constant values.

\begin{definition}
    \label{def:height-function}
    The \emph{height function} of the singular web cone is
    the continuous, piecewise smooth function
    \[
        \kappa_\tau:
        N_\tau^{\mathrm{sing}}\longrightarrow[0,1]
    \]
    determined by
    \[
        \kappa_\tau|_P=\sigma_{Z_\tau},
        \qquad
        \kappa_\tau|_{\overline R_i(\tau)}=i,
    \]
    and, for each inactive fragment $Y\in\mathcal F$,
    \[
        \kappa_\tau|_{P_Y}=i
        \quad\text{if}\quad
        Y\subset D_i\sqcup D_{\mathrm{cl},i},
        \qquad i=0,1.
    \]
    These prescriptions agree along all identified boundaries,
    so the height is well defined and continuous on the glued
    source. It need not descend to a function on the image
    of the cone in $X$: distinct source points with the same
    image may have different heights.
\end{definition}

On each active front, the deformation form is
$d\sigma_{Z_\tau}$, while on every stationary piece it
vanishes. Thus, on each smooth piece of the singular
web cone, the deformation form is $d\kappa_\tau.$
The height therefore supplies a continuous primitive
for the piecewise deformation.

The singular web cone has corners where the fronts and
inactive fragments meet the residues.
The following result rounds these corners coherently
in time, preserving the prescribed evolution away from
the corners and controlling the volume introduced by
the rounding.

The construction uses Weinstein charts relative to
the Lagrangian boundary conditions.
To preserve exactness, the rounding is performed on
the Lagrangian suspension, in charts preserving the
time coordinate.
The proof will be given in the forthcoming companion
paper~\cite{KapotaWebGeometry}.

Volumes are computed on the sources with the metric
induced by $g$, counting distinct sheets with multiplicity.
In particular, $\operatorname{Vol}_g(C_\tau^{\mathrm{sing}})$
is the sum of the volumes of its source pieces.

\begin{prop}[Coherent smoothing of the web cone]
    \label{prop:web-cone-smoothing}
    Let $\mathcal W_I$ be admissible web data on $I$,
    and let $J\subset I$ be a compact subinterval.
    For every $\eta>0$, there exist a closed manifold $N$,
    a smooth family of graded Lagrangian immersions
    \[
        \iota_\tau:N\longrightarrow X,
        \qquad \tau\in J,
    \]
    and a smooth family of functions
    \[
        \widetilde\kappa_\tau:N\longrightarrow\mathbb R
    \]
    with the following properties.
    \begin{enumerate}
        \item \emph{Smoothing and locality.}
              Each $C_\tau=[\iota_\tau]$ is obtained by rounding
              the corners of $C_\tau^{\mathrm{sing}}$.
              The roundings can be confined to arbitrarily small
              source neighborhoods of the corners, uniformly
              for $\tau\in J$.
              Outside these neighborhoods, the immersion and height
              remain unchanged, and the grading agrees with the
              prescribed grading on each piece.

        \item \emph{Exact evolution.}
              The family satisfies
              \begin{equation}
                  \label{eq:smoothed-web-cone-evolution}
                  \iota_\tau^*
                  \left(
                  \iota_{\partial_\tau\iota_\tau}\omega
                  \right)
                  =
                  d\widetilde\kappa_\tau.
              \end{equation}
              Thus $(C_\tau)_{\tau\in J}$ is an exact deformation
              with local Hamiltonian
              $H_\tau=\widetilde\kappa_\tau$.

        \item \emph{Stationary fragments.}
              The roundings along inactive fragments can
              be chosen independent of $\tau$.
              On a nonclosed fragment with boundary label
              $i$, and on its rounded attachments to the
              residue, one has
              \[
                  \widetilde\kappa_\tau=i.
              \]
              Closed fragments remain unchanged, with
              height equal to their label.

        \item \emph{Volume control.}
              The total volume of the new rounded collars
              is less than $\eta$, uniformly for $\tau\in J$.
              Consequently,
              \[
                  \operatorname{Vol}_g(C_\tau)
                  \leq
                  \operatorname{Vol}_g(C_\tau^{\mathrm{sing}})
                  +\eta,
                  \qquad \tau\in J.
              \]
    \end{enumerate}
    We call $C_\tau$ a \emph{web cone} and
    $\widetilde\kappa_\tau$ its
    \emph{smoothed height function}.
\end{prop}

\begin{remark}
    The smoothed immersions may pass through
    self-tangencies or, for curves, triple intersections.
    The family remains smooth as a family of immersions,
    and the exactness identity
    \eqref{eq:smoothed-web-cone-evolution}
    holds also at these times.
    Brane transport through isolated self-tangencies
    is treated in Subsection~\ref{subsec:self-tangencies}.
    Triple-point transitions in the two-torus
    are treated in
    Proposition~\ref{prop:elliptic-exact-continuation}.
\end{remark}

\subsection{Brane structures on the web cone}
\label{subsec:brane-structures}

The preceding construction associates a family of smoothed
web cones to admissible web data on a regular interval.
We now consider brane structures on these Lagrangians.
After recalling the necessary Fukaya-category conventions
in Sections~\ref{subsubsec:floer-conventions}
and~\ref{subsubsec:derived-category},
we describe how a brane structure chosen at one time can
be transported along the exact deformation determined
by the height function.

\subsubsection{Floer-theoretic conventions}
\label{subsubsec:floer-conventions}
Our arguments use the Floer-theoretic structures
and properties stated below, without relying on
the details of the underlying constructions or the choice of cochain models.
This subsection fixes our conventions and terminology
and provides references for Floer models satisfying these properties.

Throughout the paper, fix an algebraically closed
field $\Bbbk$ of characteristic zero.
We work over the universal Novikov field
\[
    \Lambda
    =
    \left\{
    \sum_j a_jT^{\lambda_j}
    \ \middle|\
    a_j\in\Bbbk,\quad
    \lambda_j\in\mathbb R,\quad
    \lambda_j\longrightarrow+\infty
    \right\},
\]
where finite sums are included and $T$ has degree zero.
The valuation of a nonzero series is its smallest exponent
after collecting terms, and $\operatorname{val}(0)=+\infty$.
We write
\[
    \Lambda_0=\{a:\operatorname{val}(a)\geq0\},
    \qquad
    \Lambda_+=\{a:\operatorname{val}(a)>0\}.
\]
The units of $\Lambda_0$ are precisely the elements of
valuation zero.

We use rank-one $\Bbbk$-local systems, that is,
locally constant rank-one $\Bbbk$-vector bundles,
as in \cite[Definition~2.18]{Joyce-stability}.
For an immersed Lagrangian $L=[\iota:M\to X]$,
the local system is defined on its source $M$.

We work with graded Lagrangian immersions with
compact, possibly disconnected sources without
boundary, equipped with rank-one $\Bbbk$-local
systems and relative spin structures.
The relative spin structures are taken with
respect to a fixed background class in
$H^2(X;\mathbb Z/2)$ and are used to orient
the Floer moduli spaces; see
\cite[Section~5]{AkahoJoyce}.
For two such immersions with local systems, write
\[
    CF^*((L_0,E_0),(L_1,E_1))
\]
for their Floer cochain space over $\Lambda$.
We assume that these spaces and their Floer
operations are defined with compatible choices
of data for the configurations used below.

When the two immersions intersect transversely,
the generators are indexed by pairs
$q=(p_0,p_1)$ of source points satisfying
$\iota_0(p_0)=\iota_1(p_1)$.
Their coefficient lines incorporate
$\operatorname{Hom}_{\Bbbk}((E_0)_{p_0},(E_1)_{p_1})$
and the corresponding orientation lines,
with scalars extended to $\Lambda$;
see \cite{Seidel,AkahoJoyce}. We use the following convention for their degrees.
Let $\beta_1,\ldots,\beta_n\in(0,\pi)$ be the
positive characteristic angles from the
$L_0$-branch to the $L_1$-branch.
Thus, in suitable identification of the tangent space with $\bb{C}^n$,
the
two tangent planes are
\[
    \mathbb R^n,
    \qquad
    e^{i\beta_1}\mathbb R
    \times\cdots\times
    e^{i\beta_n}\mathbb R.
\]
Writing $\theta_{L_i}$ for the grading functions,
the degree of $q$ is
\begin{equation}
    \label{eq:floer-generator-degree}
    \deg(q)
    =
    \frac{
        \theta_{L_0}(p_0)-\theta_{L_1}(p_1)
        +\sum_{j=1}^n\beta_j
    }{\pi}
    \in\mathbb Z;
\end{equation}
see \cite[Definition~2.20]{Joyce-stability}.
In particular, shifting the source grading
by $[1]$ increases the degree by one.

For an immersion $L=[\iota:M\to X]$ with
transverse self-intersections, the underlying
graded space of its self-Floer model decomposes as
\begin{equation}
    \label{eq:self-floer-cochain-decomposition}
    CF^*((L,E),(L,E))
    =
    CF^*_{\mathrm{ch}}(L,E)
    \oplus
    CF^*_{\mathrm{sw}}(L,E).
\end{equation}
The first summand consists of ordinary cochains
on $M$ with coefficients in $\operatorname{End}(E)$,
extended to $\Lambda$.
The second consists of the ordered
self-intersection generators
$q=(q_-,q_+)$, where
\[
    q_-\neq q_+,
    \qquad
    \iota(q_-)=\iota(q_+).
\]
Such a generator goes from the $q_-$-branch
to the $q_+$-branch, with degree and coefficient
line as above; see \cite{AkahoJoyce}.
The Floer operations need not preserve these
two summands separately.

For clean intersections, the isolated generators
are replaced by cochains on the source fiber product, with the corresponding coefficient systems
and grading shifts;
see \cite[Section~3]{Fukaya-immersed}.
For self-Floer cochains, the diagonal component
of the source fiber product contributes to
$CF^*_{\mathrm{ch}}$, while its other components
contribute to $CF^*_{\mathrm{sw}}$.

Suppressing local systems in the notation,
the Floer operations are $\Lambda$-multilinear maps
\[
    \mu^k:
    CF^*(L_{k-1},L_k)
    \otimes_{\Lambda}\cdots\otimes_{\Lambda}
    CF^*(L_0,L_1)
    \longrightarrow
    CF^*(L_0,L_k),
    \qquad k\geq1,
\]
of degree $2-k$, together with a curvature term
\[
    \mu^0_{(L,E)}
    \in CF^2((L,E),(L,E)).
\]
They satisfy the curved $A_\infty$ relations
and are constructed from moduli spaces of
holomorphic polygons using compatible
perturbation and orientation data.
A polygon contributing to $\mu^k$ has $k$
ordered inputs and one output.
Its contribution is weighted by local-system
parallel transport and by
$T^{\int u^*\omega}$.
The curvature term $\mu^0$ comes from discs
with no inputs. For these constructions, see
\cite{FOOO,AkahoJoyce,AFOOO}
and, for clean intersections,
\cite[Sections~3.3--3.4]{Fukaya-immersed}.

In the presence of curvature, $\mu^1$ need
not square to zero.
After deformation by bounding cochains,
introduced in Definition~\ref{def:immersed-lagrangian-brane} below, the unary operation is a
differential and the binary operation induces
composition on Floer cohomology.
The higher operations express associativity
up to homotopy and its higher coherence relations.
For the $A_\infty$ formalism, see \cite{Seidel}.

The valuation of a Floer cochain is the infimum
of the valuations of its nonzero coefficients
in the chosen cochain model, with $\operatorname{val}(0)=+\infty$.
The cochain spaces are complete for the
corresponding Novikov topology, and the operations
satisfy
\begin{equation}
    \label{eq:floer-valuation-bound}
    \operatorname{val}\bigl(
    \mu^k(x_k,\ldots,x_1)
    \bigr)
    \geq
    \sum_{j=1}^k\operatorname{val}(x_j),
    \qquad k\geq1.
\end{equation}

\begin{definition}[Branes]
    \label{def:immersed-lagrangian-brane}
    An \emph{immersed Lagrangian brane}
    is a triple
    \[
        \mathcal L=(L,E,b),
    \]
    where $L=[\iota:M\to X]$ is a Lagrangian
    immersion with compact, possibly disconnected
    source without boundary, equipped with a
    grading and relative spin structure as above,
    $E$ is a rank-one $\Bbbk$-local system on $M$,
    and
    \[
        b\in CF^1((L,E),(L,E))
    \]
    has strictly positive valuation and satisfies
    the Maurer--Cartan equation
    \begin{equation}
        \label{eq:immersed-maurer-cartan}
        \sum_{k\geq0}\mu^k(b,\ldots,b)=0.
    \end{equation}
    We call $b$ a \emph{bounding cochain}.
    The grading and relative spin structure are suppressed
    from the notation $(L,E,b)$.
\end{definition}

By~\eqref{eq:floer-valuation-bound}, the condition
$\operatorname{val}(b)>0$ ensures convergence
of the Maurer--Cartan series and of the series
obtained by inserting copies of $b$ into the
Floer operations;
see \cite[Section~3.6]{AkahoJoyce}.

\subsubsection{Morphisms, shifts, and the derived category}
\label{subsubsec:derived-category}
For two branes $\mathcal L_i=(L_i,E_i,b_i)$,
the deformed differential on
$CF^*((L_0,E_0),(L_1,E_1))$ is
\begin{equation}
    \label{eq:deformed-floer-differential}
    \mu^1_{b_1,b_0}(x)
    =
    \sum_{r,s\geq0}
    \mu^{r+s+1}
    \bigl(
    \underbrace{b_1,\ldots,b_1}_{r},
    x,
    \underbrace{b_0,\ldots,b_0}_{s}
    \bigr).
\end{equation}
The series converges by the positive valuations
of $b_0$ and $b_1$.
Their Maurer--Cartan equations and the
$A_\infty$ relations imply
$(\mu^1_{b_1,b_0})^2=0$;
see \cite[Section~3.6]{AkahoJoyce}.
We denote this complex by
$CF^*(\mathcal L_0,\mathcal L_1)$ and its
cohomology by $HF^*(\mathcal L_0,\mathcal L_1)$.
A \emph{closed degree-zero morphism} is a cochain
\[
    f\in CF^0(\mathcal L_0,\mathcal L_1),
    \qquad
    \mu^1_{b_1,b_0}(f)=0.
\]
Inserting the bounding cochains into all the
Floer operations similarly defines higher
operations $\mu^k_{b_k,\ldots,b_0}$.
These give an uncurved, $\Lambda$-linear,
$\mathbb Z$-graded $A_\infty$ category
$\mathcal F(X)$, whose objects are the
immersed Lagrangian branes
defined above and whose morphism complexes
are $CF^*(\mathcal L_0,\mathcal L_1)$.
We use a cohomologically unital model, so
the induced products on Floer cohomology
admit identity classes.
For the immersed categorical constructions,
see \cite{AFOOO}
and \cite[Section~3.4]{Fukaya-immersed}.

The shift $[1]$ adds $\pi$ to the grading of a Lagrangian and
reverses the orientation.
When referring only to the underlying graded
immersion, we also write $L[r]$ for the same
immersion with grading shifted by $r\pi$
and the corresponding orientation.
Our convention is
\begin{equation}
    \label{eq:floer-shift-convention}
    CF^k(\mathcal L_0[r],\mathcal L_1[s])
    \cong
    CF^{k+s-r}(\mathcal L_0,\mathcal L_1).
\end{equation}
In particular, a degree-zero morphism from
$\mathcal L_0$ to $\mathcal L_1$ becomes
degree one from $\mathcal L_0[1]$ to
$\mathcal L_1$.
We define the \emph{derived Fukaya category} by
\[
    \mathcal D
    :=
    H^0(\operatorname{Tw}\mathcal F(X)).
\]
Here $\operatorname{Tw}\mathcal F(X)$ is the
$A_\infty$ category of finite twisted complexes.
Its objects are finite formal direct sums
of shifted branes, equipped with degree-one
twisting matrices satisfying the
Maurer--Cartan equation and strictly triangular
with respect to an ordering of the summands.
Strict triangularity makes the sums involving
twisting entries finite, so these entries
need not have positive valuation.

Taking $H^0$ means passing to degree-zero
cohomology of the morphism complexes.
The resulting category $\mathcal D$ is
triangulated, with shifts and cones supplied
by the twisted-complex construction.
Viewing branes as one-term twisted complexes
gives natural identifications
\[
    \operatorname{Hom}_{\mathcal D}
    (\mathcal L_0,\mathcal L_1[k])
    \cong
    HF^k(\mathcal L_0,\mathcal L_1).
\]
Thus a closed degree-zero cochain $f$ determines
a categorical morphism $[f]$.
We write $\simeq$ for isomorphism in $\mathcal D$.
For these categorical constructions, see
\cite[Sections~2--3]{Seidel}.

\subsubsection{Disconnected sources}
\label{subsubsec:disconnected-sources}

The mapping-cone construction in
Subsection~\ref{subsec:initial-cone}
uses the following description of Floer
cochains on a disconnected source.
Let
\[
    L=[\iota:M_A\sqcup M_B\to X],
\]
and write
\[
    A=[\iota|_{M_A}:M_A\to X],
    \qquad
    B=[\iota|_{M_B}:M_B\to X].
\]
Thus $L=A\cup B$ denotes the union with
disjoint source.
The images of $A$ and $B$ may intersect in $X$,
so the self-Floer space of $L$ includes
cochains between the two source parts,
as well as the self-Floer cochains of each part.
Suppressing the local systems and their
restrictions to $M_A$ and $M_B$, we write
the resulting direct sum decomposition of graded modules as
\[
    CF^*(L,L)
    \cong
    \begin{pmatrix}
        CF^*(A,A) & CF^*(B,A) \\
        CF^*(A,B) & CF^*(B,B)
    \end{pmatrix}.
\]
Rows indicate targets and columns indicate sources.
For example, the upper-right block consists of cochains
from $B$ to $A$.
The diagonal blocks contain the full self-Floer spaces
of $A$ and $B$, including their internal self-intersection
generators.

We choose Floer models compatible with this decomposition.
The operations follow the same rule for matching indices
as matrix multiplication: the target part of one input
must agree with the source part of the next.
Geometrically, this reflects the fact that each boundary
segment of a counted polygon lifts to a single source part.
For example, an input from $A$ to $B$ followed by an input
from $B$ to $A$ contributes to the $A$-to-$A$ block. For the disjoint-source construction in the
transverse case, see
\cite[Section~13.2]{AkahoJoyce}.
For the cochain-level decomposition and its
compatibility with the operations in the
clean-intersection setting, see
\cite[Section~3.4 and Remark~3.48]{Fukaya-immersed}.

A bounding cochain on $L$ consists of diagonal
cochains $b_A,b_B$ together with possible
off-diagonal entries coupling the two source parts.
If the off-diagonal entries vanish, the
Maurer--Cartan equation separates into the
equations for $b_A$ and $b_B$.
The resulting brane is then isomorphic in
$\mathcal D$ to the direct sum of the
corresponding branes on $A$ and $B$.

\subsubsection{Regular transport}

We return to admissible web data $\mathcal W_I$ on
$I=(a,b)$.
Choose a compact subinterval $J\subset I$ and a coherent
smoothing
\[
    C_\tau=[\iota_\tau:N\longrightarrow X],
    \qquad \tau\in J,
\]
as in Proposition~\ref{prop:web-cone-smoothing}.
The resulting exact deformation is generated by the
local Hamiltonian $\widetilde\kappa_\tau.$

\begin{prop}[Regular transport]
    \label{prop:web-cone-brane-transport}
    Fix $\tau_0\in J$ such that $C_{\tau_0}$
    has only transverse double points, and choose
    a brane
    \[
        \mathcal C_{\tau_0}
        =
        (C_{\tau_0},E_{\tau_0},b_{\tau_0}),
        \qquad
        \operatorname{val}(b_{\tau_0})>0.
    \]
    Its brane structure admits transport to nearby
    web cones, giving branes
    \[
        \mathcal C_\tau=(C_\tau,E_\tau,b_\tau)
        \simeq\mathcal C_{\tau_0}
    \]
    in $\mathcal D$.
    Regular transport can be continued along $J$
    as long as the smoothed web cones have only
    transverse double points and $b_\tau$ has
    strictly positive valuation.

    The local system is transported using the fixed
    source $N$.
    With compatible choices of Floer data, the ordinary
    source-cochain part of $b_\tau$ is unchanged.
    For a smoothly varying ordered degree-one
    self-intersection generator
    \[
        q_\tau=(q_-(\tau),q_+(\tau)),
    \]
    its coefficient satisfies
    \begin{equation}
        \label{eq:bounding-cochain-parallel-transport}
        b_{\tau,q_\tau}
        =
        T^{-\int_{\tau_0}^{\tau}
                \left(
                \widetilde\kappa_s(q_+(s))
                -
                \widetilde\kappa_s(q_-(s))
                \right)\,ds}
        b_{\tau_0,q_{\tau_0}}.
    \end{equation}
    Here the coefficient spaces are identified by
    parallel transport in $E_{\tau_0}$ along the two
    source paths, and orientation generators are
    chosen continuously.
\end{prop}

\begin{remark}
    \label{rem:exact-brane-transport}
    The same conclusion holds for an exact Lagrangian
    deformation on a fixed source, generated by a
    local Hamiltonian $H_s$, provided that all
    intersections remain transverse double points,
    apart from fixed clean intersections between
    closed source components that remain stationary.
    In Equation~\eqref{eq:bounding-cochain-parallel-transport},
    one uses $H_s$ in place of $\widetilde\kappa_s$.
\end{remark}

We call this \emph{regular transport} of the
brane structure along the web cones.
The proof will be given in~\cite{KapotaWebFloer}.
The Novikov rescaling compensates for the
change in polygon areas under the deformation,
as in the area-comparison argument of
\cite[Example~13.16]{AkahoJoyce}.
The categorical part of the proof produces
degree-zero Floer cohomology classes in both
directions between $\mathcal C_{\tau_0}$ and
$\mathcal C_\tau$ whose compositions are the
identity classes.
These establish the stated isomorphism
in $\mathcal D$.

In particular, the valuation of each nonzero
self-intersection coefficient evolves according to
\begin{equation}
    \label{eq:web-cone-valuation-transport}
    \operatorname{val}(b_{\tau,q_\tau})
    =
    \operatorname{val}(b_{\tau_0,q_{\tau_0}})
    -
    \int_{\tau_0}^{\tau}
    \left(
    \widetilde\kappa_s(q_+(s))
    -
    \widetilde\kappa_s(q_-(s))
    \right)\,ds.
\end{equation}
Thus, if the two preimages remain at heights $0$ and $1$,
respectively, the valuation decreases at unit rate.
The height function therefore determines when a
transported coefficient can reach valuation zero,
even while the geometric evolution remains smooth. This can occur even when the web cone is geometrically
stationary, since the height may take different constant
values on distinct source components

The transversality assumption above concerns the smoothed
web cones.
Self-tangencies of these cones may occur during a regular
interval of front evolution;
transport across such times is treated in the following Subsubsection.

\subsubsection{Transport through self-tangencies}
\label{subsec:self-tangencies}

During a regular front evolution, the smoothed web cone may develop a self-tangency at which two transverse
double points appear or disappear.
The source remains unchanged, but the change in Floer
generators requires a correction to the bounding cochain.
The following result gives continuation through an
isolated self-tangency when the incoming bounding
cochains have valuation uniformly bounded away from zero.

\begin{prop}[Transport through an isolated self-tangency]
    \label{prop:self-tangency-transport}
    Let
    \[
        \mathcal C_t=(C_t,E_t,b_t),
        \qquad -\epsilon<t<0,
    \]
    be branes related by regular transport.
    Suppose that the underlying immersions extend to
    a coherently smoothed family of web cones
    \[
        C_t=[\iota_t:N\longrightarrow X],
        \qquad -\epsilon<t<\epsilon,
    \]
    with fixed closed source $N$ and exact deformation
    generated by $H_t=\widetilde\kappa_t$.
    Assume that this family has only transverse double
    points except for a single quadratic self-tangency
    at $t=0$, in the local form of
    \cite[Definition~2.7(a)]{PalmerWoodwardMaslov}.
    The tangency occurs in a neighborhood meeting only
    the two participating local source sheets.

    Suppose that there is $\eta>0$ such that
    \[
        \operatorname{val}(b_t)\geq\eta
    \]
    for all negative $t$ sufficiently close to zero.
    After shrinking $\epsilon$, the brane structures
    admit a continuation to branes
    \[
        \mathcal C_t=(C_t,E_t,b_t),
        \qquad 0<t<\epsilon,
    \]
    with $\operatorname{val}(b_t)>0$, related by regular
    transport on the outgoing interval, and satisfying
    \[
        \mathcal C_{t_-}\simeq\mathcal C_{t_+}
        \qquad
        \text{for }
        -\epsilon<t_-<0<t_+<\epsilon
    \]
    in $\mathcal D$.
    The local system is identified throughout using
    the fixed source $N$.
    The statement applies to both the birth and the
    death of a pair of transverse double points.
\end{prop}

The proof, in the case of the two-torus, will be given
in~\cite{KapotaWebFloer}.
It adapts the self-tangency analysis of
\cite[Sections~8--9]{PalmerWoodwardMaslov}.

\section{Initiation, compatible continuation, and termination}
	\label{sec:singular-web-evolution}

The preceding section constructs web cones and transports
their brane structures along regular intervals.
We now describe evolutions starting from the mapping cone
of a morphism and continuing through geometric
degenerations and obstruction walls.
After initialization, an evolution alternates between
regular transport and continuation at events.
Each complete continuation supplies outgoing web-cone
branes with bounding cochains of strictly positive
valuation, representing the same object as the incoming
branes.
If all active fronts disappear, the geometry becomes
stationary, but the brane data continue to evolve.
A later obstruction wall may require another resolution
and create new active fronts.
Terminality means that this stationary continuation
encounters no further obstruction wall.

Subsection~\ref{subsec:initial-cone} specifies the
initial cone data and defines initialization,
as used in Conjecture~\ref{conj:intro-initialization}.
Subsection~\ref{subsec:web-continuation} defines
compatible continuation and brane transfer as used in
Conjecture~\ref{conj:intro-continuation}.  In Subsection~\ref{subsec:terminal-triangle}, we define
  terminality and construct the exact
  triangle~\eqref{eq:terminal-web-triangle} from a terminal
  evolution.
The local models for the operations at initialization
and subsequent events are described in
Section~\ref{sec:singular-events}.
In Subsection~\ref{subsec:terminal-triangle}, we construct
the exact triangle~\eqref{eq:terminal-web-triangle}
from a terminal evolution.
Finally, Subsection~\ref{subsec:Dirichlet-energy}
gives energy and volume estimates, including a
bound on the total active time.

\subsection{The initial cone}
\label{subsec:initial-cone}
Let
\[
    \mathcal L_i=(L_i,E_i,b_i),
    \qquad i=0,1,
\]
be branes with
$\operatorname{val}(b_i)>0$, where
$L_i=[\iota_i:M_i\to X]$.
Assume that each component of $L_i$ is either positive,
with grading taking values in $(-\pi/2,\pi/2)$,
or closed imaginary special, with grading $\pi/2$.
Assume also that the disjoint-source union of
$L_0$ and $L_1$ has clean self-intersections.
Let
\[
    f\in CF^0(\mathcal L_0,\mathcal L_1)
\]
be closed.
By \eqref{eq:floer-shift-convention}, $f$ is
a closed degree-one morphism from
$\mathcal L_0[1]$ to $\mathcal L_1$.

\begin{definition}
    The \emph{categorical cone} of $f$ is the twisted complex
    \[
        \mathcal C_f
        :=
        \operatorname{Cone}(f)
        :=
        \left(
        \mathcal L_0[1]\oplus\mathcal L_1,\,
        \delta_f
        \right),
        \qquad
        \delta_f=
        \begin{pmatrix}
            0 & 0 \\
            f & 0
        \end{pmatrix}.
    \]
    The bounding cochains of the input branes are
    already incorporated into the two summands.
    The matrix $\delta_f$ specifies the additional
    twisting morphism between them.
\end{definition}

\begin{remark}
    \label{rem:cone-scaling}
    Multiplying $f$ by a nonzero scalar does not change
    the isomorphism class of its cone:
    \[
        \operatorname{Cone}(cf)\simeq\operatorname{Cone}(f),
        \qquad c\in\Lambda^\times.
    \]
    Indeed,
    \[
        (c\,1_{\mathcal L_1})\circ[f]
        =
        [cf]\circ 1_{\mathcal L_0}.
    \]
    Thus the identity on $\mathcal L_0$ and the
    automorphism $c\,1_{\mathcal L_1}$ give a commutative
    square between $[f]$ and $[cf]$, with invertible
    vertical maps, and hence an isomorphism of their cones.
\end{remark}
For $f\neq0$, we henceforth replace $f$ by
$T^{-\operatorname{val}(f)}f$ and continue
to denote the resulting cochain by $f$.
Thus $\operatorname{val}(f)=0$.
The notation $\mathcal C_f$ refers to the
cone of this normalized cochain, which is
isomorphic to the original cone.
No normalization is needed when $f=0$.

\begin{definition}[Initial cone]
    \label{def:initial-cone}
    The \emph{initial cone data} associated to $f$ are
    \[
        C_f=L_0[1]\cup L_1,
        \qquad
        E_f=E_0\sqcup E_1,
        \qquad
        B_f=
        \begin{pmatrix}
            b_{0} & 0 \\
            f & b_1
        \end{pmatrix},
    \]
    where the union has source $M_0\sqcup M_1$
    and $b_{0}$ denotes the bounding cochain
    of $\mathcal L_0[1]$.
    With this source decomposition understood,
    we write
    \[
        \mathcal C_f=(C_f,E_f,B_f)
    \]
    for the categorical cone $\operatorname{Cone}(f)$.
    The diagonal entries of $B_f$ record the bounding
    cochains of the summands, and $f$ is the
    additional twisting morphism.
\end{definition}

The following lemma verifies that the Floer
operations can be computed directly by inserting
the total cochain $B_f$.

\begin{lemma}[Initial cone data]
    \label{lem:geometric-cone}
    The cochain $B_f$ satisfies the
    Maurer--Cartan equation.
    All series obtained by inserting $B_f$
    into Floer operations with a fixed number
    of external inputs converge in the Novikov
    topology.
    Under the source decomposition
    $M_0\sqcup M_1$, the resulting operations
    agree with the twisted-complex operations
    on $\mathcal C_f$.
    In particular, the initial cone
    fits into the exact triangle
    \begin{equation}
        \label{eq:initial-cone-triangle}
                \mathcal L_0
        \xrightarrow{[f]}
        \mathcal L_1
        \longrightarrow
        \mathcal C_f
        \longrightarrow
        \mathcal L_0[1].
    \end{equation}
\end{lemma}

\begin{proof}
    By compatibility with the source decomposition,
    a nonzero contribution to a Floer operation
    requires a composable sequence of matrix entries.
    The entry $f$ passes from the $M_0$-source to
    the $M_1$-source, and $B_f$ has no entry in
    the reverse direction.
    Consequently, a composable sequence of entries
    of $B_f$ contains at most one occurrence of $f$.
    Choose $\delta>0$ such that
    \[
        \operatorname{val}(b_{0[1]})\geq\delta,
        \qquad
        \operatorname{val}(b_1)\geq\delta.
    \]
    Since $\operatorname{val}(f)\geq0$, every
    contribution with $k\geq1$ entries of $B_f$
    has valuation at least $(k-1)\delta$.
    Thus the Maurer--Cartan series converges.
    More generally, in an operation with $r$ fixed
    external inputs, there are $r+1$ blocks of
    inserted entries of $B_f$.
    Each block contains at most one occurrence
    of $f$, so the valuations again tend to
    $+\infty$ as the number of insertions grows.
    This proves convergence of the other
    insertion series.

    We now examine the matrix components of the
    Maurer--Cartan equation.
    The diagonal components are the
    Maurer--Cartan equations for $b_{0[1]}$ and $b_1$,
    which vanish because the input branes are
    unobstructed.
    The upper-right component vanishes because
    there is no entry from the $M_1$-source to
    the $M_0$-source.
    The lower-left component is
    \[
        \sum_{r,s\geq0}
        \mu^{r+s+1}
        \bigl(
        \underbrace{b_1,\ldots,b_1}_{r},
        f,
        \underbrace{b_{0[1]},\ldots,b_{0[1]}}_{s}
        \bigr)
        =
        \mu^1_{b_1,b_{0[1]}}(f)
        =
        0,
    \]
    where $f$ is regarded as a closed degree-one
    cochain after shifting $\mathcal L_0$.
    Hence $B_f$ satisfies the Maurer--Cartan equation.

    The twisted-complex operations insert copies
    of $f$ into the operations already deformed
    by $b_{0[1]}$ and $b_1$.
    Expanding those deformed operations gives
    exactly the composable insertions of $B_f$.
    The convergence established above permits
    this regrouping.
    Thus the operations agree, and the displayed
    triangle is the cone triangle of $\mathcal C_f$.
\end{proof}

\begin{definition}[Initial web configuration]
    \label{def:initial-web-configuration}
    Write
    \[
        L_i=P_i\sqcup S_i,
        \qquad i=0,1,
    \]
    for the decomposition by source components,
    where $P_i$ consists of the positive components
    and $S_i$ of the closed imaginary-special components.
    The \emph{initial web configuration} is the
    configuration at $\tau=0$ of the constant
    web data of Definition~\ref{def:regular-web-data}
    with boundary conditions $P_0,P_1$, empty
    active front and incoming boundary traces,
    and inactive fragments consisting of the
    components of $S_i$, each carrying label $i$.

    These data are admissible in the sense of
    Definition~\ref{def:admissible-web}.
    Their \emph{initial web cone} and its height
    function are
    \[
        C_0=C_f,
        \qquad
        \kappa_0|_{M_0}=0,
        \qquad
        \kappa_0|_{M_1}=1.
    \]
    We associate to this configuration the cone
    object $\mathcal C_0:=\mathcal C_f\in\mathcal D$.
\end{definition}

For $f\neq0$, one has $\operatorname{val}(B_f)=0$,
so the initial cone data do not define an
unobstructed brane.
Initialization supplies web-cone
branes representing the same object.

\begin{definition}[Initialization]
    \label{def:web-initialization}
    \label{def:initial-brane-comparison}
    An \emph{initialization of the web evolution of $f$}
    consists of admissible web data
    $\mathcal W_I$ on $I=(0,\epsilon)$,
    $\epsilon>0$, in the sense of
    Definition~\ref{def:admissible-web},
    with empty incoming boundary traces.
    These data are obtained from the initial web
    configuration of
    Definition~\ref{def:initial-web-configuration}
    by operations described in
    Section~\ref{sec:singular-events}
    performed at $\tau=0$.
    All source pieces not involved in these
    operations are retained with their immersions
    and labels.
    The initialization also includes a
    brane $\mathcal C_{\tau_0}$ on a smoothing of
    the web cone at some $\tau_0\in I$, whose
    regular brane transport is defined throughout
    $I$, together with a chosen isomorphism
    \[
        \Phi_{\tau_0}:
        \mathcal C_0
        \xrightarrow{\sim}
        \mathcal C_{\tau_0}
        \qquad\text{in }\mathcal D.
    \]
    We call $\Phi_{\tau_0}$ an
    \emph{initial brane comparison}. Regular transport extends this comparison
throughout $I$.
\end{definition}

When $f=0$, the initial cone data already define
a brane.
Each diagonal block is supported on a source
part of constant height, so the configuration
can remain stationary with constant brane data.
This gives an initialization with empty active
front and the identity as initial brane comparison.

Conjecture~\ref{conj:intro-initialization}
asserts that every closed degree-zero cochain
representing a morphism in $\mathcal A$
admits an initialization.

\subsection{Compatible continuation and brane transfer}
\label{subsec:web-continuation}

We now assemble regular web data into a web evolution
and specify the brane comparisons across events.
Here \emph{events} include geometric degenerations,
obstruction walls, and changes in the intersection
pattern of the smoothed web cone requiring separate
brane transport, such as isolated self-tangencies.
The local operations at geometric degenerations and
obstruction walls are described in
Section~\ref{sec:singular-events}.

\subsubsection{Geometric continuation and special Lagrangian webs}
\label{subsubsec:geometric-web-continuation}
The local models prescribe changes to the participating
fronts and fragments.
Compatible continuation requires these changes to
match the remaining pieces and preserve the regions
already swept out on the boundary sources.

\begin{definition}[Compatible continuation]
    \label{def:surgery-compatible}
    Let $\mathcal W^\pm$ be admissible web data on
    \[
        I^-=(\tau_*-\epsilon,\tau_*),
        \qquad
        I^+=(\tau_*,\tau_*+\epsilon),
    \]
    with boundary conditions
    \[
        L_i^\pm=[\iota_i^\pm:M_i^\pm\longrightarrow X],
        \qquad i=0,1.
    \]
    Write $P^\pm$ for the disjoint unions of their
    front and fragment sources, with immersions
    $F_\tau^\pm$ and boundary lifts
    $\beta_{i,\tau}^\pm$.
    Fix local surgeries or closed-fragment replacements
    at $\tau_*$, with specified attaching maps,
    allowing also changes of $M_i^-$.
    The constructions of
    Section~\ref{sec:singular-events} provide examples.
    Denote the removed and inserted source regions by
    $A^\pm\subset P^\pm$, respectively, and set
    \[
        U^\pm=P^\pm\setminus A^\pm.
    \]
    The regions $A^\pm$ are closed and may contain
    whole source components.
    These operations give rise to a
    \emph{compatible continuation} from
    $\mathcal W^-$ to $\mathcal W^+$
    if there exist orientation-preserving
    diffeomorphisms
    \[
        \phi:U^-\xrightarrow{\sim}U^+,
        \qquad
        \psi_i:V_i^-\xrightarrow{\sim}V_i^+,
    \]
    where $V_i^\pm\subset M_i^\pm$ are open,
    satisfying the following conditions.

    \begin{enumerate}
        \item \emph{Fronts and fragments.}
              The map $\phi$ preserves boundary labels:
              \[
                  \phi(U^-\cap\partial_iP^-)
                  =
                  U^+\cap\partial_iP^+.
              \]
              It sends each retained closed fragment
              to a closed fragment with the same label.
              The restrictions $F_\tau^\pm|_{U^\pm}$
              have one-sided limits $F_*^\pm$ in
              $C^\infty_{\mathrm{loc}}(U^\pm)$,
              satisfying
              \[
                  F_*^+\circ\phi=F_*^-.
              \]

        \item \emph{Boundary sources and lifts.}
              The boundary immersions satisfy
              \[
                  \iota_i^+\circ\psi_i
                  =
                  \iota_i^-|_{V_i^-}.
              \]
              The restrictions of $\beta_{i,\tau}^\pm$
              to $U^\pm\cap\partial_iP^\pm$
              have one-sided limits
              $\beta_{i,*}^\pm$ in
              $C^\infty_{\mathrm{loc}}
              (U^\pm\cap\partial_iP^\pm)$.
              We require
              \[
                  \beta_{i,*}^-
                  (U^-\cap\partial_iP^-)
                  \subset V_i^-,
              \]
              and, on $U^-\cap\partial_iP^-$,
              \[
                  \beta_{i,*}^+\circ\phi
                  =
                  \psi_i\circ\beta_{i,*}^-.
              \]

        \item \emph{Accumulated boundary traces.}
              Set
              \[
                  B_i(\tau_*^-)
                  :=
                  \operatorname{Int}_{M_i^-}
                  \overline{
                      \bigcup_{\tau\in I^-}B_i^-(\tau)
                  }.
              \]
              We require
              \[
                  \overline{B_i(\tau_*^-)}
                  \subset V_i^-,
              \]
              and
              \begin{equation}
                  \label{eq:continuation-trace-matching}
                  B_i^{+,\mathrm{in}}
                  =
                  \psi_i\bigl(B_i(\tau_*^-)\bigr).
              \end{equation}
    \end{enumerate}
\end{definition}

\begin{remark}[Auxiliary perturbations]
    \label{rem:auxiliary-perturbations}
    At a regular time, we may make an arbitrarily
$C^\infty$-small Hamiltonian perturbation of the strictly
    positive boundary representatives, fixing the
    current fronts and fragments, their labels
    and boundary lifts, and the accumulated swept
    regions under the natural source identifications.
    In particular, admissibility is preserved. These perturbations are used only to prepare
    the next event, with at most one
    perturbation between successive events.

    We require the induced deformation of the
    coherently smoothed web cones to be exact.
    When brane structures are included, they are
    transported along this deformation with
    bounding cochains of strictly positive valuation.

    Regular evolution then resumes with the perturbed
    boundary representatives and the accumulated
    swept regions as incoming traces.
    Subsequent local operations are applied to
    the resulting data.
\end{remark}

\begin{definition}[Special Lagrangian web]
    \label{def:special-lagrangian-web}
    A \emph{special Lagrangian web} $\mathcal W$
    on an open interval $J\subset\mathbb R$ consists
    of a closed discrete set $E\subset J$ of
    \emph{event times}, admissible web data
    $\mathcal W_I$ on each connected component
    $I$ of $J\setminus E$, and chosen compatible
    continuations at the times in $E$.
    We write $\mathcal I(J,E)$ for these connected
    components and call the times in $J\setminus E$
    \emph{regular times}.
\end{definition}

\subsubsection{Brane transfer}
\label{subsubsec:brane-decorated-webs}

We now equip special Lagrangian webs with brane
structures and comparison isomorphisms.
On regular intervals, we use the transport of
Proposition~\ref{prop:web-cone-brane-transport},
with the extension of
Remark~\ref{rem:exact-brane-transport}
for retained clean intersections between
    closed source components that remain stationary.

\begin{definition}[Brane transfer]
    \label{def:brane-surgery-transfer}
    \label{def:brane-decorated-web-evolution}
    Fix a compatible continuation at $\tau_*$,
    with incoming web-cone branes
    $\mathcal C_\tau^-$ on $I^-$ related by
    regular transport.

    A \emph{brane transfer} across this continuation
    consists of brane structures
    $\mathcal C_\tau^+$ on coherent smoothings
    of the outgoing web cones for $\tau\in I^+$,
    related by regular transport, together with
    a chosen isomorphism
    \[
        \Psi_*:
        \mathcal C_{\tau_-}^-
        \xrightarrow{\sim}
        \mathcal C_{\tau_+}^+
        \qquad\text{in }\mathcal D
    \]
    at nearby regular times $\tau_\pm\in I^\pm$.
    Comparisons at other regular times are obtained
    by composing $\Psi_*$ with regular transport
    on the two sides.
    If the continuation includes a closed-component
    reduction, its algebraic data, as specified in
    Definition~\ref{def:closed-component-reduction},
    are included in the transfer.

    A special Lagrangian web equipped with
    regularly transported branes
    on coherent smoothings of its web cones and
    chosen brane transfers at every event
    is called a \emph{brane-decorated web evolution}.
\end{definition}

The brane comparison is taken after the complete
continuation, including all operations performed
at the same time.
No brane of strictly positive valuation is required
at the event time itself.
For an isolated self-tangency,
the comparison is supplied by
Proposition~\ref{prop:self-tangency-transport}.
The existence of compatible continuations and brane
transfers is part of
Conjecture~\ref{conj:intro-continuation}.

\begin{definition}[Special Lagrangian web for $f$]
    \label{def:f-admissible-web-evolution}
    A \emph{special Lagrangian web for $f$} is a
    brane-decorated web evolution on $(0,T)$,
    $0<T\leq\infty$, extending a given
    initialization in the sense of
    Definition~\ref{def:web-initialization}.
\end{definition}

\begin{lemma}
    \label{lem:web-cone-comparison}
    For every regular time $\tau$ of a special
    Lagrangian web for $f$, the chosen initialization
    and transport data determine an isomorphism
    \[
        \Phi_\tau:
        \operatorname{Cone}(f)
        \xrightarrow{\sim}
        \mathcal C_\tau
        \qquad\text{in }\mathcal D.
    \]
    These isomorphisms are compatible with regular
    transport, brane transfer across events, and
    transport along auxiliary perturbations.
\end{lemma}

\begin{proof}
    Compose the initial brane comparison with
    regular transport, the intervening brane
    transfers, and transport along auxiliary
    perturbations.
    This is a finite composition, since the
    event set is closed and discrete in
    $(0,T)$ and both the initial comparison time
    and $\tau$ lie in this interval.
\end{proof}

\subsection{Termination and the terminal triangle}
\label{subsec:terminal-triangle}
We first describe stationary brane continuation
and define terminality.
We then construct the branes and exact triangle of
Proposition~\ref{prop:terminal-web-triangle}
from a terminal web.

\subsubsection{Stationary continuation and termination}
\label{subsubsec:web-terminality}

The disappearance of all active fronts need not end
the brane evolution.

\begin{definition}[Geometric exhaustion]
    \label{def:geometrically-exhausted}
    A web is \emph{geometrically exhausted} at a
    time $\tau_*$ if its active front is empty
    after all operations at that time have
    been completed.
\end{definition}

At geometric exhaustion, the height function of
Definition~\ref{def:height-function} is constant,
equal to $0$ or $1$, on each connected component
of the glued cone source.
A smoothing preserving these constants has
$d\widetilde\kappa=0$, so the prescribed deformation
keeps the immersion stationary.

The brane transport still depends on the differences
between the component heights.
By Equation~\eqref{eq:web-cone-valuation-transport},
a nonzero bounding-cochain coefficient from a
height-zero branch to a height-one branch decreases
in valuation at unit speed.
If such a coefficient is present, stationary
continuation reaches an obstruction wall after
a finite waiting time.
Continuing past this wall requires a compatible
continuation and a brane transfer, and may create
new active fronts.

\begin{definition}[Stationary brane continuation]
    \label{def:stationary-brane-continuation}
    Suppose that a brane-decorated web evolution
    is geometrically exhausted at $\tau_*$.
    If $\tau_*$ is regular, set $\tau_0=\tau_*$;
    otherwise, choose $\tau_0$ in the immediately
    outgoing regular interval.
    We also allow $\tau_*=\tau_0=0$ when the initial
    configuration is geometrically exhausted and
    its cone cochain has strictly positive valuation.

    Use a smoothing preserving the constant heights,
    and let
    \[
        \mathcal C_{\tau_0}
        =
        (C_{\tau_0},E_{\tau_0},b_{\tau_0}),
        \qquad
        \operatorname{val}(b_{\tau_0})>0,
    \]
    be the corresponding brane.
    The \emph{stationary brane continuation}
    keeps the immersion and auxiliary Floer
    data fixed,
    \[
        C_s=C_{\tau_0},
        \qquad s\geq\tau_0,
    \]
    and transports the brane data by regular
    transport with Hamiltonian
    \[
        H=\widetilde\kappa_{\tau_0},
    \]
    for as long as the bounding cochain has
    strictly positive valuation.
    The local system is identified throughout
    using the fixed source.
\end{definition}

\begin{definition}
    \label{def:brane-terminal-web}
    A brane-decorated web evolution is
    \emph{terminal from $\tau_*$} if it is
    geometrically exhausted there and its
    stationary brane continuation is defined
    without encountering an obstruction wall
    for every finite $s\geq\tau_0$.
    Here $\tau_0$ is chosen as in
    Definition~\ref{def:stationary-brane-continuation}.
\end{definition}

\begin{remark}
For $f=0$, the stationary initialization described in
  Subsection~\ref{subsec:initial-cone} is terminal from
  time zero.
\end{remark}

Conjecture~\ref{conj:intro-termination} asserts
that, after generic Hamiltonian perturbations
as in Conjecture~\ref{conj:intro-continuation},
the continued special Lagrangian web for $f$
becomes terminal in finite time, after finitely
many geometric degenerations and obstruction walls.

\subsubsection{The terminal Lagrangians $K$ and $K^*$}
\label{subsec:terminal-smoothing}

Let $\mathcal W$ be a special Lagrangian web for $f$
that is terminal from time $\tau_*$.
We use the stationary brane $\mathcal C_{\tau_0}$
chosen in
Definition~\ref{def:stationary-brane-continuation}.
All geometric data below refer to its fixed terminal
configuration.
Write
\[
    L_i=[\iota_i:M_i\longrightarrow X],
    \qquad i=0,1,
\]
for the terminal boundary conditions, after any
absorptions and auxiliary perturbations, and let $R_i\subset M_i$ denote their residues.
We regard $\overline{R_i}$ as a Lagrangian piece
through $\iota_i|_{\overline{R_i}}$.
Let $D_i$ and $D_{\mathrm{cl},i}$ denote the inactive
fragments of
Definition~\ref{def:inactive-fragment-notation}.
Since no active fronts remain, the residues and
these fragments constitute the entire terminal
web cone.

\begin{definition}[Singular terminal Lagrangians]
    \label{def:terminal-branes}
    The \emph{singular terminal Lagrangians} are
    \[
        K_{\mathrm{sing}}
        :=
        \bigl(\overline{R_0}\cup D_0[-1]\bigr)
        \sqcup D_{\mathrm{cl},0},
        \qquad
        K_{\mathrm{sing}}^*
        :=
        \bigl(\overline{R_1}\cup D_1\bigr)
        \sqcup D_{\mathrm{cl},1}.
    \]
    The nonclosed fragments are glued to the residues
    along their chosen boundary lifts, as in
    Definition~\ref{def:singular-web-cone}.
    The closed fragments are adjoined as separate
    source components.
    With the inherited gradings, the singular terminal
    web cone is
    \[
        C_{\tau_0}^{\mathrm{sing}}
        =
        K_{\mathrm{sing}}[1]\sqcup K_{\mathrm{sing}}^*.
    \]
    Here $\sqcup$ denotes a disjoint union of sources.
\end{definition}

Write the chosen stationary smoothing as
\[
    C_{\tau_0}=\widehat K[1]\sqcup\widehat K^*,
\]
where $\widehat K$ and $\widehat K^*$ smooth the
corresponding singular terminal pieces.
Let $N_K$ and $N_{K^*}$ denote their sources.
The heights used in the stationary brane continuation
are
\[
    \widetilde\kappa|_{N_K}=0,
    \qquad
    \widetilde\kappa|_{N_{K^*}}=1.
\]
The bounding cochain may still couple these two
source parts; we examine it in the next subsubsection.

We next deform the stationary smoothing to obtain
the phase conditions required for objects of
$\mathcal A$.
The following geometric statement will be proved
in~\cite{KapotaWebGeometry}, using positive
perturbations and corner smoothing.

\begin{lemma}[Positive terminal representatives]
    \label{lem:terminal-smoothing}
    The terminal corner roundings may be chosen so
    that exact perturbation and smoothing of
    \[
        \overline{R_0}\cup D_0[-1],
        \qquad
        \overline{R_1}\cup D_1,
    \]
    give strictly positive Lagrangian immersions
    $K_{\mathrm{pos}}$ and $K_{\mathrm{pos}}^*$,
    respectively.
    Set
    \[
        K:=K_{\mathrm{pos}}\sqcup D_{\mathrm{cl},0},
        \qquad
        K^*:=K_{\mathrm{pos}}^*\sqcup D_{\mathrm{cl},1}.
    \]
    There is an exact Lagrangian deformation on the
    fixed source $N_K\sqcup N_{K^*}$ from
    \[
        \widehat K[1]\sqcup\widehat K^*
        \quad\text{to}\quad
        K[1]\sqcup K^*,
    \]
    preserving the two source parts and their gradings.
    The closed fragments remain fixed, with phase
    $\pi/2$ in the unshifted Lagrangians $K$ and $K^*$.
    The perturbations and roundings, together with
    the deformation primitives, may be chosen
    arbitrarily small.
\end{lemma}

We call $K$ and $K^*$ the \emph{smooth terminal
    Lagrangians}.
Each connected component is either strictly positive
or closed imaginary special of phase $\pi/2$.

\subsubsection{The terminal bounding cochain and exact triangle}
\label{subsubsec:terminal-triangle-proof}

We first use terminality to determine the block structure
of the bounding cochain on the stationary representative
\[
    \mathcal C_{\tau_0}
    =
    \bigl(\widehat K[1]\sqcup\widehat K^*,
    E_{\tau_0},b_{\tau_0}\bigr).
\]
We then transport this structure to the positive
representatives of Lemma~\ref{lem:terminal-smoothing}
and obtain the terminal triangle.

\begin{lemma}[Vanishing of the terminal forward block]
    \label{lem:terminal-forward-block}
    The component of $b_{\tau_0}$ directed from
    $\widehat K[1]$ to $\widehat K^*$ vanishes.
    Thus, with respect to the two source parts,
    \[
        b_{\tau_0}
        =
        \begin{pmatrix}
            b_{\widehat K} & \widehat\alpha   \\
            0                 & b_{\widehat K^*}
        \end{pmatrix},
        \qquad
        \widehat\alpha
        \in CF^1(\widehat K^*,\widehat K[1]),
    \]
    with the inherited local systems.
    Each diagonal block contains all cochain data
    within its source part, including terms between
    distinct connected components.
\end{lemma}

\begin{proof}
    Write $\beta_s$ for the block of the transported
    bounding cochain directed from $\widehat K[1]$
    to $\widehat K^*$.
    The Hamiltonian has constant values
    $0$ and $1$ on these source parts.
    Thus stationary transport gives
    \[
        \beta_s
        =
        T^{-(s-\tau_0)}\beta_{\tau_0}.
    \]
    This rescaling applies to the entire cochain
    block, including cochains on clean-intersection
    loci.
    It is compatible with the Floer operations
    because their source indices compose.
    If $\beta_{\tau_0}\neq0$, then
    \[
        \operatorname{val}(\beta_s)
        =
        \operatorname{val}(\beta_{\tau_0})
        -(s-\tau_0)
    \]
    would cease to be strictly positive at a finite
    time.
    This contradicts terminality.
    Hence $\beta_{\tau_0}=0$.
\end{proof}

\begin{proof}[Proof of Proposition~\ref{prop:terminal-web-triangle}]
    Choose the stationary representative so that its
    intersections are transverse double points,
    apart from retained clean intersections between
    fixed closed components.
    Choose the exact deformation of
    Lemma~\ref{lem:terminal-smoothing} sufficiently
    small to preserve this intersection pattern
    and to keep the transported bounding cochain
    of strictly positive valuation.

    Apply Proposition~\ref{prop:web-cone-brane-transport},
    using the deformation primitive in place of
    $\widetilde\kappa$, with the extension of
    Remark~\ref{rem:exact-brane-transport}.
    This gives a brane
    \[
        \mathcal C
        =
        \bigl(K[1]\sqcup K^*,E,b\bigr)
        \simeq\mathcal C_{\tau_0},
        \qquad
        \operatorname{val}(b)>0.
    \]

    The deformation preserves the two source parts.
    Transport rescales the tracked intersection
    coefficients and the retained clean-intersection
    blocks by Novikov factors.
    It therefore preserves the vanishing block of
    Lemma~\ref{lem:terminal-forward-block}, so
    \[
        b=
        \begin{pmatrix}
            b_{K} & \alpha  \\
            0        & b_{K^*}
        \end{pmatrix}.
    \]
    Restrict the local system to the two source parts:
    \[
        E_K:=E|_{N_K},
        \qquad
        E_{K^*}:=E|_{N_{K^*}}.
    \]
    By source compatibility of the Floer operations,
    the diagonal components of the Maurer--Cartan
    equation for $b$ are the Maurer--Cartan equations
    for $b_{K}$ and $b_{K^*}$.
    The remaining component says that $\alpha$ is
    closed for the differential deformed by these
    diagonal bounding cochains.
    Undoing the shift on the first source part gives
    a bounding cochain $b_K$ on $(K,E_K)$.
    We therefore obtain branes
    \[
        \mathcal K=(K,E_K,b_K),
        \qquad
        \mathcal K^*=(K^*,E_{K^*},b_{K^*}),
    \]
    with bounding cochains of strictly positive valuation.
    Their components have the required phases by
    Lemma~\ref{lem:terminal-smoothing}, so
    $\mathcal K,\mathcal K^*\in\mathcal A$.  The closed cochain
    \[
        \alpha\in CF^1(\mathcal K^*,\mathcal K[1])
    \]
    gives the twisted-complex presentation
    \[
        \mathcal C
        \simeq
        \left(
        \mathcal K[1]\oplus\mathcal K^*,
        \begin{pmatrix}
            0 & \alpha \\
            0 & 0
        \end{pmatrix}
        \right),
    \]
    where the diagonal bounding cochains are included
    in the two branes.
    The degree-zero morphism
    \[
        r:\mathcal K^*\longrightarrow\mathcal K[2]
    \]
    corresponding to $\alpha$ therefore satisfies
    \[
        \mathcal C\simeq\operatorname{Cone}(r)[-1].
    \]
    Finally, Lemma~\ref{lem:web-cone-comparison}
    identifies the stationary brane with the initial
    cone object, giving
    \[
        \mathcal C
        \simeq\mathcal C_{\tau_0}
        \simeq\operatorname{Cone}(f).
    \]
    When $\tau_0=0$, the second identification instead
    comes directly from
    Lemma~\ref{lem:geometric-cone}.
    The exact triangle associated to this twisted
    complex therefore becomes
    \[
        \mathcal K[1]
        \longrightarrow
        \operatorname{Cone}(f)
        \longrightarrow
        \mathcal K^*
        \xrightarrow{\,r\,}
        \mathcal K[2],
    \]
    as required.
\end{proof}

The initial and terminal triangles determine morphisms
between the terminal branes and the original input
branes.
These will be used in the categorical applications
below.

\begin{definition}[Complementary morphisms]
    \label{def:complementary-morphisms}
    Write the initial cone triangle and the terminal
    triangle of Proposition~\ref{prop:terminal-web-triangle}
    as
    \[
        \begin{gathered}
            \mathcal L_0\xrightarrow{[f]}\mathcal L_1
            \xrightarrow{u}\operatorname{Cone}(f)
            \xrightarrow{v}\mathcal L_0[1],\\
            \mathcal K[1]\xrightarrow{a}\operatorname{Cone}(f)
            \xrightarrow{p}\mathcal K^*
            \longrightarrow\mathcal K[2].
        \end{gathered}
    \]
    The \emph{complementary morphisms} are the compositions
    in $\mathcal D$
    \[
        \begin{aligned}
            \kappa_0 & :=(v\circ a)[-1]:
            \mathcal K\longrightarrow\mathcal L_0, \\
            \kappa_1 & :=p\circ u:
            \mathcal L_1\longrightarrow\mathcal K^*.
        \end{aligned}
    \]
\end{definition}

\subsection{Dirichlet energy and finite active time}
\label{subsec:Dirichlet-energy}

Let $\mathcal W$ be a special Lagrangian web for $f$ on
$(0,T)$, with event set $E$ and boundary conditions
\[
    L_i(\tau)
    =[\iota_{i,\tau}:M_i(\tau)\longrightarrow X].
\]
The boundary conditions are fixed on each regular
interval. Apart from the auxiliary perturbations of
Remark~\ref{rem:auxiliary-perturbations},
the boundary conditions change only through absorption.
All integrals over immersed Lagrangians are taken over
their sources, with source multiplicities.
For the estimates below, assume that, for each
$i=0,1$ and
$\eta\in\{\operatorname{Re}\Omega,\operatorname{Im}\Omega\}$,
the swept integral
\[
    \int_{B_i(\tau)}\iota_{i,\tau}^*\eta
\]
has finite, matching incoming and outgoing limits
at each event, under the boundary-source
identifications of
Definition~\ref{def:surgery-compatible},
and tends to zero as $\tau\to0^+$.

\begin{definition}
    For an active or inactive component $Y=[F_Y:P_Y\to X]$, let
    $g_Y=F_Y^*g$ and $\rho_Y=\rho\circ F_Y$, so that
    $F_Y^*\operatorname{Im}\Omega
        =\rho_Y\,d\operatorname{vol}_{g_Y}$.
    If $Y$ is active, its \emph{weighted Dirichlet energy} is
    \[
        \mathcal E(Y)
        :=\int_{P_Y}\rho_Y|d\sigma_Y|_{g_Y}^2
        \,d\operatorname{vol}_{g_Y},
    \]
    where $\sigma_Y$ is the fundamental harmonic, with
    boundary values $0$ and $1$ on the corresponding
    labelled boundary components.
\end{definition}

\begin{lemma}[Swept real period and Dirichlet energy]
    \label{lem:web-dirichlet-energy}
    At every regular time, the two swept real integrals
    agree:
    \[
        \int_{B_0(\tau)}
        \iota_{0,\tau}^*\operatorname{Re}\Omega
        =
        \int_{B_1(\tau)}
        \iota_{1,\tau}^*\operatorname{Re}\Omega
        =:r(\tau).
    \]
    Their common value satisfies
    \begin{equation}
        r'(\tau)
        =\sum_{Y\subset Z_\tau}\mathcal E(Y),
        \label{eq:web-dirichlet-energy}
    \end{equation}
    where the sum is over the active components, and
    \begin{equation}
        0\leq r(\tau)\leq r_{\max},
        \qquad
        r_{\max}
        :=\min_{i=0,1}
        \int_{L_i^{\mathrm{in}}}\operatorname{Re}\Omega.
        \label{eq:web-R-bound}
    \end{equation}
    Here $L_i^{\mathrm{in}}=P_i$ are the initial positive
    boundary conditions.
\end{lemma}

\begin{proof}
    Denote the two swept integrals by $r_i(\tau)$.
    On a regular interval, choose parametrizations
    $F_{Y,\tau}$ and write $v_Y=\partial_\tau F_{Y,\tau}$.
    The evolution equation and the pointwise identity
    in \cite[Lemma~4.5(a)]{SY-cylinders} give
    \[
        F_{Y,\tau}^*(\iota_{v_Y}\omega)=d\sigma_Y,
        \qquad
        F_{Y,\tau}^*(\iota_{v_Y}\operatorname{Re}\Omega)
        =\rho_Y *_Yd\sigma_Y.
    \]
    The boundary orientation of $B_0(\tau)$ is opposite
    to that of $\partial_0P_Y$, while that of
    $B_1(\tau)$ agrees with that of $\partial_1P_Y$.
    The inactive fragments are stationary. The
    moving-boundary formula therefore gives
    \[
        r_0'(\tau)
        =-\sum_Y\int_{\partial_0P_Y}\rho_Y *_Yd\sigma_Y,
        \qquad
        r_1'(\tau)
        =\sum_Y\int_{\partial_1P_Y}\rho_Y *_Yd\sigma_Y.
    \]
    Writing $*_Y$ for the Hodge star of $g_Y$,
    weighted harmonicity gives
    $d(\rho_Y *_Yd\sigma_Y)=0$.
    Together with the boundary values
    $\sigma_Y|_{\partial_iP_Y}=i$,
    Stokes' theorem therefore gives
    \[
        -\int_{\partial_0P_Y}\rho_Y *_Yd\sigma_Y
        =\int_{\partial_1P_Y}\rho_Y *_Yd\sigma_Y
        =\int_{P_Y}\rho_Y\,d\sigma_Y\wedge *_Yd\sigma_Y
        =\mathcal E(Y).
    \]
    Thus $r_0'=r_1'=\sum_Y\mathcal E(Y)$ on each regular
    interval. Since both integrals start at zero and
    have matching limits at events, they agree
    throughout the evolution.

    Positivity gives
    $0\leq r_i(\tau)\leq
        \int_{L_i(\tau)}\operatorname{Re}\Omega$.
    Auxiliary perturbations preserve the boundary real
    periods, as does absorption: for the
    absorbed closed fragment $S$,
    \[
        \int_{L_0^+}\operatorname{Re}\Omega
        =\int_{L_0^-}\operatorname{Re}\Omega
        -\int_S\operatorname{Re}\Omega
        =\int_{L_0^-}\operatorname{Re}\Omega.
    \]
    This proves~\eqref{eq:web-R-bound}.
\end{proof}

\begin{corollary}[Finite active time]
    \label{cor:finite-active-time}
    Suppose that every active component occurring in
    the evolution satisfies
    \begin{equation}
        \mathcal E(Y)\geq\epsilon
        \label{eq:uniform-dirichlet-gap}
    \end{equation}
    for some $\epsilon>0$.
    Then the total active time satisfies
    \[
        T_{\mathrm{act}}
        :=\int_{\{\tau\in(0,T)\setminus E:
            Z_\tau\neq\varnothing\}}d\tau
        \leq\frac{r_{\max}}{\epsilon}.
    \]
\end{corollary}

\begin{proof}
    Whenever the active front is nonempty,
    \eqref{eq:web-dirichlet-energy} gives
    $r'(\tau)\geq\epsilon$.
    Integrating over the regular active intervals
    and using continuity and monotonicity of $r$ gives
    $\epsilon T_{\mathrm{act}}\leq r_{\max}$.
\end{proof}

To deduce finite-time termination, one must also
control stationary waiting periods and the number
of events.

We also record a bound on the total special volume.
Let $S_i$ denote the initially supplied closed
imaginary-special components from
Definition~\ref{def:initial-web-configuration}.
Let $Z_\tau^{\mathrm{tot}}$ be the disjoint-source
union of the active front and all retained
inactive fragments, each with its
imaginary-special orientation.
Write
\[
    V(\tau)
    :=\int_{Z_\tau^{\mathrm{tot}}}\operatorname{Im}\Omega,
    \qquad
    V_{\mathrm{in}}
    :=\int_{S_0}\operatorname{Im}\Omega
    +\int_{S_1}\operatorname{Im}\Omega,
\]
and set
\[
    a_i(\tau)
    :=\int_{B_i(\tau)}
    \iota_{i,\tau}^*\operatorname{Im}\Omega.
\]
Since all components of $Z_\tau^{\mathrm{tot}}$
are imaginary special,
\[
    V(\tau)
    =
    \sum_{Y\subset Z_\tau^{\mathrm{tot}}}
    \int_{P_Y}\rho_Y\,d\operatorname{vol}_{g_Y}.
\]
Thus $V(\tau)$ is their total weighted volume,
counted with source multiplicity.

\begin{lemma}[Total special volume]
    \label{lem:web-volume-bound}
    On every regular interval,
    $V+a_0-a_1$ is constant.
    If
    \[
        V(0^+)\leq V_{\mathrm{in}},
        \qquad
        V(e^+)\leq V(e^-)
        \quad\text{for every }e\in E,
    \]
    then at every regular time
    \begin{equation}
        \begin{aligned}
            V(\tau)
             & \leq V_{\mathrm{in}}-a_0(\tau)+a_1(\tau) \\
             & \leq V_{\mathrm{in}}
            +\sum_{i=0,1}\int_{M_i(\tau)}
            \left|\iota_{i,\tau}^*
            \operatorname{Im}\Omega\right|.
        \end{aligned}
        \label{eq:web-volume-bound}
    \end{equation}
\end{lemma}

\begin{proof}
    Since $d\operatorname{Im}\Omega=0$ and the inactive
    fragments are stationary, Cartan's formula and
    Stokes' theorem give, on a regular interval,
    \[
        V'(\tau)
        =\sum_Y\int_{\partial P_Y}
        F_{Y,\tau}^*
        (\iota_{v_Y}\operatorname{Im}\Omega)
        =a_1'(\tau)-a_0'(\tau).
    \]
    Thus $V+a_0-a_1$ is constant there. In particular,
    $V$ has finite one-sided limits wherever the swept
    integrals do. Their continuity and the assumed
    inequalities for $V$ imply
    $V+a_0-a_1\leq V_{\mathrm{in}}$ throughout the
    evolution. The second inequality follows by
    bounding each trace integral by the corresponding
    absolute integral over the full boundary source.
\end{proof}

The initial and event volume inequalities in
Lemma~\ref{lem:web-volume-bound} must be verified
for the chosen continuations.
When the other pieces have matching limiting
special volumes, absorption of a closed fragment
$S$ decreases $V$ by its special volume $\int_S\operatorname{Im}\Omega.$
At events involving closed-component reductions,
the effect of those replacements on $V$ must
also be included when checking the event inequality.

If the chosen evolution contains only finitely
many absorptions, its boundary conditions take
only finitely many values.
Consequently, \eqref{eq:web-volume-bound} gives
a bound uniform in $\tau$ along that evolution.

\section{Local models for singular events}
\label{sec:singular-events}

We describe the local operations used to continue
web evolutions through geometric degenerations
and obstruction walls.
The geometric models specify the source modifications
and the required behavior of the immersions and
boundary lifts.
Their compatibility with the full web data and
the brane comparisons across events are specified
in Definitions~\ref{def:surgery-compatible}
and~\ref{def:brane-surgery-transfer}.

Subsection~\ref{subsec:local-continuation-conventions}
introduces the common source and convergence
conventions and formulates the local geometric
continuation conjecture.
Subsection~\ref{subsec:geometric-degenerations}
describes geometric degenerations and their
continuation by surgery.
Subsection~\ref{subsec:obstruction-walls}
describes obstruction walls and their resolutions
by surgery or reduction of closed inactive components.

\subsection{Conventions for local continuation}
\label{subsec:local-continuation-conventions}

We first specify the notation and conventions shared
by the local constructions below.
Fix an endpoint $\tau_*$ of an incoming regular interval.
Write
\[
    Y_\tau^-
    =[F_\tau^-:P^-\longrightarrow X,
    \beta_{0,\tau}^-,\beta_{1,\tau}^-],
    \qquad \tau<\tau_*,
\]
for the incoming fronts and, where relevant,
inactive fragments.
The outgoing family of fronts and inactive fragments is denoted by $Y_\tau^+$,
with source $P^+$, for $\tau>\tau_*$.
The Lagrangian boundary conditions
\[
    L_j=[\iota_j:M_j\longrightarrow X],
    \qquad j=0,1,
\]
remain fixed.
On either side, the immersions are imaginary special
and the boundary lifts are embeddings satisfying
\[
    F_\tau^\pm|_{\partial_jP^\pm}
    =
    \iota_j\circ\beta_{j,\tau}^\pm.
\]
Each local model specifies a surgery on the source,
replacing a closed region $A^-\subset P^-$
by a region $A^+\subset P^+$ with prescribed
attaching maps.
Set
\[
    U^\pm=P^\pm\setminus A^\pm,
\]
and denote the resulting orientation-preserving
diffeomorphism by
\[
    \phi:U^-\xrightarrow{\sim}U^+.
\]
This map preserves the boundary labels.
New boundary portions receive the labels specified
by the individual construction.
Write
\[
    \jmath^-:\widehat A^-\hookrightarrow P^-,
    \qquad A^-=\jmath^-(\widehat A^-),
\]
for the incoming surgery embedding.
We choose source parameterizations for which
$\jmath^-$ is independent of $\tau$.
Thus the degeneration is described by the
time-dependent maps $F_\tau^-$ and
$\beta_{j,\tau}^-$ on a fixed source with a fixed
surgery region.

Outgoing components carrying both boundary labels
are active; the remaining components are inactive.
A nonclosed inactive component carries the label
of its boundary.
The label of a closed inactive component is
part of the continuation data.

\begin{definition}[Common limiting source]
    \label{def:source-replacement-limit}
    A \emph{common limiting source} for the source
    surgery $P^-\rightsquigarrow P^+$ is a compact
    Hausdorff space $P_*$ with a distinguished point
    $p_*$ and continuous surjections
    \[
        \pi^\pm:P^\pm\longrightarrow P_*,
    \]
    called the \emph{collapse maps}.
    Their collapsing sets
    \[
        C^\pm:=(\pi^\pm)^{-1}(p_*)
    \]
    must satisfy
    \begin{equation}
    \label{eq:collapse-compatability}
                C^\pm\subset A^\pm,
        \qquad
        \pi^+\circ\phi=\pi^-|_{U^-}.
    \end{equation}
    We require $P_*\setminus\{p_*\}$ to be a smooth
    oriented manifold with boundary and the restrictions
    \[
        \pi^\pm:
        P^\pm\setminus C^\pm
        \longrightarrow P_*\setminus\{p_*\}
    \]
    to be orientation-preserving diffeomorphisms.
    The labelled boundary images must agree:
    \[
        \Gamma_{j,*}
        :=
        \pi^-(\partial_jP^-\setminus C^-)
        =
        \pi^+(\partial_jP^+\setminus C^+),
        \qquad j=0,1.
    \]
    In the local models below, we take
    $P_*=P^-/C^-$, with $\pi^-$ the quotient map and
    the smooth structure induced by
    $P^-\setminus C^-$ away from $p_*$.
    The compatibility equation~\eqref{eq:collapse-compatability} determines $\pi^+$
    on $U^+$; its restriction to $A^+$ is specified
    in each model.
\end{definition}

\begin{definition}[Limiting front]
    \label{def:limiting-front}
    A \emph{limiting front} $Y_*$ on a common limiting
    source consists of a continuous map
    \[
        F_*:P_*\longrightarrow X,
        \qquad q:=F_*(p_*),
    \]
    which is an imaginary-special immersion away
    from $p_*$, together with boundary lift embeddings
    \[
        \beta_{j,*}:\Gamma_{j,*}\longrightarrow M_j,
        \qquad
        F_*|_{\Gamma_{j,*}}
        =\iota_j\circ\beta_{j,*},
        \qquad j=0,1.
    \]
\end{definition}

\begin{definition}[Convergence through surgery]
    \label{def:convergence-through-surgery}
    Let $Y_*$ be a limiting front.
    For either choice of sign, the family $Y_\tau^\pm$
    \emph{converges to $Y_*$ through the surgery}
    if, as $\tau\to\tau_*^\pm$,
    \[
        \sup_{x\in C^\pm}
        d_g\bigl(F_\tau^\pm(x),q\bigr)
        \longrightarrow 0,
    \]
    where $d_g$ is the distance induced by the
    K\"ahler metric $g$, and
    \[
        \begin{aligned}
            F_\tau^\pm
            &\longrightarrow F_*\circ\pi^\pm
            &&\text{in }
            C^\infty_{\mathrm{loc}}
            (P^\pm\setminus C^\pm),\\
            \beta_{j,\tau}^\pm
            &\longrightarrow\beta_{j,*}\circ\pi^\pm
            &&\text{in }
            C^\infty_{\mathrm{loc}}
            (\partial_jP^\pm\setminus C^\pm).
        \end{aligned}
    \]
\end{definition}

Convergence of the incoming and outgoing families
to the same limiting front gives the matching
equations for the immersions and boundary lifts
on $U^\pm$ in
Definition~\ref{def:surgery-compatible},
with $\psi_j$ the identity on the fixed boundary
sources.
The accumulated boundary traces must additionally
satisfy \eqref{eq:continuation-trace-matching}
when the local construction is incorporated
into the full web.

\begin{definition}[Inactive continuation]
    \label{def:inactive-continuation}
    Let $Y_*$ be a limiting front.
    An \emph{inactive continuation} on an inactive
    component $Q\subset P^+$ is a stationary family
    $Y_{\tau,Q}^+=Y_Q$ of imaginary-special fragments,
    with immersion and boundary lifts specified
    by the chosen local model.
    In the cases where $C^+\cap Q$ is a single point, we denote this point by $c_Q$ and require a homeomorphism
    \[
        \chi_Q:Q\longrightarrow\pi^+(Q),
        \qquad \chi_Q(c_Q)=p_*,
    \]
    agreeing with $\pi^+|_Q$ outside $Q \cap A^+$
    and restricting to an orientation-preserving
    diffeomorphism preserving boundary labels
    away from $c_Q$, such that
    \[
        F_Q^{\mathrm{cap}}:=F_*\circ\chi_Q
    \]
    is a smooth imaginary-special immersion.
    The maps $\beta_{j,*}\circ\chi_Q$ on
    $\partial_jQ\setminus\{c_Q\}$ must extend
    to embedded boundary lifts
    $\beta_{j,Q}^{\mathrm{cap}}$.
    If $c_Q\in\partial_jQ$, we require
$\beta_{j,Q}^{\mathrm{cap}}(c_Q)$ to be the
preimage of $F_Q^{\mathrm{cap}}(c_Q)$ under $\iota_j$.
    These maps represent the stationary fragment
    $Y_Q=Y_Q^{\mathrm{cap}}$.
    We use $\chi_Q$ as the outgoing comparison map
    in Definition~\ref{def:convergence-through-surgery}.
\end{definition}

\begin{definition}[Equivalence of outgoing continuations]
    \label{def:outgoing-continuation-equivalence}
    Fix an incoming family, a source surgery,
    and a limiting front.
    Two outgoing continuations $Y_\tau^+$ and
    $\widetilde Y_\tau^+$ on $P^+$ are
    \emph{equivalent} if, after shortening their
    common interval $(\tau_*,\tau_*+\epsilon)$,
    there is a smooth family of orientation-preserving
    diffeomorphisms
    \[
        \eta_\tau:P^+\longrightarrow P^+
    \]
    preserving the boundary labels and the labels
    of closed fragments, such that
    \[
        \widetilde F_\tau^+
        =
        F_\tau^+\circ\eta_\tau,
        \qquad
        \widetilde\beta_{j,\tau}^+
        =
        \beta_{j,\tau}^+
        \circ\eta_\tau|_{\partial_jP^+}.
    \]
    The evolution parameter $\tau$ is unchanged.
\end{definition}

\begin{conjecture}[Local geometric continuation]
    \label{conj:local-geometric-continuation}
    Fix a local model from
    Sections~\ref{subsec:geometric-degenerations},~\ref{subsubsec:puncture}
    or~\ref{subsubsec:interior-zero-surgery}.
    Suppose that the hypotheses of the model hold
    and the incoming family converges through the surgery
    to a limiting front $Y_*$, with
    \[
        F_*^{-1}(q)=\{p_*\}.
    \]
    Then there exist $\epsilon>0$ and an outgoing
    family $Y_\tau^+$ on $P^+$, defined for
    $\tau_*<\tau<\tau_*+\epsilon$,
    satisfying the geometric requirements of the
    model and the following conditions.
    \begin{enumerate}
        \item \emph{Active continuation.}
              On each active component $Q$,
              \[
                  \frac{d}{d\tau}Y_{\tau,Q}^+
                  =
                  d\sigma_{Y_{\tau,Q}^+},
              \]
              and as $\tau \to \tau_*$ the family converges through surgery
              to $Y_*$ in the sense of Definition~\ref{def:convergence-through-surgery}.

        \item \emph{Inactive continuation.}
              On each inactive component, the family
              is an inactive continuation in the sense
              of Definition~\ref{def:inactive-continuation}.
    \end{enumerate}
    This continuation is unique up to the equivalence
    of Definition~\ref{def:outgoing-continuation-equivalence}.
\end{conjecture}

Establishing the conjecture in general requires
gluing arguments beyond
the scope of this paper, which we plan to address in future work.
The limiting-source framework above is used for
the geometric degenerations, puncture, and
transverse interior $0$-surgery models.
The remaining obstruction-wall resolutions
are described separately in
Section~\ref{subsec:obstruction-walls}.
Compatibility with the full web data and brane
transfer are specified in
Definitions~\ref{def:surgery-compatible}
and~\ref{def:brane-surgery-transfer}.

\subsection{Geometric degenerations}
\label{subsec:geometric-degenerations}

We describe local surgeries on the fronts and their
effects on the assembled web cone.
The models below specify the source modifications,
the hypotheses on the incoming degeneration, and
the required properties of the outgoing geometry.
Conjecture~\ref{conj:local-geometric-continuation}
asserts the existence of these geometric
continuations, unique up to the equivalence of
Definition~\ref{def:outgoing-continuation-equivalence}.

\subsubsection{Boundary degenerations}
We first describe a boundary degeneration in which
a sphere in the boundary of the front source
contracts to a point on a chosen sheet of a
Lagrangian boundary condition.
We also consider the reverse source surgery,
obtained by interchanging the removed and
inserted pieces; see
Remark~\ref{rem:inverse-boundary-surgery}.
Write $D^r$ for the closed $r$-disk, with
$\partial D^0=\emptyset$.
For $r\geq1$, let
\[
    H_D^r=D^r\cap\{x_1\geq0\},
    \qquad
    H_S^{r-1}=\partial D^r\cap\{x_1\geq0\},
    \qquad
    H_F^{r-1}=D^r\cap\{x_1=0\}
\]
denote the half-disk and its hemispherical and flat boundary components.

\begin{definition}[Boundary degeneration and surgery]
    \label{def:boundary-degeneration}
    \label{def:boundary-surgery}
    Fix $i \in \{0,1\}$, choose $q_i\in M_i$ and a neighborhood
    $V_i\subset M_i$ on which $\iota_i$ is an
    embedding, and put $q=\iota_i(q_i)$.

    A \emph{boundary degeneration of index $k$}
    along this sheet, where $0\leq k\leq n-1$,
    is specified by an embedding
    \[
        \jmath^-:
        S^k\times H_D^{n-k}\longrightarrow P^-
    \]
    such that $ \jmath^-(S^k \times H_F^{r-1}) \subset \partial_iP^-$
    and whose other points lie in
    $\operatorname{int}(P^-)$.
    The incoming surgery region and its
    \emph{vanishing sphere} are
    \[
        A^-:=\jmath^-(S^k\times H_D^{n-k}),
        \qquad
        C^-:=\jmath^-(S^k\times\{0\}).
    \]
    We require
    \[
        \beta_{i,\tau}^-|_{C^-}
        \longrightarrow q_i
    \]
    uniformly as $\tau\to\tau_*^-$.
    The associated \emph{boundary surgery} produces
    $P^+$ by replacing
    \[
        S^k\times H_D^{n-k}
        \quad\rightsquigarrow\quad
        D^{k+1}\times H_S^{n-k-1},
    \]
    using $\jmath^-$ to identify the common attaching
    region $S^k\times H_S^{n-k-1}$.
        Denote the embedding of the inserted piece into
    $P^+$ by $\jmath^+$, and set
    \[
        \begin{aligned}
            A^+
            &:=\jmath^+
            \bigl(D^{k+1}\times H_S^{n-k-1}\bigr),\\
            C^+
            &:=\jmath^+
            \bigl(\{0\}\times H_S^{n-k-1}\bigr),\\
            \Gamma^+
            &:=\jmath^+
            \bigl(D^{k+1}\times\partial H_S^{n-k-1}\bigr).
        \end{aligned}
    \]
    The new boundary portion $\Gamma^+$ receives
    label $i$, and its outgoing boundary lift
    must satisfy
    \[
        \beta_{i,\tau}^+(\Gamma^+)\subset V_i.
    \]
\end{definition}

\begin{remark}
    For $k=n-1$, boundary surgery caps a boundary sphere:
    \[
        S^{n-1}\times H_D^1
        \quad\rightsquigarrow\quad
        D^n.
    \]
    In this case $\Gamma^+=\emptyset$ and $C^+$ consists
    of a single point.
\end{remark}

\begin{figure}[htbp]
    \centering
    \includegraphics[width=\textwidth]{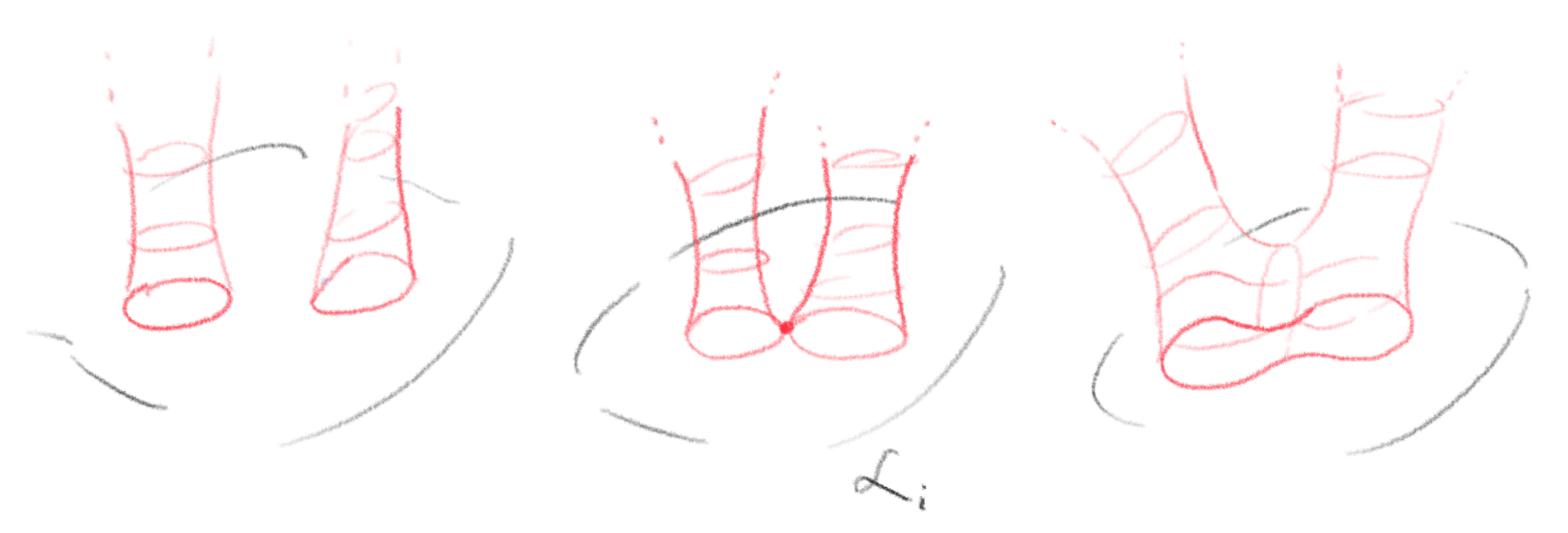}
    \caption{
        Boundary degeneration and surgery of index zero,
        shown before, at, and after the event.
    }
    \label{fig:boundary-surgery}
\end{figure}

\begin{definition}[Collapse maps for boundary surgery]
    \label{def:boundary-degeneration-limit}
    On the inserted region $A^+$, define $\pi^+$ by
    \[
        \pi^+\bigl(\jmath^+(\rho u,v)\bigr)
        :=
        \pi^-\bigl(\jmath^-(u,\rho v)\bigr),
        \qquad
        u\in S^k,\quad
        v\in H_S^{n-k-1},\quad
        0<\rho\leq1,
    \]
    together with $\pi^+(C^+)=p_*$.

    We choose gluing collars extending the coordinates
    $(u,v,\rho)$ across the attaching locus $\rho=1$,
    with $\rho>1$ parametrizing the unchanged part
    of the source.
    In these coordinates, the displayed formula joins
    smoothly to $\pi^-|_{U^-}\circ\phi^{-1}$ on $U^+$.
    These maps define a common limiting source in
    the sense of
    Definition~\ref{def:source-replacement-limit}.
\end{definition}

\begin{remark}[Inverse boundary surgery]
    \label{rem:inverse-boundary-surgery}
    As in Remark~\ref{rem:maslov-n-webbing}, exchanging
    the boundary labels reverses the prescribed evolution.
    This symmetry motivates allowing the source
    replacement of
    Definition~\ref{def:boundary-surgery}
    in both directions.
    We call the reverse operation an
    \emph{inverse boundary surgery}.
    In particular, the inverse of the capping
    operation $k=n-1$ is puncture, described
    separately in
    Definition~\ref{def:web-front-puncture}.
\end{remark}

For $k=n-1$, an inactive component $Q$ containing
$C^+$ has a single cap center in
$\operatorname{int}(Q)$.
For inverse boundary surgery with $k=0$, each
inactive component $Q$ produced by splitting
has a single cap center in $\partial_iQ$.
These are the capping cases of
Definition~\ref{def:inactive-continuation}.

\begin{remark}[Handle transfer at a boundary degeneration]
    \label{rem:boundary-degeneration-handle-transfer}
    For the boundary surgery in the direction displayed
    in Definition~\ref{def:boundary-surgery}, the expected
    local model transfers a handle
    \[
        H=D^{k+1}\times D^{n-k-1}
    \]
    from the residue to the front.
    The swept region expands through this handle,
    removing it from the residue, while boundary surgery
    adds the corresponding handle to the front.
    Figure~\ref{fig:boundary-surgery} illustrates the
    index-zero case.

    Since the web cone is assembled by gluing the front
    to the residues along the chosen boundary lifts,
    the handle is supplied by the residue before the
    event and by the front afterward.
    In this model, the assembled source therefore
    retains its diffeomorphism type after rounding
    corners; only its decomposition into front and
    residue changes.
\end{remark}

\subsubsection{Interior degenerations}

Since fronts are allowed to be immersed, an interior
self-intersection does not by itself require surgery.
A geometric interior degeneration concerns a failure
of smooth immersed continuation on the original source.
We distinguish such degenerations from the interior
$0$-surgeries used to resolve obstruction walls.

\begin{definition}[Interior surgery and degeneration]
    \label{def:interior-degeneration}
    \label{def:interior-surgery}
    An \emph{interior surgery of index $k$}, where
    $0\leq k\leq n-1$, is specified by an embedding
    \[
        \jmath^-:
        S^k\times D^{n-k}
        \longrightarrow\operatorname{int}(P^-).
    \]
    Set
    \[
        A^-:=\jmath^-(S^k\times D^{n-k}),
        \qquad
        C^-:=\jmath^-(S^k\times\{0\}).
    \]
    The surgery produces $P^+$ by replacing
    \[
        S^k\times D^{n-k}
        \quad\rightsquigarrow\quad
        D^{k+1}\times S^{n-k-1},
    \]
    using $\jmath^-$ to identify the common attaching
    region $S^k\times S^{n-k-1}$.
    Denote the embedding of the inserted piece into
    $P^+$ by $\jmath^+$, and set
    \[
        \begin{aligned}
            A^+
            &:=\jmath^+
            \bigl(D^{k+1}\times S^{n-k-1}\bigr),\\
            C^+
            &:=\jmath^+
            \bigl(\{0\}\times S^{n-k-1}\bigr).
        \end{aligned}
    \]
    For $1\leq k\leq n-1$, the incoming family has
    an \emph{interior degeneration of index $k$}
    at $q\in X$ along $\jmath^-$ if
    \[
        F_\tau^-|_{C^-}\longrightarrow q
    \]
    uniformly as $\tau\to\tau_*^-$,
    and it has no smooth immersed limit on $P^-$,
    even after source reparameterization.
    We call $C^-$ its \emph{vanishing sphere}.
\end{definition}

Here ``interior'' refers to the source.
The limit point $q$ may lie on the ambient image
of a boundary condition.
The front is attached to the residues along its
source boundary, whereas the surgery region lies
in its interior.
Consequently, the same source surgery is performed
on the assembled web cone.

\begin{figure}[htbp]
    \centering
    \includegraphics[scale=0.12]
    {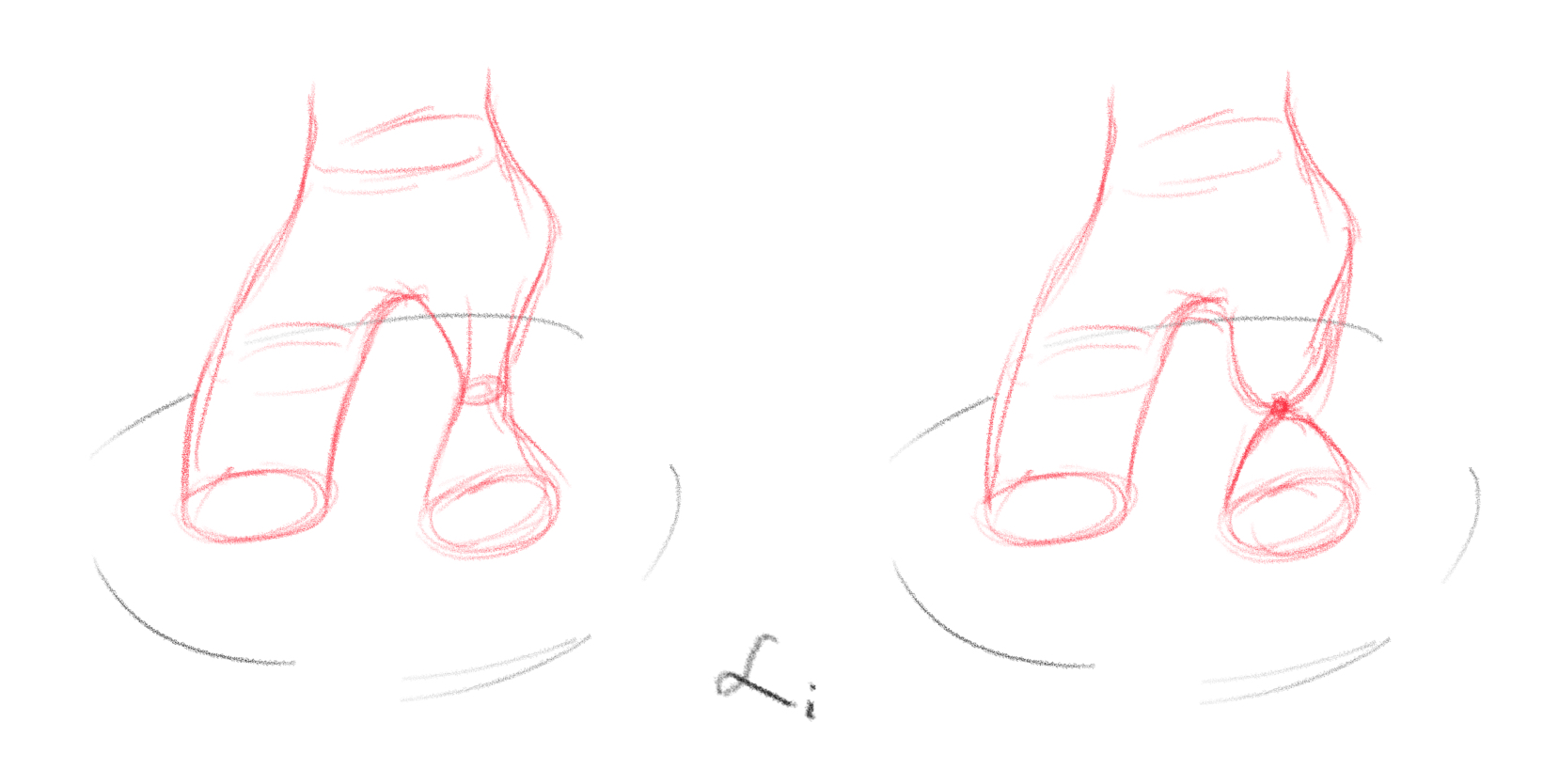}
    \caption{
        Interior degeneration of index $k=1$
        for a two-dimensional front, shown in red.
        Left: a vanishing circle lies in the
        interior of the right leg.
        Right: the circle has collapsed to a point.
        The corresponding interior anti-surgery
        cuts along this circle and caps the two
        resulting ends.
        The detached leg becomes an inactive
        fragment whose remaining boundary carries
        label $i$ and lies on the displayed sheet
        of $L_i$.
    }
    \label{fig:interior-degeneration}
\end{figure}

For $k=0$, the surgery replaces
\[
    D^n\sqcup D^n
    \quad\rightsquigarrow\quad
    D^1\times S^{n-1}.
\]
This is the source surgery used at selected
transverse intersections to resolve obstruction
walls in Subsection~\ref{subsec:obstruction-walls}.
For $k=n-1\geq1$, the surgery replaces
\[
    S^{n-1}\times D^1
    \quad\rightsquigarrow\quad
    D^n\sqcup D^n.
\]
We call this \emph{interior anti-surgery}.
It cuts the source along the vanishing sphere
and caps the two resulting ends; see Figure~\ref{fig:interior-degeneration}.

\begin{definition}[Collapse maps for interior surgery]
    \label{def:interior-degeneration-limit}
    On the inserted region $A^+$, define $\pi^+$ by
    \[
        \pi^+\bigl(\jmath^+(\rho u,v)\bigr)
        :=
        \pi^-\bigl(\jmath^-(u,\rho v)\bigr),
        \qquad
        u\in S^k,\quad
        v\in S^{n-k-1},\quad
        0<\rho\leq1,
    \]
    together with $\pi^+(C^+)=p_*$.
    With gluing collars chosen as in
    Definition~\ref{def:boundary-degeneration-limit},
    these maps define a common limiting source in
    the sense of
    Definition~\ref{def:source-replacement-limit}.
\end{definition}

\begin{remark}
    The inverse of an interior surgery of index $k$
    is an interior surgery of index $n-k-1$,
    after exchanging the product factors.
    Thus, for $1\leq k\leq n-2$, both directions
    are included among the interior degeneration
    models, subject to the hypotheses on the
    incoming evolution.
    The inverse of interior $(n-1)$-surgery is
    the interior $0$-surgery used to resolve
    obstruction walls in
    Section~\ref{subsubsec:interior-zero-surgery}.
\end{remark}

An inactive outgoing component $Q$ can occur only
when $k=n-1$ and the surgery splits $P^-$.
It then has a single cap center in
$\operatorname{int}(Q)$, giving a capping case of
Definition~\ref{def:inactive-continuation}.

If the vanishing sphere separates the $0$-labelled
boundaries from the $1$-labelled boundaries, the
two outgoing components are inactive fragments
with opposite labels, as illustrated at the bottom
of Figure~\ref{fig:milnor-web}.
A closed inactive fragment may also occur alongside
an active component.

\begin{remark}
    \label{rem:interior-higher-index-models}
    Haug's Lagrangian antisurgery
    \cite[Definition~2.1 and Theorem~1.1]
    {Haug-antisurgery}
    realizes the underlying source surgery
    when an isotropic surgery disk with
    compatible framing is available.
    The resulting Lagrangian immersion has
    a new transverse double point.
    The imaginary-special equation is not
    part of this construction.
\end{remark}

\subsubsection{Degenerations at intersections of the boundary conditions}

The following relative surgery was proposed in
\cite[\S1.3.1]{SY-webbing}.
For immersed boundary conditions, the two participating
sheets are part of the local data.

\begin{definition}[Intersection degeneration and surgery]
    \label{def:intersection-degeneration}
    \label{def:intersection-surgery}
    Choose $q_i\in M_i$, $i=0,1$, with
    \[
        \iota_0(q_0)=\iota_1(q_1)=q,
    \]
    and neighborhoods $U_i$ of $q_i$ on which
    $\iota_i$ is an embedding.
    Suppose that the selected sheets intersect
    transversely at $q$ with Maslov index
    \[
        m=m(q_0,q_1;L_0,L_1),
        \qquad
        1\leq m\leq n-1.
    \]
    An \emph{intersection degeneration} at this pair
    is specified by an embedding
    \[
        \jmath^-:
        S^{m-1}\times[0,1]\times D^{n-m}
        \longrightarrow P^-
    \]
    mapping $S^{m-1}\times\{i\}\times D^{n-m}$
    into $\partial_iP^-$ for $i=0,1$, and all
    other points into $\operatorname{int}(P^-)$.
    The incoming surgery region and its
    \emph{vanishing cylinder} are
    \[
        \begin{aligned}
            A^-
            &:=\jmath^-
            \bigl(S^{m-1}\times[0,1]\times D^{n-m}\bigr),\\
            C^-
            &:=\jmath^-
            \bigl(S^{m-1}\times[0,1]\times\{0\}\bigr).
        \end{aligned}
    \]
    We require
    \[
        F_\tau^-|_{C^-}\longrightarrow q,
        \qquad
        \beta_{i,\tau}^-|_{C^-\cap\partial_iP^-}
        \longrightarrow q_i,
        \qquad i=0,1,
    \]
    uniformly as $\tau\to\tau_*^-$.
    The associated \emph{intersection surgery}
    produces $P^+$ by replacing
    \[
        S^{m-1}\times[0,1]\times D^{n-m}
        \quad\rightsquigarrow\quad
        D^m\times[0,1]\times S^{n-m-1},
    \]
    using $\jmath^-$ to identify the common attaching
    region $S^{m-1}\times[0,1]\times S^{n-m-1}$.
    Denote the embedding of the inserted piece into
    $P^+$ by $\jmath^+$, and set
    \[
        \begin{aligned}
            A^+
            &:=\jmath^+
            \bigl(D^m\times[0,1]\times S^{n-m-1}\bigr),\\
            C^+
            &:=\jmath^+
            \bigl(\{0\}\times[0,1]\times S^{n-m-1}\bigr),\\
            \Gamma_i^+
            &:=\jmath^+
            \bigl(D^m\times\{i\}\times S^{n-m-1}\bigr),
            \qquad i=0,1.
        \end{aligned}
    \]
    The new boundary portion $\Gamma_i^+$ receives
    label $i$, and its outgoing boundary lift
    must satisfy
    \[
        \beta_{i,\tau}^+(\Gamma_i^+)\subset U_i,
        \qquad i=0,1.
    \]
\end{definition}

\begin{figure}[htbp]
    \centering
    \includegraphics[width=\textwidth]{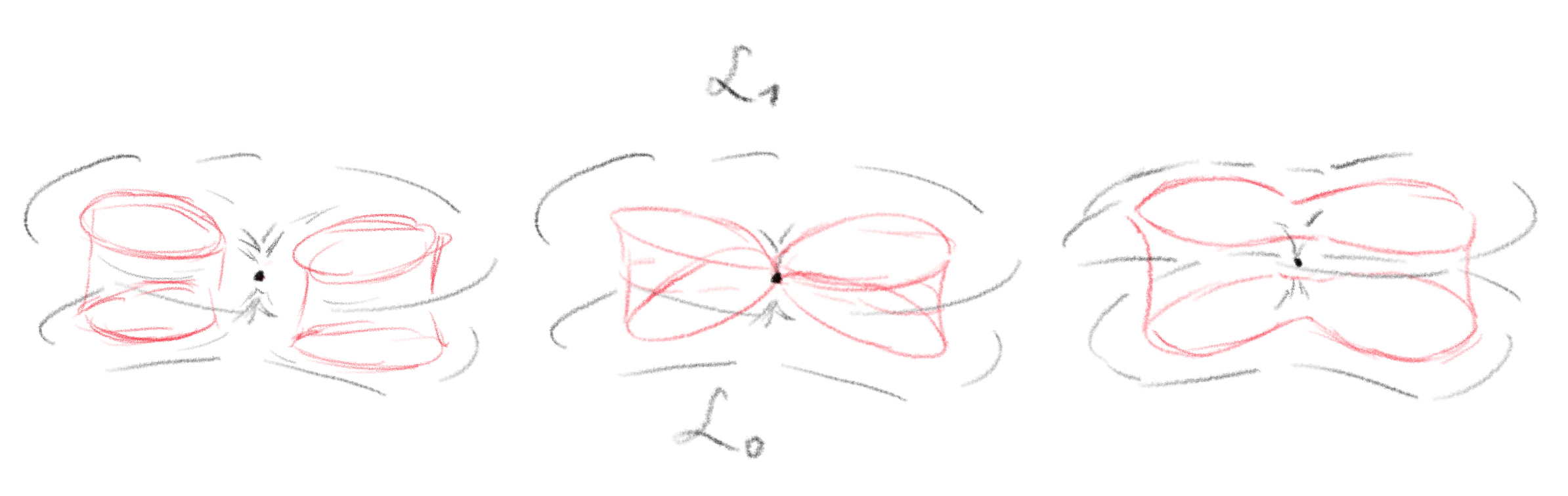}
    \caption{
        Degeneration of the front at an intersection
        of the boundary conditions.
        From left to right, the picture illustrates
        a Maslov-index-one degeneration and its surgery.
        With the ordering of the boundary conditions
        exchanged, the reverse sequence illustrates
        a Maslov-index-$(n-1)$ degeneration.
        The fronts are shown in red.
    }
    \label{fig:intersection-surgery}
\end{figure}

\begin{definition}[Collapse maps for intersection surgery]
    \label{def:intersection-degeneration-limit}
    On the inserted region $A^+$, define $\pi^+$ by
    \[
        \pi^+\bigl(\jmath^+(\rho u,t,v)\bigr)
        :=
        \pi^-\bigl(\jmath^-(u,t,\rho v)\bigr),
    \]
    where
    \[
        u\in S^{m-1},\qquad
        v\in S^{n-m-1},\qquad
        t\in[0,1],\qquad
        0<\rho\leq1,
    \]
    together with $\pi^+(C^+)=p_*$.
    With gluing collars chosen as in
    Definition~\ref{def:boundary-degeneration-limit},
    these maps define a common limiting source in
    the sense of
    Definition~\ref{def:source-replacement-limit}.
\end{definition}

\begin{remark}
    Exchanging the ordered boundary conditions gives
    \[
        m(q_1,q_0;L_1,L_0)=n-m.
    \]
    After replacing $t$ by $1-t$ and exchanging
    the product factors, the intersection model
    with Maslov index $n-m$ reverses the source
    replacement above.
    Thus both directions are included among the
    models with $1\leq m\leq n-1$.
\end{remark}

Every participating outgoing component contains
a component of $A^+$ and therefore meets both
$\Gamma_0^+$ and $\Gamma_1^+$.
Thus all such components are active, although
$P^+$ need not be connected.

\begin{remark}[Web-cone surgery at a Lagrangian intersection]
    \label{rem:intersection-degeneration-web-cone-surgery}
    The residue pieces
    adjoining the two ends of the vanishing cylinder
    $S^{m-1}\times[0,1]$ cap it to form a vanishing
    $m$-sphere in the web-cone source.
    Its neighborhood is diffeomorphic to
    $S^m\times D^{n-m}$.
    Across the event, these residue pieces become swept,
    while the outgoing front supplies
    $D^m\times[0,1]\times S^{n-m-1}$.
    After rounding corners, the assembled source
    undergoes the interior source surgery of
    index $m$
    \[
        S^m\times D^{n-m}
        \quad\rightsquigarrow\quad
        D^{m+1}\times S^{n-m-1}.
    \]
    Figure~\ref{fig:intersection-surgery} illustrates
    the front surgery for $m=1$.
\end{remark}

The endpoint indices $m=0$ and $m=n$ are
treated separately.
For $m=n$, we use the terminating webbing of
Remark~\ref{rem:maslov-n-webbing}, which
contracts to $q$.
At the collapse time, its components are
removed from the active fronts.
For $m=0$, Theorem~\ref{thm:SY-start-web},
applied to the selected embedded sheets,
supplies a new webbing.
Subsection~\ref{subsec:obstruction-walls}
describes this operation and its use at
initialization and in resolving obstruction
walls.

\subsection{Obstruction walls and their resolutions}
\label{subsec:obstruction-walls}

Even when the underlying web evolves regularly,
the transported bounding cochain may lose
strict positivity of valuation.
For transverse self-intersection generators,
the valuation evolves according to
\eqref{eq:web-cone-valuation-transport}.

\begin{definition}[Obstruction wall]
    \label{def:obstruction-wall}
    Suppose that
    \[
        \mathcal C_\tau=(C_\tau,E_\tau,b_\tau),
        \qquad a<\tau<\tau_*,
    \]
    is an incoming family of smoothed web-cone
    branes transported along the prescribed
    height deformation, with
    $\operatorname{val}(b_\tau)>0$.
    We call $\tau_*$ an \emph{obstruction wall} if
    \[
        \lim_{\tau\to\tau_*^-}
        \operatorname{val}(b_\tau)=0.
    \]
\end{definition}

In the transverse setting, such a wall occurs
when a nonzero coefficient $b_{\tau,q_\tau}$
of $b_\tau$ has valuation tending to zero.
Several coefficients may do so simultaneously.
Strict positivity is required by our brane
convention in
Definition~\ref{def:immersed-lagrangian-brane}.

For the transverse configurations classified in
Lemma~\ref{lem:possible-transverse-walls},
the proposed resolution is a graded Lagrangian
$0$-surgery of the web cone, following the
surgery proposal of
\cite[Section~3.4]{Joyce-stability}.
Its underlying source surgery is the
$k=0$ case of
Definition~\ref{def:interior-surgery}.
We describe four realizations of this surgery
in the web geometry.
The selected branches are assumed to remain
transverse at the wall.

For certain walls supported on clean intersections between closed components,
we also allow the interior $0$-surgery model
described in
Section~\ref{subsubsec:interior-zero-surgery}
and the reductions of closed inactive components
described in
Section~\ref{subsubsec:closed-component-reduction}.

At an event time $\tau_*$, let $M_i^-$ be the
incoming boundary sources and write
\[
    B_i^-:=B_i(\tau_*^-),
    \qquad
    R_i^-:=M_i^-\setminus\overline{B_i^-},
    \qquad i=0,1,
\]
where $B_i(\tau_*^-)$ is the limiting incoming
trace defined in
Definition~\ref{def:surgery-compatible}.

\subsubsection{Classification of transverse walls}

The configurations of the two branches of a
degree-one generator are constrained by the
prescribed heights and gradings of the
web-cone pieces.

\begin{lemma}[Possible transverse obstruction walls]
    \label{lem:possible-transverse-walls}
    Let $q_\tau$ be a smoothly varying ordered
    degree-one self-intersection generator whose
    nonzero coefficient $b_{\tau,q_\tau}$ in
    $b_\tau$ has valuation tending to zero as
    $\tau\to\tau_*^-$.
    Suppose that its branches belong to geometric
    web-cone pieces, intersect transversely, and
    remain away from the smoothing locus throughout
    the incoming regular interval.
    The branches fall into one of the following
    configurations:
    \begin{enumerate}
        \item The generator runs from a branch of
              $R_0[1]$ to a branch of $R_1$.
              The corresponding intersection of the
              unshifted boundary sheets has Maslov
              index zero.

        \item One branch lies in a residue and the
              other in the interior of the source
              of an active front or a nonclosed
              inactive fragment.

        \item Both branches lie in the interiors
              of sources of active fronts or
              inactive fragments, including
              closed fragments.

        \item The generator runs from a branch of
              $R_0[1]$ to a closed fragment in
              $D_{\mathrm{cl},1}$.
    \end{enumerate}
\end{lemma}

\begin{proof}
    We use the fragment notation of
    Definition~\ref{def:inactive-fragment-notation}
    and the height function of
    Definition~\ref{def:height-function},
    suppressing the time argument in the
    residue notation.
    Away from the smoothing locus,
    $\widetilde\kappa_\tau=\kappa_\tau$, and
    \[
        \begin{aligned}
            \kappa_\tau|_{R_0[1]\cup D_0
            \cup D_{\mathrm{cl},0}[1]}
             & =0, \\
            \kappa_\tau|_{R_1\cup D_1
                    \cup D_{\mathrm{cl},1}}
             & =1.
        \end{aligned}
    \]
    On the interior of each active front source,
    one has $0<\kappa_\tau<1$.

    By \eqref{eq:web-cone-valuation-transport},
    a coefficient is nondecreasing in valuation
    if its target height is at most its source
    height throughout the incoming interval.
    Such a coefficient cannot have valuation
    tending to zero from a strictly positive
    initial value.

    A generator from $D_{\mathrm{cl},0}[1]$
    to a higher-height piece has source phase
    $3\pi/2$ and target phase
    $\theta_+\leq\pi/2$.
    Writing $\beta_1,\ldots,\beta_n\in(0,\pi)$
    for the characteristic angles,
    \eqref{eq:floer-generator-degree} gives
    \[
        \deg q_\tau
        =
        \frac{
            3\pi/2-\theta_+
            +\sum_{j=1}^n\beta_j
        }{\pi}
        >1.
    \]
    Thus this configuration cannot contribute
    a degree-one generator reaching the wall.

    If both branches lie in residues, their
    constant heights leave only
    $R_0[1]\to R_1$, giving case~(1).
    By \eqref{eq:floer-shift-convention},
    shifting the source by $[1]$ increases
    the degree by one, so the corresponding
    unshifted Maslov index is zero.

    If one branch lies in a residue and the
    other in a closed fragment, the height
    constraint and the preceding degree
    calculation leave only
    $R_0[1]\to D_{\mathrm{cl},1}$,
    giving case~(4).

    All remaining configurations fall under
cases~(2) and~(3).
\end{proof}

\begin{example}[A wall on coincident closed components]
\label{ex:cancellation-wall}
    Let $\mathcal S$ be a closed
    imaginary-special brane.
    Consider the stationary disjoint-source
    union of two copies, carrying labels $0$
    and $1$, with diagonal brane data
    $\mathcal S[1]$ and $\mathcal S$.
    Choose a degree-zero cocycle $e$ of
    valuation zero representing
    $1_{\mathcal S}$.
    After the source shift, $e$ has degree one
    from $\mathcal S[1]$ to $\mathcal S$.

    Suppose that the off-diagonal entry at
    time $\tau_0$ is $T^\lambda e$, with
    $\lambda>0$.
    The stationary rescaling determined by
    the heights $0$ and $1$ gives the entry
    \[
        T^{\lambda-(\tau-\tau_0)}e.
    \]
    Its valuation reaches zero at
    $\tau_*=\tau_0+\lambda$.
    Algebraically, this pair represents the
    cone of an invertible scalar multiple
    of $1_{\mathcal S}$, and hence the zero
    object in $\mathcal D$.
    The pair can therefore be cancelled,
    as an instance of
    Definition~\ref{def:closed-component-reduction}.
\end{example}

The proposed resolutions of cases~(1)--(4) are,
respectively, \emph{webbing}, \emph{puncture},
\emph{interior $0$-surgery}, and \emph{absorption}.
We describe these operations in the following
Sections.

\subsubsection{Webbing}
\label{subsubsec:webbing}

Consider case~(1) of
Lemma~\ref{lem:possible-transverse-walls}.
The selected degree-one generator runs from
$R_0[1]$ to $R_1$ and is specified by
\[
    x_i\in R_i^-\subset M_i,
    \qquad
    \iota_0(x_0)=\iota_1(x_1)=q.
\]
Choose neighborhoods $U_i\subset R_i^-$ of
$x_i$ on which $\iota_i$ is an embedding.
Theorem~\ref{thm:SY-start-web}, applied to
these selected sheets, supplies an imaginary
special Lagrangian webbing emanating from $q$.
We use the evolution parameter prescribed by
\eqref{eq:fundamental-harmonic-evolution},
with collapse time $\tau_*$.
The new boundary lifts converge to $x_i$ as
$\tau\to\tau_*^+$.
We call this operation
\emph{webbing at $(x_0,x_1)$}.

When the pre-existing fronts evolve regularly
through $\tau_*$, they continue by their
regular evolution, while the inactive
fragments remain stationary.

The new boundary spheres sweep small disks
out of the two residues, and the new cylinder
joins the resulting residue boundaries.
After rounding corners, this realizes interior
$0$-surgery on the web-cone source, as in
Definition~\ref{def:interior-surgery}.

The same operation is used at initialization
(Definition~\ref{def:initial-brane-comparison})
for valuation-zero coefficients of $B_f$
supported at transverse intersections of the
positive input components.
Initialization is not an obstruction wall,
since it has no incoming regular interval.

\subsubsection{Puncture}
\label{subsubsec:puncture}

Consider case~(2) of
Lemma~\ref{lem:possible-transverse-walls}.
The proposed resolution punctures the
participating front or nonclosed inactive
fragment, creating a new labelled boundary
on the selected residue sheet.

\begin{definition}[Web-front puncture]
    \label{def:web-front-puncture}
    Let
    \[
        Y_*=
        [F_*:P^-\longrightarrow X,
        \beta_{0,*},\beta_{1,*}]
    \]
    be the smooth incoming limit of the
    participating component.
    Fix $i\in\{0,1\}$ and choose
    \[
        p\in\operatorname{int}(P^-),
        \qquad
        r\in R_i^-,
        \qquad
        F_*(p)=\iota_i(r)=q,
    \]
    so that the selected front and residue
    sheets intersect transversely at $q$.
    Choose a neighborhood $U_r\subset R_i^-$
    of $r$ on which $\iota_i$ is an embedding,
    and an embedding
    \[
        \jmath^-:D^n\longrightarrow
        \operatorname{int}(P^-),
        \qquad \jmath^-(0)=p.
    \]
    Set
    \[
        A^-:=\jmath^-(D^n),
        \qquad C^-:=\{p\}.
    \]
    A \emph{web-front puncture at $(p,r)$}
    is the inverse of the boundary surgery
    of index $n-1$ in
    Definition~\ref{def:boundary-surgery}.
    It produces $P^+$ by replacing
    \[
        D^n
        \quad\rightsquigarrow\quad
        S^{n-1}\times[0,1],
    \]
    with attaching map $(u,0)\sim\jmath^-(u)$.
    Denote the embedding of the inserted piece
    into $P^+$ by $\jmath^+$, and set
    \[
        \begin{aligned}
            A^+
            &:=\jmath^+(S^{n-1}\times[0,1]),\\
            C^+
            &:=\jmath^+(S^{n-1}\times\{1\}).
        \end{aligned}
    \]
    The new boundary $C^+$ receives label $i$.
    Each outgoing immersion $F_\tau^+$ must be
    embedded on a neighborhood of $C^+$, and
    \[
        \beta_{i,\tau}^+(C^+)=\partial D_\tau
    \]
    for an embedded closed disk
    $D_\tau\subset U_r$ containing $r$.
    As $\tau\to\tau_*^+$, these disks must
    eventually lie in every neighborhood of $r$.
\end{definition}

\begin{figure}[htbp]
    \centering
    \includegraphics[scale=0.14]
    {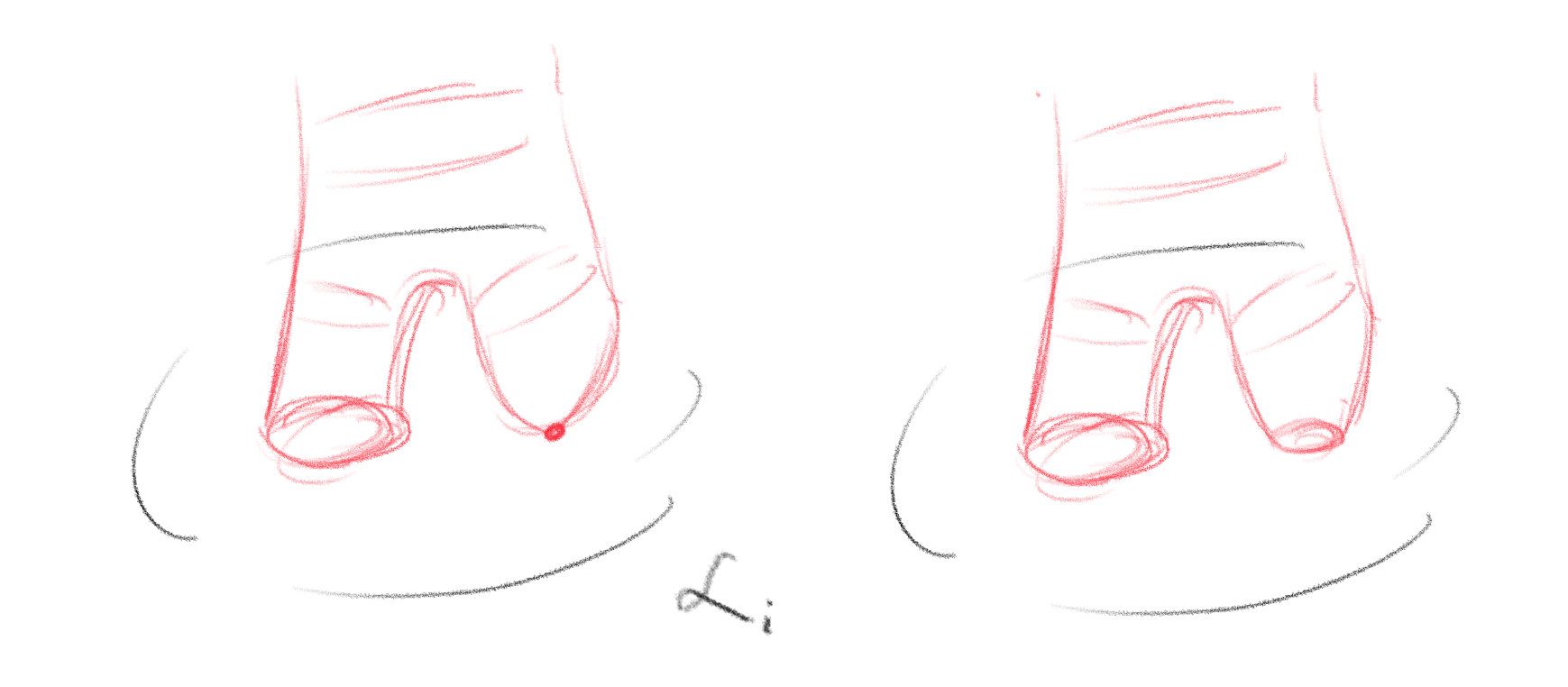}
    \caption{
        Puncture of a two-dimensional front,
        shown in red.
        Left: an interior point of the smooth
        incoming limit meets the selected residue
        sheet of $L_i$ at the marked point.
        Right: puncture creates a new boundary
        circle carrying label $i$.
        Its boundary lift is $\partial D_\tau$,
        where the disks $D_\tau$ shrink to the
        selected boundary source point as
        $\tau\to\tau_*^+$.
        The pre-existing boundary is retained.
    }
    \label{fig:web-front-puncture}
\end{figure}

\begin{definition}[Collapse maps for puncture]
    \label{def:puncture-collapse-maps}
    Since $C^-=\{p\}$, identify $P_*=P^-$,
    with $p_*=p$ and
    $\pi^-=\operatorname{id}_{P^-}$.
    On the inserted region $A^+$, define $\pi^+$ by
    \[
        \pi^+\bigl(\jmath^+(u,t)\bigr)
        :=
        \jmath^-\bigl((1-t)u\bigr),
        \qquad
        u\in S^{n-1},\quad 0\leq t\leq1.
    \]
    With gluing collars chosen as in
    Definition~\ref{def:boundary-degeneration-limit},
    using $\rho=1-t$, these maps define a common
    limiting source in the sense of
    Definition~\ref{def:source-replacement-limit}.
\end{definition}

For an incoming inactive fragment, the selected
residue must carry its missing boundary label:
otherwise the participating branches have equal
heights and cannot produce an obstruction wall.
Puncturing such a fragment therefore produces
active components.

Since $D_\tau$ shrinks to $r$, puncture adds
no open region to the boundary trace at the event.
Gluing the new boundary to the residue boundary
and rounding the seam replaces the two source
disks about $p$ and $r$ by a cylinder.
Thus puncture realizes interior $0$-surgery
of the web-cone source, as in
Definition~\ref{def:interior-surgery}.

\begin{remark}[Puncture for curves]
    When $n=1$, puncturing an active interval
    produces one active interval and one inactive
    interval.
    The latter is a boundary capping case of
    Definition~\ref{def:inactive-continuation}.
    Here $D_\tau$ is an interval with endpoint $r$.
    This configuration occurs in
    Definition~\ref{def:elliptic-puncture}
    and is illustrated in
    Figure~\ref{fig:elliptic-puncture}.
\end{remark}

\subsubsection{Interior $0$-surgery}
\label{subsubsec:interior-zero-surgery}
We describe case~(3) of
Lemma~\ref{lem:possible-transverse-walls}
when the participating source components are nonclosed.

Let $P^-$ be their disjoint union.
Assume that the incoming immersions and boundary
lifts extend smoothly to $\tau_*$, and denote
the limiting immersion and prescribed height by
$F_*^-:P^-\to X$ and $\kappa_*^-$.
Let $p_{\mathrm s},p_{\mathrm t}\in
\operatorname{int}(P^-)$ be the limiting source
and target points of the selected degree-one
generator, with common image $q$.
Assume that the selected branches remain transverse
at $q$ and that
\[
    \kappa_*^-(p_{\mathrm s})
    <
    \kappa_*^-(p_{\mathrm t}).
\]
Use the surgery of
Definition~\ref{def:interior-surgery}
with $k=0$ in disjoint source disks centered
at $p_{\mathrm s}$ and $p_{\mathrm t}$,
and the collapse maps of
Definition~\ref{def:interior-degeneration-limit}.
Thus $C^-=\{p_{\mathrm s},p_{\mathrm t}\}$,
and the smooth incoming limit descends to a
limiting front on the common limiting source.

Transversality forces $n\geq2$, so the outgoing
source $P^+$ is connected.
If an incoming component is active, it retains
both boundary labels.
Otherwise, the strict height inequality forces
the participating inactive components to have
opposite labels.
Thus $P^+$ is active.

Since the surgery is disjoint from the source
boundary, it induces the same $0$-surgery on
the assembled web-cone source.

\begin{remark}[Polterovich surgery and Lawlor necks]
    \label{rem:interior-neck-model}
    The Lagrangian handle construction of
    \cite[Section~1.3]{Polterovich-surgery}
    realizes the underlying source surgery.
    Both participating branches have grading
    $\pi/2$, so the degree-one condition in
    \eqref{eq:floer-generator-degree} gives
    \[
        \sum_{j=1}^n\beta_j=\pi,
    \]
    where $\beta_j\in(0,\pi)$ are the characteristic
    angles from the source branch to the target
    branch.
    This is Lawlor's angle condition
    \cite{Lawlor-angle-criterion}.
    For $n\geq3$, the corresponding imaginary-special
    tangent-space models and desingularization
    results are described in
    \cite[Section~9]{Joyce-SL-survey}.
    The continuation conjecture additionally
    requires a realization in $X$ with the
    prescribed boundary conditions, convergence,
    and normalized evolution.
\end{remark}

In the elliptic-curve setting,
Definition~\ref{def:elliptic-interior-surgery-closed}
also describes an interior $0$-surgery at a clean
intersection between a label-zero inactive interval
and a closed label-one fragment.
The outgoing interval is stationary.
Its geometric continuation and brane transfer
are treated in Section~\ref{sec:elliptic-curve}.

\subsubsection{Absorption}
\label{subsubsec:absorption}

Consider case~(4) of
Lemma~\ref{lem:possible-transverse-walls}.
Absorption incorporates the selected closed
label-one fragment into the boundary condition
$L_0^-$.

\begin{definition}[Absorption]
    \label{def:absorption}
    Let $S=[\iota_S:M_S\to X]$ be the selected
    component of $D_{\mathrm{cl},1}^-$.
    Choose
    \[
        r\in R_0^-,
        \qquad s\in M_S,
        \qquad
        \iota_0^-(r)=\iota_S(s)=q,
    \]
    and assume that the ordered intersection
    from $L_0^-[1]$ to $S$ is transverse
    of degree one.

    An \emph{absorption of $S$ at $(r,s)$}
    consists of the graded surgery
    \[
        \widehat L_0=L_0^-\#_q S[-1],
    \]
    followed by an exact deformation to a
    strictly positive graded immersion
    $L_0^+=[\iota_0^+:M_0^+\to X]$,
    as in Lemma~\ref{lem:positive-surgery}.
    The surgery disk about $r$ lies in $R_0^-$.
    Throughout the construction, the immersion
    is unchanged on an open neighborhood
    $V_0^-\subset M_0^-$ of $\overline{B_0^-}$.
    The resulting orientation-preserving
    diffeomorphism onto an open subset
    $V_0^+\subset M_0^+$ satisfies
    \[
        \psi_0:V_0^-\xrightarrow{\sim}V_0^+,
        \qquad
        \iota_0^+\circ\psi_0
        =
        \iota_0^-|_{V_0^-}.
    \]
\end{definition}

\begin{remark}[Outgoing web data after absorption]
    \label{rem:absorption-boundary-conditions}
    Replace the boundary pair by $(L_0^+,L_1)$,
    remove $S$ from the closed fragments,
    and retain all other fronts and fragments
    with their labels.
    On the front and fragment sources, this gives
    $A^-=M_S$, $A^+=\emptyset$, and $\phi$ equal
    to the identity on the retained sources.

    Transport the retained $0$-labelled boundary
lifts by $\psi_0$ and take
$\psi_1=\operatorname{id}_{M_1}$.
    The incoming traces are
    \[
        B_{0,\mathrm{in}}^+=\psi_0(B_0^-),
        \qquad
        B_{1,\mathrm{in}}^+=B_1^-.
    \]
    The initial residue
    $R_0^+:=M_0^+\setminus
    \overline{B_{0,\mathrm{in}}^+}$
    contains the punctured source of $S$
    and the surgery neck.

If the retained fronts admit a regular
continuation for the old boundary pair,
the same immersions and transported boundary
lifts give one for the new pair on a sufficiently
short outgoing interval.
The retained inactive fragments remain stationary.
\end{remark}

Shifting the surgery by $[1]$ gives
\[
    \widehat L_0[1]=L_0^-[1]\#_q S.
\]
Since $r$ lies in the residue, absorption
therefore realizes $0$-surgery between the
web-cone branches $R_0^-[1]$ and $S$,
followed by an exact deformation; see Figure~\ref{fig:elliptic-stationary-absorption}.

The following lemma supplies the geometric
construction, including its support and volume
control.
Its proof will be given in~\cite{KapotaWebGeometry}.
\begin{lemma}[Positive surgery]
    \label{lem:positive-surgery}
    Let $L=[\iota:M\to X]$ be a graded
    Lagrangian immersion with closed,
    possibly disconnected source and
    \[
        -\frac{\pi}{2}<\theta_L<\frac{\pi}{2}.
    \]
    Let $S=[\iota_S:M_S\to X]$ be an
    imaginary-special Lagrangian immersion
    with closed connected source and grading
    $\theta_S=\pi/2$.
    Suppose that the selected branches at
    \[
        q=\iota(p)=\iota_S(z)
    \]
    intersect transversely and admit the
    graded surgery $L\#_q S[-1]$.

    For every $\eta>0$ and source neighborhoods
    $p\in U\subset M$ and $z\in V\subset M_S$,
    there exist a graded surgery representative
    \[
        \widehat L=[G_0:N\longrightarrow X]
    \]
    and a smooth exact Lagrangian deformation
    \[
        G_t:N\longrightarrow X,
        \qquad 0\leq t\leq1,
    \]
    with the following properties.
    The initial surgery modifies the immersions
    only in source disks contained in $U$ and
    $V$.
    The deformation is stationary on the
    unchanged copy of $M\setminus U$.
    Its final immersion
    $L^+=[G_1:N\to X]$ satisfies
    \[
        -\frac{\pi}{2}
        <\theta_{L^+}
        <\frac{\pi}{2},
        \qquad
        \operatorname{Vol}_g(L^+)
        \leq
        \operatorname{Vol}_g(L)
        +\operatorname{Vol}_g(S)+\eta.
    \]
    The surgery grading extends those of
    $L$ and $S[-1]$ and is continued along
    the deformation.
    Volumes use the fixed K\"ahler metric
    $g$ and are counted on the sources,
    including multiplicity.
\end{lemma}

\subsubsection{Reduction of closed inactive components}
\label{subsubsec:closed-component-reduction}

Valuation-zero terms can also occur between coincident closed components, as in Example~\ref{ex:cancellation-wall}.
We allow these terms to be resolved by changing the
presentation of the closed components and cancelling
identity cones.

\begin{definition}[Reduction of closed inactive components]
    \label{def:closed-component-reduction}
    Let
    \[
        (\widehat C,\widehat E,\widehat b)
    \]
    be provisional smoothed cone data used
    during initialization or continuation
    across an event, where $\widehat b$ is
    a degree-one cochain with
    \[
        \operatorname{val}(\widehat b)\geq0.
    \]
    Suppose that the remaining valuation-zero
    terms occur between closed inactive
    components, while the terms involving
    the other web-cone pieces have strictly
    positive valuation.

    A \emph{reduction of closed inactive components}
    replaces selected closed inactive components
    and their local systems by closed
    imaginary-special components with rank-one
    $\Bbbk$-local systems and labels $0$ or $1$.

    At leading order, obtained by reducing
    the nonnegative Floer data modulo
    $\Lambda_+$, the replacement is given
    by label-preserving changes of presentation
    of the closed part, followed by cancellation
    of identity cones.
    The reduction includes a bounding cochain
    $b_+$ giving a brane
    \[
        \mathcal C_+=(C_+,E_+,b_+),
        \qquad
        \operatorname{val}(b_+)>0.
    \]
    The nonclosed geometric pieces and their
    local systems, the boundary conditions,
    the boundary traces, and the boundary
    lifts are unchanged.

    The reduction must be compatible with
    a brane transfer for the complete
    continuation, in the sense of
    Definition~\ref{def:brane-surgery-transfer}.
    At initialization, the required categorical
    comparison is instead the initial brane
    comparison of
    Definition~\ref{def:initial-brane-comparison}.
\end{definition}

The comparison concerns the full cone brane.
The closed components need not form a direct summand:
positive-valuation terms coupling them to the remaining components are incorporated into the outgoing bounding cochain.

The required convergence arguments for the
auxiliary data and the construction of the
brane comparisons will be given
in~\cite{KapotaWebFloer}.
For the two-torus,
Lemma~\ref{lem:elliptic-closed-cancellation}
states the corresponding existence result
under its hypotheses and controls the
resulting closed-component multiplicities.

\section{Categorical consequences}
\label{sec:categorical-consequences}

\subsection{The abelian subcategory of positive branes}
\label{subsec:positive-branes-abelian}

We work with the derived Fukaya category $\mathcal D$
and the Floer conventions of
Subsection~\ref{subsec:brane-structures}.

\begin{definition}
    \label{def:positive-brane-subcategory}
    Let $\mathcal A\subset\mathcal D$ be the full
    subcategory whose objects are isomorphic to finite
    direct sums of immersed branes
    $\mathcal L=(L,E,b)$ with $L=P\sqcup S$ and
    \begin{equation}
        \label{eq:positive-brane-window}
        -\frac{\pi}{2}<\theta_P<\frac{\pi}{2},
        \qquad
        \theta_S\equiv\frac{\pi}{2},
        \qquad
        \operatorname{val}(b)>0.
    \end{equation}
    Here $P$ is a union of strictly positive components,
    $S$ is a union of closed imaginary-special components,
    and $\sqcup$ denotes a disjoint union of sources.
    Either part may be empty.
    The gradings are the displayed real lifts.
    The empty direct sum is included.
\end{definition}

By construction, $\mathcal A$ is closed under
isomorphisms and finite direct sums.

\begin{lemma}[{cf.\ \cite[Remark~A.2]{YangLi}}]
    \label{lem:positive-negative-vanishing}
    For $\mathcal L_0,\mathcal L_1\in\mathcal A$ and
    $k<0$, one has
    \begin{equation}
        \label{eq:positive-negative-vanishing}
        \operatorname{Hom}_{\mathcal D}
        (\mathcal L_0,\mathcal L_1[k])=0.
    \end{equation}
\end{lemma}

\begin{proof}
    First consider two defining branes.
    Their sources are compact, and the phase $-\pi/2$
    is excluded.
    Their phases therefore lie in a common compact
    interval of length strictly less than $\pi$.
    Sufficiently small generic Hamiltonian perturbations
    used to compute the mixed Floer complex preserve
    this strict phase-width bound.

    At a transverse generator $q=(q_-,q_+)$,
    Equation~\eqref{eq:floer-generator-degree} gives
    \[
        \deg(q)
        =\frac{1}{\pi}
        \left(
        \theta_{L_0}(q_-)-\theta_{L_1}(q_+)
        +\sum_{j=1}^n\beta_j
        \right)>-1,
        \qquad 0<\beta_j<\pi.
    \]
    The degree is integral, hence nonnegative.
    Bounding-cochain deformation changes the differential
    but not the degrees of the generators.
    Thus the Floer complex has no negative-degree terms.
    The claim follows for finite direct sums and
    isomorphic objects.
\end{proof}

The following extension property will also be used
to construct a bounded heart.

\begin{lemma}[Closure under extensions]
    \label{lem:positive-extension-closure}
    The subcategory $\mathcal A$ is closed under
    extensions: in every distinguished triangle
    \[
        A\longrightarrow C\longrightarrow B
        \xrightarrow{a}A[1]
    \]
    with $A,B\in\mathcal A$, one has $C\in\mathcal A$.
\end{lemma}

\begin{proof}
    Finite direct sums of defining branes can be
    realized by disjoint-source unions with
    block-diagonal bounding cochains, as in
    Section~\ref{subsubsec:disconnected-sources}.
    We may therefore choose brane representatives
    \[
        \mathcal L_A=(L_A,E_A,b_A),
        \qquad
        \mathcal L_B=(L_B,E_B,b_B)
    \]
    for $A$ and $B$, respectively, satisfying
    \eqref{eq:positive-brane-window}.
    Under these identifications, represent $a$ by
    a closed cochain
    \[
        \alpha\in CF^1(\mathcal L_B,\mathcal L_A).
    \]
    Choose $\lambda\in\mathbb R$ such that
    $\operatorname{val}(T^\lambda\alpha)>0$.
    On the disjoint-source union $L_A\cup L_B$,
    with local system $E_A\sqcup E_B$, set
    \[
        b_\lambda=
        \begin{pmatrix}
            b_A & T^\lambda\alpha \\
            0   & b_B
        \end{pmatrix}.
    \]
    By source compatibility of the Floer operations,
    its Maurer--Cartan equation consists of the two
    diagonal equations and the closedness equation
    for $T^\lambda\alpha$, as in
    Lemma~\ref{lem:geometric-cone}.
    Thus $b_\lambda$ is a bounding cochain of strictly
    positive valuation.

    The gradings on $L_A$ and $L_B$ are unchanged,
    so the resulting brane $\mathcal C_\lambda$
    belongs to $\mathcal A$.
    Its triangular presentation identifies it with
    $\operatorname{Cone}(T^\lambda a)[-1]$.
    By the scalar-rescaling observation in
    Remark~\ref{rem:cone-scaling},
    \[
        \mathcal C_\lambda
        \simeq\operatorname{Cone}(T^\lambda a)[-1]
        \simeq\operatorname{Cone}(a)[-1]
        \simeq C.
    \]
    Hence $C\in\mathcal A$.
\end{proof}

We now use the terminal triangle to identify
the geometric kernels and cokernels.

\begin{prop}[Geometric kernels and cokernels]
    \label{prop:geometric-kernel-cokernel}
    Let $\mathcal L_0,\mathcal L_1$ be defining
    geometric branes in $\mathcal A$, and let
    $f\in CF^0(\mathcal L_0,\mathcal L_1)$ be closed.
    Suppose that $f$ admits a terminal special Lagrangian
    web evolution to which
    Proposition~\ref{prop:terminal-web-triangle} applies.
    Then the terminal branes $\mathcal K,\mathcal K^*$
    constructed in
    Section~\ref{subsubsec:terminal-triangle-proof}
    belong to $\mathcal A$, and the complementary
    morphisms of
    Definition~\ref{def:complementary-morphisms},
    \[
        \kappa_0:\mathcal K\longrightarrow\mathcal L_0,
        \qquad
        \kappa_1:\mathcal L_1\longrightarrow\mathcal K^*,
    \]
    are a kernel and a cokernel of $[f]$ in
    $\mathcal A$, respectively.
\end{prop}

\begin{proof}
    We follow the argument of
    \cite[Proposition~1.2.2]{BBD}.
    The terminal representatives have the componentwise
    phases of \eqref{eq:positive-brane-window} by
    Lemma~\ref{lem:terminal-smoothing}.
    Their brane structures, constructed in the proof of
    Proposition~\ref{prop:terminal-web-triangle},
    have rank-one local systems and bounding cochains
    of strictly positive valuation.
    Thus $\mathcal K,\mathcal K^*\in\mathcal A$.
 
    Fix $A\in\mathcal A$, and take all Hom spaces
    below in $\mathcal D$.
    Write $\mathcal C_f=\operatorname{Cone}(f)$. Apply $\operatorname{Hom}(A,-)$ to the terminal
    triangle~\eqref{eq:terminal-web-triangle}
    shifted by $[-1]$, and
    $\operatorname{Hom}(-,A)$ to the unshifted
    terminal triangle.
    The long exact sequences give
    \[
        \begin{gathered}
            0\longrightarrow
            \operatorname{Hom}(A,\mathcal K)
            \longrightarrow
            \operatorname{Hom}(A,\mathcal C_f[-1])
            \longrightarrow0,\\
            0\longrightarrow
            \operatorname{Hom}(\mathcal K^*,A)
            \longrightarrow
            \operatorname{Hom}(\mathcal C_f,A)
            \longrightarrow0.
        \end{gathered}
    \]
    Indeed, the adjacent terms vanish by
    Lemma~\ref{lem:positive-negative-vanishing}:
    they are $\operatorname{Hom}(A,\mathcal K^*[-j])$
    and
    \[
        \operatorname{Hom}(\mathcal K[j],A)
        \cong\operatorname{Hom}(\mathcal K,A[-j]),
        \qquad j=1,2.
    \]
    Using these natural isomorphisms in the long exact
    sequences of the initial cone triangle~\eqref{eq:initial-cone-triangle} gives
    \[
        \begin{gathered}
            0\longrightarrow
            \operatorname{Hom}(A,\mathcal K)
            \xrightarrow{\kappa_0\circ-}
            \operatorname{Hom}(A,\mathcal L_0)
            \xrightarrow{[f]\circ-}
            \operatorname{Hom}(A,\mathcal L_1),\\
            0\longrightarrow
            \operatorname{Hom}(\mathcal K^*,A)
            \xrightarrow{-\circ\kappa_1}
            \operatorname{Hom}(\mathcal L_1,A)
            \xrightarrow{-\circ[f]}
            \operatorname{Hom}(\mathcal L_0,A).
        \end{gathered}
    \]
    The initial zeros use
    $\operatorname{Hom}(A,\mathcal L_1[-1])=0$ and
    $\operatorname{Hom}(\mathcal L_0[1],A)=0$,
    respectively.
    The maps are composition with $\kappa_0$ and
    $\kappa_1$ by
    Definition~\ref{def:complementary-morphisms}.
    Since $\mathcal A$ is full, these exact sequences
    express the kernel and cokernel universal
    properties.
\end{proof}

In the presence of negative-degree vanishing
\eqref{eq:positive-negative-vanishing}, a full abelian
subcategory is called \emph{admissible} if every short
exact sequence in it extends to a distinguished triangle
in $\mathcal D$; see \cite[Definition~1.2.5]{BBD}.

\begin{proof}[Proof of Theorem~\ref{thm:intro-abelianity}]
    Every object of $\mathcal A$ can be represented
    by a disjoint-source union of defining branes.
    By Lemma~\ref{lem:positive-negative-vanishing}
    and closure under finite direct sums, the
    criterion of \cite[Proposition~1.2.4]{BBD}
    reduces abelianity and admissibility to the
    following condition: the cone of every morphism
    in $\mathcal A$ fits into a distinguished triangle
    of the form~\eqref{eq:terminal-web-triangle},
    with $\mathcal K,\mathcal K^*\in\mathcal A$.

    For the zero morphism, the split cone
    $\mathcal L_0[1]\oplus\mathcal L_1$ gives such
    a triangle with
    $\mathcal K=\mathcal L_0$ and
    $\mathcal K^*=\mathcal L_1$.
    The kernel and cokernel maps are the identity maps.

    For a nonzero morphism, choose a closed
    degree-zero Floer cochain $f$ representing it.
    Conjectures~\ref{conj:intro-initialization},~\ref{conj:intro-continuation}
    and~\ref{conj:intro-termination}
    supply a terminal special Lagrangian web evolution.
    Proposition~\ref{prop:terminal-web-triangle}
    gives the required triangle, and
    Proposition~\ref{prop:geometric-kernel-cokernel}
    identifies its terminal branes and complementary
    morphisms as the kernel and cokernel.
    BBD's criterion therefore proves the theorem.
\end{proof}

\subsection{Harder--Narasimhan filtrations}
\label{subsec:harder-narasimhan}

Throughout this subsection, assume
Conjectures~\ref{conj:intro-initialization},~\ref{conj:intro-continuation}
and~\ref{conj:intro-termination}.
Then $\mathcal A$ is abelian by
Theorem~\ref{thm:intro-abelianity}.

If $X$ is noncompact, assume that the metric
\[
    \widetilde g=\rho^{2/n}g
\]
is complete, has bounded sectional curvature,
and has injectivity radius bounded away from zero.
Here $g$ and $\rho$ are the K\"ahler metric and
Calabi--Yau density recalled in
Subsection~\ref{subsec:reminder-webbing}.
These conditions are automatic when $X$ is compact.
The form $\operatorname{Im}\Omega$ calibrates
imaginary-special Lagrangians for $\widetilde g$:
for such an immersion $F:P\to X$,
\[
    F^*\operatorname{Im}\Omega
    =d\operatorname{vol}_{F^*\widetilde g}
    =(\rho\circ F)\,d\operatorname{vol}_{F^*g}.
\]
In particular, these immersions are minimal
for $\widetilde g$.

For a geometric brane $\mathcal L=(L,E,b)$, with
$L=[\iota:N\to X]$, recall the central charge
\begin{equation}
    \label{eq:geometric-central-charge}
    Z(\mathcal L):=\int_L\Omega
    =\int_N\iota^*\Omega.
\end{equation}
For an abelian or triangulated category $\mathcal D$,
we write $K_0(\mathcal D)$ for its Grothendieck group,
generated by isomorphism classes $[E]$ with relations
$[E]=[E']+[E'']$ for each short exact sequence
$0\to E'\to E\to E''\to0$, or distinguished triangle
$E'\to E\to E''\to E'[1]$, respectively.

\begin{remark}[The geometric central charge]
    \label{rem:homological-compatibility}
    We assume that the standard filtered open--closed
    construction is available in our Fukaya-category setting.
    Its zero-energy term, together with the Hochschild
    Chern character, yields a homomorphism
    \[
        Z:K_0(\mathcal D)\longrightarrow\mathbb C,
        \qquad
        Z([\mathcal L])=\int_L\Omega.
    \]
    See \cite[Theorem~16.1(1)]{AFOOO} and
    \cite[\S3.1 and Remark~3.2]{YangLi}.
    We use its restrictions to $\mathcal A$ and the
    rotated categories $\mathcal A_c$.
\end{remark}

\begin{definition}
    \label{def:stability-function}
    A homomorphism $W:K_0(\mathcal A)\to\mathbb C$ is a
    \emph{stability function} if
    \[
        W(A)\in\mathbb H_{\mathrm R}
        :=\left\{
        r e^{i\phi}:r>0,\ -\pi/2<\phi\leq\pi/2
        \right\}
    \]
    for every nonzero $A\in\mathcal A$.
    Its \emph{phase} is
    \[
        \phi_W(A):=\arg W(A)\in(-\pi/2,\pi/2].
    \]
\end{definition}

\begin{lemma}
    \label{lem:central-charge-stability}
    The central charge $Z$ is a stability function on
    $\mathcal A$.
\end{lemma}

\begin{proof}
    For a defining representative $L=P\sqcup S$ as in
    Definition~\ref{def:positive-brane-subcategory},
    the identity
    \[
        \iota^*\Omega
        =(\rho\circ\iota)e^{i\theta_L}
        \,d\operatorname{vol}_{\iota^*g}
    \]
    gives $\operatorname{Re}Z(\mathcal L)>0$ if
    $P\neq\emptyset$.
    If $P=\emptyset$ and $S\neq\emptyset$, then
    \[
        Z(\mathcal L)
        =i\int_S\rho\,d\operatorname{vol}_g
        \in i\mathbb R_{>0}.
    \]
    The same conclusion holds for every nonempty finite
    disjoint-source union of defining representatives.
    Isomorphism invariance therefore gives
    \[
        Z(A)\in
        \{r e^{i\phi}:r>0,\ -\pi/2<\phi\leq\pi/2\}
    \]
    for every nonzero $A\in\mathcal A$.
\end{proof}

We use the definitions of semistability, destabilising
subobjects and quotients, and Harder--Narasimhan
filtrations given in Definition~\ref{def:intro-semistability-HN}.
Additivity gives the see-saw property: in a short exact
sequence $0\to\mathcal M\to\mathcal L\to\mathcal Q\to0$
with $\mathcal M,\mathcal Q\neq0$, one has
\begin{equation}
    \label{eq:destabilising-seesaw}
    \phi_Z(\mathcal M)>\phi_Z(\mathcal L)
    \quad\Longleftrightarrow\quad
    \phi_Z(\mathcal Q)<\phi_Z(\mathcal L).
\end{equation}
We now prove the geometric boundedness statement from
the introduction. All volumes below are total source
volumes, including multiplicity.

\begin{lemma}[Compact containment of terminal quotients]
    \label{lem:terminal-quotient-containment}
Fix a defining representative
$\mathcal L=(L,E,b)$ with $L=P\sqcup S$,
a compact neighbourhood $U\subset X$ of its image,
and a constant $V\geq0$.
There exists a compact subset $Y\subset X$,
depending only on $U$ and $V$,
with the following property.

    Let $\mathcal L\twoheadrightarrow\mathcal Q$
    be a quotient in $\mathcal A$, and consider
    a terminal evolution of its kernel inclusion
whose target representatives remain in $U$.
    If its terminal special fragments satisfy
    \[
        \operatorname{Vol}_{\widetilde g}(D_1)
        +\operatorname{Vol}_{\widetilde g}
        (D_{\mathrm{cl},1})
        \leq V,
    \]
    then the singular terminal representative
    \[
        K^*_{\mathrm{sing}}
        =(\overline{R_1}\cup D_1)
        \sqcup D_{\mathrm{cl},1}
    \]
    has image contained in $\operatorname{int}(Y)$.
\end{lemma}

\begin{proof}
Let $\Sigma\subset X$ be the image
of $K^*_{\mathrm{sing}}$.
The terminal residue and the boundaries of all
nonclosed fragments lie in $U$.

    We first show that every connected component
    of $\Sigma$ meets $U$.
    Suppose otherwise that a component $\Sigma'$
    is disjoint from $U$.
    It consists entirely of closed
    imaginary-special fragments and is separated
    by a positive distance from
    $U$ and the remaining support.
    Choose the terminal roundings and perturbations
    sufficiently small to preserve these separations.
    The closed fragments themselves remain fixed.

    There are then no mixed Floer cochains between
    the source parts supported on $\Sigma'$ and
    its complement.
    By Section~\ref{subsubsec:disconnected-sources}
    and Proposition~\ref{prop:geometric-kernel-cokernel},
    the quotient splits as
    \[
        \mathcal Q\simeq\mathcal Q'\oplus\mathcal Q'',
    \]
    where $\mathcal Q'\in\mathcal A$ is supported
    on $\Sigma'$.
    Its positive imaginary period implies
    $\mathcal Q'\neq0$.
    Disjointness from $U$ gives
    \[
        \operatorname{Hom}_{\mathcal D}
        (\mathcal L,\mathcal Q')=0.
    \]
    Hence the projection $\mathcal Q\to\mathcal Q'$
    becomes zero after composition with
    $\mathcal L\twoheadrightarrow\mathcal Q$.
    This contradicts the epimorphism property
    in the full subcategory $\mathcal A$.

    The special fragments are minimal for
    $\widetilde g$, and their boundaries lie
    in $U$.
    Uniform monotonicity supplies constants
    $r_0,v_0>0$, depending only on the geometric
    bounds on $\widetilde g$, such that
    \[
        \operatorname{Vol}_{\widetilde g}
        \bigl(
        (D_1\sqcup D_{\mathrm{cl},1})
        \cap B^{\widetilde g}_{r_0}(x)
        \bigr)
        \geq v_0
    \]
    whenever $x\in\Sigma$ and
    $\operatorname{dist}_{\widetilde g}
        (x,U)\geq2r_0$.
    Volume inside an ambient ball is counted
    on its full inverse image in the fragment
    sources.
    Fix $x\in\Sigma$ and put
    \[
        d:=\operatorname{dist}_{\widetilde g}
        (x,U),
        \qquad
        m:=\left\lfloor\frac{d}{3r_0}\right\rfloor.
    \]
    The connected component containing $x$
    meets $U$.
    The continuous distance function therefore
    takes every value between $0$ and $d$
    on that component.
    Choose points $x_j$ at distance $3jr_0$
    from $U$, for $1\leq j\leq m$.
    Their radius-$r_0$ balls are pairwise disjoint
    and avoid $U$.
    Monotonicity and the assumed volume bound give
    \[
        m v_0\leq V,
        \qquad
        d<R_V:=3r_0\left(1+\frac{V}{v_0}\right).
    \]
    Since $U$ is compact and
    $\widetilde g$ is complete, the Hopf--Rinow
    theorem implies that
    \[
        Y:=
        \left\{
        x\in X:
        \operatorname{dist}_{\widetilde g}
        (x,U)
        \leq R_V+1
        \right\}
    \]
    is compact.
    The preceding distance bound gives
    $\Sigma\subset\operatorname{int}(Y)$.
\end{proof}

\begin{proof}[Proof of
        Proposition~\ref{prop:intro-destabilising-volume-bound}]
    Fix a defining representative
    $\mathcal L=(L,E,b)$, with $L=P\sqcup S$.
    Finite direct sums are represented by
    disjoint-source unions.
    If $\operatorname{Re}Z(\mathcal L)=0$, additivity
    and nonnegativity of real parts imply that every
    nonzero quotient has phase $\pi/2$.
    There are then no destabilising quotients,
    and the assertion is vacuous.
    Hence suppose that
    \[
        r_{\mathcal L}:=\operatorname{Re}Z(\mathcal L)>0,
        \qquad
        \mu_{\mathcal L}
        :=\frac{\operatorname{Im}Z(\mathcal L)}
        {r_{\mathcal L}}.
    \]
    Fix a compact neighborhood $U$ of the image of $L$.
By $C^1$-continuity on the compact source,
choose the target perturbations with sufficiently
small total $C^1$-size that every perturbed strictly
positive part $P'$ lies in $U$ and satisfies
\[
    \begin{gathered}
        \operatorname{Vol}_g(P')
        \leq \operatorname{Vol}_g(P)+1,\\
        \int_{P'}\rho|\sin\theta_{P'}|\,
        d\operatorname{vol}_g
        \leq
        \int_P\rho|\sin\theta_P|\,
        d\operatorname{vol}_g+1.
    \end{gathered}
\]
    Let $\mathcal L\twoheadrightarrow\mathcal Q$
    be a destabilising quotient, and let
    $\mathcal M\hookrightarrow\mathcal L$
    be its kernel inclusion.
    Choose a defining representative of $\mathcal M$
    and a terminal web evolution for this inclusion.
    By Proposition~\ref{prop:geometric-kernel-cokernel},
    the terminal brane $\mathcal K^*$ represents
    $\mathcal Q$.
    Its singular representative is
    \[
        K^*_{\mathrm{sing}}
        =(\overline{R_1}\cup D_1)\sqcup D_{\mathrm{cl},1}.
    \]
Write $P'$ for the terminal strictly positive
target representative.
Then $R_1$ is a region in its source,
and Hamiltonian invariance gives
\[
    \int_{P'}\operatorname{Re}\Omega=r_{\mathcal L}.
\]
    All fragments in $D_1$ and $D_{\mathrm{cl},1}$
    carry their imaginary-special orientations.
    Put
    \[
        \begin{aligned}
            r & :=\int_{R_1}\operatorname{Re}\Omega,
              & a                                    & :=\int_{R_1}\operatorname{Im}\Omega, \\
            s & :=\int_{D_1}\operatorname{Im}\Omega
            +\int_{D_{\mathrm{cl},1}}
            \operatorname{Im}\Omega.
        \end{aligned}
    \]
    The terminal smoothing preserves the homology
    class of the singular representative.
    Homological compatibility therefore gives
    \begin{equation}
        \label{eq:destabilising-quotient-period}
        Z(\mathcal Q)=Z(\mathcal K^*)=r+i(a+s).
    \end{equation}
    Moreover,
    \[
        0\leq r\leq r_{\mathcal L},
        \qquad
        s\geq0,
        \qquad
        |a|\leq
        \int_P\rho|\sin\theta_P|\,
        d\operatorname{vol}_g +1 .
    \]
    Since $\mathcal Q$ is destabilising,
    \[
        \phi_Z(\mathcal Q)
        <\phi_Z(\mathcal L)<\pi/2.
    \]
    Hence $r>0$, and the phase inequality becomes
    \[
        a+s<\mu_{\mathcal L}r.
    \]
    Consequently,
    \begin{equation}
        \label{eq:destabilising-special-period-bound}
        \begin{split}
            0\leq s
             & <\mu_{\mathcal L}r-a \\
             & \leq
            \max\{0,\operatorname{Im}Z(\mathcal L)\}
            +\int_P\rho|\sin\theta_P|\,
            d\operatorname{vol}_g +1
            =:A_{\mathcal L}.
        \end{split}
    \end{equation}
    This bound depends only on the fixed
    representative of $\mathcal L$.

    The calibration identity gives
    \[
        \operatorname{Vol}_{\widetilde g}(D_1)
        +\operatorname{Vol}_{\widetilde g}
        (D_{\mathrm{cl},1})
        =s<A_{\mathcal L}.
    \]
    Lemma~\ref{lem:terminal-quotient-containment},
    applied with this $U$ and $V=A_{\mathcal L}$, therefore
    supplies a compact subset $Y_{\mathcal L}\subset X$,
    depending only on the fixed representative
    of $\mathcal L$, such that
    \[
        K^*_{\mathrm{sing}}
        \subset\operatorname{int}(Y_{\mathcal L}).
    \]
    Put
    \[
        c_{\mathcal L}:=\min_{Y_{\mathcal L}}\rho>0.
    \]
    The special fragments therefore satisfy
    \[
        \operatorname{Vol}_g(D_1)
        +\operatorname{Vol}_g(D_{\mathrm{cl},1})
        \leq c_{\mathcal L}^{-1}s.
    \]
    Using Proposition~\ref{prop:web-cone-smoothing}
    and Lemma~\ref{lem:terminal-smoothing},
    choose the terminal roundings and perturbations
    inside $Y_{\mathcal L}$, with total
    $g$-volume error at most $1$.
    The resulting brane representative
    $\mathcal K^*\simeq\mathcal Q$ satisfies
    $K^*\subset Y_{\mathcal L}$ and
    \begin{equation}
        \label{eq:destabilising-geometric-bound}
        \begin{split}
            \operatorname{Vol}_g(K^*)
             & \leq
            \operatorname{Vol}_g(\overline{R_1})
            +\operatorname{Vol}_g(D_1)
            +\operatorname{Vol}_g(D_{\mathrm{cl},1})
            +1      \\
             & \leq
            \operatorname{Vol}_g(P)
            +c_{\mathcal L}^{-1}A_{\mathcal L}
            +2
            =:C_{\mathcal L}.
        \end{split}
    \end{equation}
    Both $Y_{\mathcal L}$ and $C_{\mathcal L}$
    depend only on the fixed representative
    of $\mathcal L$, proving the proposition.
\end{proof}

\begin{lemma}[{cf.\ \cite[discussion following
                    Corollary~5.14]{YangLi}}]
    \label{lem:bounded-volume-finite-homology}
    For every compact subset $Y\subset X$ and
    every $C>0$, only finitely many classes in
    $H_n(X;\mathbb Z)$ are represented by closed
    oriented immersed $n$-manifolds whose images
    lie in $Y$ and whose total source $g$-volume
    is at most $C$.
\end{lemma}

\begin{proof}
    Choose a compact smooth codimension-zero
    submanifold $W\subset X$, possibly with boundary,
    such that $Y\subset\operatorname{int}(W)$.
    Its integral homology is finitely generated:
    \[
        H_n(W;\mathbb Z)\cong\mathbb Z^r\oplus T,
    \]
    with $T$ finite.
    Let $\lambda_1,\ldots,\lambda_r$ be the
    integer coordinate homomorphisms on the free part.
    By the universal coefficient theorem and
    de Rham's theorem, choose smooth closed
    $n$-forms $\eta_j$ on $W$ representing these
    homomorphisms over $\mathbb R$.

    Regard each immersion in the statement as
    a map $\iota:N\to W$.
    Then
    \[
        \bigl|\lambda_j(\iota_*[N])\bigr|
        =
        \left|\int_N\iota^*\eta_j\right|
        \leq C\sup_W|\eta_j|_g.
    \]
    Each integer coordinate therefore belongs
    to a fixed finite set.
    Since $T$ is finite, only finitely many
    classes in $H_n(W;\mathbb Z)$ occur.
    Their images under the inclusion-induced map
    \[
        H_n(W;\mathbb Z)\longrightarrow H_n(X;\mathbb Z)
    \]
    also form a finite set.
\end{proof}

\begin{corollary}
    \label{cor:finite-destabilising-phases}
    Under the hypotheses of
    Proposition~\ref{prop:intro-destabilising-volume-bound},
    the central charges, and hence the phases,
    of the destabilising subobjects and quotients of
    each fixed nonzero object form finite sets.
\end{corollary}

\begin{proof}
    Fix a nonzero object $\mathcal L\in\mathcal A$.
    Proposition~\ref{prop:intro-destabilising-volume-bound}
    supplies geometric representatives of its
    destabilising quotients contained in a common
    compact subset $Y_{\mathcal L}\subset X$,
    with uniformly bounded total source volume.
    By Lemma~\ref{lem:bounded-volume-finite-homology},
    their homology classes belong to a finite set.
    Since $\Omega$ is closed and its periods are
    invariant under isomorphism, the central charges
    of these quotients also form a finite set.

    If $\mathcal M\subset\mathcal L$ is a
    destabilising subobject, then
    $\mathcal Q=\mathcal L/\mathcal M$ is nonzero
    and destabilising by
    \eqref{eq:destabilising-seesaw}.
    Additivity gives
    \[
        Z(\mathcal M)
        =Z(\mathcal L)-Z(\mathcal Q).
    \]
    Thus the central charges of the destabilising
    subobjects form a finite set as well.
    The same conclusion follows for their phases.
\end{proof}

\begin{proof}[Proof of
        Theorem~\ref{thm:intro-positive-branes-HN}]
    By \cite[Proposition~2.4]{Bridgeland-stability},
    applied to $iZ$ with the corresponding increasing
    normalization of phases, it suffices to exclude
    infinite subobject chains with strictly increasing
    phases and infinite quotient chains with strictly
    decreasing phases.

    In a chain of nonzero subobjects
    \[
        \cdots\subset A_2\subset A_1\subset A_0,
        \qquad
        \phi_Z(A_{j+1})>\phi_Z(A_j),
    \]
    every $A_j$, $j\geq1$, destabilises the fixed object
    $A_0$.
    Its phase belongs to the finite set supplied by
    Corollary~\ref{cor:finite-destabilising-phases},
    which rules out such a chain.
    Likewise, in a chain of nonzero quotients
    \[
        B_0\twoheadrightarrow B_1\twoheadrightarrow B_2
        \twoheadrightarrow\cdots,
        \qquad
        \phi_Z(B_{j+1})<\phi_Z(B_j),
    \]
    every $B_j$, $j\geq1$, is a destabilising quotient
    of $B_0$.
    Finiteness of their phases rules out this chain as
    well. The criterion therefore gives the
    Harder--Narasimhan property.
\end{proof}

\subsection{Stability conditions on the derived Fukaya category}
\label{subsec:stability-conditions}

Throughout this subsection, assume the geometric
hypotheses of Subsection~\ref{subsec:harder-narasimhan},
the homological compatibility of
Remark~\ref{rem:homological-compatibility},
and Conjectures~\ref{conj:intro-initialization},~\ref{conj:intro-continuation}
and~\ref{conj:intro-termination}.
Thus $\mathcal A\subset\mathcal D$ is an admissible
abelian subcategory, closed under extensions, and
$Z$ is a stability function with the
Harder--Narasimhan property.
We first construct a stability condition on the
subcategory generated by $\mathcal A$.
We then discuss generation of $\mathcal D$ and
compare the rotated positive-brane categories
to prove Theorem~\ref{thm:intro-special-representatives}.

\subsubsection{The generated heart}

\begin{definition}
    For a collection $S$ of objects of $\mathcal D$,
    write $\langle S\rangle_{\mathrm{ext}}$ for
    the strictly full subcategory obtained by
    taking finite extensions, including the
    zero object.
    Set
    \begin{equation}
        \label{eq:BBD-generated-subcategories}
        \begin{aligned}
            D^b_{\mathcal A}
             & :=
            \bigl\langle
            \mathcal A[m]:m\in\mathbb Z
            \bigr\rangle_{\mathrm{ext}}, \\
            D^{b,\leq0}_{\mathcal A}
             & :=
            \bigl\langle
            \mathcal A[m]:m\geq0
            \bigr\rangle_{\mathrm{ext}}, \\
            D^{b,\geq0}_{\mathcal A}
             & :=
            \bigl\langle
            \mathcal A[m]:m\leq0
            \bigr\rangle_{\mathrm{ext}}.
        \end{aligned}
    \end{equation}
    Here $D^b_{\mathcal A}$ denotes an extension
    closure inside $\mathcal D$.
    An equivalence with the abstract derived
    category $D^b(\mathcal A)$ is not asserted.
\end{definition}

\begin{definition}
    A strictly full triangulated subcategory is
    \emph{thick} if it is closed under direct
    summands in the ambient category.
    We write $\operatorname{thick}_{\mathcal D}(S)$
    for the smallest such subcategory containing $S$.
\end{definition}

\begin{prop}
    \label{prop:BBD-generated-heart}
    The category $D^b_{\mathcal A}$ is a thick
    triangulated subcategory of $\mathcal D$, and
    \[
        \bigl(D^{b,\leq0}_{\mathcal A},
        D^{b,\geq0}_{\mathcal A}\bigr)
    \]
    is a bounded $t$-structure with heart $\mathcal A$.
    Moreover,
    \begin{equation}
        \label{eq:extension-envelope-is-thick}
        D^b_{\mathcal A}
        =\operatorname{thick}_{\mathcal D}(\mathcal A).
    \end{equation}
\end{prop}

\begin{proof}
    Admissibility and extension closure, established in
    Theorem~\ref{thm:intro-abelianity} and
    Lemma~\ref{lem:positive-extension-closure},
    allow us to apply
    \cite[Proposition~1.3.13]{BBD}.
    This gives the triangulated subcategory and the
    stated $t$-structure with heart $\mathcal A$.
    It is bounded because every object of
    $D^b_{\mathcal A}$ is a finite extension of shifts
    of objects of $\mathcal A$.

    The theorem of Le--Chen~\cite{LeChen} then implies
    that $D^b_{\mathcal A}$ is idempotent complete.
    If an object of $D^b_{\mathcal A}$ has a direct
    summand in $\mathcal D$, its projector belongs to
    $D^b_{\mathcal A}$ by fullness and splits there.
    Uniqueness of the image of a split projector
    identifies that summand with an object of
    $D^b_{\mathcal A}$.
    Thus this subcategory is thick.
    Equality~\eqref{eq:extension-envelope-is-thick}
    follows: one inclusion uses thickness, and the other
    uses closure under shifts and extensions.
\end{proof}

Write
\[
    H^j_{\mathcal A}:D^b_{\mathcal A}
    \longrightarrow\mathcal A
\]
for the cohomology functors of this $t$-structure.
Every object $E\in D^b_{\mathcal A}$ has only
finitely many nonzero cohomology objects and
a truncation filtration whose factors are
$H^j_{\mathcal A}(E)[-j]$.
Our conventions are cohomological:
for $A\in\mathcal A$, the object $A[m]$ has
cohomology concentrated in degree $-m$.

\subsubsection{The induced stability condition}

Bridgeland's definition extends stability from
an abelian heart to its ambient triangulated
category by including all shifts of semistable
objects and their Harder--Narasimhan filtrations.

\begin{definition}[{\cite[Definition~1.1]{Bridgeland-stability}}]
    \label{def:stability-condition}
    A \emph{stability condition} on a triangulated
    category $\mathcal T$ consists of a homomorphism
    $W:K_0(\mathcal T)\to\mathbb C$ and strictly full
    additive subcategories $\mathcal P(t)\subset\mathcal T$,
    for $t\in\mathbb R$, satisfying:
    \begin{enumerate}
        \item For $0\neq F\in\mathcal P(t)$,
              $W(F)=m(F)e^{i\pi t}$ with $m(F)>0$.
        \item $\mathcal P(t+1)=\mathcal P(t)[1]$.
        \item If $t_1>t_2$ and $F_i\in\mathcal P(t_i)$,
              then $\operatorname{Hom}_{\mathcal T}(F_1,F_2)=0$.
        \item Every nonzero $E\in\mathcal T$ admits
              distinguished triangles
              \[
                  E_{j-1}\longrightarrow E_j
                  \longrightarrow F_j
                  \longrightarrow E_{j-1}[1],
                  \qquad j=1,\ldots,k,
              \]
              with $E_0=0$, $E_k=E$, and
              $0\neq F_j\in\mathcal P(t_j)$ for
              $t_1>\cdots>t_k$.
    \end{enumerate}
    The homomorphism $W$ is the \emph{central charge},
    and $\mathcal P$ is the \emph{slicing}.
    The nonzero objects of $\mathcal P(t)$ are called
    \emph{semistable of phase $t$}. The \emph{stable} objects are the simple objects of
  the categories $\mathcal P(t)$.
\end{definition}

To pass to this phase convention, put
\begin{equation}
    \label{eq:Bridgeland-phase-normalization}
    \widehat Z:=iZ,
    \qquad
    \varphi(A):=
    \frac{\phi_Z(A)+\pi/2}{\pi}\in(0,1],
\end{equation}
for $0\neq A\in\mathcal A$.
This increasing change of phase preserves semistability
and Harder--Narasimhan filtrations.

\begin{prop}
    \label{prop:stability-on-generated-category}
    The central charge $\widehat Z$ extends uniquely
    to $K_0(D^b_{\mathcal A})$ and, together with
    the heart $\mathcal A$, determines a stability
    condition on $D^b_{\mathcal A}$ in the sense
    of Definition~\ref{def:stability-condition}.
    For $0<t\leq1$, its slicing is
    \begin{equation}
        \label{eq:positive-brane-slicing}
        \mathcal P(t)
        =\{0\}\cup\bigl\{
        A\in\mathcal A \setminus\{0\}:
        A\text{ is }Z\text{-semistable},\
        \varphi(A)=t
        \bigr\},
    \end{equation}
    extended to all real phases by
    $\mathcal P(t+m)=\mathcal P(t)[m]$ for
    $m\in\mathbb Z$.
    In particular, $\mathcal P((0,1])=\mathcal A$,
    where, for an interval $I$, $\mathcal P(I)$ denotes
    the extension closure of the subcategories
    $\mathcal P(t)$ with $t\in I$.
\end{prop}

\begin{proof}
    Since $\mathcal A$ is a bounded heart, its inclusion
    induces an isomorphism
    \[
        K_0(\mathcal A)\xrightarrow{\sim}
        K_0(D^b_{\mathcal A}),
    \]
    whose inverse sends $[E]$ to
    $\sum_j(-1)^j[H^j_{\mathcal A}(E)]$.
    Thus the unique extension is
    \begin{equation}
        \label{eq:central-charge-extension}
        \widehat Z(E)
        =i\sum_j(-1)^j
        Z\bigl(H^j_{\mathcal A}(E)\bigr),
    \end{equation}
    with a finite sum by boundedness.
    The construction of
    \cite[Proposition~5.3]{Bridgeland-stability}
    gives the stated stability condition.
    Its filtrations are obtained by refining the
    cohomology filtrations by the Harder--Narasimhan
    filtrations in $\mathcal A$.
\end{proof}

\begin{remark}
    \label{rem:imaginary-special-semistable}
    Every nonzero object of $\mathcal A$ with
    normalized phase $1$ is semistable, since
    no subobject has larger phase.
    Hence
    \[
        \mathcal P(1)
        =\{0\}\cup
        \bigl\{
        A\in\mathcal A\setminus\{0\}:
        \operatorname{Re}Z(A)=0
        \bigr\}.
    \]
    By the proof of
    Lemma~\ref{lem:central-charge-stability},
    these are precisely the objects admitting
    defining representatives with only closed
    imaginary-special source components,
    together with the zero object.
\end{remark}

\subsubsection{Generation}
\label{subsubsec:stability-generation}
Recall that a collection $S$ \emph{split-generates}
$\mathcal D$ if
\[
    \operatorname{thick}_{\mathcal D}(S)
    =\mathcal D.
\]
The following criterion extends the stability
condition to all of $\mathcal D$.

\begin{proof}[Proof of Theorem~\ref{prop:intro-stability-from-generation}]
    By Proposition~\ref{prop:BBD-generated-heart},
    $D^b_{\mathcal A}$ is thick in $\mathcal D$.
    If it contains a split-generating collection $S$,
    then
    \[
        \mathcal D
        =\operatorname{thick}_{\mathcal D}(S)
        \subset D^b_{\mathcal A}\subset\mathcal D.
    \]
    Hence $D^b_{\mathcal A}=\mathcal D$.
    The remaining assertions follow from
    Propositions~\ref{prop:BBD-generated-heart}
    and~\ref{prop:stability-on-generated-category}.
\end{proof}

\begin{question}[Generation]
    \label{conj:intro-generation}
    Does the subcategory $\mathcal A$ split-generate
    $\mathcal D$?
\end{question}

\begin{remark}[The two-torus]
    \label{rem:elliptic-generation-compatibility}
    On an two-torus, one can choose two
    strictly positive embedded affine circles
    whose homology classes form an integral
    basis of $H_1(X;\mathbb Z)$.
    Equipped with trivial rank-one local
    systems and zero bounding cochains,
    the corresponding branes split-generate
    the embedded Fukaya category
    by~\cite[Proposition~6.1]{AbouzaidSmithFourTorus}.
    Applying this result to the immersed category
    $\mathcal D$ used here requires a comparison
    with the embedded framework.
    If the same circles split-generate
    $\mathcal D$, they lie in $\mathcal A$, and
    Theorem~\ref{prop:intro-stability-from-generation}
    applies.
    This comparison is not undertaken here.
\end{remark}

\subsubsection{Special representatives and stability}

The Thomas--Yau picture relates stability to
the existence of special Lagrangian representatives
\cite{ThomasYau,Joyce-stability}.
We now prove the correspondence stated in
Theorem~\ref{thm:intro-special-representatives}
by comparing the positive-brane categories
at different phases.
For $c\in\mathbb R$, put
\[
    \Omega_c=e^{i(\pi/2-c)}\Omega,
\]
and let $\mathcal A_c$ be the positive-brane
subcategory defined using $\Omega_c$.
In terms of the original grading, its defining
representatives have source components with
phases in $(c-\pi,c)$, together with closed
special components of phase $c$.
In particular,
\[
    \mathcal A_{\pi/2}=\mathcal A,
    \qquad
    \mathcal A_{c+m\pi}=\mathcal A_c[m],
    \qquad m\in\mathbb Z.
\]
Rotation leaves the density $\rho$ and the
metric $\widetilde g$ unchanged.
We use the common central charge $Z$ supplied
by Remark~\ref{rem:homological-compatibility},
computed using the original form $\Omega$.

\begin{proof}[Proof of
        Theorem~\ref{thm:intro-special-representatives}]
    The generation hypothesis and
    Propositions~\ref{prop:BBD-generated-heart}
    and~\ref{prop:stability-on-generated-category}
    give a stability condition on $\mathcal D$
    with heart $\mathcal A$ and central charge
    $\widehat Z=iZ$.

    We first prove the correspondence for
    $c\in(-\pi/2,\pi/2]$.
    A semistable object of geometric phase $c$
    lies in $\mathcal A$.
    An object admitting a special brane
    representative of phase $c$ also lies
    in $\mathcal A$.
    For $c=\pi/2$, the correspondence is
    Remark~\ref{rem:imaginary-special-semistable}.

    Suppose that $-\pi/2<c<\pi/2$, and write
    $\mathcal B=\mathcal A_c$.
    By the hypotheses and
    Proposition~\ref{prop:BBD-generated-heart},
    both $\mathcal A$ and $\mathcal B$ are
    bounded hearts in $\mathcal D$.
    Applying the Maslov-index estimate from
    Lemma~\ref{lem:positive-negative-vanishing}
    to their phase windows gives
    \[
        \operatorname{Hom}_{\mathcal D}
        (\mathcal A,\mathcal B[k])=0
        \quad(k<0),
        \qquad
        \operatorname{Hom}_{\mathcal D}
        (\mathcal B,\mathcal A[k])=0
        \quad(k<-1).
    \]
    Write $\mathcal D_{\mathcal A}^{\leq j}$
    and $\mathcal D_{\mathcal B}^{\leq j}$
    for the corresponding $t$-structure subcategories.
    Boundedness and the preceding vanishings imply
    \[
        \mathcal D_{\mathcal A}^{\leq0}
        \subset
        \mathcal D_{\mathcal B}^{\leq0}
        \subset
        \mathcal D_{\mathcal A}^{\leq1}.
    \]
    By \cite[Lemma~1.1.2]{Polishchuk-constant},
    \[
        \bigl(
        \mathcal A\cap\mathcal B[1],
        \mathcal A\cap\mathcal B
        \bigr)
    \]
    is a torsion pair in $\mathcal A$.
    Thus every $E\in\mathcal A$ admits a
    short exact sequence in $\mathcal A$,
    \[
        0\longrightarrow T_E
        \longrightarrow E
        \longrightarrow F_E
        \longrightarrow0,
    \]
    with
    \[
        T_E\in\mathcal A\cap\mathcal B[1],
        \qquad
        F_E\in\mathcal A\cap\mathcal B.
    \]
    Compatibility of the central charge
    allows us to compare the phase windows:
    \[
        T_E\neq0
        \ \Longrightarrow\
        \phi_Z(T_E)>c,
        \qquad
        F_E\neq0
        \ \Longrightarrow\
        \phi_Z(F_E)\leq c.
    \]
    If $E$ is semistable of phase $c$, the
    subobject $T_E$ must vanish.
    Hence $E\in\mathcal B$.
    Since
    \[
        \operatorname{Re}
        \bigl(e^{i(\pi/2-c)}Z(E)\bigr)=0,
    \]
    a defining representative in $\mathcal B$
    cannot contain a component strictly
    positive with respect to $\Omega_c$:
    each such component contributes positively
    to this real part.
    All its source components are therefore
    special of phase $c$.

    Conversely, suppose that $E$ admits a
    special brane representative of phase $c$.
    Then $E\in\mathcal A\cap\mathcal B$.
    Let $M\subset E$ be a nonzero subobject
    in $\mathcal A$, and apply the preceding
    decomposition to $M$.
    Since $T_M\in\mathcal B[1]$ and
    $E\in\mathcal B$,
    \[
        \operatorname{Hom}_{\mathcal D}(T_M,E)=0.
    \]
    The composite of the monomorphisms
    $T_M\hookrightarrow M\hookrightarrow E$
    therefore forces $T_M=0$.
    Thus $M\in\mathcal A\cap\mathcal B$, and
    $\phi_Z(M)\leq c$.
    This proves semistability.

    Finally, for arbitrary $c\in\mathbb R$,
    choose $m\in\mathbb Z$ such that
    $c-m\pi\in(-\pi/2,\pi/2]$.
    Shifting by $[-m]$ decreases both the
    geometric phase of a semistable object
    and the phase of a special representative
    by $m\pi$.
    The general assertion follows from the
    case already proved.
\end{proof}

In particular, every stable object admits a
special brane representative.
The converse need not hold: a direct sum of
two nonzero special branes of the same phase
has a special representative but is not stable.
Representatives are understood up to isomorphism
in $\mathcal D$ and may carry nontrivial
bounding cochains.

\begin{remark}[A sufficient generation hypothesis]
    \label{rem:special-generators}
    Suppose that $\mathcal D$ admits a
    split-generating collection of special
    Lagrangian branes $G_j$, of phases
    $\theta_j$.
    For every $c$, choose integers $m_j$ with
    \[
        c-\pi<\theta_j+m_j\pi\leq c.
    \]
    Then $G_j[m_j]\in\mathcal A_c$.
    Whenever the continuation and termination
    conjectures hold for $\Omega_c$,
    Proposition~\ref{prop:BBD-generated-heart}
    therefore gives
    \[
        D^b_{\mathcal A_c}
        =
        \operatorname{thick}_{\mathcal D}
        (\mathcal A_c)
        =
        \mathcal D.
    \]
    Thus special split-generators supply the
    generation hypothesis of
    Theorem~\ref{thm:intro-special-representatives}.
    The affine circles in
    Remark~\ref{rem:elliptic-generation-compatibility}
    provide such generators in the embedded
    elliptic-curve framework, subject to the
    comparison discussed there.
\end{remark}

\section{The two-torus}

\label{sec:elliptic-curve}

We use the flat Calabi--Yau structure on
$X=\mathbb C/\Gamma$ fixed in
Subsection~\ref{subsec:intro-elliptic-curve}.
Throughout this section, we assume that the vertical
direction is periodic and normalize its primitive
period to one:
\[
    \Gamma\cap\sqrt{-1}\mathbb R
    =\sqrt{-1}\mathbb Z.
\]

\subsection{Regular evolution}
	\label{subsec:elliptic-regular-evolution}
\subsubsection{The evolution of a vertical interval}

During a regular time interval, write
\[
    L_i=[\iota_i:M_i\longrightarrow X],
    \qquad i=0,1,
\]
for the current strictly positive boundary curves,
which remain fixed on this interval.
The closed imaginary-special components supplied
with the initial data are retained as closed
fragments of the web, as in
Definition~\ref{def:initial-web-configuration}.

We choose generic positive boundary representatives
whose disjoint-source union has only transverse
double intersections.
In particular, no third boundary branch passes
through an intersection of $L_0$ and $L_1$.
When choosing these representatives, we transport
the input brane structures and the morphism $f$
by the corresponding continuation comparisons.
These initial perturbations and comparisons are
included in the companion constructions
\cite{KapotaWebGeometry,KapotaWebFloer}.
The boundary conditions then remain fixed until
an absorption or an auxiliary deformation
changes them.

Each connected active front has interval source,
with one endpoint of each label.
A connected inactive fragment is either an interval
with both endpoints of the same label or a closed
vertical curve, possibly multiply covered.
The endpoint attachments specify points of the
boundary sources $M_i$, and therefore specify
which local boundary branch each endpoint follows.

To describe the length and motion of an active
interval, we use a lift to $\mathbb C$.
Write $q:\mathbb C\to X$ for the quotient map.
A lift of an immersion $F:[0,1]\to X$ means a map
$\widetilde F:[0,1]\to\mathbb C$ satisfying
$q\circ\widetilde F=F$.
It is uniquely determined once we choose the lift
of the label-zero endpoint.
In particular, the lift of the other endpoint is
then determined by the whole interval, including
its winding around the vertical fiber.
We continue this choice as the interval evolves.

\begin{figure}[htbp]
    \centering
    \includegraphics[width=0.4\textwidth]
    {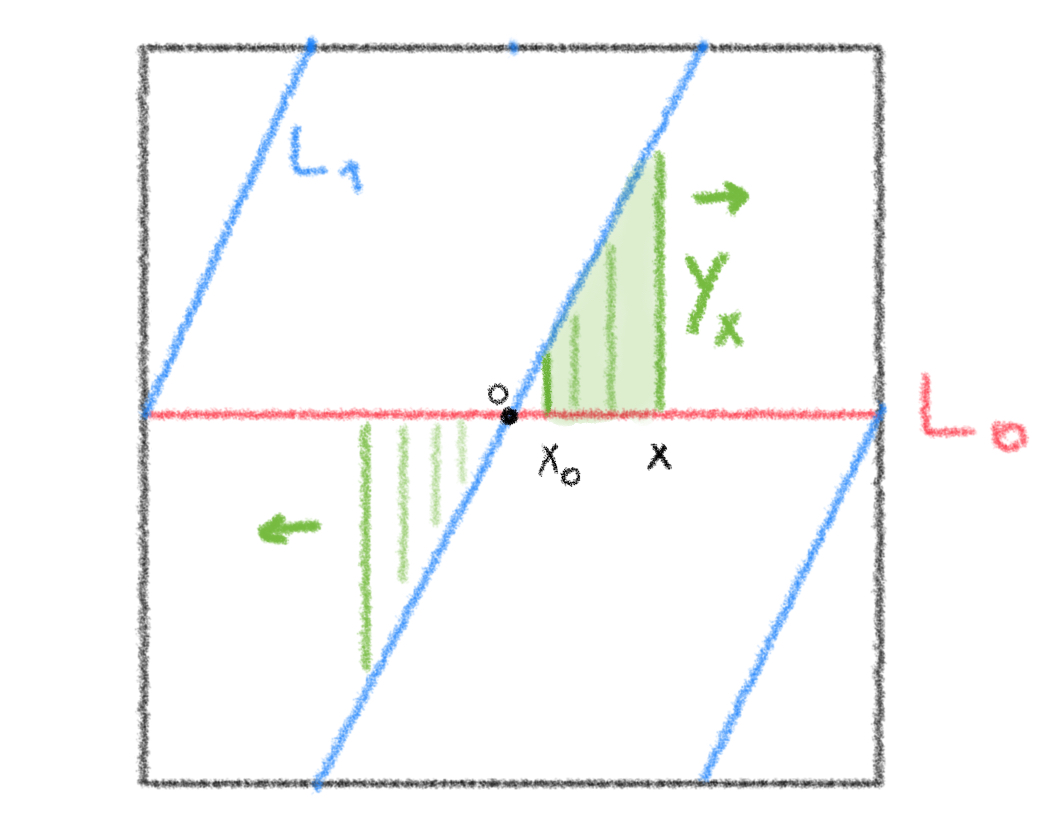}
    \caption{
        Regular evolution with boundary curves
        $L_0=\{y=0\}$ and $L_1=\{y=2x\}$ in
        $X$.
        Opposite sides of the square are identified.
        The green segments show successive positions
        of the active fronts.
        The shaded region is swept by the right-moving
        front between $x_0$ and $x$ and has area
        $A(x)=\int_{x_0}^{x}2u\,du=x^2-x_0^2$,
        equal to the elapsed evolution time. The horizontal speed is the reciprocal of the
        interval length, so shorter intervals move faster.
    }
    \label{fig:elliptic-regular-evolution}
\end{figure}

Strict positivity means that $dx$ is positive
along each oriented boundary curve.
Thus the two boundary branches followed by the
endpoints have local lifts given by graphs
\[
    y=y_i(x),
    \qquad x\in J,\quad i=0,1,
\]
where $J\subset\mathbb R$ is an interval.
Their vertical positions are fixed by the lift
of the active interval just described.
Set
\begin{equation}
    h(x):=y_1(x)-y_0(x).
    \label{eq:elliptic-h-function}
\end{equation}
The sign of $h$ records the vertical order of the
labels: label $1$ lies above label $0$ when $h>0$,
and below it when $h<0$.
For each $x$ with $h(x)\neq0$, define the immersed
vertical interval
\[
    Y_x=[F_x:[0,1]\longrightarrow X]
\]
by
\begin{equation}
    \widetilde F_x(s)
    =x+\sqrt{-1}\bigl(y_0(x)+s h(x)\bigr),
    \qquad
    F_x=q\circ\widetilde F_x.
    \label{eq:elliptic-interval-parametrization}
\end{equation}
Here $s=0$ and $s=1$ are the endpoints of labels
$0$ and $1$, attached to their respective boundary
branches.
Define its source length by
\[
    \ell(Y_x)
    :=\int_0^1
    \left|\partial_s\widetilde F_x\right|\,ds
    =|h(x)|.
\]
This length includes every traversal of the
vertical fiber, even when the image overlaps itself.

The following lemma describes the evolution in
terms of the signed area between the two boundary
graphs.
We use the Dirichlet energy $\mathcal E$ defined
in Subsection~\ref{subsec:Dirichlet-energy}.

\begin{lemma}
    \label{lem:elliptic-harmonic-flux}
    Suppose that $h$ does not vanish on $J$.
    The fundamental harmonic of $Y_x$, viewed as
    a function on its source $[0,1]$, is
    \begin{equation}
        \sigma_{Y_x}(s)=s.
        \label{eq:elliptic-fundamental-harmonic}
    \end{equation}
    Equivalently, in the lifted vertical coordinate
    $y=y_0(x)+s h(x)$, it is
    \[
        \sigma_{Y_x}(y)
        =\frac{y-y_0(x)}{h(x)}.
    \]
    Fix $x_0\in J$ and an initial time $\tau_0$,
    and set
    \[
        A(x):=\int_{x_0}^{x}h(u)\,du.
    \]
    Let $x(\tau)$ denote the horizontal coordinate
    of the evolving interval, with
    $x(\tau_0)=x_0$.
    Its regular evolution, with the prescribed
    endpoint attachments, is uniquely determined
    locally by
    \begin{equation}
        Y_\tau=Y_{x(\tau)},
        \qquad
        A(x(\tau))=\tau-\tau_0.
        \label{eq:elliptic-area-evolution}
    \end{equation}
    Its horizontal velocity and Dirichlet energy
    satisfy
    \begin{equation}
        \dot x(\tau)=\frac{1}{h(x(\tau))},
        \qquad
        |\dot x(\tau)|
        =\mathcal E(Y_\tau)
        =\frac{1}{\ell(Y_\tau)}.
        \label{eq:elliptic-speed-energy}
    \end{equation}
    Its relative flux is
    \[
        \operatorname{RelFlux}
        \bigl((Y_u)_{u\in[\tau_0,\tau]}\bigr)
        =-(\tau-\tau_0).
    \]
\end{lemma}

\begin{proof}
    The induced metric on the source of $Y_x$ is
    $h(x)^2\,ds^2$, and the Calabi--Yau density
    is one.
    Since this metric is constant along the
    interval, its harmonic functions are affine
    in $s$.
    The boundary values $0$ and $1$ therefore give
    \eqref{eq:elliptic-fundamental-harmonic}.

    Parametrize the evolving interval by
    $F_\tau=F_{x(\tau)}$, using
    \eqref{eq:elliptic-interval-parametrization}.
    In covering coordinates,
    \[
        \partial_s\widetilde F_\tau
        =\bigl(0,h(x(\tau))\bigr),
        \qquad
        \partial_\tau\widetilde F_\tau
        =\dot x(\tau)
        \bigl(1,y_0'(x(\tau))+s h'(x(\tau))\bigr).
    \]
    Since $\omega=dx\wedge dy$, we obtain
    \[
        F_\tau^*
        \bigl(\iota_{\partial_\tau F_\tau}\omega\bigr)
        =
        \omega\bigl(
        \partial_\tau\widetilde F_\tau,
        \partial_s\widetilde F_\tau
        \bigr)\,ds
        =
        h(x(\tau))\dot x(\tau)\,ds.
    \]
    The regular evolution equation
    \eqref{eq:fundamental-harmonic-evolution}
    requires this one-form on the source to equal
    $d\sigma_{Y_\tau}=ds$.
    Hence
    \[
        h(x(\tau))\dot x(\tau)=1.
    \]
    By the chain rule,
    \[
        \frac{d}{d\tau}A(x(\tau))
        =h(x(\tau))\dot x(\tau)
        =1.
    \]
    Integrating from $\tau_0$ to $\tau$, and using
    $A(x_0)=0$, gives
    \eqref{eq:elliptic-area-evolution}.
    Conversely, differentiating that identity
    recovers the evolution equation.
    Since $A'=h$ is nowhere zero, $A$ has a local
    inverse, which gives existence and uniqueness
    of $x(\tau)$.

    The velocity formula in
    \eqref{eq:elliptic-speed-energy} follows from
    $h\dot x=1$ and $\ell(Y_\tau)=|h(x(\tau))|$.
    For the energy, the source metric gives
    \[
        |d\sigma_{Y_\tau}|^2
        =\frac{1}{h(x(\tau))^2},
        \qquad
        d\operatorname{vol}_{g_{Y_\tau}}
        =|h(x(\tau))|\,ds.
    \]
    Therefore
    \[
        \mathcal E(Y_\tau)
        =\int_0^1
        \frac{|h(x(\tau))|}{h(x(\tau))^2}\,ds
        =\frac{1}{\ell(Y_\tau)}.
    \]
    Finally, the normalized deformation primitive
    is $\sigma_{Y_\tau}$, whose boundary values
    are $0$ and $1$.
    Definition~\ref{def:relative-lagrangian-flux}
    consequently gives
    \[
        \operatorname{RelFlux}
        \bigl((Y_u)_{u\in[\tau_0,\tau]}\bigr)
        =-\int_{\tau_0}^{\tau}1\,du
        =-(\tau-\tau_0).
    \]
\end{proof}

Thus an active interval moves right when its
label-one endpoint lies above its label-zero
endpoint, and left when it lies below.
In either case, its speed is the reciprocal
of its source length.

\subsubsection{Smoothing and regular brane transport}

Section~2 supplies coherent exact smoothings,
regular brane transport, and continuation through
isolated quadratic self-tangencies.
For the case of the two-torus, we need two
additional facts: the relevant families admit
generic exact perturbations, and
transport also
extends through triple-point passages.
The following companion result supplies these
facts.

\begin{prop}[Generic perturbations and brane continuation]
    \label{prop:elliptic-exact-continuation}
    Let $(C_s)_{s\in[0,1]}$ be a smooth exact
    family of graded immersed curves with fixed
    source, obtained either by coherently smoothing
    a regular web evolution or by the auxiliary
    deformation of
    Lemma~\ref{lem:elliptic-collapse-isolation}.
    Equip $C_0$ with a brane
    $(C_0,E_0,b_0)$ satisfying
    $\operatorname{val}(b_0)>0$.
    Assume that its intersections are transverse
    double points, apart from fixed clean
    intersections between stationary closed
    vertical components.
    Then the following hold.
    \begin{enumerate}
        \item \emph{Generic representatives.}
              The family admits an arbitrarily
              small exact perturbation, unchanged
              near $s=0$, retaining the stationary
              closed vertical components.
              Away from their retained clean loci,
              the perturbed family has only
              transverse double intersections,
              except at finitely many quadratic
              self-tangencies and isolated
              triple-point passages.

        \item \emph{Transport through the exceptional
                  passages.}
              The regular and self-tangency
              transport results of Section~2
              apply to these perturbed exact
              families, with the retained clean
              intersections treated in a compatible
              clean Floer model.
              Transport also extends across the
              isolated triple-point passages.
              At each exceptional passage, assume
              that the incoming bounding cochains
              have a uniform positive valuation
              bound.
              The continuation preserves the
              represented object in $\mathcal D$
              and gives outgoing bounding cochains
              of strictly positive valuation.
    \end{enumerate}
    Throughout, transport uses the source
    Hamiltonian of the chosen perturbed family,
    as in Remark~\ref{rem:exact-brane-transport}.
\end{prop}

The proof is deferred to
\cite{KapotaWebGeometry,KapotaWebFloer}.
Both the genericity argument and the continuation
through triple-point passages build on
\cite{PalmerWoodwardMaslov}; for genericity,
see in particular their Theorem~2.6.
Continuation at obstruction walls is treated in
Section~\ref{subsubsec:elliptic-initialization-continuation}.

\subsection{Geometric degenerations}

The previous subsection describes the regular
motion of a vertical interval.
We now describe the geometric degenerations that
this motion can reach, specializing the discussion
in Subsection~\ref{subsec:geometric-degenerations}
to dimension one.

\subsubsection{Boundary collisions}
The first possibility is a boundary degeneration.
In dimension one, this occurs when two endpoints
with the same label meet along the corresponding
boundary curve, exhausting the residue interval
between them.

\begin{definition}[Boundary collision]
    \label{def:elliptic-boundary-collision}
    Use the incoming and outgoing notation of
    Subsection~\ref{subsec:local-continuation-conventions}.
    Let $P^-$ be the disjoint union of the
    participating source intervals, which may
    be active or inactive, and let
    \[
        e_1,e_2\in\partial_iP^-
    \]
    be distinct endpoints with the same label.
    Each participating source interval is
    included only once.

    A \emph{boundary collision} along $L_i$
    occurs at $\tau_*$ when a residue interval
    with endpoints
    $\beta_{i,\tau}^-(e_1)$ and
    $\beta_{i,\tau}^-(e_2)$ contracts to a point
    $q_i\in M_i$.
    The associated boundary reconnection
    removes $e_1,e_2$ and joins the corresponding
    ends of the source intervals.
    The remaining endpoint labels and boundary
    lifts are retained.
    In the image, the vertical pieces concatenate.

    When the incident intervals have positive
limiting lengths, this is the $n=1$, $k=0$
case of Definition~\ref{def:boundary-surgery},
with vanishing sphere $C^-=\{e_1,e_2\}$.

    We also use this boundary reconnection
    as part of a simultaneous collapse, or
    for coincident endpoints produced by
    another operation when no intervening
    residue remains.
\end{definition}

\begin{figure}[htbp]
    \centering
    \includegraphics[width=0.3\textwidth]
    {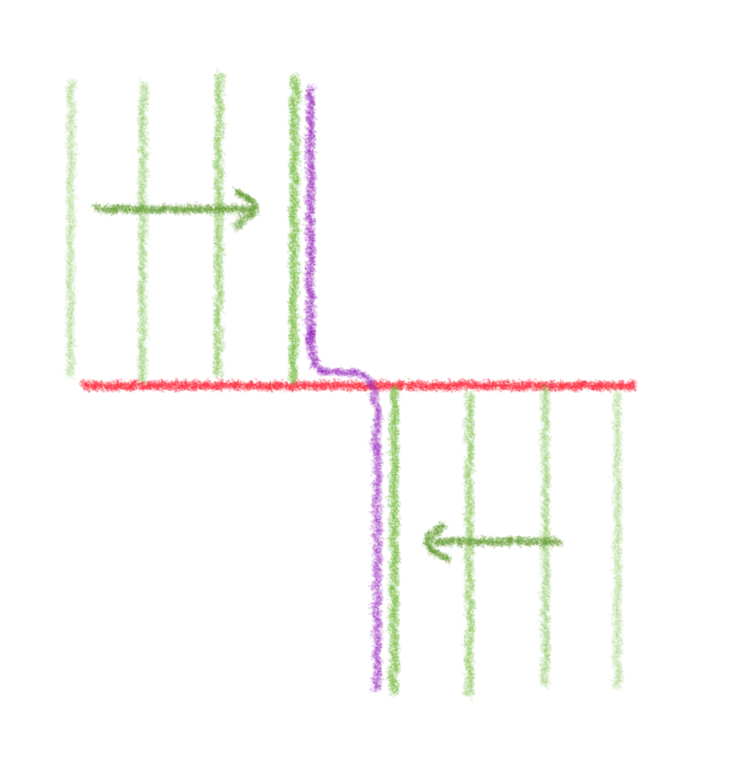}
    \caption{
        A boundary collision along $L_0$.
        The green segments show successive positions
        of two active fronts whose label-zero endpoints
        approach along the red boundary branch.
        The purple arc is a local smoothing of the
        web cone at the latest depicted time.
        As the intervening residue interval shrinks
        to a point, the fronts concatenate and the
        purple arc straightens.
        No surgery of the smoothed web cone is required.
    }
    \label{fig:elliptic-boundary-collision}
\end{figure}

\begin{remark}
    \label{rem:elliptic-boundary-collision-fragments}
    Joining two active intervals along $L_i$
    produces an inactive interval: its two
    remaining endpoints both have label $1-i$.
    If a sequence of reconnections removes all
    endpoints of a component, the result is a
    closed vertical fragment.
    It inherits the unshifted imaginary-special
    grading of the fronts and receives label $1$, hence belongs to $K^*$ by Definition~\ref{def:terminal-branes}.
\end{remark}

\begin{remark}
    \label{rem:elliptic-boundary-collision-lifts}
    A boundary collision is defined on the
    boundary source $M_i$.
    Equality of endpoint images in $X$ is
    insufficient: the attachment points must
    meet on the same labelled boundary source,
    with the intervening residue shrinking
    to that point.
\end{remark}

At a boundary collision with positive
limiting lengths, concatenation preserves
the smooth source of the web cone.
After rounding the approaching corners,
the collision is the straightening of
a smooth arc; no surgery is required.
The following companion result supplies the
smoothing and brane transfer across a boundary
collision; its proof is deferred to
\cite{KapotaWebGeometry,KapotaWebFloer}.

\begin{prop}[Continuation at a boundary collision]
    \label{prop:elliptic-boundary-continuation}
    Let a special Lagrangian web for $f$ approach a
    boundary collision at time $\tau_*$.
    Concatenate the incident intervals as in
    Definition~\ref{def:elliptic-boundary-collision}.
    Assume that each concatenated component
    contains an incoming interval of positive
    limiting length.

    The incoming and outgoing web cones admit
    smoothings that form a smooth exact family
    of immersed curves across $\tau_*$. The swept real and imaginary integrals
    have matching one-sided limits.

    If the incoming bounding cochains have
    a uniform positive valuation bound,
    the comparison gives outgoing
    branes with strictly positive valuation,
    representing the same object in $\mathcal D$.

If coefficients reach valuation zero at
$\tau_*$, the continuation supplies provisional
smoothed cone data on the concatenated
configuration, with cochain of nonnegative
valuation, for the subsequent obstruction-wall
resolutions.
The construction is compatible with brane
transfer across the complete event in the sense
of Definition~\ref{def:brane-surgery-transfer}.
\end{prop}

\subsubsection{Paired collapses}

Another geometric degeneration concerns intersections of
the two boundary curves $L_0$ and $L_1$.
A transverse intersection of strictly positive
curves has Maslov index zero or one.
At index zero, a new webbing may be initiated;
we discuss this together with obstruction walls
in Section~\ref{subsubsec:elliptic-geometric-resolutions}.
At index one, an active interval may shrink to
the intersection, which we call a \emph{collapse}. A \emph{paired collapse} is the one-dimensional
terminating webbing of
Remark~\ref{rem:maslov-n-webbing}.
It belongs to the endpoint case $m=n$,
treated separately in
Subsection~\ref{subsec:geometric-degenerations}.

\begin{remark}
    \label{rem:elliptic-maslov-slope}
    At a transverse intersection of two positive
    boundary branches $y=y_i(x)$, their tangent
    angles are
    \[
        \theta_i=\arctan y_i'(x_*)
        \in(-\pi/2,\pi/2).
    \]
    With their natural positive gradings, the
    Maslov-index convention of
    \cite[Definition~3.12]{SY-webbing} gives index
    zero when $\theta_1>\theta_0$ and index one
    when $\theta_1<\theta_0$.
    Since $\arctan$ is strictly increasing, writing
    $h=y_1-y_0$ gives
    \begin{equation}
        m(q_0,q_1;L_0,L_1)
        =
        \begin{cases}
            0, & h'(x_*)>0, \\
            1, & h'(x_*)<0.
        \end{cases}
        \label{eq:elliptic-maslov-slope}
    \end{equation}
\end{remark}

\begin{definition}[Collapse]
    \label{def:elliptic-collapse}
    Fix source points $q_i\in M_i$ whose boundary
    branches intersect transversely with
    \[
        \iota_0(q_0)=\iota_1(q_1)=q,
        \qquad
        m(q_0,q_1;L_0,L_1)=1.
    \]
    Let
    \[
        Y_\tau=[F_\tau:P\longrightarrow X]
    \]
    be a connected active interval, and let $e_i$
    be its endpoint of label $i$.
    We say that $Y_\tau$ undergoes a
    \emph{collapse} at $(q_0,q_1)$ at time $\tau_*$
    if
    \[
        \lim_{\tau\nearrow\tau_*}\ell(Y_\tau)=0,
        \qquad
        \lim_{\tau\nearrow\tau_*}
        \beta_{i,\tau}(e_i)=q_i,
        \quad i=0,1.
    \]
    Thus the whole source interval shrinks to
    the intersection, with each endpoint
    approaching its specified boundary source point.
\end{definition}

\begin{remark}
    \label{rem:elliptic-componentwise-collapse}
    A local webbing at a Maslov-index-zero
    intersection has source $S^0\times[0,1]$,
    so it starts with two active intervals.
    Subsequent reconnections may change their
    attachments, and the intervals born together
    need not collapse together.
    We therefore define collapse for an
    individual active interval.

    Lemma~\ref{lem:elliptic-paired-collapse}
    shows that a spatially isolated collapse
    nevertheless involves two intervals,
    approaching from opposite horizontal sides.
    These intervals need not have been born
    in the same webbing.
    Under the hypotheses of
    Lemma~\ref{lem:elliptic-collapse-isolation},
    simultaneous arrivals of other active fronts
    can be separated from the collapse by a
    Hamiltonian change of the boundary
    representatives.
\end{remark}

\begin{figure}[htbp]
    \centering
    \includegraphics[width=0.8\textwidth]
    {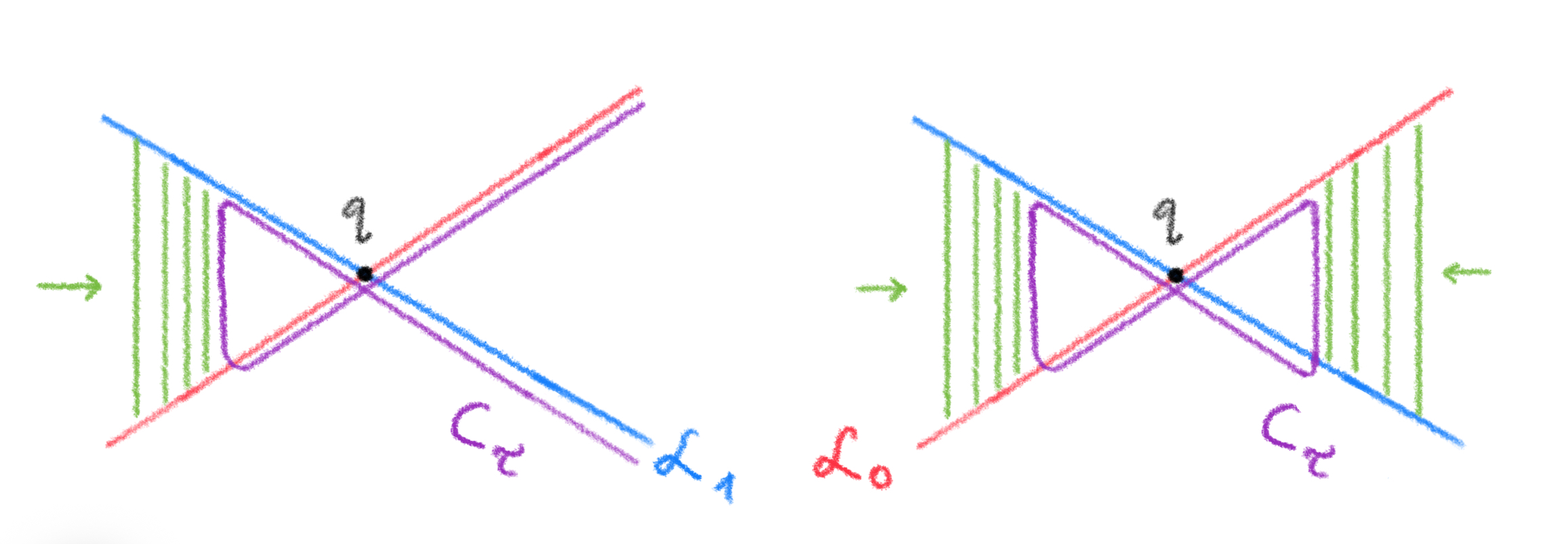}
    \caption{
        Unpaired and paired collapse at a
        Maslov-index-one intersection $q$.
        The green segments show successive positions
        of the active fronts, and the purple curves
        show local smoothings of the web cone. The smoothed web cone is drawn slightly displaced
        from the underlying Lagrangians for clarity.
        Left: the hypothetical unpaired collapse
        considered in the proof of
        Lemma~\ref{lem:elliptic-paired-collapse},
        producing a shrinking teardrop at $q$.
        Right: two fronts approach from opposite
        sides and collapse together; their smoothing
        forms a shrinking figure-eight component.
    }
    \label{fig:elliptic-paired-collapse}
\end{figure}

\begin{lemma}[Collapse isolated from other active fronts is paired]
    \label{lem:elliptic-paired-collapse}
    Let a special Lagrangian web for $f$ evolve regularly on
    $(\tau_*-\epsilon,\tau_*)$, with fixed strictly
    positive boundary conditions
    \[
        L_i=[\iota_i:M_i\longrightarrow X],
        \qquad i=0,1.
    \]
    Let $q$ be a transverse Maslov-index-one
    intersection of $L_0$ and $L_1$, with specified
    source points $q_i\in M_i$:
    \[
        \iota_0(q_0)=\iota_1(q_1)=q,
        \qquad
        m(q_0,q_1;L_0,L_1)=1.
    \]
    Suppose that an active interval collapses
    at this intersection at time $\tau_*$.
    Assume that there is a fixed neighborhood
    $U$ of $q$ such that, for all $\tau<\tau_*$
    sufficiently close to $\tau_*$, the web
    inside $U$ consists only of active intervals
    collapsing at this intersection, residue
    arcs on the two specified boundary branches,
    and interior portions of stationary inactive
    fragments.
    No endpoint of an inactive
    fragment lies in $U$.
    Assume also that the smoothed cone branes
    \[
        \mathcal C_\tau=(C_\tau,E_\tau,b_\tau)
    \]
    in the brane decoration can be chosen so that,
    inside $U$, the corners are rounded by embedded
    arcs introducing no additional intersections.
    The roundings are supported in neighborhoods
    of the active endpoints that avoid $q$ and
    shrink to $q$ as $\tau\nearrow\tau_*$.
    Outside these neighborhoods, the curve
    inside $U$ is unchanged.

    Then exactly two active intervals collapse
    at this intersection at time $\tau_*$,
    one approaching from each horizontal side
    of $q$.
    The residue arc joining their endpoints
    on each boundary source $M_i$ also shrinks
    to $q_i$.
    Thus the two intervals and the intervening
    residue arcs shrink together, forming
    the paired collapse model.
\end{lemma}

\begin{proof}
    Put $x(q)=0$ and write the chosen boundary
    branches as $y_i(x)$. Set
    \[
        h(x)=y_1(x)-y_0(x),
        \qquad
        A(x)=-\int_0^x h(u)\,du.
    \]
    By Remark~\ref{rem:elliptic-maslov-slope}, $h'(0)<0$. Since $h(0)=0$, the function
    $A$ is positive and strictly monotone on each
    side of zero, sufficiently close to zero.
    Along an active interval,
    \[
        \frac{d}{d\tau}A(x(\tau))
        =-h(x(\tau))\dot x(\tau)=-1.
    \]
    Thus every interval contracting at $\tau_*$
    satisfies
    \[
        A(x(\tau))=\tau_*-\tau.
    \]
    This determines its position uniquely on
    each side. Two intervals at that position
    would have the same labelled boundary lifts,
    contrary to source admissibility.
    Hence there is at most one contracting
    interval on either side.

    Suppose that only one interval contracts.
    After shrinking $U$, every stationary
    fragment meeting $U$ lies on the fiber
    $x=0$.
    These fragments are disjoint from the
    contracting interval, and its corner
    roundings can also be chosen disjoint
    from them.
    The contracting interval and the two
    residue arcs leading to $q$ therefore
    still bound an embedded disk with one
    corner at $q$.
    After rounding the other corners, its
    area $a_\tau$ satisfies
    \[
        a_\tau>0,
        \qquad
        a_\tau\longrightarrow0.
    \]
    In the cone, the ordered intersection from
    $R_0[1]$ to $R_1$ has degree
    \[
        \deg_{C_\tau}(q)
        =\deg_{L_0,L_1}(q)+1=2.
    \]
    The disk therefore contributes to the
    $q$-coefficient of $\mu^0$.
    Its boundary is the only source path inside
    $U$ joining the two preimages of $q$ without
    another branch switch.
    The Riemann mapping theorem gives the simple
    disk uniquely up to reparametrization.
    A repeated traversal would require additional
    branch switches at $q$, so it is not another
    disk with one output and no inputs.
    The local contribution is consequently
    \[
        u_\tau T^{a_\tau},
        \qquad
        \operatorname{val}(u_\tau)=0,
    \]
    where the sign and local-system transport
    give the unit $u_\tau$.

    Every other disk contributing to this coefficient
    must leave $U$.
    Choose a fixed coordinate disk $V$ with
    $q\in V$ and $\overline V\subset U$.
    Choose $\tau_0<\tau_*$ sufficiently close to
    $\tau_*$ that the contracting interval and its
    corner roundings lie inside $V$ whenever
    $\tau\in(\tau_0,\tau_*)$.
    On the annulus $U\setminus\overline V$,
    the cone consists of fixed residue arcs
    and stationary vertical pieces.
    Their images divide the annulus into
    finitely many connected regions
    $A_1,\ldots,A_m$, each of positive area.

    A holomorphic disk with output $q$ that leaves
    $U$ must meet one of these regions.
    Its boundary lies on the cone and hence avoids
    the interior of every $A_j$.
    Within each $A_j$, the image of the compactified
    disk is closed by compactness and open by the
    open mapping theorem.
    Consequently, if the disk meets $A_j$, it covers
    all of $A_j$.
    Its area, counted with multiplicity, is therefore
    at least
    \[
        \eta:=\min_{1\leq j\leq m}\operatorname{Area}(A_j)>0.
    \]
    Since the regions are fixed, this bound is
    independent of $\tau\in(\tau_0,\tau_*)$.
    We now consider the contributions with
    bounding-cochain inputs.
    Write
    \[
        b_\tau
        =b_{\tau,\mathrm{ch}}+b_{\tau,\mathrm{sw}}
    \]
    for the ordinary source-cochain and
    branch-switching parts.
    We distinguish generators by their source
    points and use fixed compatible Floer data
    for the stationary intersections.

    The two ordered generators between the
    residue branches at $q$ have degrees
    $2$ and $-1$, so neither occurs in
    $b_\tau$.
    Thus a polygon contained in $U$ with a
    branch-switching input must switch onto
    a stationary fragment.
    It must also switch back to a residue:
    the stationary fragments have no
    attachments inside $U$ and do not meet
    the contracting active interval there.

    Consider a consecutive sequence of these
    switches, starting on a residue, passing
    through stationary fragments, and ending
    on a residue.
    Every branch in this sequence has constant
    height either $0$ or $1$.
    There are at least two switches, so the
    height cannot strictly increase at every
    switch.
    At least one input generator $p$ therefore
    satisfies
    \[
        \kappa(p_+)\leq\kappa(p_-).
    \]
    By the transport formula
    \eqref{eq:web-cone-valuation-transport},
    its coefficient has nondecreasing valuation:
    \[
        \frac{d}{d\tau}
        \operatorname{val}(b_{\tau,p})
        =
        -\bigl(\kappa(p_+)-\kappa(p_-)\bigr)
        \geq0.
    \]
    At the fixed earlier time $\tau_0$, all
    nonzero coefficients have strictly positive
    valuation.
    Since there are finitely many relevant
    generators, there is a constant $\nu>0$
    such that every nonzero local contribution
    with a branch-switching input has valuation
    at least $\nu$.
    Here we also use the nonnegative areas
    and nonnegative valuations of all other
    input coefficients.

    A polygon contained in $U$ with only
    ordinary source-cochain inputs follows
    the same source path as the teardrop.
    Its area is therefore $a_\tau$.
    Each contribution with at least one such
    input has valuation strictly greater than
    $a_\tau$, since
    \[
        \operatorname{val}
        (b_{\tau,\mathrm{ch}})>0.
    \]
    Finally, every polygon leaving $U$ has
    area at least $\eta$.
    Consequently,
    \[
        [\mu^0_{b_\tau}]_q
        =
        T^{a_\tau}(u_\tau+v_\tau)+R_\tau,
    \]
    where
    \[
        \operatorname{val}(u_\tau)=0,
        \qquad
        \operatorname{val}(v_\tau)>0,
        \qquad
        \operatorname{val}(R_\tau)
        \geq\min\{\eta,\nu\}.
    \]
    Here $v_\tau$ collects the local
    source-cochain insertions, while $R_\tau$
    collects the contributions using stationary
    fragments or leaving $U$.
    For each $\tau<\tau_*$, these sums converge
    because $\operatorname{val}(b_\tau)>0$.

    For $\tau$ sufficiently close to $\tau_*$,
    \[
        a_\tau<\min\{\eta,\nu\}.
    \]
    The displayed curvature coefficient then
    has valuation $a_\tau$ and is nonzero,
    contradicting the Maurer--Cartan equation.

    A second interval must therefore contract
    at the same source pair.
    By uniqueness on each side, it approaches
    from the opposite side.
    The residue arcs between the two intervals
    shrink to $q_i$ on each boundary source,
    giving the paired collapse.
\end{proof}

\begin{remark}[Continuation through a paired collapse]
    \label{rem:elliptic-geometric-continuation}
    Suppose that the shrinking figure-eight
    component obtained by smoothing a paired
    collapse is disjoint from the remaining
    cone at incoming regular times sufficiently
    close to the collapse.

    There are then no branch-switching generators
    between this component and the remaining cone.
    The bounding cochain splits accordingly,
    so the contracting component carries a
    brane structure that is a direct summand
    of the incoming cone brane.

    Since it is Hamiltonian displaceable, this
    summand represents the zero object in
    $\mathcal D$; see
    \cite[Section~3.7]{Joyce-stability}.
    Deleting it therefore preserves the
    represented object and leaves the remaining
    brane data unchanged.
\end{remark}

\subsubsection{Isolation of collapses}

\begin{figure}[htbp]
    \centering
    \includegraphics[width=\textwidth]{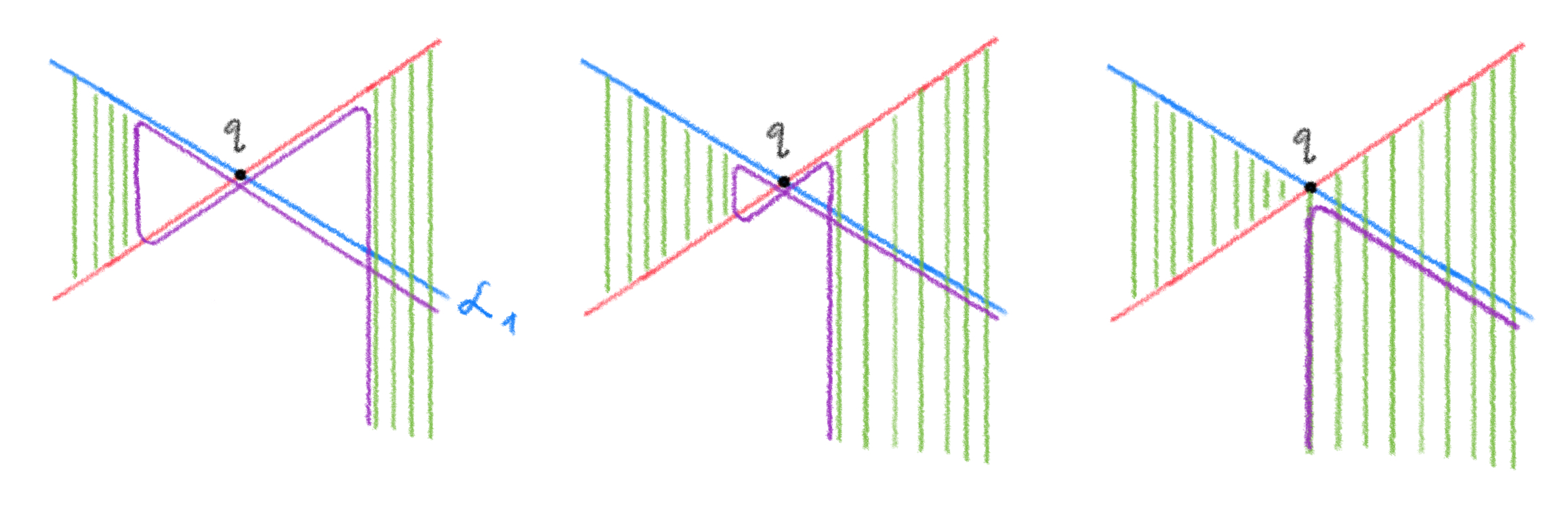}
    \caption{
        Simultaneous arrival and collapse at $q$,
        with time increasing from left to right.
        The boundary curves $L_0$ and $L_1$ are
        shown in red and blue, and the green
        segments indicate successive front positions.
        An active front contracts to $q$ from the
        left, while another active front approaches
        from the right with positive limiting length.
        They reach $q$ at the same time.
        The figure-eight portion of the smoothed
        web cone shrinks to $q$, and concatenation
        at the simultaneous boundary collision
        leaves an inactive fragment.
        The smoothed web cone is shown in purple,
        slightly displaced for clarity.
    }
    \label{fig:elliptic-simultaneous-collapse}
\end{figure}

We now explain how a small Hamiltonian change
of the boundary representatives can isolate
a collapse from other fronts and inactive fragments.
Figure~\ref{fig:elliptic-simultaneous-collapse}
illustrates a simultaneous arrival and collapse
that the perturbation will separate.

We first distinguish collapses meeting an inactive
endpoint, which are continued by boundary
concatenation.

\begin{remark}
    \label{rem:elliptic-collapse-inactive-endpoints}
    Suppose that an active interval collapses at
    $q$, and its label-$i$ endpoint approaches
    the attachment point $q_i\in M_i$ of an
    inactive interval.
    The intervening residue interval shrinks,
    giving a simultaneous boundary collision
    along $L_i$.
    Here $q_i$ lies on the boundary of the
    residue, so this is not a crossing of
    two residue interiors.

    At the inactive attachment, the singular
    cone turns from the residue into the
    stationary fragment.
    Rounding this corner opens the apparent
    teardrop, so the obstruction argument of
    Lemma~\ref{lem:elliptic-paired-collapse}
    does not apply.

    We concatenate the incident intervals
    before discarding the zero-length active
    piece.
    The resulting interval has positive
    limiting length and remaining endpoint
    labels $i$ and $1-i$, so it is active.
    Its smoothing and brane continuation are
    supplied by
    Proposition~\ref{prop:elliptic-boundary-continuation}.
\end{remark}

\begin{lemma}[Isolation of collapses]
    \label{lem:elliptic-collapse-isolation}
    Let a special Lagrangian web for $f$ evolve regularly
    up to time $\tau_*$, with fixed strictly
    positive boundary conditions.
    Suppose that active intervals collapse at
    a transverse Maslov-index-one intersection $q$.

    Assume that the only boundary branches through
    $q$ are the two residue branches at this
    intersection, and that every other active
    interval whose limiting image contains $q$
    has positive limiting source length.
    Stationary inactive fragments may pass
    through $q$, but none has an endpoint there.

    At a regular time $\tau_0<\tau_*$ sufficiently
    close to $\tau_*$, an arbitrarily small
    Hamiltonian change of the boundary
    representatives, with the brane transported
    accordingly, gives the following alternative:
    either a boundary collision or obstruction
    wall occurs before $\tau_*$, or the collapse
    at the moved intersection is spatially
    isolated and hence paired.

    At $\tau_0$, the change fixes all vertical pieces
    and their attachments and leaves the swept
    traces unchanged.
    In particular, the closed components and
    the number of labelled endpoints are unchanged.
\end{lemma}

\begin{proof}
    Put $q=(0,0)$ in local coordinates.
    Choose a small neighborhood $U$ of $q$
    containing no active interval or boundary
    attachment at $\tau_0$, and meeting the boundary
    curves only along the two residue branches
    through $q$.
    Shrink $U$ so that every stationary piece
    meeting it lies on $x=0$, and so that it
    avoids every limiting endpoint distinct
    from $q$.

    Choose a smooth Hamiltonian $H$ supported
    in $U$ and equal to $y$ near $q$.
    With the convention
    $\iota_{X_H}\omega=dH$, its time-$s$ flow
    moves $q$ to
    \[
        q_s=(s,0).
    \]
    Apply this flow to both boundary curves.
    For sufficiently small $s$, they remain
    strictly positive.
    On the web cone, apply the flow only to
    the residue arcs, leaving all vertical
    pieces fixed.
    Since the support avoids their attachments,
    this gives a smooth deformation of the
    smoothed cone.
    Write $F_s$ for its immersion.
    The function
    \[
        K_s=
        \begin{cases}
            H\circ F_s, & \text{on the residue portions}, \\
            0,          & \text{on the fixed pieces}
        \end{cases}
    \]
    is smooth on the cone source, since both
    expressions vanish near the attachments.
    Moreover,
    \[
        F_s^*(\iota_{\partial_sF_s}\omega)=dK_s,
    \]
    so the deformation is exact.
    It fixes the vertical pieces, their
    attachments, and the swept traces.

    Apply
    Proposition~\ref{prop:elliptic-exact-continuation}
    to transport the brane along this exact
    deformation.
    This includes the comparison needed when
    stationary pieces remain fixed.
    Choose the deformation sufficiently small
    that the transported bounding cochain
    retains strictly positive valuation at $\tau_0$.

    For each contracting interval, its remaining
    collapse time is the area bounded by that
    interval and the two residue arcs leading
    to $q$, by
    Lemma~\ref{lem:elliptic-harmonic-flux}.
    The Hamiltonian flow fixes the interval
    and preserves this area.
    Its collapse time therefore remains
    $\tau_*$, as long as no earlier operation
    changes its attachments.
    The geometric effect of this perturbation is
    illustrated in
    Figure~\ref{fig:elliptic-collapse-perturbation}.

    \begin{figure}[htbp]
        \centering
        \includegraphics[width=0.4\textwidth]{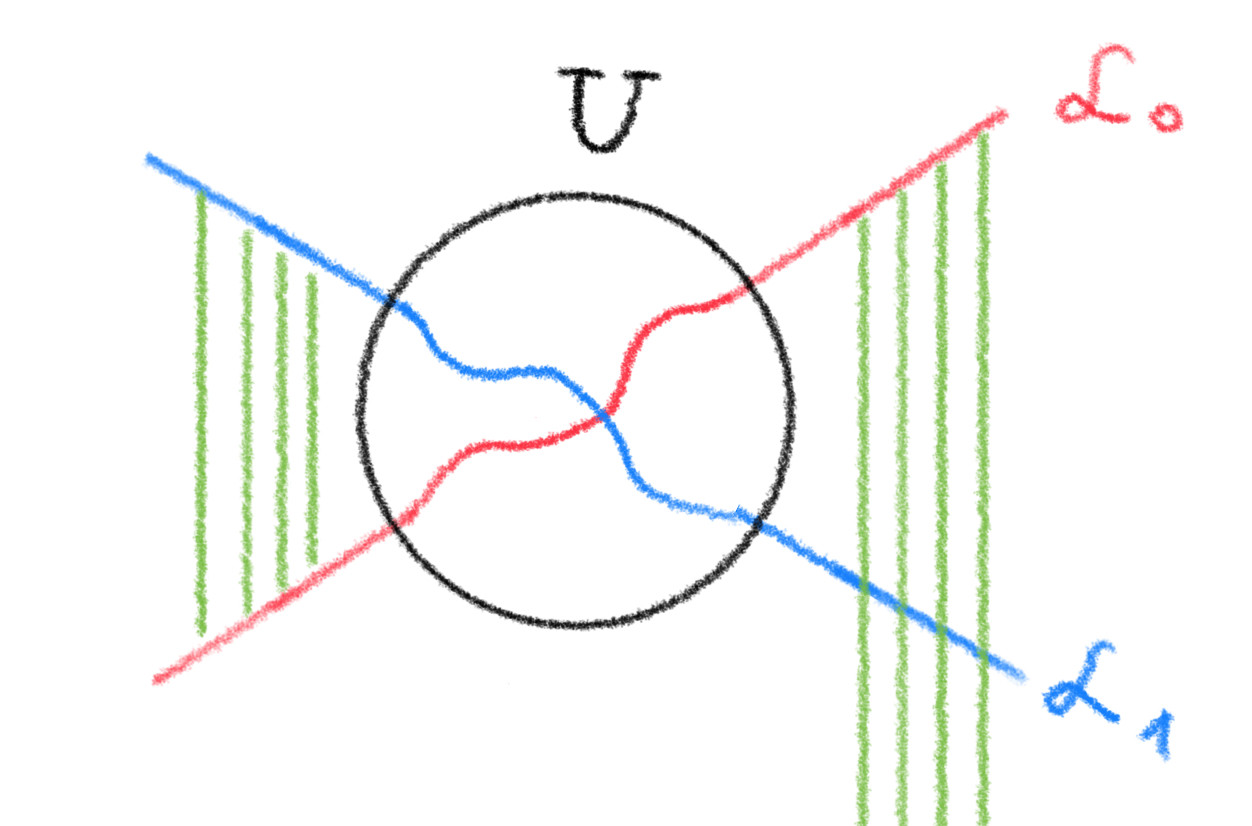}
        \caption{
            A Hamiltonian perturbation supported in $U$
            moves the intersection of $L_0$ (red) and
            $L_1$ (blue) to the right, while fixing the
            fronts at the chosen time $\tau_0$.
            The areas of the sectors cut out by the
            two curves in $U$ are preserved.
            Consequently, the remaining collapse time
            is unchanged, provided no earlier event
            occurs.
            The incoming front from the right has its
            other endpoint outside $U$ and retains
            positive length at arrival.
            Its arrival time is therefore not determined
            by these local sector areas: moving the
            intersection to the right makes it arrive
            earlier, with the first-order change
            governed by its limiting length.
            Green segments indicate successive front
            positions outside $U$.
        }
        \label{fig:elliptic-collapse-perturbation}
    \end{figure}
    Consider an additional active interval $Q$,
    with positive limiting source length $L_Q$.
    Write its lifted boundary graphs after the
    perturbation as $y_{0,s}(x),y_{1,s}(x)$, and set
    \[
        h_{Q,s}(x):=y_{1,s}(x)-y_{0,s}(x).
    \]
    Let $b_Q$ be the horizontal coordinate of
    $Q$ at time $\tau_0$.
    For each fixed perturbation parameter $s$,
    write $T_Q(s)$ for the evolution time at
    which $Q$ reaches the fiber through $q_s$,
    following the same boundary branches.
    Since $q_s=(s,0)$, its target horizontal
    coordinate is $s$.
    The regular evolution equation gives
    \begin{equation}
        \label{eq:elliptic-perturbed-arrival-time}
        T_Q(s)-\tau_0
        =\int_{b_Q}^{s}h_{Q,s}(x)\,dx,
        \qquad T_Q(0)=\tau_*.
    \end{equation}
    Use the same notation $H$ for its lift to
    the covering plane.
    At a fixed horizontal coordinate, the
    Hamiltonian deformation of a boundary graph
    satisfies
    \[
        \left.\frac{\partial y_{i,s}(x)}
        {\partial s}\right|_{s=0}
        =-\frac{d}{dx}H(x,y_{i,0}(x)).
    \]
    The Hamiltonian vanishes at the endpoints
    over $b_Q$ and at both limiting endpoints
    over $x=0$: those outside $U$ lie outside
    its support, and $H(q)=0$.
    Thus, for $i=0,1$,
    \[
        H(b_Q,y_{i,0}(b_Q))
        =H(0,y_{i,0}(0))=0.
    \]
    Differentiating
    \eqref{eq:elliptic-perturbed-arrival-time}
    and using the graph deformation formula gives
    \begin{align}
        T_Q'(0)
         & =
        h_{Q,0}(0)
        +\int_{b_Q}^{0}
        \left.\frac{\partial h_{Q,s}(x)}
        {\partial s}\right|_{s=0}\,dx \\
         & =
        h_{Q,0}(0)
        -\int_{b_Q}^{0}\frac{d}{dx}
        \bigl(
        H(x,y_{1,0}(x))-H(x,y_{0,0}(x))
        \bigr)\,dx                    \\
         & =
        h_{Q,0}(0)
        -\bigl[
            H(x,y_{1,0}(x))-H(x,y_{0,0}(x))
        \bigr]_{x=b_Q}^{x=0}          \\
         & =
        h_{Q,0}(0)
        =
        \begin{cases}
            -L_Q, & \text{if $Q$ approaches from the right}, \\
            L_Q,  & \text{if $Q$ approaches from the left}.
        \end{cases}
    \end{align}
    In particular, this derivative is nonzero.
    Since there are only finitely many such
    intervals, every sufficiently small nonzero
    $s$ separates all their arrival times from
    $\tau_*$.

    Follow the perturbed evolution until its
    first event.
    If no boundary collision or obstruction
    wall occurs before $\tau_*$, none of the
    additional active intervals meets $q_s$
    at the collapse time.
    The stationary pieces in $U$ remain on
    $x=0$, so they also avoid $q_s$ for
    $s\neq0$.
    A sufficiently small fixed neighborhood
    of $q_s$ therefore contains, for times
    sufficiently close to $\tau_*$, only the
    contracting intervals and their two
    residue branches.
    The collapse is spatially isolated, and
    Lemma~\ref{lem:elliptic-paired-collapse}
    shows that it is paired.
\end{proof}

\begin{figure}[htbp]
    \centering
    \includegraphics[width=0.7\textwidth]{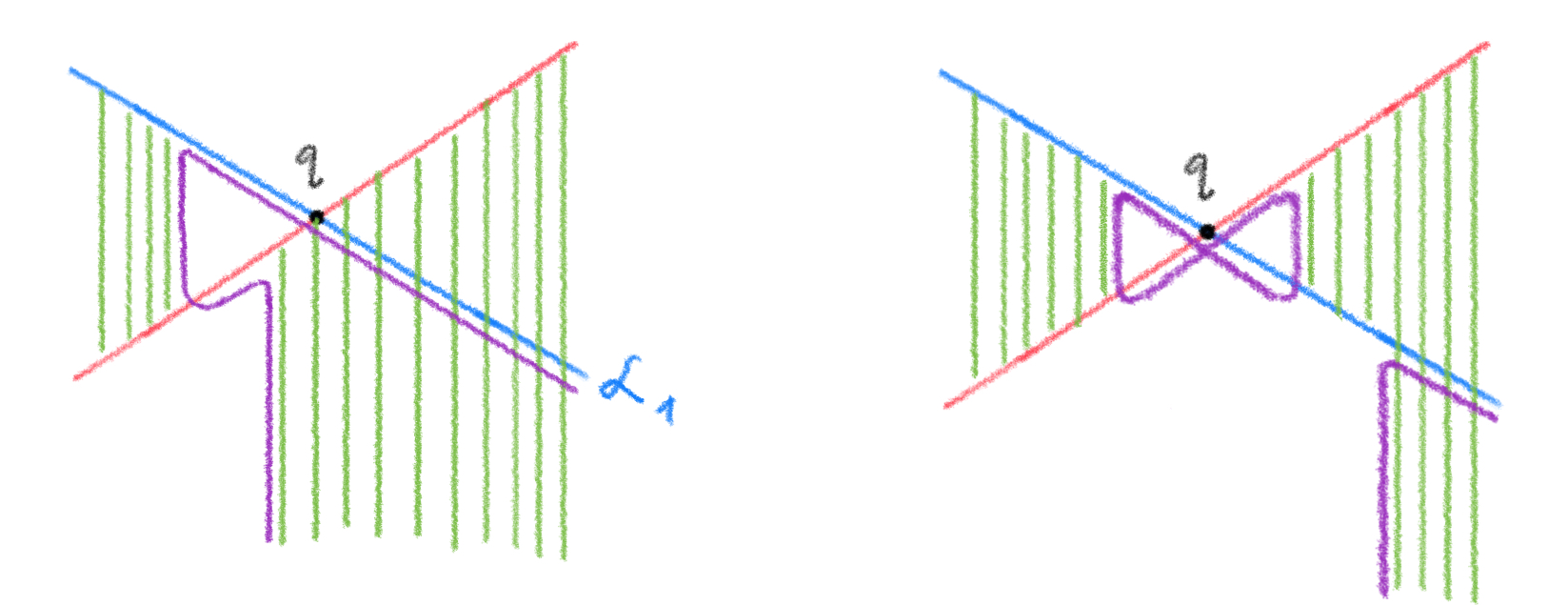}
    \caption{
        Two outcomes of perturbing the simultaneous
        arrival and collapse.
        In both panels, $q$ denotes the perturbed
        intersection.
        Left: the incoming front reaches $q$
        before the collapse time.
        The smoothed web cone passes through a
        self-tangency and subsequently undergoes
        a boundary collision, leaving an inactive
        fragment.
        The drawing shows the stage between
        these two events.
        Right: perturbing in the opposite direction
        delays the incoming front.
        An obstruction wall occurs first and is
        resolved by puncture, producing the
        configuration shown.
        One fragment remains inactive, while the
        separate figure-eight component shrinks
        to $q$ and disappears in a paired collapse.
        Green segments indicate successive front
        positions; the smoothed web cone is shown
        in purple, slightly displaced for clarity.
    }
    \label{fig:elliptic-collapse-separated}
\end{figure}

Figure~\ref{fig:elliptic-collapse-separated}
shows the two perturbed evolutions of the
configuration in
Figure~\ref{fig:elliptic-simultaneous-collapse}.

\begin{lemma}[Perturbation before the next geometric event]
    \label{lem:elliptic-finite-preparation}
    Let a special Lagrangian web for $f$ evolve regularly
    on $(a,\tau_*)$, with fixed boundary
    conditions and no obstruction walls.
    Suppose that its first boundary collision
    or collapse occurs at $\tau_*$.
    Assume the hypotheses of
    Lemma~\ref{lem:elliptic-collapse-isolation}
    at every collapse point that contains
    no inactive endpoint.

    At a regular time $\tau_0<\tau_*$ sufficiently
    close to $\tau_*$, the perturbations of
    that lemma can be made simultaneously
    at all these collapse points.

    After this perturbation, the evolution
    reaches an obstruction wall, boundary
    collision, or collapse no later than
    $\tau_*$, without another such perturbation.
    Every collapse reached at that time is
    either spatially isolated and paired,
    or part of a boundary collision with
    an inactive endpoint.
    In the latter case, concatenate the
    shrinking active interval with the
    inactive interval, obtaining an outgoing
    active interval of positive length.
\end{lemma}

\begin{proof}
    There are finitely many active intervals,
    hence finitely many collapse points at
    $\tau_*$.
    Choose disjoint neighborhoods of these
    points.
    In those containing no inactive endpoint,
    choose the Hamiltonians of
    Lemma~\ref{lem:elliptic-collapse-isolation}.
    Leave the neighborhoods containing
    inactive endpoints unchanged.

    We first show that these perturbations
    cannot create a collapse before $\tau_*$.
    Choose $\tau_0$ sufficiently close to $\tau_*$
    that each contracting interval, together
    with the two boundary arcs leading to
    its collapse point, bounds a small disk
    in a coordinate neighborhood.
    By Lemma~\ref{lem:elliptic-harmonic-flux},
    the area of this disk is its remaining
    evolution time,
    \[
        \tau_*-\tau_0.
    \]
    The common Hamiltonian deformation moves
    both boundary arcs and fixes the interval
    at $\tau_0$.
    It therefore carries this disk to the
    corresponding disk for the new boundary
    curves and preserves its area.
    Consequently, if the interval continues
    along the same boundary branches, it
    still collapses at $\tau_*$.

    Every other active interval has positive
    limiting length at $\tau_*$.
    Its equation $\dot x=1/h(x)$ therefore
    remains regular up to that time.
    By continuous dependence on the boundary
    curves, a sufficiently small perturbation
    keeps its length positive throughout
    $[\tau_0,\tau_*]$.
    Thus no interval collapses earlier.

    At each moved collapse point,
    Lemma~\ref{lem:elliptic-collapse-isolation}
    shows that the arrival time of each
    additional active interval changes with
    nonzero derivative, whereas the collapse
    time stays fixed.
    There are finitely many such arrival-time
    equalities, so the perturbation parameters
    can be chosen simultaneously to avoid
    all of them.
    The moved collapse points also avoid
    the stationary vertical pieces.

    Now resume evolution from the perturbed
    data at $\tau_0$.
    If an obstruction wall or boundary
    collision occurs before $\tau_*$,
    continue that event first.
    There can be no earlier collapse,
    by the preceding argument.

    Otherwise, the contracting intervals
    reach $\tau_*$.
    At the moved intersections, their
    collapses are spatially isolated and
    hence paired.
    At an unchanged intersection containing
    an inactive endpoint, the shrinking
    residue gives a boundary collision.
    Concatenating the incident intervals
    preserves the positive length of the
    inactive interval and gives the outgoing
    endpoint labels.
    If there were no collapses to perturb,
    the original boundary collision at
    $\tau_*$ is unchanged.
\end{proof}

\subsubsection{Classification of geometric degenerations}

\begin{lemma}
    \label{lem:elliptic-geometric-degenerations}
    Consider a regular interval of a special Lagrangian web for $f$ with fixed, transverse,
    strictly positive boundary conditions.
    Assume that the perturbations of
    Lemma~\ref{lem:elliptic-collapse-isolation}
    have been made as necessary, so that every
    collapse whose limiting point contains no
    inactive endpoint is spatially isolated.
    Then the only finite-time geometric
    degenerations of active fronts are
    boundary collisions and paired collapses.
\end{lemma}

\begin{proof}
    Consider a connected active interval $Y_\tau$
    as $\tau\nearrow\tau_*<\infty$.
    By Lemma~\ref{lem:elliptic-harmonic-flux},
    its signed height is $h(x(\tau))$ and its
    horizontal coordinate satisfies
    $\dot x=1/h(x)$.
    Consequently,
    \[
        \frac{d}{d\tau}\ell(Y_\tau)^2
        =
        \frac{d}{d\tau}h(x(\tau))^2
        =2h'(x(\tau)).
    \]
    Compactness and strict positivity of the
    boundary curves give a uniform bound on
    their graph slopes, and hence on $h'$.
    Thus the squared length has bounded derivative
    and a finite limit as $\tau\nearrow\tau_*$.
    Write $\ell_*\geq0$ for the limiting length.

    Suppose first that $\ell_*>0$.
    Then $|h(x(\tau))|$ stays bounded away from
    zero, so $\dot x=1/h(x)$ remains bounded.
    The bounded boundary slopes also bound the
    speeds of the endpoint attachments in $M_i$.
    Since $\tau_*$ is finite, these attachment
    points converge.
    Near their limits, the boundary branches
    remain smooth graphs and $h$ remains nonzero.
    The ordinary differential equation
    \[
        \dot x=\frac{1}{h(x)}
    \]
    therefore extends $x(\tau)$ smoothly through
    $\tau_*$.
    Substituting this solution into
    \eqref{eq:elliptic-interval-parametrization}
    extends the entire vertical interval.
    The parametrization remains an immersion
    because its source derivative is $(0,h)$.
    In particular, no interior degeneration occurs.

    This continuation must also satisfy
    Definition~\ref{def:admissible-web}.
    In dimension one, its embedding condition
    says that all active and inactive endpoints
    of label $i$ have distinct attachment points
    in $M_i$.
    As long as these points remain distinct,
    their order on the boundary source is
    unchanged, and they continue to bound the
    swept and residue intervals.
    Admissibility can therefore fail only when
    two such attachment points meet.
    Since the swept regions expand, the
    intervening residue interval then shrinks
    to a point.
    This is the boundary collision of
    Definition~\ref{def:elliptic-boundary-collision}.

    Suppose now that $\ell_*=0$.
    The endpoint images are joined by an interval
    whose length tends to zero.
    Write $p_i(\tau)\in M_i$ for the endpoint
    images.
    Every accumulation point $(q_0,q_1)$ of
    $(p_0(\tau),p_1(\tau))$ satisfies
    $\iota_0(q_0)=\iota_1(q_1)$, because the
    distance between the endpoint images tends
    to zero.
    There are finitely many such intersection
    pairs, by compactness and transversality.

    Choose disjoint neighborhoods of these pairs.
    For all $\tau$ sufficiently close to $\tau_*$,
    the pair $(p_0(\tau),p_1(\tau))$ lies in their union:
    otherwise, compactness would give an
    accumulation point outside the intersection set.
    By continuity, it cannot move from one
    neighborhood to another without leaving
    this union.
    It therefore stays in a single neighborhood,
    whose intersection pair is its only possible
    accumulation point.
    Thus
    \[
        (p_0(\tau),p_1(\tau))
        \longrightarrow(q_0,q_1)
        \qquad\text{as }\tau\nearrow\tau_*.
    \]
    Choose graph coordinates near this
    intersection, with intersection coordinate
    $x_*$.
    Then
    \[
        h(x_*)=0,
        \qquad h'(x_*)\neq0.
    \]
    If $h'(x_*)>0$, the derivative formula above
    would give
    \[
        \frac{d}{d\tau}\ell(Y_\tau)^2
        =2h'(x(\tau))>0
    \]
    for all $\tau$ sufficiently close to
    $\tau_*$.
    This is impossible: a positive increasing
    function cannot tend to zero.
    Hence $h'(x_*)<0$, so
    \eqref{eq:elliptic-maslov-slope} gives
    Maslov index one.
    The degeneration is therefore a collapse
    in the sense of
    Definition~\ref{def:elliptic-collapse}.

    If it reaches an inactive endpoint, the
    intervening residue shrinks and gives a
    boundary collision.
    Concatenation with the inactive interval
    has positive limiting length. The remaining endpoints have opposite labels,
    so the outgoing interval is active.
    Otherwise, the preparation makes the
    collapse spatially isolated, and
    Lemma~\ref{lem:elliptic-paired-collapse}
    shows that it is paired.
\end{proof}

\subsection{Obstruction walls and brane transfer}
\label{subsec:elliptic-obstruction-walls}

We describe the surgery operations of
Subsection~\ref{subsec:obstruction-walls}
on the two-torus.
Their relation to graded $0$-surgery on a
smoothed web cone is stated in
Proposition~\ref{prop:elliptic-surgery-smoothing}.

Superscripts $-$ and $+$ denote incoming and
outgoing data.
The same operations may be performed successively
at one wall time, while other coefficients still
have valuation zero.

\begin{definition}[Elementary wall crossing operation]
    \label{def:elliptic-elementary-wall-operation}
    An \emph{elementary wall crossing operation}
    is a webbing, puncture, interior $0$-surgery,
    or absorption, as described in
    Definitions~\ref{def:elliptic-webbing}--\ref{def:elliptic-absorption},
    performed at a selected degree-one generator
    whose coefficient in the bounding cochain
    has valuation zero.
    This includes the active-return resolution
    of Remark~\ref{rem:elliptic-overlap-puncture}:
    a puncture followed by boundary reconnection.
\end{definition}

\subsubsection{Webbing, puncture, interior surgery and absorption}
\label{subsubsec:elliptic-geometric-resolutions}

\begin{definition}[Webbing]
    \label{def:elliptic-webbing}
    In dimension one, the webbing of
    Section~\ref{subsubsec:webbing}
    creates two active intervals, one on each
    horizontal side of the selected intersection.
    Adjoining the remaining residues gives
    the two arcs of the graded crossing
    resolution; see Figure~\ref{fig:web-cone}.
\end{definition}

\begin{figure}[htbp]
    \centering
    \includegraphics[width=0.6\textwidth]
    {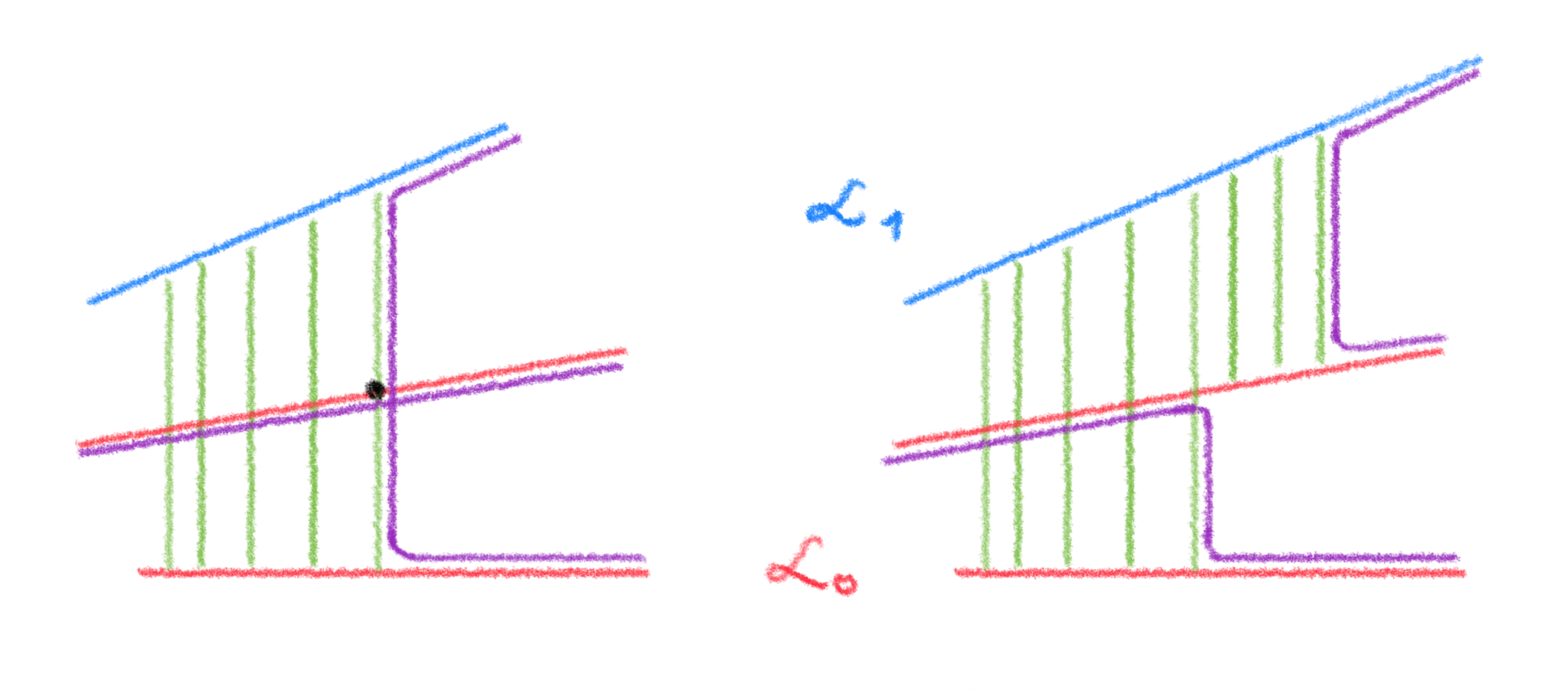}
    \caption{
        Puncture at a residue branch of $L_0$.
        Left: the green segments show successive
        positions of an active front moving to the
        right. An additional branch of $L_0$ meets
        its interior at the black dot.
        Right: puncturing at this intersection gives
        both new endpoints label $0$.
        The upper interval has endpoint labels $0,1$
        and continues moving to the right;
        the lower interval has labels $0,0$
        and remains stationary.
        The corresponding smoothed web cone,
        shown in purple and slightly displaced
        for clarity, undergoes a graded $0$-surgery.
    }
    \label{fig:elliptic-puncture}
\end{figure}

\begin{definition}[Puncture]
    \label{def:elliptic-puncture}
    A \emph{puncture} is the $n=1$ case of
    Definition~\ref{def:web-front-puncture}.
    It cuts the participating interval at the
    selected interior source point and gives
    both new endpoints the selected residue label.

    Puncturing an active interval produces
    one active and one inactive interval;
    see Figure~\ref{fig:elliptic-puncture}.
    Puncturing an inactive interval uses its
    missing boundary label and produces two
    active intervals.
\end{definition}

\begin{figure}[htbp]
    \centering
    \includegraphics[width=0.6\textwidth]
    {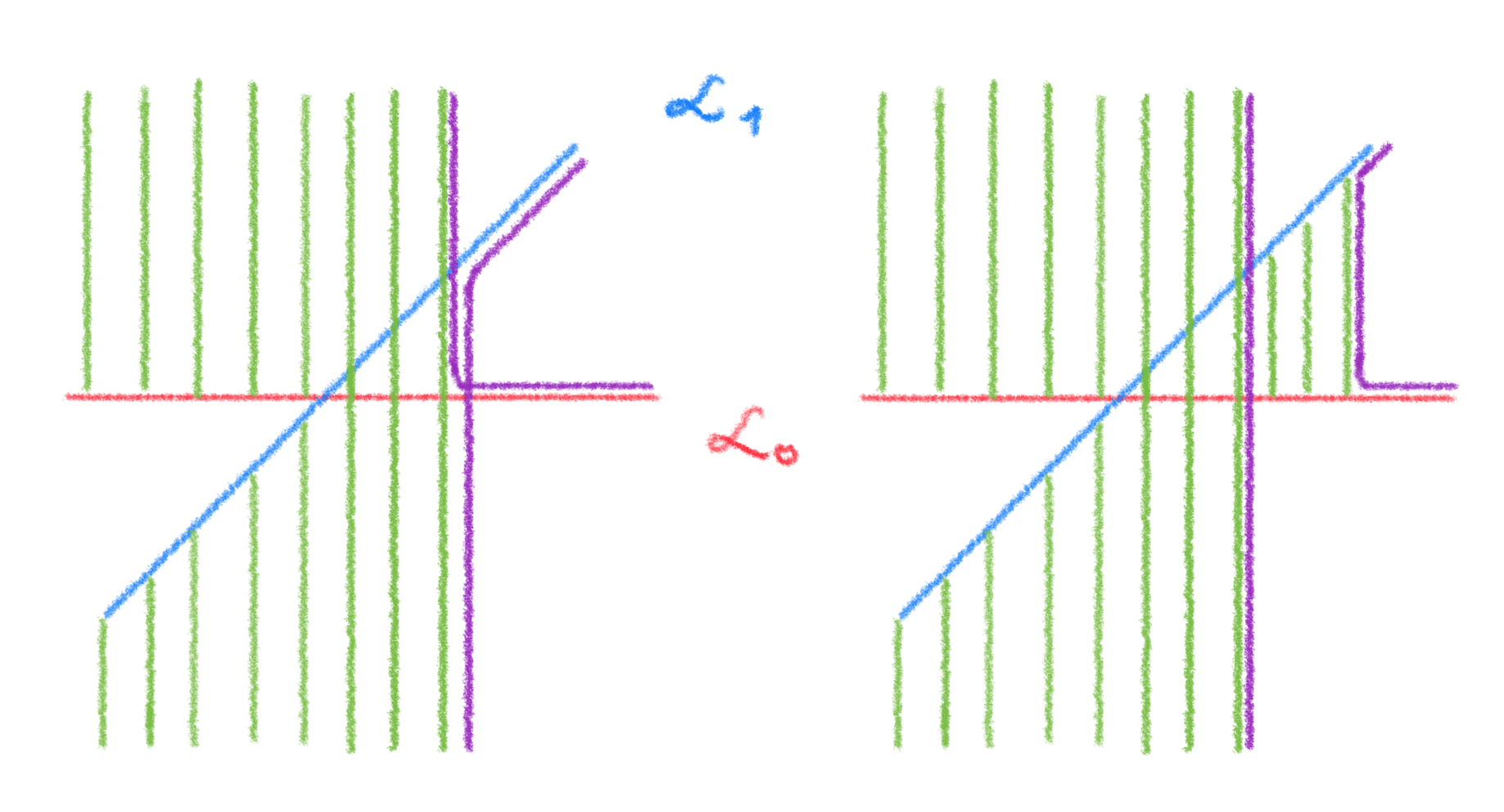}
    \caption{
        Puncture at an active return with $k=1$.
        The vertical direction is periodic with
        period one, and green segments show
        successive front positions.
        Left: a right-moving front starts on $L_0$,
        winds around the vertical fiber, and ends
        on $L_1$.
        Its endpoint images coincide at
        $L_0\cap L_1$, but their different labels
        prevent this from being a boundary collision.
        The front continues regularly, covering
        the segment between $L_0$ and $L_1$ twice.
        Right: puncture at the return to $L_0$,
        followed by boundary reconnection, produces
        a closed vertical fragment of degree one
        and label $1$.
        The remaining active interval has limiting
        length $\ell(P)-1$ and continues moving
        to the right.
        The smoothed web cone is shown in purple,
        with a small perturbation making the
        selected self-intersection transverse.
        The resolution is realized by a graded
        $0$-surgery.
    }
    \label{fig:elliptic-active-return}
\end{figure}

\begin{remark}[Puncture at an active return]
    \label{rem:elliptic-overlap-puncture}
    Suppose that the selected crossing represents
    a return of an active interval $P$ to its
    label-$i$ endpoint.
    Let the returning source point have distance
    \[
        k\in\mathbb Z_{>0},
        \qquad k<\ell(P),
    \]
    from that endpoint.

    The return meets an existing attachment
    on $\partial\overline{R_i^-}$.
    We use the source cut of
    Definition~\ref{def:web-front-puncture},
    giving both new endpoints label $i$.
    The piece adjacent to the original
    label-$i$ endpoint then has both endpoints
    lifting to the same boundary source point.
    Immediately apply the boundary reconnection
    of Definition~\ref{def:elliptic-boundary-collision}
    to close this piece.
    The cut and reconnection form one combined
    operation.

    The result is a closed label-one vertical
fragment of covering degree $k$ and an active
interval of limiting length $\ell(P)-k$;
see Figure~\ref{fig:elliptic-active-return}.
\end{remark}

The following two operations realize the
$k=0$ source surgery of
Definition~\ref{def:interior-surgery}
at clean intersections.

\begin{figure}[htbp]
    \centering
    \includegraphics[width=\textwidth]
    {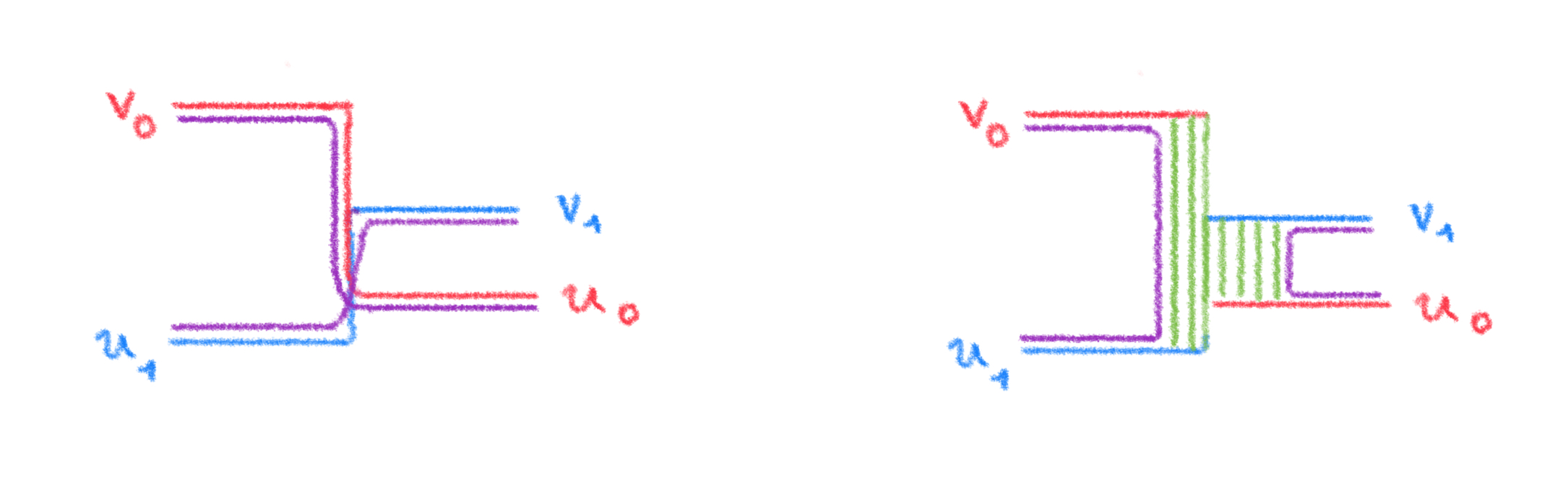}
    \caption{
        Vertical reconnection.
        Left: two overlapping nonclosed inactive
        intervals, of label $0$ in red and
        label $1$ in blue.
        The smoothed web cone, shown in purple,
        is slightly perturbed to make the selected
        crossing transverse.
        Right: reconnection produces two active
        intervals with limiting lifted ranges
        $[u_0,v_1]$ and $[u_1,v_0]$.
        Their endpoint labels, from bottom to top,
        are $(0,1)$ and $(1,0)$, so they move
        rightward and leftward, respectively.
        Green segments show their successive positions.
        On the smoothed web cone, the operation
        is a graded $0$-surgery.
    }
    \label{fig:elliptic-vertical-reconnection}
\end{figure}

\begin{definition}[Vertical reconnection]
    \label{def:elliptic-vertical-reconnection}
    Suppose that the selected crossing resolves
    an overlap of two nonclosed inactive intervals
    of opposite labels.
    Choose the participating lifts
    \[
        \{a\}\times[u_0,v_0],
        \qquad
        \{a\}\times[u_1,v_1],
    \]
    where the subscript denotes the boundary
    label, and choose the common height
    \[
        c\in(u_0,v_0)\cap(u_1,v_1)
    \]
    corresponding to the selected source pair.

    Cut both intervals at that height and join
    the lower part of each to the upper part
    of the other; see Figure~\ref{fig:elliptic-vertical-reconnection}.
    The resulting intervals have lifted ranges
    \[
        [u_0,v_1],
        \qquad
        [u_1,v_0].
    \]
    Their endpoint labels, from bottom to top,
    are $(0,1)$ and $(1,0)$, respectively.
    They therefore evolve rightward and
    leftward.
    All incoming limiting boundary lifts are retained.
    We call this operation a
    \emph{vertical reconnection}.
\end{definition}

\begin{definition}[Interior $0$-surgery with a closed fragment]
    \label{def:elliptic-interior-surgery-closed}
    Let
    \[
        Y^-=[F_*^-:P^-\longrightarrow X]
    \]
    be a label-zero inactive interval, and let
    $S=[F_S:M_S\to X]$ be a closed label-one
    fragment on the same vertical fiber.
    Choose source points
    \[
        z\in\operatorname{int}(P^-),
        \qquad
        w\in M_S,
        \qquad
        F_*^-(z)=F_S(w).
    \]
    Cut $P^-$ at $z$ and $M_S$ at $w$.
    Following the upward vertical orientation,
    insert $M_S\setminus \{w\}$ between
    the two parts of $P^-$.
    This produces an inactive interval
    $Y^+=[F^+:P^+\to X]$ with
    \[
        \ell(Y^+)=\ell(Y^-)+\ell(S).
    \]
    Its endpoints retain their label $0$,
    boundary lifts, and residue attachments.
    The selected component $S$ is removed
    from $D_{\mathrm{cl},1}$.
\end{definition}

\begin{definition}[Absorption]
    \label{def:elliptic-absorption}
    An \emph{absorption} is the $n=1$ case of
    Definition~\ref{def:absorption},
    with outgoing web data as in
    Remark~\ref{rem:absorption-boundary-conditions}.
\end{definition}

\subsubsection{Surgery and brane transfer}
At an obstruction wall, first continue boundary
collisions using
Proposition~\ref{prop:elliptic-boundary-continuation}
and remove paired collapsing components disjoint
from the remaining cone, as in
Remark~\ref{rem:elliptic-geometric-continuation}.
The following companion result relates surgery
on the smoothed cone to elementary wall crossing
operations and supplies the brane comparisons.
Its proof is deferred to
\cite{KapotaWebFloer}.

\begin{prop}[Surgery, smoothing, and brane transfer]
    \label{prop:elliptic-surgery-smoothing}
    Fix a time $\tau$ and admissible web data
    whose geometric pieces admit regular
    evolution for a short time after $\tau$.
    Suppose that the web cone carries provisional
    cone data arising during initialization or
    continuation, with rank-one source local
    systems and a degree-one cochain of
    nonnegative valuation.

    The smoothing and its Floer data can be
    chosen to give provisional data
    \[
        (C,E,b),
        \qquad \operatorname{val}(b)\geq0.
    \]
    Here $b$ satisfies the Maurer--Cartan
    equation, and its Floer insertion series
    converge.
    These data have the following properties.

    \begin{enumerate}
        \item \emph{Choosing a surgery.}
              Suppose that $b$ has a nonzero
              coefficient of valuation zero.
              Then either there is a degree-one
              generator $p$ at a transverse
              double point of $C$, with
              \[
                  \operatorname{val}(b_p)=0,
              \]
              for which item~{\rm(2)} applies,
              or all valuation-zero coefficients
              run from shifted closed label-zero
              components to closed label-one
              components on the same fibers.
              The latter case is handled by
              closed-component reduction in
              Section~\ref{subsubsec:elliptic-closed-reduction}.

        \item \emph{Surgery and smoothing.}
              For a generator $p$ chosen in
              item~{\rm(1)}, there is an elementary
              wall crossing operation on the web
              such that graded $0$-surgery at $p$,
              followed by an exact deformation,
              gives a smoothing of the resulting
              web cone.
              The graded surgery can be chosen
              without introducing additional
              self-intersections.

              The resulting web is a compatible
              continuation in the sense of
              Definition~\ref{def:surgery-compatible}.
              Its geometric pieces admit regular
              evolution for a short time after
              the same time $\tau$.

              In the absorption case, $L_0^+$ can
              be chosen so that replacing the surgery
              diagram by a diagram with residue $R_0^+$
              does not increase the number of counted
              crossings, in the sense of
              Subsection~\ref{subsec:elliptic-counting}.
              This replacement preserves the exterior
              diagram and the surviving attachments,
              pairing, and winding of the surgery diagram.

        \item \emph{The outgoing data.}
              The surgery and exact deformation
              in item~{\rm(2)} give provisional
              outgoing cone data
              \[
                  (C_+,E_+,b_+),
                  \qquad
                  \operatorname{val}(b_+)\geq0,
              \]
              with rank-one source local systems
              and the same Maurer--Cartan and
              convergence properties.
    \end{enumerate}
    If valuation-zero coefficients remain, the
construction can be repeated at the same time
$\tau$, before regular evolution resumes.
If a finite sequence of these operations
and closed-component reductions ends with
a cochain of strictly positive valuation,
the resulting brane admits the comparison of
Definition~\ref{def:brane-surgery-transfer}
at an obstruction wall, or of
Definition~\ref{def:initial-brane-comparison}
at initialization.
\end{prop}

\subsubsection{Closed-component reduction}
\label{subsubsec:elliptic-closed-reduction}

We specialize the reduction of
Definition~\ref{def:closed-component-reduction}
to closed vertical components.
Throughout this subsection, a bar denotes
reduction modulo $\Lambda_+$.

\begin{definition}[Covering multiplicity]
    Let $Q$ be a connected closed immersed
    vertical curve with parametrization
    $\iota_Q:\widetilde Q\to X$.
    Write $S_Q:=\operatorname{im}\iota_Q$
    for its supporting primitive vertical fiber.
    The \emph{covering multiplicity} $m(Q)$
    is the number of sheets of the finite covering
    $\iota_Q:\widetilde Q\to S_Q$.
\end{definition}

\begin{definition}[Closed-component multiplicities]
    \label{def:elliptic-fiber-multiplicities}
    Fix admissible web data.
    For $i=0,1$, let $\mathcal Q_i$ be the collection
    of closed vertical components carrying label $i$,
    keeping distinct source components separate even
    when their images coincide.
    Let
    \[
        \mathcal S
        :=
        \bigl\{
        S_Q:
        Q\in\mathcal Q_0\cup\mathcal Q_1
        \bigr\}
    \]
    be the finite set of their supporting primitive
    vertical fibers.
    For each $S\in\mathcal S$, define
    \[
        d_i(S)
        :=
        \sum_{\substack{
                Q\in\mathcal Q_i\\
                S_Q=S}}
        m(Q),
        \qquad i=0,1.
    \]
    Thus $d_i(S)$ counts the total number of sheets
    of closed label-$i$ components over $S$.
    Empty sums are understood to be zero.
\end{definition}

The following reduction result is proved in the
Floer-theoretic companion~\cite{KapotaWebFloer}.

\begin{lemma}[Closed-component cancellation]
    \label{lem:elliptic-closed-cancellation}
    Let
    \[
        \widehat{\mathcal C}
        =
        (\widehat C,\widehat E,\widehat b)
    \]
    be the smoothed provisional cone data at initialization or
    after the geometric operations at an event.
    Assume that the source local systems have
    rank one and that $\widehat b$ is a degree-one
    cochain of nonnegative valuation satisfying
    the Maurer--Cartan equation.
    Use the clean Floer model for coincident closed
    vertical components. Suppose that every term of
    valuation zero in $\widehat b$ runs from a shifted
    closed label-zero component to an unshifted closed
    label-one component on the same supporting fiber.

    Then there exists a reduction in the sense of
Definition~\ref{def:closed-component-reduction}
whose replacement components are primitive
copies of the original supporting fibers.
For each original supporting fiber $S$,
there is an integer $r_S\geq0$ such that
\[
    d_i^+(S)=d_i^-(S)-r_S,
    \qquad i=0,1.
\]
If $\overline{\widehat b}\neq0$, then
\[
    \sum_S r_S>0,
\]
where the sum is over the original supporting
fibers.
The reduction supplies the brane transfer
for the completed continuation, or the
initial brane comparison at initialization.
\end{lemma}

\subsection{Counting diagrams}
\label{subsec:elliptic-counting}
\subsubsection{Counting diagrams}

We want to count the crossings of the singular web cone.
At a given time, however, active intervals,
inactive intervals, and closed vertical components
may overlap along vertical segments.
To obtain a crossing count immediately after
that time, we construct an auxiliary diagram.

Fix a time $\tau$ and a specified regular geometric
evolution for a short time after $\tau$.
At a singular event, we use the outgoing evolution
after the chosen surgeries have been completed.
We construct the diagram by advancing the active
intervals a small amount in time and then
straightening the inactive intervals and residue
arcs in thin vertical strips, as described below.
The minimality of these counting diagrams
is established in
Lemma~\ref{lem:elliptic-outgoing-minimality}.

\begin{remark}
    \label{rem:elliptic-residue-orientation}
    At an endpoint of an active interval, the
    attached residue extends in the direction
    of motion. Thus both attached residue arcs
    extend rightward for a right-moving interval
    and leftward for a left-moving interval.
    For an inactive interval, the two attached
    residue arcs extend in opposite horizontal
    directions.
\end{remark}

\begin{figure}[htbp]
    \centering
    \includegraphics[scale=0.15]
    {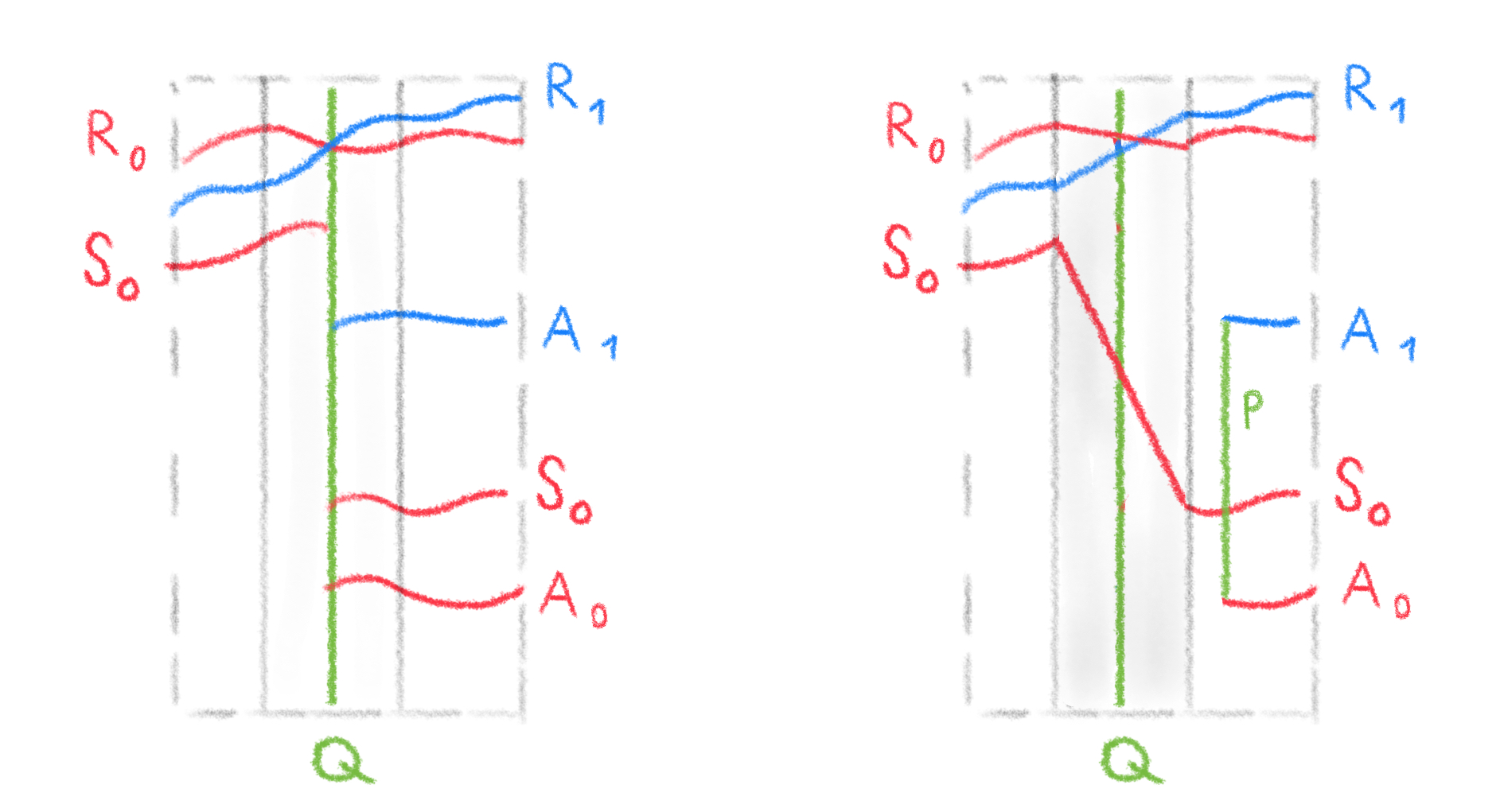}
    \caption{
        Construction of an outgoing counting diagram.
        The outer frames and shaded bands indicate
        the outer and inner strips, respectively;
        the top and bottom edges are identified.
        Left: the active interval $P$, with residue
        arms $A_0,A_1$, and an inactive fragment,
        with both residue arms marked $S_0$,
        share a supporting fiber with the closed
        component $Q$.
        The curves $R_0,R_1$ are additional
        residue branches.
        Right: $P$ has advanced to the right,
        while the stationary strands have been
        straightened inside the inner strip,
        keeping their endpoints fixed.
        The closed component $Q$ remains unchanged.
        The crossing of $P$ with the $S_0$ residue
        arm outside the inner strip also contributes
        to the crossing count.
    }
    \label{fig:elliptic-counting-diagram}
\end{figure}

\begin{definition}[Outgoing counting diagram]
    \label{def:elliptic-counting-diagram}
    We construct a single diagram on $X$ by
    performing the following steps around every
    fiber supporting a vertical piece at time $\tau$; see Figure~\ref{fig:elliptic-counting-diagram}.

    \begin{enumerate}
        \item \emph{Choose the outer strips.}
              For each distinct fiber containing
              an active interval, an inactive
              interval, or a closed vertical
              component, choose a thin vertical
              neighborhood.
              At a singular time, use the
              specified outgoing limits.

              Choose these neighborhoods pairwise
              disjoint and call them
              \emph{outer strips}.
              On $X$, their vertical coordinate
              is periodic, so they are annuli;
              their lifts to $\mathbb C$ are
              ordinary vertical strips.
              Keep the diagram outside the union
              of these strips fixed.

        \item \emph{Advance the active intervals.}
              Choose a sufficiently small $\epsilon>0$.
              Advance each active interval $P$
              from $\tau$ to $\tau+\epsilon+\delta_P$,
              following $\dot x=1/h(x)$ and its
              prescribed boundary branches.
              Choose the auxiliary times
              $\delta_P\geq0$ sufficiently small
              that the resulting fibers are distinct
              and preserve every strict horizontal
              ordering at the common time
              $\tau+\epsilon$.
              Keep every interval inside its
              original strip, and update
              all attached residues accordingly.

        \item \emph{Choose the inner strips.}
              Inside each outer strip, choose
              a narrower strip around the original
              supporting fiber that contains no
              advanced active interval.
              This is possible because the active
              intervals have moved away from that
              fiber.
              The stationary pieces remain on
              the original fiber, inside this
              \emph{inner strip}.

        \item \emph{Straighten the stationary pieces.}
              Work in a lift of each inner strip
              to $\mathbb C$.
              Each inactive interval, together
              with the portions of its attached
              residue arcs inside the strip,
              forms an arc joining the two sides
              of the strip.
              Replace this arc by the straight
              segment with the same endpoints.
              Do the same for every residue arc
              crossing the strip.

              Make identical replacements on all
              integer vertical translates and
              project back to $X$.
              Leave the closed vertical components
              and the diagram outside the inner
              strips unchanged.
    \end{enumerate}

    Denote the resulting diagram by $\cc{G}_\tau$
    and call it an \emph{outgoing counting diagram
        at time $\tau$}.
\end{definition}

The next definition counts crossings in the entire
outgoing diagram $\cc{G}_\tau$, after completing
the construction around all supporting fibers.
Although the straightening takes place only inside
the inner strips, crossings are counted throughout
$X$.
In particular, we include crossings between
advanced active intervals and residue arcs outside
the inner strips.

Corollary~\ref{cor:elliptic-count-well-defined}
shows that, when the strip
widths, time advance, and additional separations
are sufficiently small, the total count
$D(\cc{G}_\tau)$ defined below is independent
of these choices.
It depends only on the specified outgoing
evolution at $\tau$, and we then write $D(\tau)$
for this number.

\begin{definition}[Crossing count]
    \label{def:elliptic-crossing-count}
    Let $c(\cc{G}_\tau)$ be the number of transverse
    crossings between the interiors of the branches
    of $\cc{G}_\tau$, counted as unordered pairs
    of distinct source points.
    Exclude pairs supported entirely on closed
    vertical components, but include crossings
    between closed components and other branches.
    All source multiplicities are retained.
    In particular, a transverse passage through
    a closed component of covering degree $m$
    contributes $m$ crossings.

    We also assign a return count to each active
    interval.
    Write $\mathcal P_\tau$ for the active source
    intervals in $\cc{G}_\tau$.
    For $P\in\mathcal P_\tau$, let $F_P:P\to X$
    be its immersion in the diagram and let
    $e_0$ be its label-zero endpoint.
    Define its \emph{return count} by
    \[
        o_{\cc{G}_\tau}(P)
        :=
        \#\left\{
        z\in\operatorname{int}P
        \;\middle|\;
        F_P(z)=F_P(e_0)
        \right\}.
    \]
    These returns are not transverse crossings
    and are therefore not included in
    $c(\cc{G}_\tau)$.
    No return term is assigned to inactive
    intervals or closed components.
    Since the vertical period is one,
    \[
        o_{\cc{G}_\tau}(P)
        =
        \#\left\{
        k\in\mathbb Z_{>0}
        \;\middle|\;
        k<\ell_{\cc{G}_\tau}(P)
        \right\},
    \]
    where $\ell_{\cc{G}_\tau}(P)$ is the source
    length of the interval in the diagram.
    Using its label-one endpoint gives the
    same return count.
    The \emph{total crossing count} is
    \[
        D(\cc{G}_\tau)
        :=
        c(\cc{G}_\tau)
        +
        \sum_{P\in\mathcal P_\tau}
        o_{\cc{G}_\tau}(P).
    \]
\end{definition}

\subsubsection{Local outgoing minimality}

Fix an outgoing counting diagram $\cc{G}_\tau$
and one of its outer strips $A$.
Cut the web-cone source at the preimages of
$\partial A$.
We call the resulting source intervals whose
images lie in $A$ \emph{strands}.

The following lemma shows that the total count
inside $A$, including the active return terms,
is the smallest possible number of counted
transverse intersections with the following
data fixed:
the boundary endpoints
with their source identities, their pairing
by strands, and the winding of each strand.
The pairing records which two endpoints are
joined by each strand.
In the covering strip, fixing the winding means
fixing the two lifted endpoints up to simultaneous
integer vertical translation.
For closed components, we retain their
orientations, labels, and covering degrees.

\begin{lemma}[Local outgoing minimality]
    \label{lem:elliptic-outgoing-minimality}
    Fix a time $\tau$ and a specified regular
    geometric evolution immediately after $\tau$.
    The time $\tau$ may be regular or singular;
    in the latter case, use the outgoing
    configuration after the chosen surgeries.
    The boundary curves $L_0,L_1$ remain fixed
    during this outgoing evolution.

    Choose one supporting fiber at time $\tau$,
    translated to $x=0$, and let
    \[
        A=(-a,a)\times(\mathbb R/\mathbb Z)
    \]
    be its outer strip.
    Before advancing the active intervals,
    all vertical pieces in $A$ lie on $x=0$,
    using their outgoing limits if $\tau$ is
    singular.
    Choose $A$ sufficiently thin that every
    intersection between distinct residue
    branches inside $A$ lies on $x=0$, and
    the strand endpoints on $\partial A$
    have distinct images.

    Let $\cc{G}_\tau$ be an outgoing counting
    diagram constructed using this outer strip
    and sufficiently small choices in
    Definition~\ref{def:elliptic-counting-diagram}.
    Write $D_A(\cc{G}_\tau)$ for the sum of its
    crossing count inside $A$ and the return
    counts of its active intervals in $A$,
    with the conventions of
    Definition~\ref{def:elliptic-crossing-count}.
    Then the following hold.
    \begin{enumerate}
        \item \emph{Realization by transverse crossings.}
              An individual vertical interval may
              overlap itself by winding around
              its vertical fiber.
              Each active interval $P$ admits
              an arbitrarily small deformation,
              fixing its endpoints, that replaces
              these self-overlaps by exactly
              $o_{\cc{G}_\tau}(P)$ transverse
              self-intersections.
              These deformations preserve the
              crossings already counted in
              $\cc{G}_\tau$.
              After rounding the corners,
              the resulting diagram has exactly
              $D_A(\cc{G}_\tau)$ counted transverse
              crossings inside $A$.

        \item \emph{Minimality.}
              Every diagram of immersed strands
              and closed components that agrees
              with $\cc{G}_\tau$ outside $A$,
              retains the boundary endpoints,
              pairing, winding, and closed-component
              data specified above, and has
              transverse intersections for the
              branches being counted, has at least
              $D_A(\cc{G}_\tau)$ counted crossings
              inside $A$.
    \end{enumerate}
\end{lemma}

\begin{proof}
    Work in the covering strip
    $\widetilde A=(-a,a)\times\mathbb R$.
    A strand containing an active interval is
    called active; the remaining strands are
    called stationary.

    For two lifted strands, choose a sufficiently
    large rectangle in the covering strip
    containing both strands.
    Their boundary endpoints \emph{alternate}
    if, on traversing the boundary of the rectangle,
    they occur alternately from the two strands.
    Alternating endpoints force an intersection.
    This gives a lower bound for every pair of
    lifted strands, including distinct translates
    of the same strand.
    We show that the outgoing diagram attains
    these bounds simultaneously.

    A stationary strand joins the two sides
    of the strip.
    Its portion in the inner strip is straight,
    and its residue arms have no intersections
    outside that strip.
    Thus two lifted stationary strands intersect
    once precisely when their boundary endpoints
    alternate, and otherwise do not intersect.

    An active strand lies entirely on its outgoing
    side of the inner strip; see Remark~\ref{rem:elliptic-residue-orientation}.
    Both residue arms extend outward from its
    vertical interval.
    A stationary strand has only one residue arm
    on this side.
    Since residue branches do not intersect there,
    the two strands can meet only where that arm
    crosses the active interval, and hence at
    most once.
    Active strands moving in opposite directions
    lie on opposite sides of the inner strip
    and are disjoint.

    Consider next two active strands moving in
    the same direction.
    First consider their intervals at the common
    time $t=\tau+\epsilon$, before the additional
    separations in the counting-diagram construction.
    Use the outward coordinate $r=x$ for
    right-moving intervals and $r=-x$ for
    left-moving intervals.
    For an active interval $P$, let $\ell_P(r)$
    be the length between its boundary branches
    at coordinate $r$.
    By Lemma~\ref{lem:elliptic-harmonic-flux},
    the elapsed time needed to reach $r$ is
    \[
        T_P(r):=\int_0^r\ell_P(u)\,du,
        \qquad
        T_P(r_P(t))=t-\tau.
    \]
    The order of the boundary graphs is fixed
    throughout this half-strip.
    Consequently, if the interval between the
    boundary branches of $Q$ is strictly contained
    in that of $P$ at one positive coordinate $r$,
    the same containment holds throughout the
    half-strip.
    It follows that
    \[
        \ell_Q(r)<\ell_P(r),
        \qquad
        T_Q(r)<T_P(r),
        \qquad
        r_Q(t)>r_P(t).
    \]
    Thus the shorter, nested interval lies
    farther outward at the common time $t$.

    Only the farther-out vertical interval can
    meet the residue arms of the other strand.
    If it met both arms, it would contain the
    interval between the other strand's boundary
    branches at that horizontal coordinate.
    The preceding comparison would then require
    it to lie closer to the original fiber,
    a contradiction.
    Hence two same-direction active strands
    also meet at most once.

    If their horizontal positions coincide at
    the common time, strict nesting is excluded
    by the same argument.
    Their lifted intervals are therefore disjoint
    or have alternating endpoint heights.
    Either sufficiently small horizontal
    separation gives the same mutual count:
    zero in the first case and one in the second.

    The additional advances preserve every
    horizontal ordering that was already strict.
    For intervals sharing a fiber, either
    sufficiently small separation gives the
    same mutual count, as just shown.
    Thus these intersection bounds also hold
    in $\cc{G}_\tau$.

    It remains to resolve the self-overlaps
    of each active interval.
    Consider a right-moving interval at $x=b>0$,
    with lifted endpoints $u<v$.
    Replace its vertical segment by
    \[
        x=b-\rho(y-u)(v-y),
        \qquad u\leq y\leq v,
    \]
    for sufficiently small $\rho>0$.
    Its interior lies just behind the original
    interval, while its attached residues lie
    in $x\geq b$.
    The deformed arc and its translate by a
    positive integer $k$ intersect at exactly
    one point,
    \[
        y=\frac{u+v+k}{2},
    \]
    precisely when $k<v-u$.
    These intersections are transverse.
    Counting positive $k$ counts each unordered
    source pair once and gives exactly
    $o_{\cc{G}_\tau}(P)$.
    For a left-moving interval, reflect this
    construction horizontally.

    Choose these deformations and the subsequent
    corner roundings sufficiently small to avoid
    the inner strip and preserve all previously
    counted crossings.
    This can be done simultaneously, since only
    finitely many strands and relevant translates
    occur.
    The resulting diagram therefore has exactly
    $D_A(\cc{G}_\tau)$ counted transverse crossings,
    proving the first assertion.

    After these deformations, every pair of
    distinct lifted strands has at most one
    intersection.
    Such an intersection occurs precisely when
    their boundary endpoints alternate:
    alternating endpoints force an intersection,
    while intersection parity excludes a single
    crossing for nonalternating endpoints.
    Thus all the pairwise lower bounds for
    strands are attained.

    Finally, a closed component of covering
    degree $m$ meets each stationary through-strand
    exactly $m$ times.
    This realizes their absolute algebraic
    intersection lower bound.
    It is disjoint from every outgoing active
    strand, realizing the lower bound zero
    for those pairs.

    Summing over unordered source pairs and
    relative vertical translates shows that
    every competing diagram has at least
    $D_A(\cc{G}_\tau)$ counted crossings.
    This proves minimality.
\end{proof}

\begin{corollary}[Independence of choices]
    \label{cor:elliptic-count-well-defined}
    Fix a time $\tau$ and a specified regular
    outgoing evolution.
    The number $D(\cc{G}_\tau)$ is independent
    of the sufficiently small choices of outer
    strips, time advance, inner strips, and
    auxiliary separations in
    Definition~\ref{def:elliptic-counting-diagram}.
    We denote this number by $D(\tau)$.
\end{corollary}

\begin{proof}
    With the outer strips fixed, all these
    constructions have the same exterior diagram,
    boundary endpoints, pairing, winding, and
    closed-component data.
    By Lemma~\ref{lem:elliptic-outgoing-minimality},
    their counts equal the minimum for these data.

    For two choices of outer strips, choose
    smaller strips contained in both.
    Taking the time advances and inner strips
    sufficiently small gives a counting diagram
    allowed by both constructions.
    Their counts therefore agree.
\end{proof}

\begin{corollary}[Crossing count after surgery]
    \label{cor:elliptic-source-operation-count}
    Let $D^-$ be the count before a surgery
    and $D^+$ the count for a specified regular
    outgoing evolution.
    Suppose that surgery on a transverse diagram
    realizing $D^-$ produces an admissible
    competitor for the outgoing counting diagram,
    in the sense of
    Lemma~\ref{lem:elliptic-outgoing-minimality},
    removing at least $q\geq0$ counted crossings
    and introducing no new ones.
    Then
    \[
        D^+\leq D^- -q.
    \]
    In particular, $D$ strictly decreases if $q>0$.
\end{corollary}

\begin{proof}
    The modified diagram has at most $D^- -q$
    counted crossings.
    By Lemma~\ref{lem:elliptic-outgoing-minimality},
    the outgoing count is no greater than
    the count of this competitor.
\end{proof}

\subsubsection{Regular evolution and geometric degenerations}

We first compare the crossing counts during
regular evolution and across the geometric
continuations described above.

\begin{corollary}[Regular evolution]
    \label{cor:elliptic-crossing-monotonicity}
    Suppose that the fronts  evolve regularly on an interval
    $I$, with fixed boundary conditions.
    Then $D$ is non-increasing on $I$:
    \[
        D(t)\leq D(s),
        \qquad s<t,\quad s,t\in I.
    \]
\end{corollary}

\begin{proof}
    Fix $\tau\in I$ and choose outer strips around
    its supporting fibers.
    For nearby times, the exterior diagram,
    boundary endpoints, pairing, winding, and
    closed-component data remain fixed.

    Advancing the constructions at $\tau$ and
    at a nearby later time to the same time
    gives the same evolved pieces.
    Choose the auxiliary separations and inactive
    straightenings identically.
    The remaining residue straightenings preserve
    crossings, so $D$ is constant immediately
    to the right of $\tau$.
    At a nearby earlier time $t<\tau$, construct
    a counting diagram with time advance smaller
    than $\tau-t$.
    Resolving its active self-overlaps gives
    a competitor with $D(t)$ crossings for the
    outgoing minimum at $\tau$.
    Lemma~\ref{lem:elliptic-outgoing-minimality}
    therefore gives
    \[
        D(\tau)\leq D(t).
    \]
    Fix $s<t$ in $I$.
    Take the maximal subinterval in $[s,t]$ starting
    at $s$ on which $D(u)\leq D(s)$ for $u$ in the interval,
    and denote its
    right endpoint by $\tau$.
    Right local constancy gives $\tau>s$,
    and the comparison from the left gives
    $D(\tau)\leq D(s)$.
    If $\tau<t$, right local constancy extends
    this inequality beyond $\tau$, contradicting
    maximality.
    Hence $\tau=t$, and $D(t)\leq D(s)$.
\end{proof}

\begin{remark}[Crossing count under isolation perturbations]
    \label{rem:elliptic-isolation-count}
    We check that the perturbation of
    Lemma~\ref{lem:elliptic-collapse-isolation}
    preserves $D(\tau_0)$.
    Choose the counting advances small enough
    that the advanced active intervals avoid
    the perturbation support.
    All crossings involving them, and their
    return counts, are unchanged.

    Choose the strip boundaries away from residue
    intersections, and the perturbation sufficiently
    small that no such intersection passes through
    a strip boundary and the strand endpoints
    retain their order.
    Outside the inner strips, residue intersections
    are preserved by the common Hamiltonian
    deformation.
    Inside them, the straightened arcs have the
    same crossing counts because these depend
    only on the endpoint order.
    The same through-strands cross each closed
    component, with unchanged covering
    multiplicities.
    Thus $D(\tau_0)$ is unchanged.
\end{remark}

\begin{lemma}[Crossing count at geometric degenerations]
    \label{lem:elliptic-geometric-count}
    Continue boundary collisions by concatenation
    and spatially isolated paired collapses by
    deleting the contracting components.
    Then
    \[
        D^+\leq D^-.
    \]
    Each boundary collision removes two labelled
    endpoints, and each paired collapse removes
    four.
\end{lemma}

\begin{proof}
    At a boundary collision, the shrinking
    residue arc and its incident intervals
    have the same exterior connections and
    winding as their concatenation.
    The incoming resolved diagram therefore
    has the outgoing connections
    and no additional crossings.
    This comparison also applies when one
    incident interval has zero limiting length:
    keep its incoming source connections until
    concatenation, and then contract that
    zero-length portion.

    At a spatially isolated paired collapse,
    delete the contracting component from the
    incoming resolved diagram.
    Again, no crossing is introduced.
    Lemma~\ref{lem:elliptic-outgoing-minimality}
    gives $D^+\leq D^-$ in both cases.
    The endpoint losses follow from the
    definitions.
\end{proof}

\subsubsection{Obstruction reductions and the combined count}

We next compare the counts under elementary
wall crossing operations and closed-component
reductions.

\begin{lemma}[Decrease under an elementary wall crossing]
    \label{lem:elliptic-wall-crossings}
    Every elementary wall crossing operation
    chosen in
    Proposition~\ref{prop:elliptic-surgery-smoothing}
    satisfies
    \[
        D^+\leq D^- -1.
    \]
    It introduces at most four new labelled
    endpoints and creates no closed label-zero
    component.
\end{lemma}

\begin{proof}
    By Corollary~\ref{cor:elliptic-source-operation-count},
    it suffices to construct an outgoing comparison
    diagram with at most $D^- -1$ counted crossings.
    Start with an incoming diagram realizing $D^-$
    by transverse crossings, as in
    Lemma~\ref{lem:elliptic-outgoing-minimality}.

    \emph{Webbing, puncture, and vertical reconnection.}
    Each operation reconnects two strands at a
    crossing.
    For vertical reconnection, these are the
    straightened strands of the two inactive
    intervals.
    Resolve the crossing in a small disk,
    joining the branches as prescribed by the
    operation.
    Outside this disk, the paths are unchanged,
    so the resulting strands have the prescribed
    connections and winding.
    One counted crossing is removed and none
    is added.

    \emph{Active returns.}
    A return of order $k$ is represented by a
    self-crossing whose two source points are
    separated by $k$ vertical turns.
    Resolving this crossing splits off those
    $k$ turns as a closed component and leaves
    the shorter active strand with its original
    endpoints.
    This is precisely the result of puncture
    followed by boundary reconnection.
    The selected crossing disappears, and any
    further crossings lying entirely on the
    new closed component are no longer counted.

    \emph{Interior surgery with a closed fragment
        and boundary absorption.}
    Let $m$ be the covering degree of the
    participating closed component.
    In the comparison diagram, replace it by
    $m$ nearby disjoint primitive circles,
    also disjoint from the remaining closed
    components.
    This preserves its counted contribution:
    each through-strand still meets $m$ closed
    sheets.

    Join these circles successively to the
    participating inactive strand or residue
    strand, resolving one crossing at each join.
    The result inserts $m$ full turns into that
    strand, with the orientation prescribed by
    the operation.
    It therefore has the required connections
    and winding, while removing $m\geq1$
    counted crossings and introducing none.
    For boundary absorption, the subsequent
    replacement by the positive representative
    $L_0^+$ preserves these comparison data
    and does not increase the count, by
    Proposition~\ref{prop:elliptic-surgery-smoothing}.

    In every case, we have constructed the
    required outgoing comparison diagram.
    Corollary~\ref{cor:elliptic-source-operation-count}
    therefore gives $D^+\leq D^- -1$.

    Finally, webbing adds four labelled endpoints
    and ordinary puncture adds two.
    The remaining operations add none overall,
    since an active-return resolution creates
    two endpoints and then removes two by
    boundary reconnection.
    Every newly created closed fragment has
    label $1$ by definition.
\end{proof}

\begin{definition}[Closed label-zero multiplicity]
    \label{def:elliptic-closed-multiplicity}
    Let $\mathcal G$ be a counting diagram.
    Write $\mathcal Q_i$, $\mathcal S$, and $d_i(S)$
    for the closed-component data of its underlying
    web, as in
    Definition~\ref{def:elliptic-fiber-multiplicities}.
    The closed label-zero multiplicity of
    $\mathcal G$ is
    \[
        M_0(\mathcal G)
        :=
        \sum_{Q\in\mathcal Q_0}m(Q)
        =
        \sum_{S\in\mathcal S}d_0(S).
    \]
    Empty sums are understood to be zero.
\end{definition}

\begin{lemma}[Counts under closed-component reduction]
    \label{lem:elliptic-closed-reduction-count}
    Apply the reduction of
    Lemma~\ref{lem:elliptic-closed-cancellation},
    and let $r_S$ be the integers in its statement.
    Then
    \[
        D^+\leq D^-,
        \qquad
        M_0^+=M_0^--\sum_S r_S.
    \]
    The sum is over the original supporting fibers.
    If the incoming cochain is nonzero modulo
    $\Lambda_+$, then
    \[
        M_0^+\leq M_0^--1.
    \]
\end{lemma}

\begin{proof}
    Each passage of a nonclosed strand through
    a supporting fiber $S$ contributes
    $d_0(S)+d_1(S)$ counted crossings.
    By Lemma~\ref{lem:elliptic-closed-cancellation},
    each multiplicity decreases by $r_S\geq0$,
    so these contributions cannot increase.
    The nonclosed pieces are unchanged, hence
    so are their mutual crossings and
    active-return terms.
    Since closed--closed intersections are
    excluded from $D$, we obtain $D^+\leq D^-$.

    Moreover,
    \[
        M_0^+
        =\sum_S d_0^+(S)
        =\sum_S\bigl(d_0^-(S)-r_S\bigr)
        =M_0^--\sum_S r_S.
    \]
    If the incoming cochain is nonzero modulo
    $\Lambda_+$, the same lemma gives
    $\sum_S r_S\geq1$, proving the final assertion.
\end{proof}

\begin{corollary}[Count across a sequence of operations]
    \label{cor:elliptic-group-count}
    Consider a finite sequence of operations at
    initialization or at an event time.
    First continue any boundary collisions and
    delete any spatially isolated paired collapses,
    as in Lemma~\ref{lem:elliptic-geometric-count}.
    Then perform elementary wall crossing operations
    chosen by
    Proposition~\ref{prop:elliptic-surgery-smoothing}
    and closed-component reductions supplied by
    Lemma~\ref{lem:elliptic-closed-cancellation}.

    Let $r$ be the number of elementary wall
    crossing operations plus the number of
    closed-component reductions whose incoming
    cochain is nonzero modulo $\Lambda_+$.
    Each such reduction is counted once,
    regardless of how many identity pairs it cancels.
    Write $D^\pm$ and $M_0^\pm$ for the incoming
    and outgoing counts, and set
    \[
        N^\pm:=D^\pm+M_0^\pm.
    \]
    Let $E^\pm$ be the numbers of labelled endpoints
    of the active and inactive intervals, counted
    on their sources.
    Then
    \[
        N^+\leq N^- -r,
        \qquad
        E^+\leq E^-+4r.
    \]
    If at least one boundary collision or paired
    collapse is continued and $r=0$, then
    \[
        E^+\leq E^- -2.
    \]
\end{corollary}

\begin{proof}
    Geometric continuations do not increase $D$,
    by Lemma~\ref{lem:elliptic-geometric-count},
    and leave $M_0$ unchanged: any new closed
    fragment has label $1$, by
    Remark~\ref{rem:elliptic-boundary-collision-fragments}.
    Thus they do not increase $N$.
    Each boundary collision or paired collapse
    removes at least two labelled endpoints.

    By Lemma~\ref{lem:elliptic-wall-crossings},
    each elementary wall crossing decreases $D$
    by at least one and adds at most four endpoints.
    Since it does not increase $M_0$, it decreases
    $N$ by at least one.
    By Lemma~\ref{lem:elliptic-closed-reduction-count},
    closed-component reduction leaves $E$
    unchanged and does not increase $N$.
    It decreases $N$ by at least one whenever
    its incoming cochain is nonzero modulo
    $\Lambda_+$.

    Summing these estimates over the sequence
    gives the first two inequalities.
    If $r=0$, no operation adds endpoints,
    so any boundary collision or paired collapse
    gives the final inequality.
\end{proof}

\subsection{Continuation and termination}
\label{subsec:elliptic-continuation-termination}

\subsubsection{Initialization and continuation}
\label{subsubsec:elliptic-initialization-continuation}

We first show that the reductions required at a
fixed time can be completed after finitely many
operations. Set $N:=D+M_0$.

\begin{lemma}[Finite reduction at a fixed time]
    \label{lem:elliptic-finite-reduction}
    Fix a time $\tau$ and provisional cone data
    arising during initialization or continuation,
    as in
    Proposition~\ref{prop:elliptic-surgery-smoothing}.
    Assume that the underlying web admits regular
    geometric evolution for a short time after
    $\tau$.
    Let $N_{\mathrm{start}}$ be the initial
    value of $N$.

    After at most $N_{\mathrm{start}}$ obstruction
    reductions, the resulting data define an
    unobstructed brane with strictly positive
    bounding-cochain valuation, and regular
    brane transport can resume.
    The completed procedure supplies the brane
    transfer of
    Definition~\ref{def:brane-surgery-transfer},
    or, at initialization, the comparison of
    Definition~\ref{def:initial-brane-comparison}.
\end{lemma}

\begin{proof}
    Whenever a valuation-zero coefficient remains,
    Proposition~\ref{prop:elliptic-surgery-smoothing}
    supplies an elementary wall crossing operation,
    or Lemma~\ref{lem:elliptic-closed-cancellation}
    supplies a closed-component reduction.
    The resulting data again satisfy the
    hypotheses needed to continue this procedure.

    Each such reduction decreases $N$ by at
    least one, by
    Corollary~\ref{cor:elliptic-group-count}.
    Thus, after $k$ reductions, its current
    value $N_k$ satisfies
    \[
        0\leq N_k\leq N_{\mathrm{start}}-k.
    \]
    More than $N_{\mathrm{start}}$ reductions
    are therefore impossible.
    Since another reduction is available whenever
    a valuation-zero coefficient remains, the
    procedure must reach data with no such
    coefficients after at most $N_{\mathrm{start}}$ reductions.
    The outgoing
    cochain consequently has strictly positive
    valuation.
\end{proof}

We now apply this reduction procedure at
initialization and after the geometric operations
at an event.

\begin{proof}[Proof of
        Theorem~\ref{thm:intro-elliptic-continuation}]
    At $\tau=0$,
    Proposition~\ref{prop:elliptic-surgery-smoothing}
    applies to the initial cone data.
    Lemma~\ref{lem:elliptic-finite-reduction}
    resolves them in finitely many steps,
    giving an outgoing brane of strictly positive
    valuation and its comparison with
    $\operatorname{Cone}(f)$.
    Regular evolution and brane transport
    therefore give a special Lagrangian web for $f$
    on some interval $(0,\epsilon)$.

    Follow regular evolution toward the next
    event, using
    Lemma~\ref{lem:elliptic-harmonic-flux}
    and
    Proposition~\ref{prop:elliptic-exact-continuation}.
    If the event is geometric, apply
    Lemma~\ref{lem:elliptic-finite-preparation}
    when isolation is needed.
    The prepared evolution reaches an event
    without further preparation, possibly an
    earlier obstruction wall.
    Any collapse reached there is either
    spatially isolated and paired, or part of
    a boundary collision with a positive-length
    inactive interval.

    At the event, concatenate the boundary
    collisions and delete the spatially isolated
    paired collapsing components, using
    Proposition~\ref{prop:elliptic-boundary-continuation}
    and
    Remark~\ref{rem:elliptic-geometric-continuation}.
    These geometric operations are finite,
    since each removes labelled endpoints.
    They leave provisional cone data whose
    geometric pieces admit regular outgoing
    evolution.

    Apply Lemma~\ref{lem:elliptic-finite-reduction}
    to resolve any remaining valuation-zero
    terms in finitely many steps.
    The completed event supplies the comparison
    between incoming and outgoing regular branes,
    and regular evolution and brane transport
    resume.
    This includes simultaneous geometric
    degenerations and obstruction walls.

    The initial and event comparisons give the
    required special Lagrangian web for $f$.
    Lemma~\ref{lem:web-cone-comparison}
    then identifies every regular cone brane
    with $\operatorname{Cone}(f)$.
\end{proof}

\subsubsection{Finiteness of event times}
\label{subsubsec:elliptic-finiteness}

Let $N_{\mathrm{in}}$ and $E_{\mathrm{in}}$
be the values of $N$ and the number of labelled
endpoints before initialization.
For the initial configuration of
Definition~\ref{def:initial-web-configuration},
we have $E_{\mathrm{in}}=0$.

\begin{prop}[Bounds on reductions and event times]
    \label{prop:elliptic-event-bounds}
    Consider a finite initial segment of the
    construction in the proof of
    Theorem~\ref{thm:intro-elliptic-continuation},
    ending at a regular time or after completing
    all operations at its final event.

    Let $R$ count the obstruction reductions,
    including those at initialization.
    Among subsequent geometric degenerations
    and obstruction walls, let $W$ count the
    event times involving at least one obstruction
    reduction, and let $G$ count the remaining
    event times.
    Then
    \begin{equation}
        \label{eq:elliptic-event-bounds}
        W\leq R\leq N_{\mathrm{in}},
        \qquad
        G\leq \frac{E_{\mathrm{in}}}{2}+2R.
    \end{equation}
    In particular,
    \[
        G+W
        \leq \frac{E_{\mathrm{in}}}{2}
        +3N_{\mathrm{in}}.
    \]
\end{prop}

\begin{proof}
    Regular evolution and isolation perturbations
    leave $M_0$ and the endpoint count unchanged.
    They do not increase $D$, by
    Corollary~\ref{cor:elliptic-crossing-monotonicity}
    and Remark~\ref{rem:elliptic-isolation-count}.
    Thus $N$ does not increase between events.

    Write $N_{\mathrm{out}}$ and $E_{\mathrm{out}}$
    for the final counts.
    Summing the estimates of
    Corollary~\ref{cor:elliptic-group-count}
    over initialization and the subsequent events
    gives
    \[
        N_{\mathrm{out}}+R\leq N_{\mathrm{in}},
        \qquad
        E_{\mathrm{out}}+2G
        \leq E_{\mathrm{in}}+4R.
    \]
    Indeed, every obstruction reduction decreases
    $N$ by at least one and adds at most four
    endpoints, while each event counted by $G$
    removes at least two endpoints.

    Since $N_{\mathrm{out}},E_{\mathrm{out}}\geq0$
    and every event counted by $W$ contains a
    reduction, these inequalities give all
    the stated bounds.
\end{proof}

\subsubsection{A bound on active time}
\label{subsubsec:elliptic-active-time}

The bound on $N$ controls the lengths of active
intervals. Together with
\eqref{eq:elliptic-speed-energy}, this gives
a uniform positive lower bound for their
Dirichlet energy.

\begin{prop}[Bound on active time]
    \label{prop:elliptic-active-time-bound}
    Every finite initial segment of the construction
    in the proof of
    Theorem~\ref{thm:intro-elliptic-continuation}
    satisfies
    \begin{equation}
        T_{\mathrm{act}}
        \leq (N_{\mathrm{in}}+1)r_{\max},
        \label{eq:elliptic-active-time-bound}
    \end{equation}
    where $T_{\mathrm{act}}$ is its total active
    time and $r_{\max}$ is defined in
    Equation~\eqref{eq:web-R-bound}.
\end{prop}

\begin{proof}
    Fix a regular time $\tau$ and an active
    interval $P$.
    Proposition~\ref{prop:elliptic-event-bounds},
    applied to the initial segment ending at
    $\tau$, gives
    \[
        D(\tau)\leq N(\tau)\leq N_{\mathrm{in}}.
    \]
    In an outgoing counting diagram $\mathcal G$
    at $\tau$, the definition of the return
    count gives
    \[
        \ell_{\mathcal G}(P)
        \leq o_{\mathcal G}(P)+1
        \leq D(\tau)+1
        \leq N_{\mathrm{in}}+1.
    \]
    Letting the time advances and auxiliary
    separations tend to zero gives the same
    bound for the actual interval:
    \[
        \ell(P)\leq N_{\mathrm{in}}+1.
    \]
    By Lemma~\ref{lem:elliptic-harmonic-flux},
    \[
        \mathcal E(P)
        =\frac{1}{\ell(P)}
        \geq\frac{1}{N_{\mathrm{in}}+1}.
    \]
    The isolation perturbations leave the swept
    real integral unchanged and preserve the
    boundary real periods, hence $r_{\max}$.
    Thus the integration argument of
    Corollary~\ref{cor:finite-active-time}
    applies with
    $\epsilon=(N_{\mathrm{in}}+1)^{-1}$,
    proving \eqref{eq:elliptic-active-time-bound}.
\end{proof}

Since the estimate holds for every finite initial
segment, it also bounds the total active time
over the entire domain of the construction.

\subsubsection{Finite-time termination}
\label{subsubsec:elliptic-termination}

It remains to control the time during which
the active front is empty but the evolution
is not yet terminal.
We use the stationary continuation of
Definition~\ref{def:stationary-brane-continuation}.

\begin{proof}[Proof of
        Theorem~\ref{thm:intro-elliptic-termination}]
    Follow the construction of
    Theorem~\ref{thm:intro-elliptic-continuation},
    using stationary brane continuation whenever
    the active front is empty.
    Proposition~\ref{prop:elliptic-event-bounds}
    gives finitely many geometric degenerations
    and obstruction walls.
    There are also finitely many isolation
    perturbations, by
    Lemma~\ref{lem:elliptic-finite-preparation}.
    In particular, these times cannot accumulate.

    As explained in
    Section~\ref{subsubsec:web-terminality},
    every nonterminal stationary interval reaches
    an obstruction wall after a finite waiting
    time.
    Since there are finitely many walls, the
    total nonterminal stationary time is finite.
    The total active time is finite by
    Proposition~\ref{prop:elliptic-active-time-bound}.
    Thus the construction cannot remain
    nonterminal for an infinite amount of time.
    Nor can it stop at a finite nonterminal time:
    regular evolution extends until a geometric
    degeneration classified by
    Lemma~\ref{lem:elliptic-geometric-degenerations}
    or an obstruction wall, and
    Theorem~\ref{thm:intro-elliptic-continuation}
    supplies continuation across that event,
    including simultaneous occurrences.

    Therefore the construction becomes terminal
    in finite time, after finitely many geometric
    degenerations and obstruction walls.
    The resulting evolution is a special Lagrangian web for $f$
    by construction.
\end{proof}

\subsection{Examples}
\label{subsec:elliptic-examples}

\subsubsection{Puncture and boundary collisions}

The following degree-three example extends
Example~\ref{ex:intro-degree-two} by showing how
varying the morphism coefficients can replace
webbing by puncture, or allow a boundary
collision to occur before the obstruction wall
is reached.

\begin{example}
    \label{ex:elliptic-degree-three}
    Consider the square torus
    \[
        X=\mathbb R^2/\mathbb Z^2,
        \qquad
        L_0=\{y=0\},
        \qquad
        L_1=\{y=3x\},
    \]
    with the standard Calabi--Yau structure.
    Equip the curves with their positive gradings,
    trivial rank-one local systems, and zero bounding cochains.
    Denote the resulting branes by $\mathcal L_0$
    and $\mathcal L_1$.
    In the fundamental domain $[-1/2,1/2)^2$,
    their three intersections
    \[
        p_0=(0,0),
        \qquad
        p_-=\left(-\frac13,0\right),
        \qquad
        p_+=\left(\frac13,0\right)
    \]
    have Maslov index zero.
    Consider the morphism
    \[
        f=p_0-T^{l_--\frac16}p_-
        -T^{l_+-\frac16}p_+
        \in HF^0(\mathcal L_0,\mathcal L_1),
        \qquad
        \frac16<l_+<\frac13,
        \qquad
        l_+<l_-<\frac12.
    \]
    The shifts by $\frac16$ make $l_\pm$ agree
    with the valuations of the corresponding
    theta coefficients on the mirror curve.
    We fix $l_+$ and vary $l_-$.
    Write $F_c=\{x=c\}$ for the vertical circle
    at horizontal coordinate $c$.

    The coefficient at $p_0$ has valuation zero.
    Following the initialization procedure of
    Section~\ref{subsubsec:elliptic-initialization-continuation},
    we resolve it by initiating webbing at $p_0$
    at time $\tau=0$.
    Each of its two fronts sweeps out a triangle
    between the lifted lines $y=0$ and $y=3x$.
    By Lemma~\ref{lem:elliptic-harmonic-flux},
    the swept area equals the elapsed time:
    \[
        \tau=\int_0^x3u\,du=\frac32x^2.
    \]
    While $p_\pm$ remains an intersection of the
    unswept regions,
    Equation~\eqref{eq:web-cone-valuation-transport}
    gives
    \[
        \operatorname{val}(b_{\tau,p_\pm})
        =l_\pm-\frac16-\tau.
    \]
    Thus the coefficient at $p_+$ reaches
    valuation zero at $\tau=l_+-\frac16$,
    before the right-moving front from $p_0$
    reaches that intersection.
    We resolve this wall by the webbing of
    Definition~\ref{def:elliptic-webbing}.
    The fronts moving towards one another from
    $p_0$ and $p_+$ meet where their swept areas
    differ by their initiation times:
    \[
        \frac32x^2
        -\frac32\left(\frac13-x\right)^2
        =l_+-\frac16,
        \qquad x=l_+.
    \]
    Boundary surgery joins them into $F_{l_+}$.
    This occurs in all three cases below.
    Varying $l_-$ changes whether webbing begins
    at $p_-$ and how the remaining fronts meet.

    Suppose first that
    \[
        l_+<l_-<\frac13.
    \]
    By the same valuation-transport formula,
    the coefficient at $p_-$ reaches valuation
    zero at time
    $\tau=l_--\frac16<\frac16$.
    The left-moving front from $p_0$ has not yet
    reached $x=-\frac13$, so $p_-$ still lies
    in the unswept regions.
    We therefore initiate a third webbing there,
    again using Definition~\ref{def:elliptic-webbing}.
    Its right-moving front meets the left-moving
    front from $p_0$, while its left-moving front
    meets the right-moving front from $p_+$,
    across the boundary of the fundamental domain.
    Boundary surgery at these collisions produces
    the remaining two circles.
    The terminal configuration is
    \[
        F_{l_+},
        \qquad
        F_{-l_-},
        \qquad
        F_{-\frac12+l_--l_+}.
    \]
    This is the same mechanism as in
    Example~\ref{ex:intro-degree-two}:
    each intersection initiates its own webbing
    before being reached by another front.

    \begin{figure}[htbp]
        \centering
        \includegraphics[width=\textwidth]{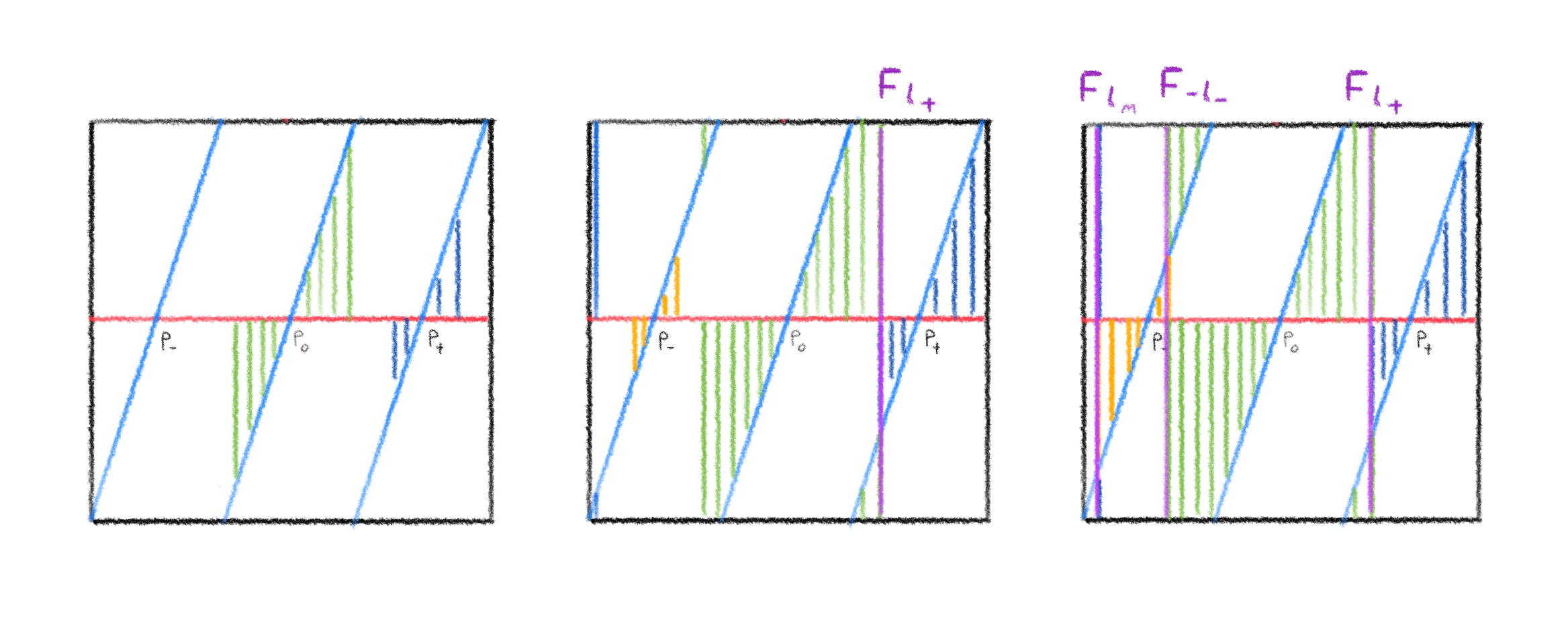}
        \caption{
            Three stages of the example
            in the regime $l_+<l_-<\frac13$.
            The webbings initiated at $p_0$, $p_+$,
            and $p_-$ are shown in green, blue,
            and orange, respectively.
            Closed vertical fibers are shown in purple.
            Left: webbing has begun at $p_0$ and $p_+$,
            but their fronts have not yet met.
            Middle: boundary collision has produced
            $F_{l_+}$, and webbing has begun at $p_-$.
            Right: the remaining fronts have met in
            pairs, producing the terminal configuration
            $F_{l_m}$, $F_{-l_-}$,
            and $F_{l_+}$, listed from left to right, where $l_m =-\frac12+l_--l_+$.
        }
        \label{fig:elliptic-degree-three-first-regime}
    \end{figure}

    Now suppose that
    \[
        \frac13<l_-<\frac14+\frac{l_+}{2}.
    \]
    At time $\tau=\frac16$, the left-moving front
    from $p_0$ reaches $x=-\frac13$.
    The coefficient associated with $p_-$ still
    has positive valuation
    \[
        \left(l_--\frac16\right)-\frac16
        =l_--\frac13>0.
    \]
    Thus no webbing is initiated there.
    Instead, the existing front passes through
    the intersection and continues to the left.
    Its length is now greater than one, so it
    winds around the vertical fiber and returns
    to its boundary endpoint.

    To follow the corresponding coefficient,
    lift the front to the interval $[3x,0]$.
    Its height function is
    $\kappa(y)=y/(3x)$, so the return at $y=-1$
    has height $-1/(3x)$.
    Writing $v$ for the coefficient's valuation,
    Equation~\eqref{eq:web-cone-valuation-transport}
    and Lemma~\ref{lem:elliptic-harmonic-flux}
    give
    \[
        \frac{dv}{d\tau}
        =\frac1{3x}
        =\frac{dx}{d\tau}.
    \]
    Since $v=l_--\frac13$ at $x=-\frac13$,
    it follows that
    \[
        v=l_-+x.
    \]
    The coefficient therefore reaches valuation
    zero when the front reaches $x=-l_-$.

    This wall is resolved by the active-return
    puncture of
    Remark~\ref{rem:elliptic-overlap-puncture}.
    Puncture followed by boundary reconnection
    splits off $F_{-l_-}$ and leaves an active
    front shorter by one vertical period.
    That front continues to the left and meets
    the remaining front from $p_+$, producing
    $F_{-\frac12+l_--l_+}$.
    The terminal fibers therefore have the same
    formulas as in the first case, but
    $F_{-l_-}$ now forms by puncture after the
    front has passed through $p_-$.

    \begin{figure}[htbp]
        \centering
        \includegraphics[width=0.8\textwidth]{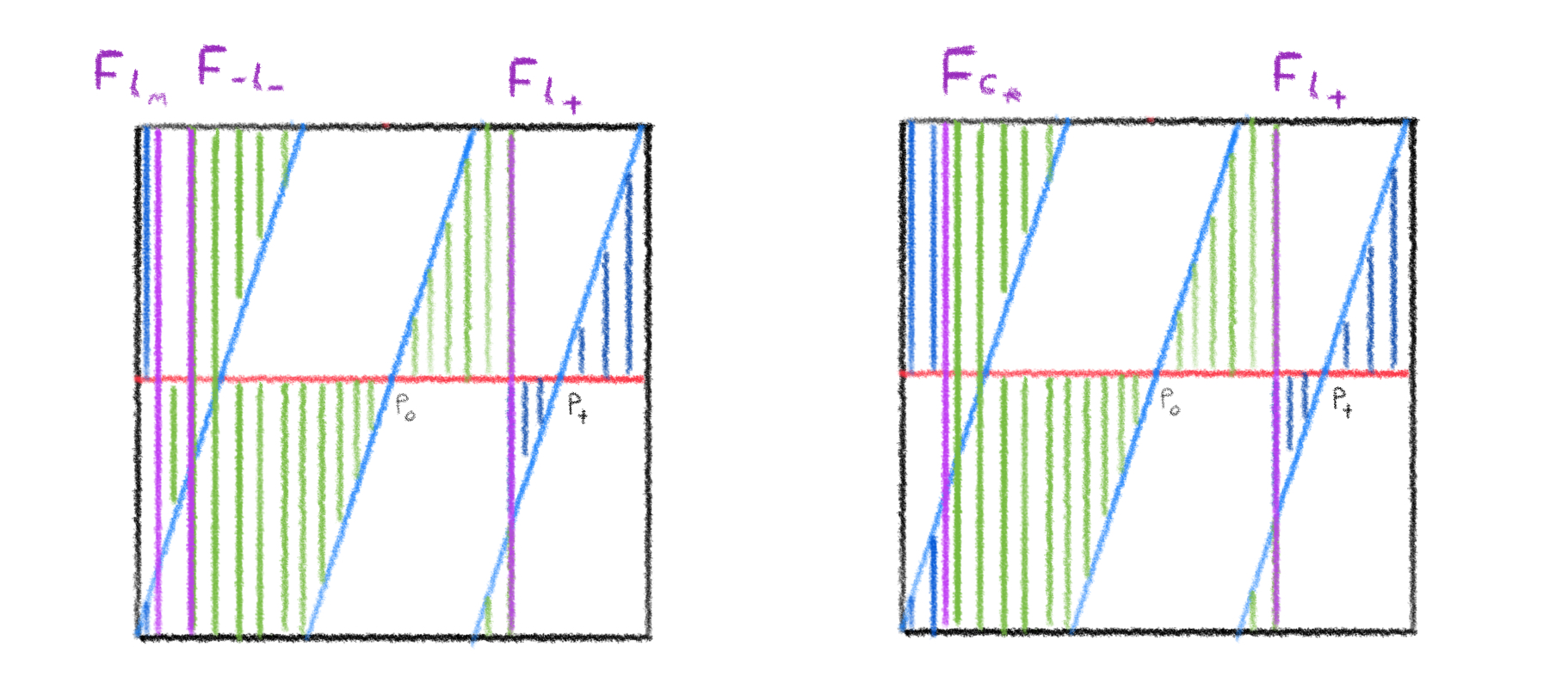}
        \caption{
            Terminal configurations in the two regimes
            where the left-moving front from $p_0$
            passes through $p_-$ without initiating
            webbing there.
            Left: for
            $\frac13<l_-<\frac14+\frac{l_+}{2}$,
            puncture followed by reconnection splits
            off $F_{-l_-}$.
            The remaining fronts subsequently meet,
            producing $F_{l_m}$, where
            $l_m=-\frac12+l_--l_+$.
            Right: for
            $\frac14+\frac{l_+}{2}<l_-<\frac12$,
            boundary collision occurs before the obstruction
            wall is reached, producing a double cover
            of $F_{c_*}$, where
            $c_*=-\frac14-\frac{l_+}{2}$.
            No puncture occurs.
            In both cases, $F_{l_+}$ has already formed.
            The webbings from $p_0$ and $p_+$ are shown
            in green and blue, respectively;
            terminal fibers are shown in purple.
        }
        \label{fig:elliptic-degree-three-puncture-collision}
    \end{figure}

    Finally, suppose that
    \[
        \frac14+\frac{l_+}{2}<l_-<\frac12.
    \]
    The left-moving front again passes through
    $p_-$.
    Using the lift of $p_+$ at $x=-\frac23$,
    its collision with the remaining front
    from $p_+$ is determined by
    \[
        \frac32x^2
        -\frac32\left(x+\frac23\right)^2
        =l_+-\frac16,
        \qquad
        x=c_*:=-\frac14-\frac{l_+}{2}.
    \]
    At this collision, the coefficient that
    could have caused puncture has valuation
    \[
        l_-+c_*
        =l_--\frac14-\frac{l_+}{2}>0.
    \]
    Thus the boundary collision occurs before
    the obstruction wall.
    The two fronts have total length two,
    and boundary surgery joins them into a
    double cover of $F_{c_*}$.
    No puncture occurs.

    The positive-valuation contribution associated
    with $p_-$ is retained in the transferred
    brane data.
    At this point the residues are exhausted,
    and every remaining component has label $1$.
    The stationary height is therefore identically
    $1$, so all height differences in
    Equation~\eqref{eq:web-cone-valuation-transport}
    vanish.
    The remaining coefficients stay constant,
    and no later obstruction wall occurs.

    The brane on the double cover is equivalent
    to two rank-one branes on $F_{c_*}$ with
    distinct local systems, whose holonomies
    incorporate the transported bounding cochain.
    Together with the rank-one brane on $F_{l_+}$,
    these form the terminal brane.
    Thus varying $l_-$ in this regime changes
    the local systems while leaving the
    geometric support fixed.

    In all three cases, the terminal components
    have label $1$.
    Hence the kernel of $f$ is zero, and these
    branes represent its cokernel.

    On the mirror elliptic curve, $f$ corresponds
    to the section
    \[
        s=\theta_0-T^{l_-}\theta_-
        -T^{l_+}\theta_+
    \]
    of $\mathcal O(3p)$, where $p$ is fixed by
    the mirror normalization.
    Here $\theta_0,\theta_-,\theta_+$ are the theta-basis
  sections normalized by $p_0\leftrightarrow\theta_0$
  and $p_\pm\leftrightarrow T^{1/6}\theta_\pm$.
    The associated morphism
    \[
        \mathcal O\xrightarrow{s}\mathcal O(3p)
    \]
    is injective, with cokernel supported on
    the zero divisor of $s$.
    In the three strict regimes considered here,
    this divisor consists of three distinct points.
    In the first two cases, their valuations are
    \[
        l_+,\qquad -l_-,
        \qquad -\frac12+l_--l_+,
    \]
    recovering the three terminal fibers.
    In the third case, two zeros have the same
    valuation $c_*$ but remain distinct points
    of the mirror curve.
    Their unit coordinates specify the two
    different local systems on $F_{c_*}$.
    Thus the mirror cokernel retains the
    dependence on $l_-$ even when the terminal
    supporting fibers no longer move.
\end{example}

\subsubsection{Stationary geometry and absorption}

The following example shows that geometric
exhaustion need not end the brane evolution.
Morphisms with the same coefficient valuations
reach the same stationary configuration, but
their units determine whether a later absorption
removes the prospective cokernel.

\begin{example}
    \label{ex:elliptic-stationary-absorption}
    On the square torus
    $X=\mathbb R^2/\mathbb Z^2$, let $L_0$ and
    $L_1$ be the images of the lines
    \[
        y=\frac{x}{2},
        \qquad
        y=2x,
    \]
    respectively.
    Equip them with their positive gradings,
    trivial rank-one local systems, and zero
    bounding cochains.
    Their homology classes are
    \[
        [L_0]=(2,1),
        \qquad
        [L_1]=(1,2),
    \]
    and their three intersections
    \[
        p_0=(0,0),
        \qquad
        p_+=\left(\frac13,-\frac13\right),
        \qquad
        p_-=\left(-\frac13,\frac13\right)
    \]
    have Maslov index zero.
    Fix
    \[
        0<\lambda<\frac1{12}
    \]
    and consider
    \[
        f_u=p_0-T^\lambda p_+
        -T^{\frac13-\lambda}u\,p_-,
        \qquad
        u\in1+\Lambda_+.
    \]
    We vary only the unit $u$, leaving all
    coefficient valuations unchanged.

    Webbing begins at $p_0$ at time $\tau=0$.
    Since the difference of the two slopes is
    $\frac32$, its fronts satisfy
    \[
        \tau=\int_0^x\frac32s\,ds
        =\frac34x^2.
    \]
    At time $\tau=\lambda$, the coefficient
    at $p_+$ reaches valuation zero and a
    second webbing begins there.

    The right-moving front from $p_0$ and
    the left-moving front from $p_+$ meet at
    \[
        \frac34x^2
        -\frac34\left(\frac13-x\right)^2
        =\lambda,
        \qquad
        x=d:=\frac16+2\lambda.
    \]
    Their endpoints on $L_1$ coincide, but
    their endpoints on $L_0$ lie on its two
    different sheets over $x=d$.
    Boundary surgery therefore joins them
    into an inactive interval $J$ of length
    $\frac12$, with both endpoints labeled $0$.
    It remains attached to the two residue
    branches of $L_0$.

    \begin{figure}[htbp]
        \centering
        \includegraphics[width=0.8\textwidth]{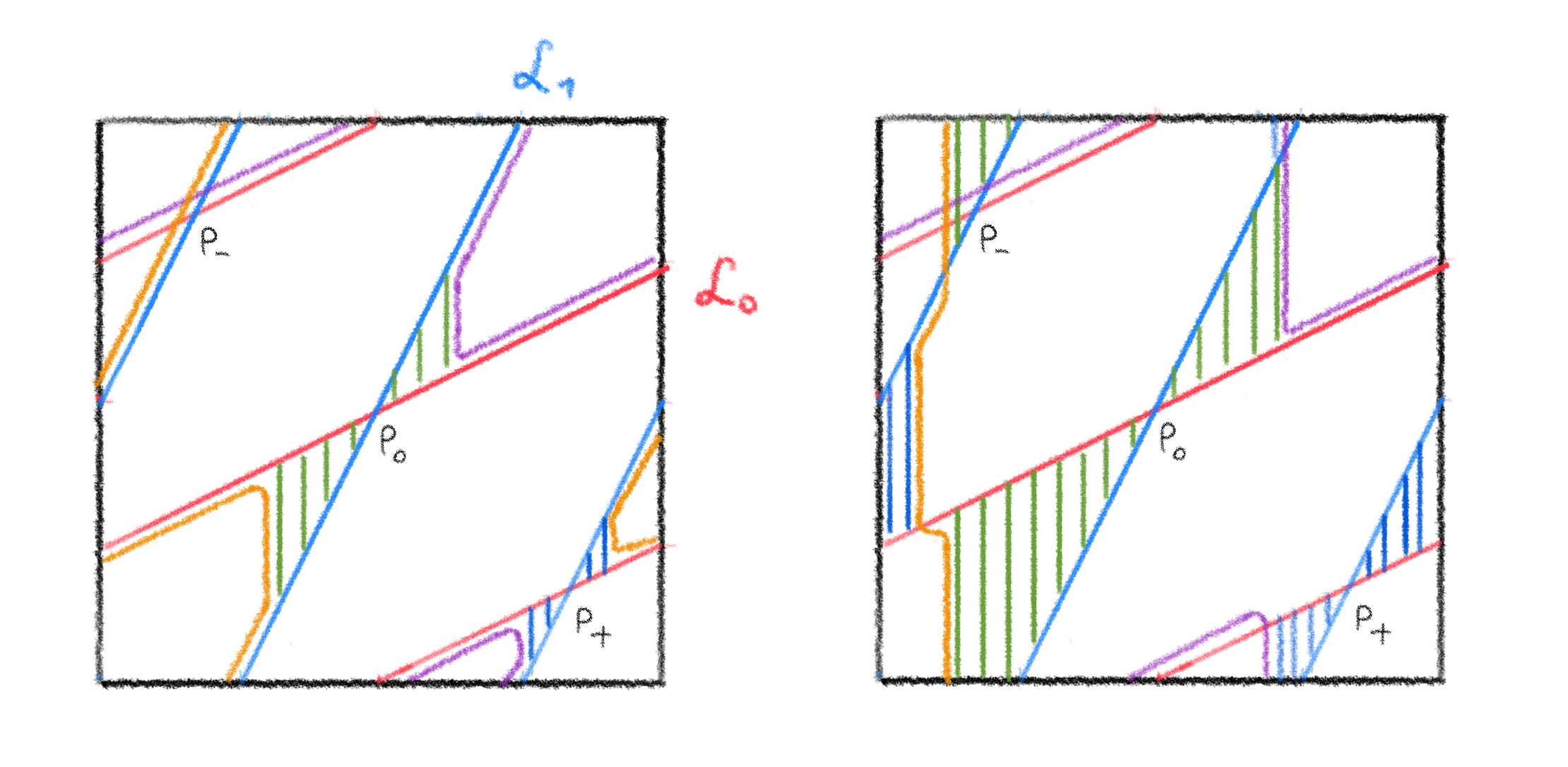}
        \caption{
            Webbing and the first boundary collision in
            Example~\ref{ex:elliptic-stationary-absorption}.
            Left: webbings have been initiated at $p_0$
            and $p_+$.
            The web-cone has two connected components,
            shown in purple and orange.
            Right: the right-moving front from $p_0$
            meets the left-moving front from $p_+$.
            Boundary surgery produces an inactive
            interval with both endpoints on $L_0$,
            shown in purple.
            The web-cones are slightly displaced
            for clarity.
        }
        \label{fig:elliptic-absorption-first-collision}
    \end{figure}

    The left-moving front from $p_0$ reaches
    $p_-$ at time $\tau=\frac1{12}$.
    No webbing is initiated there.
    With the standard local brane transfer,
    the remaining coefficient stays of
    positive valuation as this front passes
    through $p_-$ and continues towards the
    remaining front from $p_+$.
    Using the lift of $p_+$ at $x=-\frac23$,
    these two fronts meet at
    \[
        \frac34x^2
        -\frac34\left(x+\frac23\right)^2
        =\lambda,
        \qquad
        x=c:=-\frac13-\lambda.
    \]
    Their lengths add to one, and boundary
    surgery joins them into the closed
    label-one circle $S=F_c$.

    No active fronts remain.
    The residue $R_1$ is empty, while $J$
    remains attached to $R_0$.
    The curve
    \[
        K_{\mathrm{sing}}= \overline R_0\cup J[-1]
    \]
    admits a positive smoothing $K$.
    Since
    \[
        -[K]+[S]=[L_1]-[L_0]=(-1,1),
    \]
    we have $[K]=(1,0)$.
    Thus the stationary cone has geometric
    components
    \[
        K[1]\cup S.
    \]
    They meet once, at the intersection of
    $S$ with the remaining residue branch
    of $L_0$.

    \begin{figure}[htbp]
        \centering
        \includegraphics[width=0.8\textwidth]{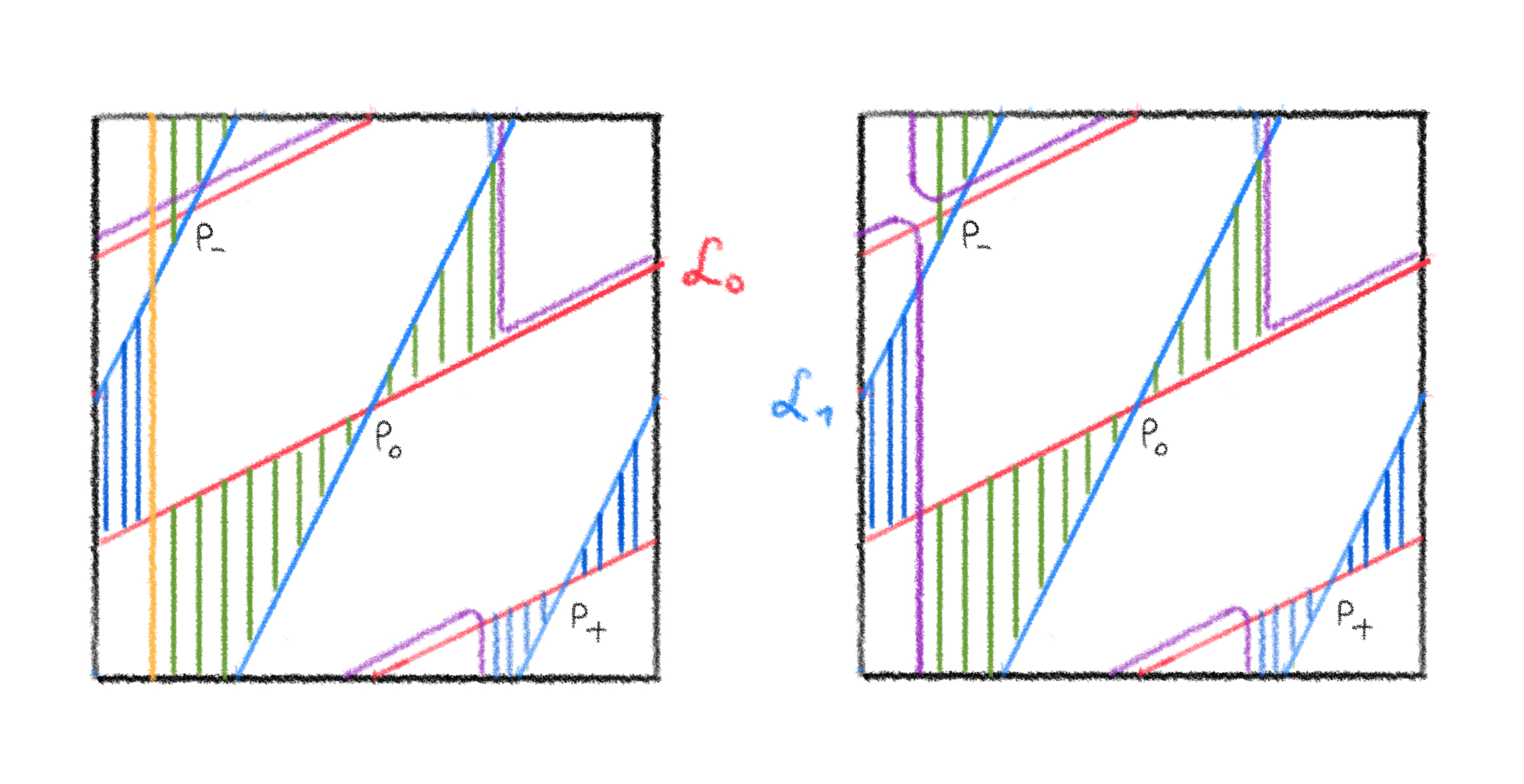}
        \caption{
            Stationary continuation and absorption in
            Example~\ref{ex:elliptic-stationary-absorption}.
            Left: boundary surgery at the collision
            of the left-moving front from $p_0$
            and the right-moving front from $p_+$
            produces a closed imaginary-special
            fragment, shown in orange.
            No active fronts remain, and $R_1$ is empty.
            The purple component consists of the
            remaining residue of $L_0$ joined by
            the inactive interval.
            When the coefficient coupling the two
            components vanishes, this configuration
            is terminal, with kernel of class $(1,0)$
            and cokernel of class $(0,1)$.
            Right: when that coefficient is nonzero,
            stationary continuation reaches an
            obstruction wall.
            Absorption joins the two components,
            yielding a single kernel of class $(1,-1)$
            and zero cokernel.
            The web-cones are slightly displaced
            for clarity.
        }
        \label{fig:elliptic-stationary-absorption}
    \end{figure}
    The geometry alone does not yet determine
    the kernel and cokernel.
    Let $\beta(u)$ denote the coefficient
    at this crossing in the transported
    bounding cochain, at a regular time
    $\tau_*$ just after the final collision.
    Although its valuation is positive,
    the coefficient need not vanish.
    Its dependence on $u$ distinguishes
    two outcomes.

    For a distinguished unit
    $u_*\in1+\Lambda_+$, the coefficient
    $\beta(u_*)$ vanishes exactly.
    The two components are then uncoupled,
    and stationary continuation produces
    no further obstruction wall.
    The configuration is terminal:
    $K$ represents the kernel, of class
    $(1,0)$, and $S$ represents the cokernel,
    of class $(0,1)$.

    For a general choice of $u$, the coefficient
    $\beta(u)$ is nonzero.
    The stationary height is $0$ on $K[1]$
    and $1$ on $S$, so
    Equation~\eqref{eq:web-cone-valuation-transport}
    gives
    \[
        \operatorname{val}(b_{\tau,q})
        =
        \operatorname{val}(\beta(u))
        -(\tau-\tau_*),
    \]
    where $q$ is the surviving crossing.
    The geometry remains fixed, but this
    coefficient reaches valuation zero after
    a finite waiting time.

    We resolve the resulting wall by absorption
    into the boundary, as in
    Definition~\ref{def:elliptic-absorption}.
    The closed component $S$ is absorbed into
    the residue of $L_0$.
    On the positive smoothing, this replaces
    $K$ by a curve $K'$ obtained by joining
    it to $S[-1]$.
    Consequently,
    \[
        [K']=[K]-[S]=(1,-1).
    \]
    No closed label-one component remains,
    and the resulting configuration is terminal.
    It represents a kernel of class $(1,-1)$
    and zero cokernel.

    The mirror description explains the
    exceptional cancellation.
    The input branes correspond to a rank-two
    bundle $E$ of degree one and a line bundle
    $M$ of degree two.
    Locally, the mirror morphism
    \[
        f_u^\vee:E\longrightarrow M
    \]
    is a row of two functions.
    It fails to be surjective precisely where
    those functions vanish simultaneously.
    The theta-function calculation gives a
    distinguished unit
    \[
        u_*=1-T^{3\lambda}+r,
        \qquad
        \operatorname{val}(r)>3\lambda,
    \]
    for which the two functions have a common
    zero $z_*$ lying over $x=c$.
    In this case,
    \[
        \operatorname{coker}(f_{u_*}^\vee)
        \cong\mathcal O_{z_*},
        \qquad
        \ker(f_{u_*}^\vee)
        \cong\det(E)\otimes M^{-1}(z_*).
    \]
    The kernel has degree zero.
    For a general unit, including $u=1$,
    there is no common zero, and
    \[
        \operatorname{coker}(f_u^\vee)=0,
        \qquad
        \ker(f_u^\vee)
        \cong\det(E)\otimes M^{-1},
    \]
    of degree $-1$.

    In particular, $u=u_*$ and $u=1$ give
    morphisms with the same valuations and
    the same leading coefficients.
    Their initial geometric evolution agrees,
    but only the exact cancellation for $u=u_*$
    allows the closed component to survive
    as a cokernel.
\end{example}

\section{An example on an abelian surface}
\label{sec:abelian-surface}

We consider two product tori meeting at a single
point.
The mirror cokernel is an extension of two
sheaves, whose geometric counterparts are
two tori meeting once.
This gives a simple prediction for the
terminal cokernel: a connected surface of
genus two, with a prescribed intersection
pattern with the input tori.
We first describe these predictions and then
compare them with a proposed web evolution.

\subsection{The morphism and its mirror}
\label{subsec:abelian-product-morphism}

Let $E_j=\mathbb R^2/\mathbb Z^2$, with coordinates
$z_j=x_j+iy_j$, and equip $X=E_1\times E_2$ with
\[
    \omega=dx_1\wedge dy_1+dx_2\wedge dy_2,
    \qquad
    \Omega=e^{-i\alpha}dz_1\wedge dz_2,
    \qquad
    \alpha=\arctan(1/2).
\]
Write $L_k^j\subset E_j$ for the image of
$y_j=kx_j$.
Equip these curves with their standard gradings,
trivial rank-one local systems over the complex
Novikov field $\Lambda$, and zero bounding
cochains, and denote the resulting branes
by $\cc{L}_k^j$.
Consider
\[
    \cc{L}_0=\cc{L}_0^1\boxtimes\cc{L}_0^2,
    \qquad
    \cc{L}_1=\cc{L}_1^1\boxtimes\cc{L}_1^2.
\]
Each slope-one factor has angle $\pi/4$,
so the product tori have phases $-\alpha$
and $\pi/2-\alpha$.
Both are strictly positive.
The phase of $\Omega$ is chosen so that
the two input tori have equal real periods:
\[
    \int_{L_1}\Omega-\int_{L_0}\Omega
    =e^{-i\alpha}(2i-1)
    =i\sqrt5.
\]
This gives an evolution that
exhausts both residues and ends in a closed
imaginary-special surface.

The input tori meet only at
$P=((0,0),(0,0))$, which has Maslov index zero.
Every nonzero morphism is a scalar multiple
of this generator, so we normalize
$f=P\in HF^0(\cc{L}_0,\cc{L}_1)$.

Let $\check E_1,\check E_2$ be the mirror
elliptic curves and put
$\check X=\check E_1\times\check E_2$.
The input branes correspond to
$\cc{O}_{\check X}$ and
\[
    \cc{G}
    =\cc{O}_{\check E_1}(o_1)
    \boxtimes\cc{O}_{\check E_2}(o_2),
\]
where $o_j$ is fixed by the degree-one
mirror normalization.
Let $s_j$ be a nonzero section of
$\cc{O}_{\check E_j}(o_j)$, with its single
zero at $o_j$.
The mirror morphism is the product section
\[
    f^\vee=s_1\boxtimes s_2
\]
Its zero divisor is
\[
    D=V+H,
    \qquad
    V=\{o_1\}\times\check E_2,
    \qquad
    H=\check E_1\times\{o_2\}.
\]
The two components meet transversely at
$(o_1,o_2)$.
To describe the cokernel, put
\[
    \cc{V}
    =\cc{O}_{o_1}\boxtimes\cc{O}_{\check E_2},
    \qquad
    \cc{H}
    =\cc{O}_{\check E_1}(o_1)
    \boxtimes\cc{O}_{o_2}.
\]
We use ordinary skyscraper sheaves here,
choosing bases in the one-dimensional
line-bundle fibers at the indicated points.

\begin{lemma}[The mirror cokernel]
    \label{lem:abelian-mirror-cokernel}
    \label{lem:abelian-cokernel-extension}
    The mirror morphism is injective, with
    cokernel $\cc{C}=\cc{O}_D(D)$:
    \[
        0\longrightarrow\cc{O}_{\check X}
        \xrightarrow{f^\vee}\cc{G}
        \longrightarrow\cc{C}
        \longrightarrow0.
    \]
    Here $\cc{O}_D(D)$ denotes the restriction
    of $\cc{O}_{\check X}(D)$ to $D$.
    Moreover, there is a non-split exact sequence
    \begin{equation}
        \label{eq:abelian-cokernel-extension}
        0\longrightarrow\cc{V}
        \longrightarrow\cc{C}
        \longrightarrow\cc{H}
        \longrightarrow0.
    \end{equation}
\end{lemma}

\begin{proof}
    The product section has divisor $D$, so
    $\cc{G}\cong\cc{O}_{\check X}(D)$.
    This gives the first exact sequence.

    Factor the morphism through
    $\cc{B}=\cc{O}_{\check E_1}(o_1)
        \boxtimes\cc{O}_{\check E_2}$,
    multiplying first by $s_1$ and then by $s_2$.
    The successive inclusions
    $\cc{O}_{\check X}\subset\cc{B}\subset\cc{G}$
    have quotients
    $\cc{B}/\cc{O}_{\check X}\cong\cc{V}$
    and $\cc{G}/\cc{B}\cong\cc{H}$.
    Taking the image of $\cc{B}$ in
    $\cc{C}=\cc{G}/\cc{O}_{\check X}$
    gives \eqref{eq:abelian-cokernel-extension}.

    At the crossing $(o_1,o_2)$, the fiber of
    $\cc{C}$ is one-dimensional, since
    $\cc{C}$ is a line bundle on $D$.
    A splitting as $\cc{V}\oplus\cc{H}$
    would give a two-dimensional fiber,
    with one dimension from each branch.
    Hence the extension is non-split.
\end{proof}

\begin{figure}[htbp]
    \centering
    \includegraphics[width=0.6\textwidth]{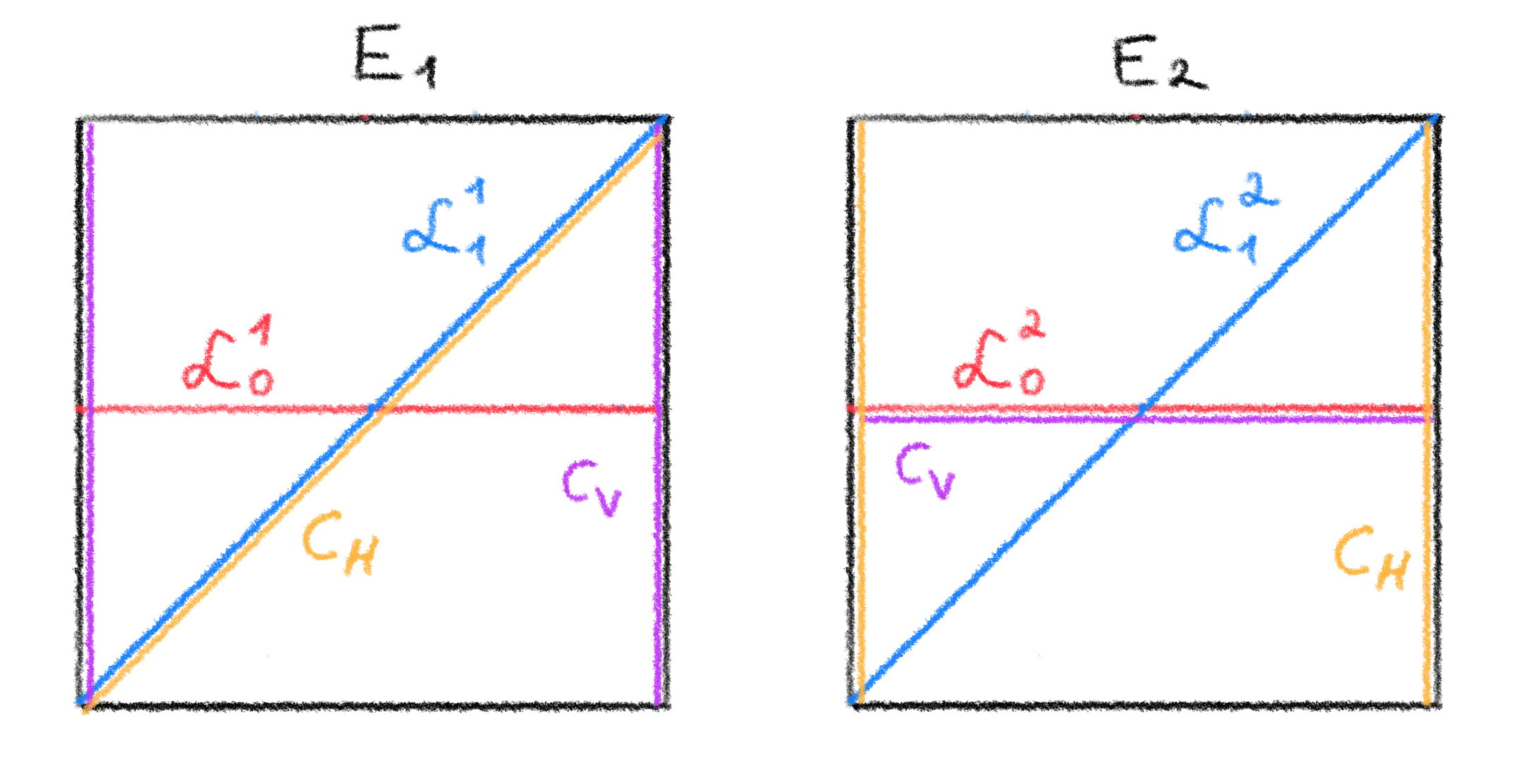}
    \caption{
        The two factors $E_1$ and $E_2$, with opposite
        edges identified.
        The curves $L_0^j$ are shown in red and
        $L_1^j$ in blue.
        The factors of
        $C_H=L_1^1\times F_{o_2}^2$ are shown in orange,
        and those of
        $C_V=F_{o_1}^1\times L_0^2$ in purple.
        Coincident curves are drawn slightly apart
        for clarity.
    }
    \label{fig:abelian-product-factors}
\end{figure}

\subsection{The predicted intersections and genus}
\label{subsec:abelian-terminal-surface}

Let $\cc{F}_{o_j}^j$ denote the vertical brane
corresponding to the skyscraper sheaf at $o_j$.
In the standard degree-one normalization,
its supporting fiber is $x_j=1/2$,
with local system determined by $o_j$.
The two sheaves in
\eqref{eq:abelian-cokernel-extension}
therefore correspond to
\[
    \cc{C}_V
    =\cc{F}_{o_1}^1\boxtimes\cc{L}_0^2,
    \qquad
    \cc{C}_H
    =\cc{L}_1^1\boxtimes\cc{F}_{o_2}^2.
\]
Their underlying tori meet transversely
at the single point
$r=((1/2,1/2),(1/2,0))$.
Under the mirror correspondence,
the nonzero extension class gives an
exact triangle in $\mathcal D$,
\[
    \cc{C}_V
    \longrightarrow\operatorname{Cone}(f)
    \longrightarrow\cc{C}_H
    \xrightarrow{e}\cc{C}_V[1],
    \qquad
    0\neq e\in HF^1(\cc{C}_H,\cc{C}_V).
\]

This suggests joining the two tori by
graded $0$-surgery at $r$.
Topologically, the result is their connected
sum, a surface of genus two.
We use this surgery picture to predict
the intersections of the terminal cokernel
with $L_0$ and $L_1$.
The proposed web evolution will seek an
imaginary-special representative of the
same derived object.

The joining point $r$ lies on neither
input torus.
Thus the surgery can be supported away
from $L_0\cup L_1$, preserving their
intersections with the two pieces.

For $L_0$, the intersection with $C_H$
is the single point
\[
    Q=((0,0),(1/2,0)),
\]
while its intersection with $C_V$ is
the circle
$\Gamma_0=\{(1/2,0)\}\times L_0^2$.
A small Hamiltonian perturbation replaces
this circle by two transverse intersections,
labelled $m_0,M_0$ according to the minimum
and maximum of a Morse function on $\Gamma_0$.

For $L_1$, the intersection with $C_V$
is the single point
\[
    R=((1/2,1/2),(0,0)),
\]
while its intersection with $C_H$ is
the circle
$\Gamma_1=L_1^1\times\{(1/2,1/2)\}$.
Perturbing this circle in the same way
gives two intersections $m_1,M_1$.

The product gradings give the following
Maslov indices, taken from the input brane
to the surgery surface:
\begin{center}
    \begin{tabular}{c|cc}
        Input      & Degree $0$ & Degree $1$ \\ \hline
        $\cc{L}_0$ & $Q,m_0$    & $M_0$      \\
        $\cc{L}_1$ & $m_1$      & $R,M_1$
    \end{tabular}
\end{center}
This is the intersection pattern we will
look for in the proposed terminal surface.

\begin{figure}[htbp]
    \centering
    \includegraphics[width=0.7\textwidth]{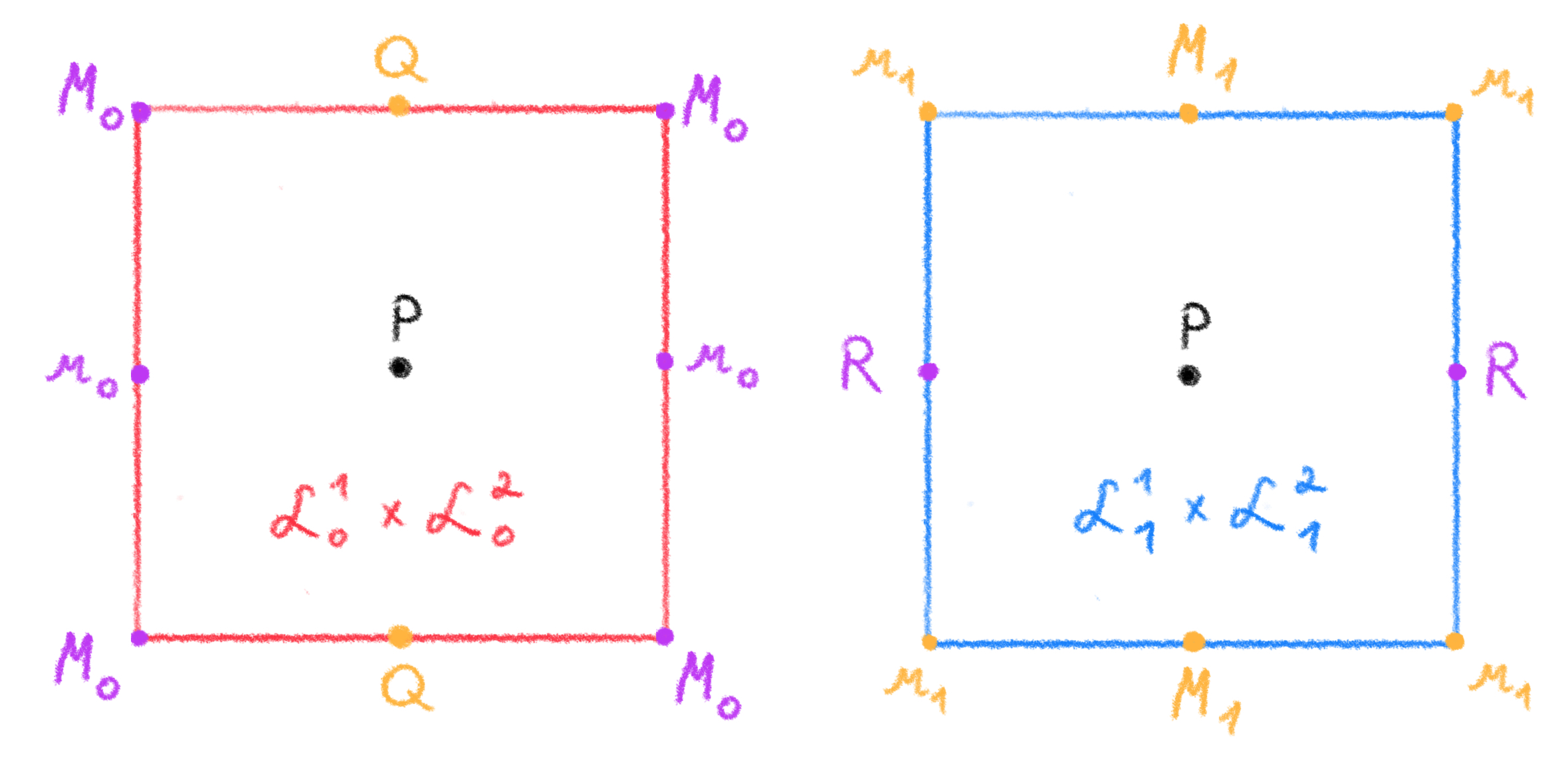}
    \caption{
        The product tori $L_0$ (red) and $L_1$ (blue),
        each represented by a square with opposite
        edges identified.
        Their unique intersection $P$ is marked
        in black in both panels.
        Intersections with $C_H$ are marked in orange,
        and those with $C_V$ in purple.
        The isolated intersections are
        $Q\in L_0\cap C_H$ and $R\in L_1\cap C_V$.
        The labels $m_0,M_0$ and $m_1,M_1$ mark
        the minima and maxima chosen on the
        intersection circles $L_0\cap C_V$ and
        $L_1\cap C_H$, respectively.
        Repeated boundary labels represent the
        same point on the torus.
    }
    \label{fig:abelian-input-intersections}
\end{figure}

The genus also follows directly from
the mirror algebra.
Write
$\chi(\cc{E},\cc{F})
    =\sum_k(-1)^k
    \dim_\Lambda\operatorname{Ext}^k(\cc{E},\cc{F})$
for the Euler pairing.

\begin{lemma}[The predicted genus]
    \label{lem:abelian-mirror-invariants}
    \label{lem:abelian-surgery-genus}
    We have
    \[
        D^2=2,
        \qquad
        \chi(\cc{C},\cc{C})=-2.
    \]
    Consequently, any connected embedded
    unobstructed rank-one brane representing
    $\cc{C}$ has genus two.
\end{lemma}

\begin{proof}
    Both components of $D$ have
    self-intersection zero, and they meet once.
    Hence $D^2=2V\cdot H=2$.
    The K\"unneth formula gives
    $\chi(\check X,\cc{O}_{\check X})=0$
    and
    $\chi(\check X,\cc{G})
        =\chi(\check X,\cc{G}^{-1})=1$.
    Using
    $[\cc{C}]=[\cc{G}]-[\cc{O}_{\check X}]$
    and additivity of the Euler pairing,
    we obtain
    $\chi(\cc{C},\cc{C})=0-1-1+0=-2$.
    For a connected embedded rank-one
    representative of genus $g$, the
    self-Floer Euler characteristic is
    $2-2g$.
    Thus $2-2g=-2$, so $g=2$.
\end{proof}

We therefore have two concrete predictions
for the web evolution: the intersection
pattern listed above, and a connected
terminal surface of genus two.

\subsection{The proposed web evolution}
\label{subsec:abelian-web-evolution}

We now describe a proposed web evolution whose
topology and intersections agree with the mirror
calculation. The initial webbing is a family of
imaginary special Lagrangian cylinders
$S^1\times[0,1]$ emerging from $P$, with one boundary
circle on each of $L_0$ and $L_1$.
These circles expand, sweeping out discs in the
two input tori, until they meet their periodic
translates.

\begin{figure}[htbp]
    \centering
    \includegraphics[width=\textwidth]{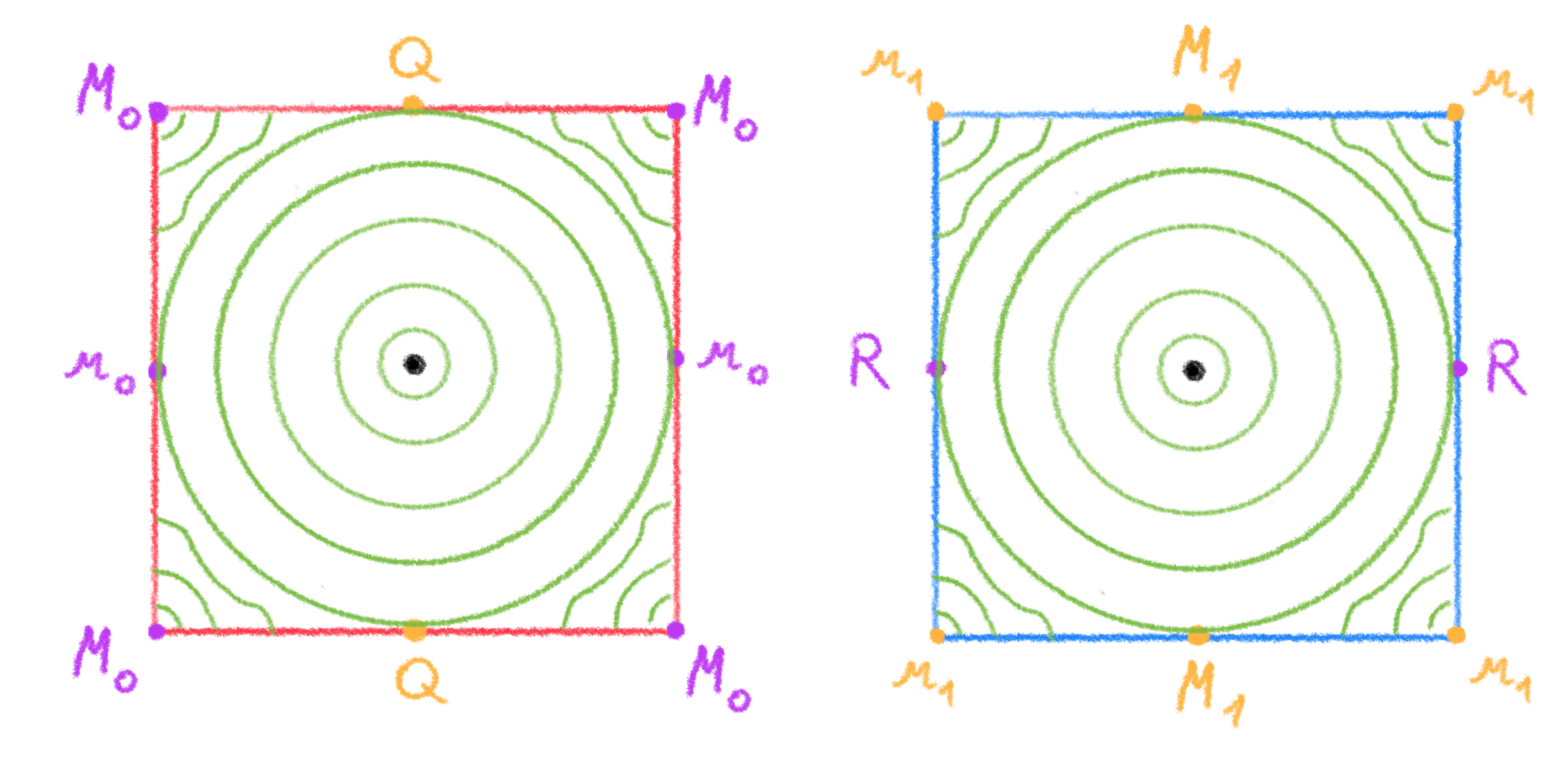}
    \caption{
        Boundary traces of the proposed web evolution
        on $L_0$ (left) and $L_1$ (right).
        Green curves indicate successive positions
        of the moving boundary.
        Starting at $P$, the boundary expands and
        undergoes periodic self-intersections at
        $Q,m_0$ on $L_0$ and at $R,M_1$ on $L_1$.
        After boundary surgery, the remaining circle
        contracts to $M_0$ and $m_1$, respectively.
        Opposite edges are identified; in particular,
        the four corners in each panel represent
        a single point.
    }
    \label{fig:abelian-boundary-traces}
\end{figure}

Figure~\ref{fig:abelian-boundary-traces} shows the
evolution from the perspective of each input torus.
In the symmetric square-torus model, the initial
boundary circles are round in the coordinates
of the fundamental squares.
On $L_0$, the first periodic self-intersections
therefore occur simultaneously at $Q$ and $m_0$.
Similarly, on $L_1$, they occur simultaneously
at $R$ and $M_1$.
The opposite edges of each square are identified,
so each contact joins two portions of the same
boundary circle.

We continue through these contacts by boundary
surgery. To understand the topology, imagine
separating the two simultaneous surgeries on
each side. The first splits one boundary circle
into two. The second joins these two circles,
adding a handle to the front and leaving one
boundary circle on that side.
Thus the two surgeries on $L_0$ add one handle,
and the two surgeries on $L_1$ add another.
Starting from a cylinder, we obtain a connected
surface of genus two with two boundary circles.

After these surgeries, the unswept region of
each input torus is a disc.
Its boundary appears in
Figure~\ref{fig:abelian-boundary-traces}
as four arcs near the corners of the fundamental
square. These arcs form a single circle after
the edges are identified.
The remaining boundary circles contract to
$M_0$ on $L_0$ and to $m_1$ on $L_1$.
The two collapses cap the boundary components
of the front, producing a closed surface of
genus two.

\begin{figure}[htbp]
    \centering
    \includegraphics[width=\textwidth]{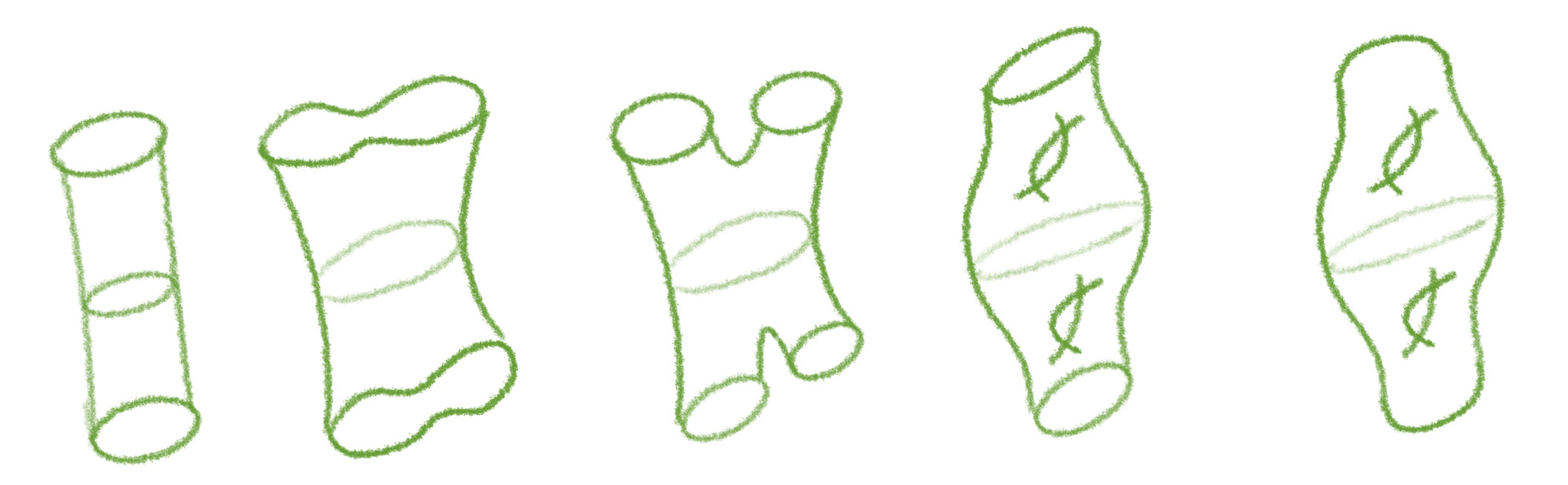}
    \caption{
        Schematic topology of the web-front,
        read from left to right.
        The first two surfaces are cylinders.
        The third shows the result of one boundary
        surgery on each side: a sphere with four
        boundary circles.
        The fourth follows the second surgery on
        each side and has genus two with two
        boundary circles.
        Capping these circles gives the closed
        genus-two surface in the final drawing.
        For illustration, the simultaneous surgeries
        are drawn in two stages.
    }
    \label{fig:abelian-front-topology}
\end{figure}

The boundary events also account for the predicted
intersections with the input tori.
The tori themselves remain fixed while the front
changes. In the proposed continuation, the new
bands and caps meet the corresponding input torus
near the surgery points, leaving intersections
after the moving boundary has passed.
On $L_0$, these intersections are $Q,m_0,M_0$;
on $L_1$, they are $R,M_1,m_1$.
This recovers the intersection pattern described
in Subsection~\ref{subsec:abelian-terminal-surface}.

At the end, both input tori have been exhausted
and no residues remain.
The closed front is inactive and carries label
$1$, so it is the candidate terminal cokernel
$K^*$; there is no kernel component.
With compatible brane transfers, it represents
$\operatorname{Cone}(f)$.
Its genus agrees with the mirror prediction of
Lemma~\ref{lem:abelian-mirror-invariants}.

\bibliographystyle{amsalpha}
\bibliography{ref}

\end{document}